\documentclass[11pt,a4paper,twoside,openright]{report}

\usepackage{color}
\definecolor{MyDarkBlue}{rgb}{0.15,0.25,0.65}

\usepackage[british]{babel}
\usepackage{fancyhdr}
\usepackage{lastpage}
\usepackage[margin=2.5cm, includehead, includefoot]{geometry}  

\usepackage{soul} 
\renewcommand{\underline}[1]{\ul{#1}}

\usepackage{tikz}
\usepackage{verbatim}

\usepackage{makeidx}

\let\fn\footnote
\renewcommand{\footnote}[1]{\linespread{1.1}\fn{#1}\linespread{1.29}}

\makeatletter
\renewcommand\@idxitem{\par\hangindent 40\p@}
\renewcommand\subitem{\@idxitem \hspace*{20\p@}}
\renewcommand\subsubitem{\@idxitem \hspace*{30\p@}}
\renewcommand\indexspace{\par \vskip 10\p@ \@plus5\p@ \@minus3\p@\relax}

\makeatletter\renewcommand{\section}{\@startsection {section}{1}{\z@}{-3.5ex plus -1ex minus -.2ex}{2.3ex plus .2ex}{\bf\mathversion{bold}\large }}

\makeatletter\renewcommand{\subsection}{\@startsection{subsection}{2}{\z@}{-3.25ex plus -1ex minus  -.2ex}{1.5ex plus .2ex}{\bf\mathversion{bold} }}

\makeatletter\renewcommand{\subsubsection}{\@startsection{subsubsection}{3}{\z@}{-3.25ex plus -1ex minus -.2ex}{1.5ex plus .2ex}{\it }}

\renewcommand{\thesection}{\arabic{chapter}.\arabic{section}.}
\renewcommand{\thesubsection}{\arabic{chapter}.\arabic{section}.\arabic{subsection}.}

\renewcommand{\@seccntformat}[1]{\@nameuse{the#1}~}

\renewcommand{\theequation}{\thechapter.\arabic{equation}}
\makeatletter \@addtoreset{equation}{chapter}

\makeatletter \@addtoreset{figure}{chapter}

\renewcommand*\l@chapter[2]{%
  \ifnum \c@tocdepth >\m@ne
    \addpenalty{-\@highpenalty}%
    \vskip 0.3em \@plus\p@   
    \setlength\@tempdima{1.2em}
    \begingroup
      \def\numberline##1{\hb@xt@\@tempdima{##1.\hfil}}
      \parindent \z@ \rightskip \@pnumwidth
      \parfillskip -\@pnumwidth
      \leavevmode 
      \advance\leftskip\@tempdima
      \hskip -\leftskip
      #1\nobreak\mdseries
      \leaders\hbox{$\m@th
        \mkern \@dotsep mu\hbox{.}\mkern \@dotsep
        mu$}\hfill
      \nobreak\hb@xt@\@pnumwidth{\hss #2}\par
      \penalty\@highpenalty
    \endgroup
  \fi}
  
\renewcommand*\l@section{\vskip 0.3em \@plus\p@ \@dottedtocline{2}{1.2em}{2em}}
\renewcommand*\l@subsection{\vskip 0.3em \@plus\p@ \@dottedtocline{3}{3.2em}{2.8em}}

\renewcommand\tableofcontents{%
    \section*{\large\contentsname
        \@mkboth{%
          \MakeUppercase\contentsname}{\MakeUppercase\contentsname}}%
       {\baselineskip=15pt plus 2pt minus 1pt
    \@starttoc{toc}}%
}

\newcommand{\appendices}{
\chapter*{Appendices}\label{appendices}
\clearpage{\pagestyle{empty}\cleardoublepage}
\setcounter{chapter}{0}
\renewcommand{\thechapter}{\Alph{chapter}}
\addcontentsline{toc}{chapter}{Appendices}
\setcounter{equation}{0}
\makeatletter
\renewcommand{\theequation}{\Alph{chapter}.\arabic{equation}}
\renewcommand{\thesection}{\Alph{chapter}.\arabic{section}.}
\renewcommand{\thesubsection}{\Alph{chapter}.\arabic{section}.\arabic{subsection}.}
\@addtoreset{equation}{chapter}
\makeatother
}

\newcommand{\clearemptydoublepage}{\newpage\phantom{}\thispagestyle{empty}\newpage}

\usepackage{amsmath,amssymb,amsthm}
\usepackage{amsfonts}
\usepackage{mathrsfs}
\usepackage{mathabx}
\usepackage{bbm}
\usepackage{bm}
\usepackage{epsfig,rotating}
\usepackage{enumerate}

\newtheorem{de}{Definition}
\numberwithin{de}{chapter}
\newtheorem{lem}[de]{Lemma}
\newtheorem{thm}[de]{Theorem}
\newtheorem{pro}[de]{Proposition}
\newtheorem{cor}[de]{Corollary}
\newtheorem{re}[de]{Remark}
\newtheorem{exa}[de]{Example}

\tikzstyle{mybox} = [draw=black!50, fill=gray!20,  thick,
    rectangle, rounded corners, inner sep=10pt, inner ysep=20pt]
\tikzstyle{fancytitle} =[fill=black!50, text=white]

\newcounter{exe}
\numberwithin{exe}{chapter}
\renewcommand{\theexe}{\thechapter.\roman{exe}}

\usepackage[Bjornstrup]{fncychap}
\usepackage{cmbright}  

\usepackage[linktocpage=true]{hyperref}
\hypersetup{
colorlinks=true,
citecolor=MyDarkBlue,
linkcolor=MyDarkBlue,
urlcolor=MyDarkBlue,
pdfauthor={Marcel Dengler}, 
pdftitle={New regularity and uniqueness results in the multidimensional Calculus of Variations}, 
breaklinks=true
}

\usepackage{lmodern}
\usepackage{graphicx}
\graphicspath{ {/home/user/documents/PhDThesisMD/} }
\usepackage{enumitem}
\usepackage[utf8]{inputenc}
\usepackage{csquotes}
\usepackage{textcomp}
\usepackage{amsfonts, ifthen, amsmath, amssymb, amsthm}
\usepackage{esint}
\usepackage{lmodern}

\usepackage{enumitem}

\usepackage{color}

\renewcommand{\[}{\begin{equation*}}
\renewcommand{\]}{\end{equation*}}
\newcommand{\grad}{\ensuremath{\nabla}}
\newcommand{\p}{\ensuremath{\partial}}
\newcommand{\ve}{\ensuremath{\varepsilon}}
\newcommand{\vp}{\ensuremath{\varphi}}
\newcommand{\al}{\ensuremath{\alpha}}
\newcommand{\be}{\ensuremath{\beta}}
\newcommand{\ga}{\ensuremath{\gamma}}
\newcommand{\ta}{\ensuremath{\tau}}
\newcommand{\ka}{\ensuremath{\kappa}}
\newcommand{\la}{\ensuremath{\lambda}}
\renewcommand{\th}{\ensuremath{\theta}}
\newcommand{\om}{\ensuremath{\omega}}

\renewcommand{\ss}{\ensuremath{\subset}}
\newcommand{\sm}{\ensuremath{\setminus}}
\newcommand{\mt}{\ensuremath{\mapsto}}
\newcommand{\ra}{\ensuremath{\rightarrow}}
\newcommand{\Llra}{\ensuremath{\Longleftrightarrow}}
\newcommand{\Lra}{\ensuremath{\Longrightarrow}}
\newcommand{\rhu}{\ensuremath{\rightharpoonup}}
\newcommand{\mb}[1]{\ensuremath{\;\mbox{#1}\;}}
\newcommand{\cd}{\ensuremath{\cdot}}
\newcommand{\ot}{\ensuremath{\otimes}}
\newcommand{\ti}{\ensuremath{\times}}

\newcommand{\xb}{\ensuremath{\bar{x}}}

\newcommand{\Id}{\ensuremath{\mbox{Id}}}

\newcommand{\wt}[1]{\widetilde{#1}}
\newcommand{\ol}[1]{\overline{#1}}

\newcommand{\cof}{\textnormal{cof}\;}
\renewcommand{\div}{\textnormal{div}\;}
\newcommand{\adj}{\textnormal{adj}}

\newcommand{\dor}{\ensuremath{\dot{r}}}
\newcommand{\doz}{\ensuremath{\dot{z}}}
\newcommand{\ddoz}{\ensuremath{\ddot{z}}}
\newcommand{\dod}{\ensuremath{\dot{d}}}

\newcommand{\ddor}{\ensuremath{\ddot{r}}}

\newcommand{\dist}{\textnormal{dist}\;}
\newcommand{\supp}{\textnormal{supp}\;}

\newcommand{\R}{\mathbb{R}}

\newcommand{\N}{\mathbb{N}}

\newcommand{\bb}[1]{\ensuremath{{\mathbb{#1}}}}

\newcommand{\D}{\Delta}
\renewcommand{\d}{\delta}
\newcommand{\A}{\mathcal{A}}
\newcommand{\F}{\mathcal{F}}

\newcommand{\W}{\mathcal{W}}
\newcommand{\La}{\mathcal{L}}

\newcommand{\Om}{\Omega}

\allowdisplaybreaks
\usepackage[nosort]{cite}  
\usepackage[bulletsep]{collref}  
\begin{document}


\begin{titlepage}

\setcounter{page}{0}
\renewcommand{\thefootnote}{\fnsymbol{footnote}}

\begin{center}

    {\huge\textbf{\mathversion{bold}New regularity and uniqueness results in the multidimensional Calculus of Variations}\par}

    \vspace{2.5cm}

    {\Large\bf Marcel Dengler}

    \vspace{2cm}

    Thesis submitted to the University of Surrey\\ for the degree of Doctor of Philosophy

    \vspace{2cm}

    {\it\large Department of Mathematics\\ University of Surrey\\ Guildford GU2 7XH, United Kingdom}

    \vspace{2.5cm}

    \href{http://www.surrey.ac.uk/}{\includegraphics[width=3cm]{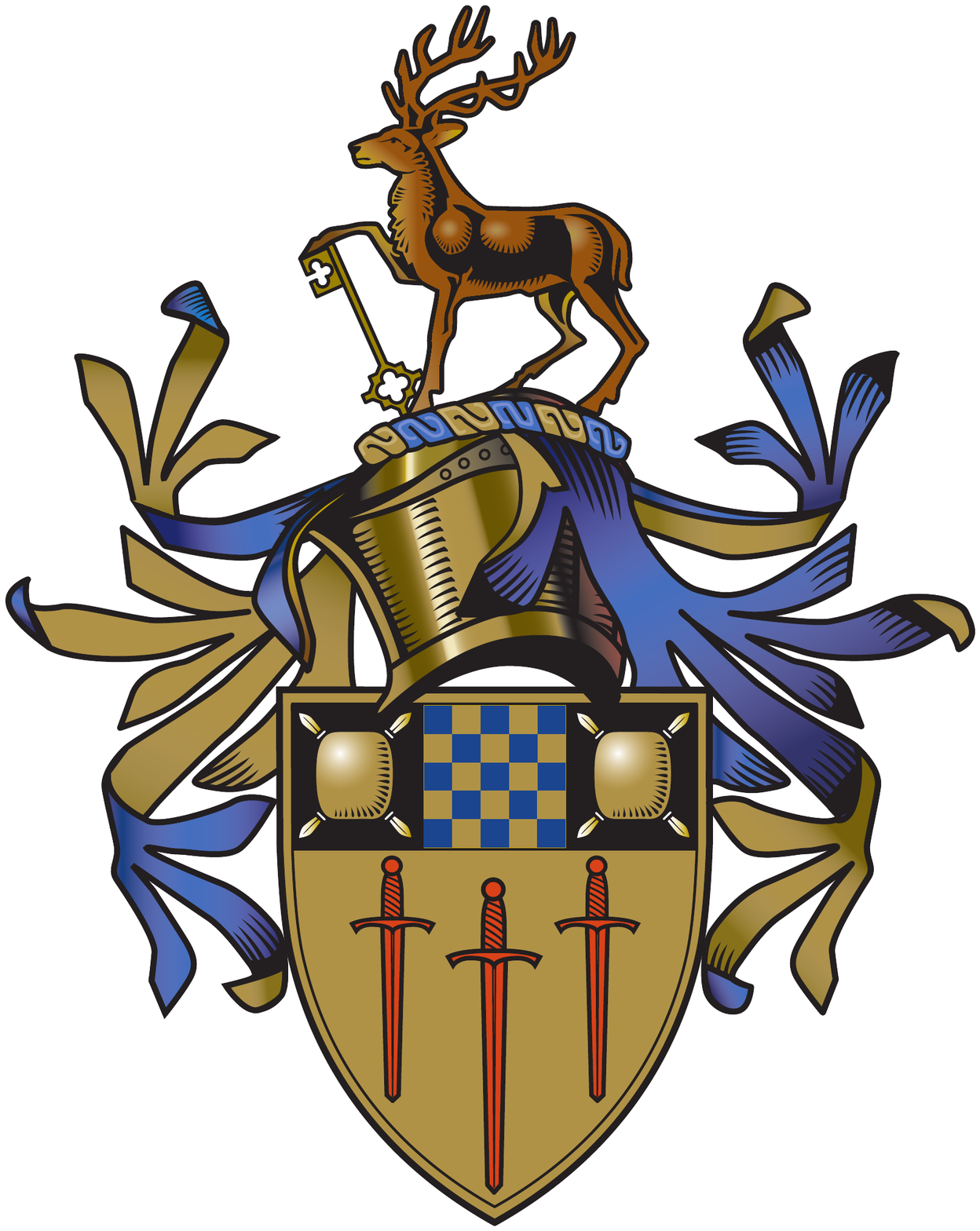}}

    \vfill
    
    Copyright \copyright\ 2021 by Marcel Dengler. All rights reserved.\\
    E-mail address: \href{mailto:m.dengler@surrey.ac.uk}{\ttfamily m.dengler@surrey.ac.uk}

\end{center}

\setcounter{footnote}{0}\renewcommand{\thefootnote}{\arabic{thefootnote}}

\end{titlepage}
\clearemptydoublepage
\pagenumbering{roman}
\section*{Scientific abstract}

In the first part of the Thesis we develop a Regularity Theory for a polyconvex functional in compressible elasticity. In particular, we consider energy minimizers/stationary points of the functional
\begin{equation}I(u)=\int\limits_\Om{\frac{1}{2}|\grad u|^2+\rho(\det\grad u)\;dx},\label{eq:SA.1.1}\end{equation}
where $\Omega\ss\R^2$ is open and bounded, $u\in W^{1,2}(\Om,\R^2)$ and $\rho:\R\ra\R_0^+$ smooth and convex with $\rho(s)=0$ for all $s\le0$ and $\rho$ becomes affine when $s$ exceeds some value $s_0>0.$ Additionally, we may impose boundary conditions.\\

The first general result we will establish is that every stationary point needs to be locally Hölder-continuous. 
Secondly, we prove that if the growth of $\rho$ is `small' s.t.\! the integrand is still uniformly convex, then all stationary points have to be in $W_{loc}^{2,2}.$ Next, a higher-order regularity result is shown. 
We show that all stationary points that are additionally of class $W_{loc}^{2,2}$ and whose Jacobian is Hölder-continous are of class $C_{loc}^{\infty}.$ In particular, these results show that all stationary points have to be smooth for $\rho'$ `small' enough. \\

The theory described above works for fairly general domains and boundary conditions. 
We specify those by restricting to the unit ball, and we consider M-covering maps, which take the unit sphere to itself, covering the image $M$ times in the process, on the boundary. Under these circumstance, we construct radial symmetric M-covering stationary points to the functional, as given in \eqref{eq:SA.1.1}, which are at least of class $C^1.$ In certain situations, depending mainly on the behaviour of $\rho$ and the stationary point itself, we are even able to guarantee maximal smoothness. \\

In the second part, we will concentrate on uniqueness questions in various situations of finite elasticity. Starting in incompressible elasticity, the central point is to show 
that for problems with uniformly convex integrands with "small" pressure, a unique global minimizer can be guaranteed. We make use of this statement by considering various examples and applications. 
One such application is the construction of a counterexample to regularity. Indeed, on the unit ball and for smooth boundary conditions we give a non-autonomous uniformly convex functional $f(x,\xi)$ depending smoothly on $\xi$ however discontinously on $x,$ where the unique global minimizer is Lipschitz but no better. Then we generalise the main result in various ways, for instance, we show that if the pressure is too large to guarantee uniqueness in the full class of admissible maps, one can still guarantee uniqueness up to the first Fourier-modes. 
Lastly, we discuss analogous statements for polyconvex integrands in compressible elasticity. \\

\vfill
\noindent
Keywords and AMS Classification Codes: Calculus of Variations, Polyconvexity, Regularity, Uniqueness
\vspace{1cm}

\clearemptydoublepage
\section*{Lay summary}

In this thesis we consider rigorous mathematical models of Material Science, focussing in particular on Elastostatics.\\

In the first part of this work, we are interested in elastic situations and the corresponding deformed configuration. For specific compressible settings, where a deformation may change the 'volume' of the material, 
we will investigate the question, how smooth does the profile of a deformed configuration need to be?\\

Secondly, it is well known that for the same setup, that is, the same material, the same domain, the same boundary data, etc., in elastic situations only one configuration is possible in other situations multiple configurations 
are possible, i.e. non-uniqueness can occur in general. Although it is easy to imagining both types of settings, it turns out to be a hard task to give criteria which guarantee that for a certain setup only one configuration is possible.    
This is exactly what we are going to focus on in the second part of this thesis. We will discuss criteria for both compressible and incompressible situations, in which only deformations that preserve the volume are possible, such that a unique configuration can be guaranteed. 
We then give various applications. Surprisingly, this leads to an incompressible situation where the configuration clearly possesses an edge, i.e. the deformed configuration possesses a rough profile. \\

\clearemptydoublepage
\section*{Acknowledgements}

This work has profited to a huge amount from my supervisor Jonathan Bevan. I can't describe how grateful I am for the many enlightening discussions, the countless insightful comments and remarks, guiding this project in the right direction.\\ 

I'm also appreciative to my cosupervisor James Grant for friendly discussions and support.\\

To Bin Cheng and Stephen Gourley, for being more than just my Confirmation Examiners, but also for great advice, feedback and guidance through the course of my PhD. We want to emphasise that this might be one of the latter works Stephen Gourley, who sadly passed away recently, has greatly influenced.\\

I want to thank all the people of the Department of Mathematics for providing a friendly welcome, when I first joined the university and creating a stimulating environment ever since. In particular, to my fellow PhD's for many inspiring discussions about various topics.\\

I'm grateful to the Engineering \& Physical Sciences Research Council (EPRSC) for funding my work.\\

Lastly, thanks to the Library-Team for providing a good service, the IT-Support for always offering help and efficient solutions and to the RDP-Team for offering very informative workshops.

\vspace{0.5cm}

\begin{flushright}
Guildford, 24.07.2021\\
Marcel Dengler
\end{flushright}
\clearemptydoublepage
\section*{Declaration of Originality}
This thesis and the work to which it refers are the results of my own efforts. Any ideas, data, images or text resulting from the work of others(whether published or unpublished) are fully identified as such within the work and attributed to their originator in the text, bibliography or in footnotes. This thesis has not been submitted in whole or in part for any other academic degree or professional qualification. I agree that the University has the right to submit my work to the plagiarism detection service TurnitinUK for originality checks. Whether or not drafts have been so-assessed, the University reserves the right to require an electronic version of the final document(as submitted) for assessment as above.

\begin{flushright}
Guildford, 24.07.2021\\
Marcel Dengler
\end{flushright}
\clearemptydoublepage
\tableofcontents 
\bigskip
\bigskip
\hrule
\bigskip
\bigskip
\clearemptydoublepage
\pagenumbering{arabic}



\pagestyle{fancy}
\renewcommand{\sectionmark}[1]{\markboth{\thesection~#1}{}} 
\fancyhead[RO,LE]{\thepage}
\fancyhead[CO,CE]{}
\fancyhead[LO,RE]{\leftmark}
\fancyfoot{}
\chapter{Introduction}
Our project originates from Solid Mechanics and Elastostatics. The standard situation in
these areas is to study deformations, that is, maps $u:\Om\ra\R^n$ which act on a set $\Om\ss\R^m,$ representing a material in a reference to configuration, and which obey certain constraints and requirements.
Of course, which deformation will arise in
which situation will depend on what material the object is made of and how it behaves. Such a situation is usually described by a given energy functional of the form
\begin{equation}I(u)=\int\limits_\Om {W(x,u(x),\grad u(x))\;dx},
\label{eq:Intro.1.1}
\end{equation}
where $\Om\ss\R^m$ a measurable set, representing the undeformed object. The deformation $u:\Om\ra\R^n$ has to be in a suitable function space, 
where suitable means that at least the integral $I$ needs to make sense. Since, $W:\Om\ti \R^{n}\ti\R^{m\ti n}\ra\ol{\R}$ depends on the gradient of $u,$ usually this will be a Sobolev space. 
Additionally, deformations may be required to obey boundary conditions.\\

The working hypothesis of Elastostatics is that an elastic body will deform in order to minimize its energy, as defined in \eqref{eq:Intro.1.1}. So the general problem we are interested in is finding solutions to
\begin{equation}\inf\limits_{u\in\A} I(u).
\label{eq:Intro.1.2}
\end{equation}
The main questions that then need to be addressed and answered include: Does a minimizer exist and if so is it unique? Do the minimizers satisfy specific properties, for instance, how regular (smooth) are those minimizers? Can one find upper and lower bounds on \eqref{eq:Intro.1.2}?
Can one find the explicit solution(s) to \eqref{eq:Intro.1.2} and compute its numerical value?
In this thesis, we address a fair share of all these types of questions.\\

The problem given in \eqref{eq:Intro.1.2} depends to a great extend on the behaviour of the integrand $W.$  While minimization problems with uniformly convex integrands are well understood, 
many problems arising from Elastostatics possess instead polyconvex integrands, making it necessary to extend the theory to such functions.\footnote{see the next section for precise definitions of these notions of Semiconvexities.}
Indeed, well-known physical models (like the neo-Hookean model, the Mooney–Rivlin model or the Ogden models, all of which describe hyperelastic materials\footnote{for more on these models see \cite{R18}.}) have polyconvex integrands of the form
\vspace{-0.35cm}
\begin{align}W(\xi)=f(\xi)+h(\det \xi)\;\mb{for all}\;\xi\in\R^{m\ti m}, \label{eq:Intro1.3}\end{align}
where $m\ge2,$ $f$ is uniformly convex, $h$ is convex, $h(d)\ra+\infty$ when $d\searrow0,$ and $h(d)=+\infty$ if $d\le0.$
The meaning of the latter properties can easily be explained. The property $h(d)\ra+\infty$ when $d\searrow0$ describes that infinite energy is required to compress a m-dimensional material element to zero m-dimensional volume. 
Property $h(d)=+\infty$ if $d\le0$ forbids any deformation where the object needs to penetrate itself. As a consequence, if one is looking for minimizers one can usually restrict the search to deformations $u\in \A,$ where $I(u)<+\infty$ s.t. $\det \grad u>0$ a.e. in $\Om.$\\

Roughly speaking, polyconvex functionals are energies for which in addition to the general size of the change of a given deformation described by $\grad u$ as measured by $\int\limits_\Om {f(\grad u(x))\;dx}$ in \eqref{eq:Intro1.3} for example, also the volume change of the deformation described by $\det\grad u$ is governed, ignoring the $(x,u)$ dependence for now.\\

We will refer in the future to models where such a condition i.e. where the integrand $W$ is allowed to take on the value $+\infty,$ is in place as a Nonlinear Elasticity (NLE) model. In contrast we will refer to models where $W<\infty$ as a Finite Elasticity (FE) model. 
Firstly, in this thesis we are only concerned with the latter type, although we will compare our work with works that address the NLE scenario. Secondly, it is important to mention that physically relevant models can only be described by NLE models.
FE models don't describe real materials. However, due to the properties of $h$ described above, NLE models are technical and difficult to handle analytically. Hence one studies FE models first and maybe some of the insights can be transferred to NLE models.\\

To any Finite Elasticity scenario, where we only consider admissible maps that are measure-preserving, we refer to as an incompressible situation. In contrast we call any FE model, where the latter constraint is not present, a compressible one. Recall, that mathematically we can represent the property of a suitable injective smooth map (for members of a Sobolev-space one might approximate) to be measure-preserving by $\det\grad u=1$ a.e. in $\Om.$ Indeed, for a measure-preserving injective Lipschitz map $u:\Om\ss\R^m\ra u(\Om)\ss\R^m$ by the injectivity, the elementary properties of integrals, and the area formula we get
$\La^m(u(\Om'))=\int\limits_{u(\Om')}{dy}=\int\limits_{\Om'}{\det\grad u(x)\;dx}$
and hence $\det\grad u=1$ a.e. in $\Om$ implies $\La^m(\Om')=\La^m(u(\Om'))$ for any $\Om'\ss\Om.$ It is standard procedure to replace the measure-preservability condition $\La^m(\Om')=\La^m(u(\Om'))$ for any $\Om'\ss\Om$ by the constraint $\det\grad u=1$ a.e.\!\! in $\Om$ even in situations, where $u$ might not be injective and therefore not necessarily measure-preserving in the above sense.\\

\section{Semiconvexities and their relations}
\label{sec:1.1}
It seems unavoidable in this area of research to give a brief overview on Semiconvexities.\footnote{For complete treatises on different types of convexities and their relations, see Dacorogna's famous book \cite{D08} or the PhD Thesis \cite{K16} of Käbisch. The latter one also introduces a whole new type of semiconvexity (n-polyconvexity).}\\

From a mathematical point of view, Semiconvexities arise when generalising the notion of convexity. Recall, that a function $f:\R^{N\ti n}\ra\ol{\R}$ is called convex if
\[f(t\xi+(1-t)\xi')\le tf(\xi)+(1-t)f(\xi')\]
holds for every $\xi,\xi'\in\R^{N\ti n}$ and $t\in[0,1].$\\

In one dimension it is well known that if $f$ is regular enough, $f\in C^2,$ then
\[f \;\mb{convex}\; \Longleftrightarrow f''(x)\ge0\;\mb{for all}\;x\in\R.\]
A similar characterisation is true in higher dimensions:
\begin{align}f \;\mb{convex}\; \Longleftrightarrow D_\xi^2f(\xi)M\cd M\ge0\;\mb{for all}\;M\in\R^{N\ti n},  \xi\in \R^{N\ti n}. \label{eq:Convexity 1.1}\end{align} 
However, this turns out to be not the most natural condition (and a fairly restrictive one, too!) in higher dimensions. A more natural quantity, originating from the 2nd-variation, gives rise to a first new kind of convexity:

\begin{de}[Rank-one convexity]
The function $f:\R^{N\ti n}\ra\ol{\R}$ is rank-one convex if
\[f(\la\xi+(1-\la)\xi')\le\la f(\xi)+(1-\la)f(\xi')\]
for all $\la\in[0,1]$ and $\xi,\xi'\in \R^{N\ti n}$ s.t. $\mbox{rank}\;(\xi-\xi')\le1.$ 
\end{de}
Obviously, $f$ convex implies that $f$ is rank-one convex, since for rank-one convexity the convexity inequality only needs to hold on `rank-1 lines'.\\

Again, assuming $f\in C^2,$ one can see that rank-one convexity agrees with the Legendre-Hadamard condition i.e.
\[f\; \mb{rank-one convex} \;\Longleftrightarrow D_\xi^2f(\xi)(a\ot b)\cdot(a\ot b)\ge0, \]

for all $\xi\in\R^{N\ti n},$ $a\in \R^N$ and $b\in\R^n.$\footnote{In components this reads

\[\sum\limits_{i,j=1}^N\sum\limits_{\al,\be=1}^n\frac{\p^2f(\xi)}{\p\xi_\al^i\p\xi_\be^j}a^ia^jb_\al b_\be\ge0.\]}
Since, convexity implies rank-one convexity, clearly the Legendre-Hadamard condition needs to be weaker then the one given in \eqref{eq:Convexity 1.1} and it surely is. An example of a rank-one convex function which is not convex is given by the determinant $\xi \mt\det\xi$ if $\xi\in\R^{N\ti N},$ see the discussion at the end of this paragraph.\\

The notion of rank-one convexity turns out to be fairly weak. The existence of a minimizer to a minimization problem with a rank-one convex integrand can in general not be guaranteed. Moreover, regularity results for general rank-1 convex functions are rare as well. They either come with strong structural assumptions or apply to small classes of functions.\\ 

Morrey, in 1952, motivated by the aforementioned lack of existence, introduced another notion:

\begin{de}[Quasiconvexity]
A Borel-measurable function $f:\R^{N\ti n}\ra \R,$ is called quasiconvex, if
\[f(\xi)\le\fint\limits_{B_N(0,1)}{f(\xi+\grad\psi(x))\;dx}\]
for all $\xi\in \R^{N\ti n}$ and $\psi\in W_0^{1,\infty}(B_N(0,1),\R^n).$
\end{de}

\begin{re}
We want to point out that quasiconvexity can only be defined for functions of the form $f:\R^{N\ti n}\ra \R,$ this agrees with the majority of the texts in this field including the standard texts \cite{D08} and \cite{R18} and the PhD Thesis \cite{K16}. 
However, sometimes the notion of quasiconvexity, it seems with a small abuse of notion, has been extend to functions with $f:\R^{N\ti n}\ra \ol{\R},$ see for instance the works \cite{BOP91MS}, \cite{Y06} and \cite{BY07}.
\end{re}

Quasiconvexity turns out to be weaker than convexity and stronger than rank-1 convexity; for a counterexample see Remark \ref{re:1.1.1.6.a}. Note that it has been defined in a completely different manner than the convexities we already introduced and the one which comes next. 
This is due to the fact that it was introduced purely to guarantee existence via the direct method of the Calculus of Variations 
of which it is a major part to show that the considered functional is weak lower semicontinuous\footnote{see \cite[\S 2.1]{R18} for the direct method and \cite[\S 5.5]{R18} for the statement that quasiconvex functionals are weak lower semicontinuous.}. 
This makes it a very hard task to compare quasiconvexity with the other notions. \\

The last important semiconvexity is due to Ball, in \cite{B77}. He introduced polyconvexity for the sake of modelling and describing phenomena in the theory of nonlinear elasticity. The notion is given by

\begin{de}[Polyconvexity] We call $f:\R^{N\ti n}\ra\ol{\R}$ polyconvex, if there exists a convex function $g:\R^{\ta(N,n)}\ra\ol{\R}$ s.t. 
\[f(\xi)=g(T(\xi))
\label{eq:A1.3}\]
where $T:\R^{N\ti n}\ra\R^{\ta(n,N)}$ represents a list of the minors of $\xi,$ i.e.
\[T(\xi)=(\xi,\adj_2\xi,\ldots,\adj_{n\wedge N}\xi),\]
\end{de}
where 
$\ta(n,N):=\sum\limits_{k=1}^{n\wedge N}\binom{N}{k}\binom{n}{k}$ and $\adj_k\xi$ denotes the matrix of all $k\ti k-$minors of the matrix $\xi\in \R^{N\ti n}.$

For $N=n=2$ this reduces to $T(\xi)=(\xi,\det\xi)$ and we call $f:\R^{2\ti 2}\ra\ol{\R}$ polyconvex if there exists a convex function $g:\R^5\ra\ol{\R}$ s.t. $f(\xi)=g(\xi,\det\xi)$ for all $\xi\in\R^{2\ti2}.$ If $N=n=3,$ $f$ will additionally depend on the cofactor of the gradient, i.e. $T(\xi)=(\xi,\det\xi,\cof\xi).$\\

Polyconvexity is weaker than convexity but stronger than quasiconvexity. Therefore, we have the following account.

\begin{thm} [Thm 5.3.(i), \cite{D08}] (i) Assume $f:\R^{N\ti n}\ra \R$ then
\[f\;\mb{convex}\;\Lra\;f  \;\mb{polyconvex}\;\Lra\;f  \;\mb{quasiconvex}\;\Lra\;f  \;\mb{rank-one convex}.\]
(ii) Assume $f:\R^{N\ti n}\ra \ol{\R}$ then
\[f\;\mb{convex}\;\Lra\;f  \;\mb{polyconvex}\;\Lra\;f  \;\mb{rank-one convex}.\]

\end{thm}
\begin{re}\label{re:1.1.1.6.a}
The opposite implications fail in general. It is straight forward to see that not every polyconvex function is convex (counterexample: $f(\xi)=\det\xi$). 
Much more effort is necessary to show quasiconvexity does not imply polyconvexity (Alibert-Dacorogna-Marcellini example, see \cite[\S 5.3.8]{D08}).  Intriguingly, it has turned out that it is very difficult to compare rank-one convexity with quasiconvexity. In \cite{M52} Morrey conjectured that rank-one convexity does not necessarily imply quasiconvexity. In the 90's a counterexample was given by Sverak, in \cite{S92}, in the cases when $N\ge3,$ $n\ge2.$ The cases $N=2,$ $n\ge2$ remain open in general. However, partial results are available on subspaces of matrices and/or specific classes of functionals, see for example \cite{M99,CM03,HKL18, GMN17, N18, VMGN21}.
\end{re}

For a novel notion of polyconvexity for functionals depending only on the symmetric part of the gradient, see the recent article \cite{BKS19}.

\section{Introduction for chapters 2 and 3}
In the chapters 2 and 3 we focus on specific situations in compressible elasticity. In fact, we will consider the special class of polyconvex integrands with $n=m=2$ of the form
\begin{equation}W(\xi)=\frac{1}{2}|\xi|^2+\rho(\det \xi), \label{eq:Intro1.4}\end{equation}
where $\xi\in\R^{2\ti2}$ and $\rho:\R\ra\R$ is a smooth convex function, s.t. $\rho(s)=0$ for $s\le0$ and if $s\ge s_0>0$ then $\rho(s)=\ga s+\ka,$ where $0<\ga<\infty$ and $\ka\ge-\ga s_0.$ Hence, $\rho$ becomes affine if $s$ exceeds $s_0$. 
The behaviour of the functional depends on the parameter $\ga:$ if $0<\ga<1$ the functional is actually uniformly convex, while if $\ga\ge1$ the functional becomes genuinely polyconvex, i.e. polyconvex but not uniformly convex.\\
In chapter 2 we allow general domains and general boundary conditions ($u_0\in L^2(\Om, \R^2),$ see section \ref{sec:2.1} for more information). In chapter 3 we then consider the case, where $\Om$ is the unit ball and $u_0$ the M-covering map on the boundary.\\

\textbf{Partial regularity:} What can we say about the regularity of the stationary points and (local) minimizers corresponding to problem \eqref{eq:1.1} and the given conditions using established theory?\\
The integrand $W(\xi)=\frac{1}{2}|\xi|^2+\rho(\det \xi)$ is polyconvex by the convexity of $\rho.$ Therefore, although we would like to apply Evans's partial regularity theorem, see \cite{Ev86}, 
we find that the condition that the second derivative of $W$ is bounded, i.e. $|\grad_\xi^2W(\xi)|\le C$ for some $C>0$ independent of $\xi$ is not satisfied. Indeed, 
\begin{equation}
\grad_\xi^2W(\xi)=\Id_{2\times2}\ot \Id_{2\times2}+\rho
'(\det \xi)\grad_\xi(\cof \xi)+\rho''(\det \xi)\cof \xi\ot \cof \xi
\label{eq:Intro1.5}
\end{equation}
where for two matrices $A,B\in\R^{2\ti2},$ $A\ot B:=(A_{i,j}B_{k,l})_{1\le i,j,k,l\le2}.$ Then the first term of \eqref{eq:Intro1.5} and $\grad_\xi(\cof \xi)$ are constant matrices, 
for the last term it holds\footnote{see Section \ref{sec:2.1} for the definition of the cofactor and the identity $|\cof\xi|=|\xi|.$} $|\cof \xi\ot \cof \xi|=|\xi|^2$ and hence it grows quadratically. 
Therefore, Evans's partial regularity theorem is not applicable. Nevertheless, there is a related result by E. Acerbi and N. Fusco \cite{AF87}, 
which includes our case and states that there exists an open subset $\Sigma\subset\Omega$ s.t. $\La^2(\Omega\setminus\Sigma)=0,$ and the derivative of any stationary point/(local) minimizer $u\in \A$ is locally Hölder continuous on $\Sigma,$
i.e. $u\in C_{loc}^{1,\mu}(\Sigma,\R^2)$ for some $\mu\in(0,1].$\\

In this context, our goal for chapter 2 and 3 is to get more information on the set $\Sigma.$\\

\textbf{Plan for chapter 2:} We start in chapter 2 by giving a precise description of the problem and introducing the notation which is used in this thesis. 
Furthermore, for any $u_0\in L^2(\Om,\R^2)$ we show the existence of a global minimizer in the class $\A:=W_{u_0}^{1,2}(\Om,\R^2)$ of the functional \eqref{eq:Intro.1.1} with integrand, as defined in \eqref{eq:Intro1.4}, 
justifying our discussion. In the chapters 2 and 3 we will only consider this functional \eqref{eq:Intro.1.1} with integrand, as defined in \eqref{eq:Intro1.4}, so minimizers and stationary points are all wrt.\; this functional.\\

In Section \ref{sec:2.2} we prove that every stationary point is locally Hölder continuous. The central part will be to establish a Caccioppoli type inequality.\\

In Section \ref{sec:2.3} we use a variant of De Maria’s argument, see \cite{DM09}, that for general domains and boundary conditions all stationary points need to be in $W_{loc}^{2,2}(\Om,\R^2)$ if the integrand is uniformly convex, i.e. $0<\ga<1.$ The main obstacle one needs to overcome to establish a Caccioppoli type inequality is that the second derivative of the integrand $W,$ as defined in \eqref{eq:Intro1.4}, is not bounded but rather grows quadratically. The choice of test function is crucial to control this behaviour. \\

In Section \ref{sec:2.4} we develop higher-order regularity. We will consider stationary points $u\in \A$ and assume additionally $u\in W_{loc}^{2,2}(\Om,\R^2)$ and that $x\mt\det\grad u(x)$ is Hölder-continuous for any $x\in\Om.$ 
We will show that such stationary points obtain slightly higher integrability, i.e. there exists $\d>0$ s.t. $u\in W_{loc}^{2,2+\d}(\Om,\R^2).$ This is done by establishing a Reverse Poincar\'{e} inequality. 
However, there are some technical difficulties to overcome, making it necessary to rely on measure theoretic and compensated compactness type arguments. \\
For $u\in W_{loc}^{2,2+\d}(\Om,\R^2),$ Schauder Theory takes over and guarantees maximal smoothness, i.e. $u\in C_{loc}^\infty(\Om,\R^2).$ \\ 
Combining the results of Section \ref{sec:2.2} and Section \ref{sec:2.3} yields that stationary points are smooth as long as the integrand is uniformly convex, i.e. $0<\ga<1.$\\

\textbf{Plan for chapter 3:} In chapter 3 we still consider the functional \eqref{eq:Intro.1.1} with integrand, as defined in \eqref{eq:Intro1.4}. 
However, we fix the domain to be the unit ball and the boundary conditions to be M-covering maps, which are maps $u_M:S^1\ra S^1$ for any $M\in \N$ which can be represented by $\th\mt e_R(M\th):=(\cos(M\th),\sin(M\th)),$ $\th\in[0,2\pi).$ 
Furthermore, we restrict the class of admissible functions $\A$ to radially symmetric M-covering maps (rsMc) $\A_r^M$. Those are maps $u\in \A$ s.t. $u$ can be represented as $u(x)=r(R)e_R(M\th)$ for any $R\in [0,1],$ any $\th\in [0,2\pi),$ 
where $r:[0,1]\ra\R$ is the radial part of $u$ that is independent of $\th,$ and satisfies $r(1)=1$, which needs to be true so that $u$ agrees with $u_M$ on the boundary.\\

This scenario has been studied in the NLE case by Bauman, Owen and Phillips in a series of innovative papers \cite{BOP91,BOP91MS} and has been further investigated by Yan and Bevan in \cite{Y06,BY07}. In chapter 3 we develop this method in compressible elasticity.\\

Chapter 3 starts again with a precise statement of the considered problem and by showing the existence of a minimizer in $\A_r^M.$ Moreover, for a stationary point $u\in \A_r^M$ we derive the elementary properties that the radial part $r$ needs to satisfy. \\
In Section \ref{sec:3.1} we show that in case of radially symmetric maps $(M=1)$ the identity is the unique (smooth) global minimizer in the full class $\A$ for any $0<\ga<\infty$.\\
In Section \ref{sec:3.2} we will discuss the classical BOP strategy and we conclude that stationary points $u\in \A_r^M$ need to be at least in $C^1$ (for all $0<\ga<\infty$). For $\ga\ge1$ this is a minor improvement from the Hölder-continuity before. In case of $0<\ga<1,$ we already know by chapter 2 that $u$ has to be smooth.\\
This maps are further investigated in Section \ref{sec:3.3} following the newer results by Yan and Bevan. One of our novel findings is that there are two different cases to consider depending on the behaviour of $\rho$. 
The first such case occurs if $\rho$ is lifting-off delayed, i.e. there exists $\tilde{s}>0$ s.t.  $\rho\equiv0$ on $[0,\tilde{s}].$ Then all rsMc.\! stationary points have to be smooth.\\
Else let $\rho$ be an immediate lift-off function, i.e. $\rho(s)>0$ for all $s>0.$ Then, if we assume that $u$ is a rsMc.\! stationary point, $u$ either has to be smooth and its radial map $r$ needs to be a delayed lift-off function
(i.e. there exists some positive $\d>0$ s.t. $r(R)=0$ for all $R\in[0,\d]$) or $u$ is of class $C^1$ and the corresponding $r$ an immediate lift-off function (i.e. $r(0)=0$ and $r>0$ on $(0,1]$). For the latter case it remains open, whether $r$
needs to be any smoother than $C^1,$ except for the case $\ga<1,$ where this is known.  We then at least give a necessary condition which needs to be satisfied if $r$ is of class $C^{1,\al}$ for some $\al\in(0,1],$ see Lemma \ref{lem:3.3.200} for details.\\

\section{Introduction for chapter 4}
\textbf{Plan for Chapter 4:}\\

In this chapter we discuss various uniqueness results in incompressible and compressible situations. The whole chapter can be seen as a contribution to John Ball's agenda \cite[\S 2.6]{BallOP}.\\

In Section \ref{sec:4.1} we add to a discussion which has started recently in \cite{BeDe20}. They consider a uniformly convex functional with the identity as boundary data and incompressibility condition. 
What we add here is that the identity needs to be the (unique) global minimizer as long as a crucial parameter the energy depends on remains small enough. \\

In the second Section \ref{sec:4.2} we introduce a small pressure criteria which when satisfied implies uniqueness of the global minimizer of any quadratic uniformly convex functional on the unit ball and suitable boundary conditions in the incompressible case. The small pressure condition is then discussed under affine boundary data and $N-$covering boundary data. \\

Surprisingly, this leads to a counterexample to regularity in incompressible elasticity, which appears to be new.\\ 
More precisely, on the unit ball $B\ss\R^2$ and for smooth boundary conditions (however, with a topological change), we construct a uniformly convex functional
which depends discontinuously on $x,$ but depends smoothly on $\grad u,$ s.t. the corresponding energy is uniquely globally minimized by a map $u:B\ra\R^2$ that is everywhere Lipschitz but not any better. See, Corollary \ref{cor:4.2.Counter.1} for details.\\

We then also give a partial uniqueness result: if the pressure is too high to guarantee uniqueness one can still guarantee uniqueness up to low-order Fourier-modes. See, Section \ref{sec:4.2.4} for details.\\
This is done by extending ideas  of J. Bevan \cite{JB14}, who recognised this in the special case of the Dirichlet energy and with the double-covering map on the boundary.\\

In the next Section \ref{sec:4.2.5} we give similar conditions in case of the $p-$Dirichlet functional.\\

Finally, in Section \ref{sec:4.3} \!we transfer some of the ideas to compressible elasticity and show a high frequency uniqueness result for the polyconvex functional discussed in the Chapters 2 and 3 see Section \ref{sec:4.3.1}, 
and once again in more general situations, this time for general $p-$growth polyconvex functionals, see Section \ref{sec:4.3.2}\\

\section{Regularity results related to Elasticity}
\subsection{Compressible elasticity}
We give an overview of some important regularity results: a complete treatment would be beyond the scope of this thesis.\\
For an overview of the field of regularity theory in the Calculus of Variations, see the classic texts by \cite{MG83}, \cite{GIU03}, or the famous article \cite{M06} by Mingione. We are not discussing results related to linear growth functionals: 
for an overview see the monograph \cite{B03} by Bildhauer. We don't list results which are concerned with material failure, like cavitation, fracture, and cracks. In particular, there will be no discussion of the spaces $BD,$ $BV,$ etc. \\

\textbf{Uniformly convex functionals:}\\
We start with the uniformly convex case. As a reminder, for a uniformly convex functional of the form \eqref{eq:Intro.1.1} with integrand $W\in C^2$ the ELE is a strongly elliptic nonlinear system (see the discussion in Section \ref{sec:1.1}). We will recall that notion. For this let $\Om \ss\R^m$ be an open and bounded domain with $\p\Om\in C^{0,1}.$ Let $A:\Om \ti\R^{m\ti n}\ra\R^{m\ti n}$ be a matrix-valued, measurable function s.t. there exist a $1\le p,q\le\infty$ s.t. $x\mt A(x,\grad u(x))\in L^q(\Om,\R^{m\ti n})$ for any $u\in W^{1,p}(\Om,\R^n).$ Then we call $u\in W^{1,p}(\Om,\R^n)$ a weak solution to a strongly elliptic nonlinear system if $u$ solves\footnote{Recall that the divergence-operator applies to any  matrix valued map $A\in C^1(\Om\times\R^{m\ti n},\R^{m\ti n})$ 'rowwise', that is 
\begin{align*}\div A(x,\grad u(x))=(\div A_{ij}(x,\grad u(x)))_{j=1,\ldots,n}=\left(\sum\limits_{i=1}^m\p_{i}A_{ij}(x,\grad u(x))\right)_{j=1,\ldots,n}.\end{align*}
}
\begin{align}\div A(x,\grad u(x))=0\; \mb{in} D'(\Om,\R^n),\label{eq:UC.1.1}\end{align}
which must be understood in the sense that
\begin{align*}\int\limits_{\Om}{A(x,\grad u(x))\cd\grad \vp(x)\;dx}=0 \;\mb{for any}\; \vp\in D(\Om,\R^n):=C_c^\infty(\Om,\R^n).\end{align*}
Note, that in this case we can allow $\vp\in  W_0^{1,q'}(\Om,\R^n)$ for any $1\le q'\le\infty,$ since $C_c^\infty$ is dense in $ W_0^{1,q'}$ for any $1\le q'\le\infty.$\\ 

Additionally, in order for the given system $A$ to be called strongly elliptic, there must be some constant $\nu>0$ s.t.
\begin{align}\grad_\xi A(x,\xi)F\cdot F\ge\nu|F|^2 \mb{for all} F\in\R^{m\ti n}\sm\{0\},\mb{for a.e.} x\in\Om, \mb{and for any} \xi\in\R^{m\ti n}.\label{eq:UC.1.1.a}\end{align}
The relation between $A$ and $W$ is then given by $A(x,\grad u(x))=\grad_\xi W(x,\grad u(x))$ for a.e. $x\in\Om,$ ignoring any possible $u$ dependence for now.\footnote{It is crucial to realise that if $W$ depends explicitly on $u$ the relevant ELE is given by \begin{align*}\p_sW(x,u,\grad u)-\div \grad_\xi W(x,u,\grad u)=0.\end{align*}.} We will focus mainly on the autonomous case, that is, if $W$ does not explicitly depend on $x$ and $u$ and we will make sure to mention it, if it becomes important to us. Many of the following results remain valid if $W$ depends on $x$ and $u$ in a suitable manner.\\

\textbf{Scalar case:} In case of $u$ being a scalar function $n=1,$ again, assuming $W\in C^2,$ then by a standard argument the ELE can be interpreted as a `linearized' strongly elliptic equation, 
where we call $u\in W^{1,2}(\Om)$ a weak solution to such an equation if for a matrix-valued, measurable function $A:\Om \ss\R^m\ra\R^{m\ti m},$  $u$ weakly solves
\begin{align}\div (A(x)\grad u(x))=0 \label{eq:UC.1.2}\end{align}
and there exists some constant $\nu>0$ s.t.
\begin{align}A(x)F\cdot F\ge\nu|F|^2 \mb{for all} F\in\R^{m}\mb{and a.e.} x\in\Om.\label{eq:UC.1.2.a}\end{align}
For these types of equations, the De Giorgi-Nash-Moser Theorem states that for $A\in L^\infty$ one gets $C_{loc}^{0,\al}-$regularity for solutions of \eqref{eq:UC.1.2}-\eqref{eq:UC.1.2.a}, see \cite[\S 14]{JJ12}.

If $A\in W^{k,\infty}$ for some $k\in\N$ then one can obtain $C_{loc}^{k,\al}-$regularity by Schauder Theory, see below.\\
Combining the De Giorgi-Nash-Moser Theorem with Schauder Theory yields, that for $n=1$ and $A$ smooth enough every weak solution needs to be maximally smooth. Hence, full regularity can be obtained.\\

Mooney, in \cite{CM20}, constructed a uniformly convex functional with a $C^1-$integrand but where the 2nd derivative of this integrand blows up on a subset of $\R^4,$ s.t. the global minimizer $u:\R^4\ra\R,$ is Lipschitz but no better. This example is limiting, since it shows that some assumptions on the regularity of $A$ or $W$ have to be made in order for the De Giorgi-Nash-Moser result to hold.\\

\textbf{Systems:} For systems ($n\ge2$), as given in \eqref{eq:UC.1.1}-\eqref{eq:UC.1.1.a}, the De Giorgi-Nash-Moser Theorem fails in general. For a list of counterexamples by De Giorgi and others, see the citations above Theorem 2 on p.2 in \cite{B03}. 
Since one can no longer rely on a general theory in the vectorial case, regularity results need to be obtained individually for each situation.\\

A classical result states that if $W\in C^2$ is uniformly convex and $|D^2W|<\infty$ then for all stationary points $u$ it holds $u\in W_{loc}^{2,2}.$ This holds for arbitrary dimensions $m,n\ge2.$ However, in particular for $m=2$ and $n=2$ combining Meyers's Theorem, see \cite{Me63}, and Schauder Theory 
one even gets maximal smoothness. These results were first obtained by Morrey, see \cite{CM38, CM66}. De Maria \cite{DM09} showed that the condition $|D^2W|<\infty$ can be dropped, for a fairly general subclass of functionals $W$. 
However, then one does not necessarily get any more regularity for free in $2\ti 2$ dimensions, since Meyers's no longer applies. Moreover, from the fact that $u\in W_{loc}^{2,2}$ one is further able to extract a partial regularity result and an estimate on the dimension of the singular set $\dim_{\mathcal{H}} \Sigma\le m-2,$ see \cite[\S 3.8]{GIU03}. \\

In \cite{MS16} Mooney and Savin construct a uniformly convex functional $W\in C^\infty(\R^{3\ti2})$ with a singular global minimizer. Note, the dimensions are optimal, since in the scalar case and $2\ti 2$ full smoothness can be guaranteed. \\

Kristensen and Mingione in \cite{KM05,KM05/2} give estimates on the dimension of the singular set for $\omega-$minima (almost minimizers) for suitable uniformly convex functionals.\\

The uniformly convex case is the case which is best understood and with the strongest results. One can then start weakening the assumptions in various ways and see what kind of results remain valid.\\

\textbf{Asymptotically convex functionals:}\\
For asymptotically convex functionals many everywhere regularity results have been established, starting with the work of K. Uhlenbeck in \cite{U77} and refined by many authors: see \cite{M06} for an overview. 
They obtained $C^{1,\al}-$regularity for local minimizers of functionals of the form $\xi\mt W(|\xi|^p),$ for $1<p<\infty$ where $W$ is asymptotically convex. 
Regularity results are available for a larger class of functionals, in which the structural assumption can be modified. However, the price to pay is that only Lipschitz regularity can be obtained, see \cite{DSV09}.\\

\textbf{Rank-1 convex functionals:}\\
We now introduce a weaker notion of ellipticity. A system is said to be elliptic if $u:\Om\ss\R^m\ra\R^n$ is a solution of the following nonlinear elliptic system: Let $A:\Om\ti\R^{m\ti n}\ra \R^{m\ti n}$ be a matrix-valued, measurable function s.t. there exist a $1\le p,q\le\infty$ s.t. $x\mt A(x,\grad u(x))\in L^q(\Om,\R^{m\ti n})$ for any $u\in W^{1,p}(\Om,\R^n).$ Then  $u\in W^{1,p}(\Om,\R^n)$ is called a weak solution to a nonlinear elliptic system if $u$ weakly solves
\begin{align}\div A(x,\grad u(x))=0,\label{eq:R1.1.1}\end{align}
and where there exists some constant $\nu>0$ s.t.
\begin{align}
 \;\grad_\xi A(x,\xi)(a\ot b)\cdot (a\ot b)\ge\nu|a|^2|b|^2 \mb{for all} a\in\R^{m},b\in\R^{n},\xi\in\R^{m\ti n}, \mb{and a.e.} x\in\Om, \label{eq:R1.1.1.a}\end{align}
where the latter is known as the Legendre-Hadamard condition. If the functional is rank-1 convex (which in particular is satisfied if the integrand is poly- or quasiconvex) the ELE is such an elliptic system. Like above, we only focus on results where $A$ only depends on $\xi.$\\

Phillips showed, in \cite{Ph02}, that every Lipschitz continuous one-homogeneous weak solution $u:\R^2\rightarrow \R^N$ of \eqref{eq:R1.1.1} and \eqref{eq:R1.1.1.a} for $A$ smooth enough is necessarily linear. This result has been extended by J. Bevan in \cite{JB10}.\\
In contrast, Bevan in \cite{BE05} constructed a $2\ti2-$dimensional irregular rank-1 functional with a 1-homogenous map as the global minimizer, which is Lipschitzian but not $C^1.$ This shows that the regularity assumption of $A$ or $W$ in Phillips result is indeed necessary. \\

In \cite{CFLM20}, Cupini et al. establish local Hölder regularity for local minimizers of rank-1 and polyconvex functionals under certain structural assumptions on the integrand. \\

Higher-order regularity for linear elliptic systems is what is known as Schauder Theory. \\
For this let $A:\Om\ss\R^{m}\ra \R^{m\ti n}$ be a matrix-valued, measurable function for which we can find a constant $\nu>0$ s.t. it holds
\begin{align*}A(x)(a\ot b)\cdot (a\ot b)\ge\nu|a|^2|b|^2 \mb{for all} a\in\R^{m},b\in\R^{n}\mb{and a.e.} x\in\Om.\end{align*} Moreover, let $f:\Om\ss\R^{m}\ra \R^{n}$ be a measurable function. Then $u\in W^{1,2}(\Om, \R^{n})$ is a weak solution of a linear elliptic system if it weakly solves 
\begin{align*}\div (A(x)\grad u(x))=f(x).\end{align*}
Then if $A\in C^0(\ol{\Om},\R^{m\ti n}),$ $f\in L^{2,\la}(\Om,\R^{n})$ for some $0<\la<n,$ then $\grad u\in L^{2,\la}(\Om,\R^{m\ti n})$ for some $0<\la<n,$ here $L^{2,\la}$ denotes a Morrey space.
If $A, f\in C^{k,\al}$ for some $k\in\N,$ $0<\al<1$ then $u\in C_{loc}^{k+1,\al},$ consult \cite[\S3.3]{MG83}.
This is of interest because it starts a process which is known as bootstrapping. If one starts with $k=0,$ once applied one gets $u\in C_{loc}^{1,\al},$ which in turn improves the regularity of $A$ and $f,$ then one gets even 
more regularity on $u\in C_{loc}^{2,\al},$ etc. This process only ends when the regularity of $A$ and $f$ is reached.\\

\textbf{Quasiconvex functionals:}\\
For general quasiconvex functionals, only partial regularity results and a higher order regularity result are known.\\

The classical works concerning partial regularity for stationary points of quasiconvex functionals by Evans \cite{Ev86}, \cite{AF87}. Since then much more research has been done, see for example \cite{KT03, SS09}. Recently, a partial regularity result in the context of $\A-$quasiconvex functionals has been established in \cite{CG20}. \\

In contrast to Morrey's result, which we previously discussed, the counterexamples by Müller and Sverak in \cite{MS03} and Kristensen and Taheri \cite{KT03} show that there are weak Lipschitz solutions to elliptic systems as given in \eqref{eq:R1.1.1} and $A$ being completely smooth that are nowhere $C^1.$ 
Indeed, they construct quasiconvex smooth integrands s.t. the (local) minimizer is everywhere Lipschitz but nowhere $C^1$. 
These results are obtained using Gromov's convex integration method. These results, in particular, give examples of stationary points where the singular set has full dimension, destroying any hope on a sharper estimate on $\dim_{\mathcal{H}}\Sigma$ below the space dimension for general stationary points. \\

Kristensen and Mingione in \cite{KM07} were able to establish such an estimate on the dimension of the singular set, strictly less than the full dimension, for Lipschitzian $\omega-$minima for suitable quasiconvex functionals. 
This shows that the 'wild' stationary points discussed before can never be global minimizers.\\  

A higher order regularity result can be found in the PhD Thesis \cite{C14}. If the integrand is strongly quasiconvex and $C^2,$ the boundary conditions are smooth enough and are small in some $L^p-$norm, 
then maximal smoothness and uniqueness can be guaranteed. Note that in Section \ref{sec:2.4} we allow more general boundary conditions and no smallness condition is assumed.\\

\textbf{Polyconvex functionals:}\\
In \cite{FH94} Fusco and Hutchinson consider a $2\ti2-$dimensional polyconvex model problem and obtain everywhere continuity for $Q-$minimzers. In \cite{CLM17} a special class of polyconvex functionals in $3\ti3-$dimension is considered and for local minimizers local boundedness is obtained.\\ 

For partial regularity results, see for instance \cite{FH94,EM01,CY07}.\\

In the spirit of the results by Müller et al., Szekelhydi \cite{Sz04} constructed a counterexample with a smooth polyconvex integrand.\\

In \cite{SS11} singular weak solutions to the energy-momentum equations in finite elasticity, are discussed. 

\subsection{Incompressible elasticity}

A partial regularity result is available in the incompressible case, in \cite{EVGA99}. They are able to establish partial regularity for strongly quasiconvex autonomous integrands in $2\ti2-$dimension for Lipschitzian-minimizer 
and with a non-degeneracy condition on the gradient. It seems to be an open question whether one can find an estimate on the dimension of the singular set that is strictly smaller then 2. 
Note that the previous result does not apply to the example we describe in Section \ref{sec:4.2.3}, since our integrand is clearly non-autonomous with a discontinuous dependency on $x$ and the constructed minimizer does not satisfy the non-degeneracy condition.\\

Higher order regularity for the special case of the Dirichtlet functional $\mathbb{D}(\xi)= |\xi|^2/2$ has been discussed in \cite{BOP92}.

\subsection{Nonlinear elasticity}
In \cite{BOP91MS} and later improved by Yan in \cite{Y06} a higher order regularity result was established, which applies to $W^{2,2}\cap C^{1}-$solutions of the equilibrium equations of $W(\xi)=f(\xi)+h(\det \xi),$ where $h$ as in \eqref{eq:Intro1.3}, satisfies some additional properties, and $f$ is of $p-$growth and quasiconvex. On the other hand, a $C^1$ but not any better counterexample is constructed, destroying the hope of developing a full regularity theory in those situations.\\ 

In \cite{JB17} Bevan considered the energy $W(\xi)=|\xi|^2/2+h(\det \xi)$ and picked a specific $h$ with properties as in \eqref{eq:Intro1.3}. He further introduced the notion of positive twist maps and obtained (local) Hölder regularity for global minimizers that possess the positive twist property. Consult this paper also for a nice overview on related regularity results.\\
 
Minimizers, in the set of twist and shear maps to nonlinear elasticity type integrands, have been considered in \cite{BK19}. \\

Partial regularity in NLE remains completely open. However, Fuchs and Reuling in \cite{FR95} obtain partial regularity for a sequence of regularised functionals converging in some sense to the actual functional.\\


\begin{re}The next paragraph discusses mainly uniqueness results. However, some of those results coincide with regularity statements and therefore might be interesting from a regularity point of view, too. \end{re}

\section{Overview of uniqueness results in Elasticity}
In this section we are concerned with the following question: do minimizers/ stationary points in any given elastic situation have to be unique? In \cite[\S 2.6]{BallOP} Ball raised awareness to these types of uniqueness questions in elasticity although Problem 8 of his paper is concerned with a very specific setting.
Simple considerations from material science make apparent that in general we can not expect uniqueness. For instance, imagine an elastic rod, being indented parallel to the direction of the stick. 
It will bend perpendicular to the direction of the stick. 
But it can bend in an arbitrary direction of the plane perpendicular to the stick. This destroys any hope of achieving uniqueness in general, even so this is just a motivating example not a rigorous argument.\\

Here we give an incomplete list of related uniqueness results.

\subsection{Compressible elasticity}

It is well known that uniformly convex functionals possess unique global minimizers, see for example \cite[\S 3.3]{Kl16}. Note this is true even under fairly weak assumptions: it is enough to allow $W(x,\xi)$ to be measurable in $x$ and uniformly convex and $C^1$ wrt. \@ $\xi.$ This can be no longer true in this generality in the incompressible case, 
as our counterexample clearly shows, at least if the integrand is non-autonomous and depends discontinuously on $x.$ Maybe one can still recover a uniqueness result in incompressible elasticity assuming some kind of continuity or smoothness in $x$. However, we don't believe this, we are convinced that one can also construct energies that depend smoothly on $x,$ by adjusting our method.\\

 Knorps and Stuart showed in \cite{KS84}, that for a strongly quasiconvex integrand defined on a star-shaped domain and subject to linear boundary data $u_0=Ax$ any $C^2$ stationary point needs to agree with $Ax$ everywhere. 
 A generalisation can be found in \cite{T03}. These results have been transferred to the incompressible and the nonlinear elasticity case, see below. \\

Note that \cite{C14} also contains a uniqueness result, guaranteeing a unique minimizer for strongly quasiconvex $C^2$ integrands and for smooth and small enough boundary conditions. \\

John in \cite{J72} and Spector and Spector in \cite{SS19} obtain uniqueness of equilibrium solutions for small enough strains and under various boundary conditions. In sharp contrast, Post and Sivaloganathan in \cite{PS97} construct multiple equilibrium solutions in finite elasticity. \\

Counterexamples to uniqueness, for strongly polyconvex functionals, have been established by Spedaro in \cite{S08}. However, these counterexamples rely highly on allowing the determinant to take on negative values, which is neither possible in the incompressible nor in the NLE stetting. 

\subsection{Incompressible elasticity}

In \cite{ST10} Shahrokhi-Dehkordi and Taheri give an analogous result to the one by Knorps and Stuart in the incompressible case.\\

Ball in \cite{B77} first discussed the following problem: consider the Dirichtlet Energy $\mathbb{D}(\xi)= |\xi|^2/2$ on the unit ball and as the boundary condition we have the double covering map given by $u_2=(\cos(2\th),\sin(2\th))$. 
This is widely known as Double Covering Problem (DCP).\\ 
 
Our result can be understood as a contribution to the DCP. We are able to give a non-autonomous uniformly convex integrand $W(x,\xi)$ depending discontinuously on $x,$ with a non-smooth global minimizer subject to smooth boundary conditions. 
The integrand although it depends on $x$ is in a sense close to the Dirichlet functional. So although we are not able to give an answer to the DCP, we were able to construct a functional "close" to the Dirichlet with a global minimizer that
is Lipschitzian but no better. In \cite{JB14} Bevan showed that $u_2$ is the unique global minimizer up to the first Fourier-mode. In \cite{BeDe21} Bevan and Deane obtained that $u_2$ is the unique global minimizer for general purely inner 
as well as general purely outer variations. Additionally, local minimality is shown for a class of variations allowing a mixture of certain inner and outer ones.\\

Another consequence of our uniqueness result is that it divides all incompressible situations in small and high pressure ones.
It turns out that the DCP is in the high pressure regime. In chapter 4 it is explained why these types of problems are much more difficult to treat than low pressure situations. 
Even so we could shed some light on the DCP and give a partial answer, determining the DCP is by no means unimportant; on the contrary it is even more intriguing, since it would be a first uniqueness or non-uniqueness result in a high pressure situation.\\

Much more research has been done that is closely related to the DCP. For example, Morris and Taheri in \cite{T09} and \cite{MT19} consider more general functionals of the form $W(x,s,\xi)=F(|x|^2,|s|^2)|\xi|^2/2$, $F\in C^2$ on the annulus $A$ 
and the set of admissible maps $\A=W_{\Id}^{1,2}(A,\R^2)$ and they show that there are countably many unique solutions, exactly one for each homotopy class.\\

On the negative side of things, in the paper \cite{BeDe20}, which we already mentioned, equal energy stationary points of an inhomogeneous uniformly convex functional $(x,\xi)\mapsto f(x,\xi)$ depending discontinuously on $x$ are constructed. 
It remains unknown for now if these stationary points are actually global minimizers.

\subsection{Nonlinear elasticity}
 In \cite{B11} Bevan extended the aforementioned results by Knorps, Stuart and Taheri to the NLE case.\\

Bevan and Yan show in \cite{BY07} that the BOP-map is the unique global minimizer in a special sub-class of admissible maps.\\

Uniqueness has been discussed very recently by Sivaloganathan and Spector in \cite{SS18}, which is closely connected to our work. They discuss a uniqueness criteria for polyconvex integrands of the form as given in \eqref{eq:Intro1.3}, which when satisfied, implies uniqueness. Moreover, various examples are discussed.\\

Uniqueness and regularity of Twist and Shear Maps have been addressed in \cite{BK19} by Bevan and Käbisch. In the latter paper a nice overview is given on further literature regarding twist and shear maps in NLE. \\


\textbf{Plan for the conclusion of the thesis:}\\ 
The main parts of this thesis can be found in the Chapter 2 to 4. It is followed by Chapter 5, which includes a small overview of our results, some additional thoughts and some intriguing open questions. 
We included some additional results in the Appendix that we believe could be helpful. The thesis is concluded by a list of references, and finally a CV is included.

\clearpage{\pagestyle{empty}\cleardoublepage} 
\chapter{Regularity results for stationary points of a polyconvex functional}
\section{Introduction and notation}
\label{sec:2.1}

For $\Om\ss \R^2$ open and bounded with $\p\Om\in C^{0,1},$ define the functional $I:W^{1,2}(\Om,\R^2)\rightarrow \R$ by
\begin{equation}
I(u):=\int\limits_{\Om}{\frac{1}{2}|\grad u|^2+\rho(\det \grad u)\; dx}
\label{eq:1.1}
\end{equation}
for all $u\in W^{1,2}(\Om,\R^2).$ The function $\rho\in C^{\infty}(\R)$ is defined by
\begin{equation}
\rho(s)=\left\{\begin{array}{ccc}
0& {\mbox{if}}& s\le0,\\
\rho_1(s)& {\mbox{if}}& 0\le s\le s_0,\\
\gamma s+\ka&{\mbox{if}}& s_0\le s,
\end{array}
\right.
\label{eq:1.2}
\end{equation}
for some constants $\gamma>0,$ $s_0\ge0$ and $\ka\ge-\ga s_0.$ 
Here $\rho_1:[0,s_0]\rightarrow \R$ is a smooth and convex function on $[0,s_0]$ satisfying the boundary conditions $\rho_1(0)=0$ and $\rho_1(s_0)=\gamma s_0+\ka$ and the connections need to be in such a way that $\rho$ is smooth everywhere.\footnote{It would be enough to assume $\rho\in C^k$ for some $k\ge2.$ The results below remain also true for more general integrands of the form
\begin{equation*}
\rho(s)=\left\{\begin{array}{ccc}
\gamma |s|+\ka&{\mbox{if}}& s\le -s_0,\\
\rho_1(s)& {\mbox{if}}& -s_0\le s\le s_0,\\
\gamma s+\ka&{\mbox{if}}& s_0\le s,
\end{array}
\right.
\end{equation*} 
where $\rho_1:[-s_0,s_0]\rightarrow \R_0^+$ is again a smooth and convex function on $[-s_0,s_0]$ satisfying the conditions $\rho_1(0)=0,$ $\rho_1(s_0)=\gamma s_0+\ka$ and $\rho_1(-s_0)=\gamma |s_0|+\ka$ with smooth connections.}  Note that $\rho$ is convex on the whole real line. Hence, the complete integrand is polyconvex. Recall, again, we call $f:\R^{2\ti2}\ra\ol{\R}$ polyconvex, if there exists a convex function $g:\R^5\ra\ol{\R}$ s.t. $f(\xi)=g(\xi,\det\xi)$ for all $\xi\in\R^{2\ti2}.$\\
The behaviour of the functional depends mainly on the parameter $\ga.$ If $\ga\ra0$ then the functional turns into the well known Dirichlet energy. In the regime $0<\ga<1$ the functional is uniformly convex, a proof of which can be found in Section \ref{Ap:A.3} If $\ga\ge1$ then the functional is genuinely polyconvex.\\
Recall that for $u\in W^{1,2}(\Om,\R^2)$ the Jacobian is in $L^1,$ i.e. $\det\grad u\in L^1(\Om)$ so $\rho(\det\grad u(\cd))\in L^1(\Om)$ and hence $I(u)$ is well defined for all $u\in W^{1,2}(\Om,\R^2).$\\

Furthermore, we introduce the set of admissible functions
\[\A_{u_0}:=\{u\in W^{1,2}(\Om,\R^2): u=u_0 \;\mb{on}\; \p\Om\}.\]
The boundary condition needs to be understood in the trace sense; for a small discussion, see, Section \ref{Ap:A.1}\\

\textbf{Notation:} Let $A,B\in \R^{m_1\times\cdots\times m_n},$ where $m_k\in\N\sm\{0\}$ for all $0\le k\le n.$ We will denote the Frobenius inner product by $A\cdot B:=\sum\limits_{i_1,\ldots, i_n}A_{i_1,\ldots, i_n}B_{i_1,\ldots, i_n}$ and the corresponding norm by $|A|:=(A\cdot A)^{\frac{1}{2}}$ for all $A,B\in \R^{m_1\times\cdots\times m_n}.$\\
For a matrix $A\in \R^{n\ti n}$ we denote the determinant by $d_A:=\det A$ (when there is no confusion, we will suppress the dependence on $A$.) and the cofactor by $\cof A:=(\adj \;A)^T.$ It is well known, that the cofactor is the derivative of the determinant, i.e. $\grad_A(\det A)=\cof A.$ Moreover, for a $2\times2-$matrix $A$ the cofactor takes the simple form 
\begin{equation}
\cof A=\begin{pmatrix}a_{22} &-a_{21}\\-a_{12} & a_{11}\end{pmatrix}
\label{eq:1.3}
\end{equation}
which is a rearrangement of the original matrix (up to signs). Note, that $|\cof A|=|A|.$\\
For two vectors $a\in \R^n,b\in \R^m$ we define the tensor product $a\ot b\in \R^{n\ti m}$ by $(a\ot b)_{i,j}:=(ab^T)_{i,j}=a_ib_j$ for all $1\le i\le n,$ $1\le j\le m.$\\
$B(x,r):=\{y\in\bb{R}^2:|y-x|<r\}$ be the open ball with center $x$ in $\R^2$ and denote its boundary by $S(x,r).$ If there is no confusion we use the abbreviations $B_r,S_r$ and in particular, we denote the unit ball with center $0$ by $B$ and its boundary by $S^1.$\\
We want $\mathcal{H}^s$ to be the $s-$dimensional Hausdorff measure for any $s\in \R_0^+,$ $\La^n$ to be the $n-$dimensional Lebesgue measure for any $n\in \N:$ As usual, we use $dx=d\La^n$ for short.\\

As a reminder, we give a short overview over the considered function spaces: For all $k\in \N\cup\{\infty\},$ $C^k$ is the class of $k-$times differentiable functions and the $k-$th derivative is continuous in particular, we allow $k=\infty$ and call all functions in $C^\infty$ `smooth'. We use $C_c^k$ for functions of class $C^k$ and with compact support. Functions of class $C^{k,\al}$ with $0<\al<1$ are again in $C^k$ but additionally, the $k-$th derivative needs to be Hölder continuous with exponent $\al.$ Moreover, we use $C^{0,1}$ for Lipschitz continuous functions.\\
The Lebesgue spaces are defined by 
\[L^p(U,\R^m,\mu):=\left\{f:U\ra \R^m:\|f\|_{L^p(U,\R^m,\mu)}<\infty\right\}\]
 for any $1\le p\le\infty,$ any measure $\mu$ and any $\mu-$measurable set $U\ss \R^n.$ Recall
 \begin{align*}
  \|f\|_{L^p(U,\R^m,\mu)}&:=\left(\int\limits_U{|f(x)|^p\;d\mu}\right)^{\frac{1}{p}} \;\mb{for any}\; 1\le p<\infty \;\mb{and}&\\
 \|f\|_{L^\infty(U,\R^m,\mu)}&:=\lim_{p\ra\infty}\left(\int\limits_U{|f(x)|^p\;d\mu}\right)^{\frac{1}{p}}.&
\end{align*}
 The Sobolev spaces are given by
\[W^{k,p}(U,\R^m,\mu):=\left\{f:U\ra \R^m:\|f\|_{W^{k,p}(U,\R^m,\mu)}=\left(\sum\limits_{l=0}^k\|\grad^lf\|_{L^p(U,\R^{m(n^l)},\mu)}^p\right)^{\frac{1}{p}}<\infty\right\},\] for all $1\le p\le\infty,$ $k\in\N$  and where the $p=\infty$ case is again thought of as taking the limit. As usual we will suppress the measure if $\mu=\La^n,$ the target space is suppressed if $m=1$ and sometimes all of it might be suppressed if there can't be any confusion. Note, $W^{0,p}=L^p$. Local versions of these spaces are indicated by adding the subscript `loc' to the spaces before ($L_{loc}^p,$ $C_{loc}^{k,\al},$ etc.).\footnote{Roughly speaking, in these local versions the definition of the space only holds for all compact subsets in the considered set. For instance, consider an open and bounded set $U\ss \R^n$ then $f\in L_{loc}^p(U,\R^m)$ (similar for $W_{loc}^{k,p}$) means that for all compact sets $K\ss\ss U,$ $\|f\|_{L^p(\Om,\R^m)}<\infty$ but in general $\|f\|_{L^p(U,\R^m)}$ might be infinity. Another example we will use is $f\in C_{loc}^{k}(U,\R^m)$ then $f$ is $k-$times continuous differentiable for every point $x\in U$ but not necessarily up to the boundary.}\\
 
Furthermore, we introduce the spaces of $W^{1,p}-$functions with zero boundary conditions defined by $W_0^{1,p}(U,\R^m,\mu):=\ol{(C_c^\infty(U,\R^m,\mu))}^{W^{1,p}(U,\R^m,\mu)},$ for any $1<p<\infty.$ By Theorem \ref{thm:Trace theorem2} we have in particular $W_0^{1,2}(\Om,\R^2)=\A_0.$\\
Finally, we use `$\ra$' for strong and `$\rhu$'  for weak convergence.\\

In the next theorem we show that there must be a global minimizer w.r.t. \@ the boundary conditions.

\begin{thm}[Existence] The functional $I$ attains its minimum in $\A_{u_0}$.
\label{thm:2.1.1Ex}
\end{thm}
\textbf{Proof:}\\ We apply the direct method of the Calculus of Variations.\footnote{For a general description, see e.g. Evans's PDE book \cite[Section 8.2.]{LE10}.} First note that $I(u)\ge0$ for all $u\in\A_{u_0}$ by the properties of $\rho.$ Therefore, there exists a minimizing sequence $(v_k)_{k\in\N}\ss \A_{u_0},$ s.t. $I(v_k)\rightarrow \inf\limits_{\A_{u_0}} I$ for $k\rightarrow \infty.$ Then there exists a convergent subsequence (without relabbeling) $v_k\rightharpoonup v\in \A_{u_0}$ for $k\rightarrow \infty,$ since $\A_{u_0}$ is closed w.r.t. weak convergence. Moreover, $I$ is weakly lower semicontinuous (wlsc.) since the norm of a Hilbert space is wlsc., i.e. $\liminf\limits_{k\rightarrow\infty}\|v_k\|_{W^{1,2}(\Om,\R^2)}\ge\|v\|_{W^{1,2}(\Om,\R^2)}$, in particular,
\begin{equation*}
\liminf\limits_{k\rightarrow\infty}\int\limits_{\Om}{\frac{1}{2}|\grad v_k|^2\; dx}\ge\int\limits_{\Om}{\frac{1}{2}|\grad v|^2\; dx}.
\label{eq:1.5}
\end{equation*}
Furthermore, $v_k\rightharpoonup v\in W^{1,2}(\Om,\R^2)$ implies $\det\grad v_k\rightharpoonup \det \grad v$ in $D'(\Om),$ where the latter convergence needs to be understood in the sense that
\begin{equation*}
\int\limits_{\Om}{(\det\grad v_k) \phi\; dx}\ra\int\limits_{\Om}{(\det\grad v) \phi\; dx}\mb{for any} \phi\in D'(\Om):=C_c^\infty(\Om).
\end{equation*}
Additionally, for any $k\in\N$ and the limit it holds that $\det\grad v_k,\det\grad v\in L^1(\Om).$ These are the requirements for Theorem \ref{thm:A.2.1} yielding the desired weak lower semicontinuity 
\[
\liminf\limits_{k\rightarrow\infty}\int\limits_{\Om}{\rho(\det \grad v_k)\; dx}\ge\int\limits_{\Om}{\rho(\det \grad v)\; dx}.
\label{eq:1.6}
\]
Hence, 
\begin{equation*}
\liminf\limits_{k\rightarrow\infty}I(v_k)\ge I(v).
\label{eq:1.7}
\end{equation*}
Together this implies 
\begin{equation*}
\inf\limits_{u\in \A_{u_0}}I(u)\le I(v)\le\liminf\limits_{k\rightarrow\infty}I(v_k)\le\inf\limits_{u\in \A_{u_0}}I(u).
\label{eq:1.8}
\end{equation*}
\vspace{0.5cm}
\begin{re}
This is slightly more involved than one would expect from similar literature. Polyconvex integrands have been studied intensively, the main difference in those papers is that the function $\rho$ takes infinity for negative values and tends to infinity, close to the origin. Therefore, if one is interested in minima one can restrict the search to functions possessing finite energy, which is only possible for maps, whose Jacobian satisfies $\det\grad u>0$ a.e., which would yield $L^1-$convergence instead of the weaker $D'$ convergence of the determinants(see for instance \cite{JB15} and \cite{BCO}). But this is not the case for this type of $\rho.$ A priori we don't have any information on $\det\grad u.$
If instead of convergence in $D'$ we would have weak convergence of the determinants in $L^1$ this would directly imply weak lower semi-continuity via the subdifferential estimate
\begin{equation*}
\int\limits_\Om{\rho(d_{\grad u_k})\;dx}\ge\int\limits_\Om{\rho(d_{\grad u})d_{\grad u}\;dx}+\int\limits_\Om{\rho'(d_{\grad u})(d_{\grad u_k}-d_{\grad u})\;dx}.
\label{eq:1.9}
\end{equation*}
Since $\rho'(d_{\grad u}(\cdot))\in L^\infty(\Om)$ the rightmost integral vanishes if $k\rightarrow\infty.$ But this is not the case.
\end{re}

\section{Hölder regularity for stationary points, for any $0<\ga<\infty.$}
\label{sec:2.2}

The Euler-Lagrange equations (ELE) in the weak form can be obtained by calculating the first variation: For this sake, take $\vp\in C_c^{\infty}(\Om,\R^2)$ and we get
\begin{align}
\partial_\ve|_{\ve=0} I(u+\ve\vp)=&\partial_\ve|_{\ve=0}\int\limits_{\Om}{\frac{1}{2}|\grad u+\ve\grad\vp|^2+\rho(\det (\grad u+\ve\vp))\; dx}&\nonumber\\
=&\int\limits_{\Om}{\partial_\ve|_{\ve=0}\left[\frac{1}{2}|\grad u|^2+\frac{\ve^2}{2}|\grad\vp|^2+\ve\grad u\cdot\grad \vp\right]\;dx}&\nonumber\\
&+\int\limits_{\Om}{\partial_\ve|_{\ve=0}\rho(\det (\grad u+\ve\vp))\; dx}& \nonumber\\
=&\int\limits_{\Om}{\grad u\cdot\grad \vp+\rho'(\det \grad u)\cof\grad u\cdot\grad\vp \; dx},&
\nonumber
\label{eq:1.10}
\end{align}
where we used Lebesgue's dominated convergence theorem.\footnote{see \cite{WZ15}, Theorem 5.36.}
Then the Euler-Lagrange Equations in the weak form are given by
\begin{equation}
\int\limits_{\Om}{(\grad u+\rho'(d)\cof \grad u) \cdot\grad\vp\; dx}=0\; {\mbox{for all}}\; \vp\in C_c^\infty(\Om,\R^2).
\label{eq:ELE1.1}
\end{equation}
By density the equations hold true even for all $\vp\in W_0^{1,2}(\Om,\R^2).$ We are now in the position to state the main theorem of this section.

\begin{thm}[Hölder continuity]
Suppose that $u\in\A_{u_0}$ and $u$ satisfies \eqref{eq:ELE1.1}. Then $u$ is locally Hölder continuous, i.e. there exists $\al\in(0,1)$ s.t. $u\in C_{loc}^{0,\al}(\Om,\R^2).$
\label{thm:2.2.1}
\end{thm}
\textbf{Proof:}\\
The fact that test functions are allowed to be of class $W_0^{1,2}$ enables one to choose $\vp=\eta^2 (u-a)$ for a suitable cut-off function $\eta\in C_c^\infty(\Om)$ and some $a\in\R,$ which we will determine later. Then the gradient of $\vp$ is given by $\grad \vp=2\eta \grad\eta\ot(u-a)+\eta^2\grad u.$ Plugging this into \eqref{eq:ELE1.1}, noting that $\cof A\cdot A=2d_A$ for all $A\in \R^{2\ti2}$ and collecting all $\eta^2-$ terms to one side and all terms containing $2\eta \grad\eta\ot(u-a)$ to the other side yields,
\begin{equation}
\int\limits_{\Om}{\eta^2|\grad u|^2+2\eta^2\rho'(d)d\; dx}=-\int\limits_{\Om}{2\eta(\grad u+\rho'(d)\cof \grad u) \cdot(\grad\eta\ot(u-a))\; dx}.
\label{eq:ELE1.2}
\end{equation}
Since $\rho'$ is smooth and monotonically increasing it holds $0\le\rho'\le\gamma$ and recalling $|\grad u|=|\cof \grad u|,$ then the RHS of \eqref{eq:ELE1.2} can be controlled by
\begin{equation}
C(\ga)\int\limits_{\Om}{|\eta||\grad u||u-a||\grad\eta|\; dx},
\label{eq:1.11}
\end{equation}
where $C(\ga)=2(1+\ga)>0$ is a positive constant.\\

The left-hand side of \eqref{eq:ELE1.2} can be estimated from below by
\begin{equation}
\int\limits_{\Om}{\eta^2|\grad u|^2\; dx}.
\label{eq:1.12}
\end{equation}
However, this requires that the second term is non-negative for nearly every point in $\Om,$ i.e. \[\rho'(\det \grad u(x))\det \grad u(x)\ge0\] for almost every $x\in \Om.$ Lemma \ref{lem:L2.2.1} shows that this is indeed true.

Collecting all of the above we arrive at an inequality of the form
\begin{equation}
\int\limits_{\Om}{\eta^2|\grad u|^2\; dx}\le C(\ga) \int\limits_{\Om}{|\eta||\grad u||u-a||\grad\eta|\; dx}.
\label{eq:1.13}
\end{equation}
From here we follow the standard Caccioppoli method. 
Applying the Cauchy-Schwarz inequality leads to 
\begin{equation}
\int\limits_{\Om}{\eta^2|\grad u|^2\; dx}\le C(\ga) \left(\int\limits_{\Om}{\eta^2|\grad u|^2\;dx}\right)^{1/2} \left(\int\limits_{\Om}{|u-a|^2|\grad\eta|^2\; dx}\right)^{1/2}
\label{eq:1.14}
\end{equation} 
and dividing by the square root of the LHS and then squaring the equation again yields
\begin{equation}
\int\limits_{\Om}{\eta^2|\grad u|^2\; dx}\le C(\ga) \int\limits_{\Om}{|u-a|^2|\grad\eta|^2\; dx}.
\label{eq:1.15}
\end{equation}
Note that from now on, the constant $C(\ga)>0$ may change from line to line as usual.\\

Let $x_0\in \Om$ and let $r>0$ s.t. $B(x_0,2r)\ss\ss \Om$ and choose $\eta\in C_c^\infty(\Om)$ to be $1$ on $B(x_0,r)$ and $0$ on $\Om\setminus B(x_0,2r).$ Furthermore, assume that there exists some constant $c>0$ s.t. $|\grad \eta|\le \frac{c}{r}$ and $\supp\grad\eta \ss B(x_0,2r)\sm B(x_0,r),$ and define $a:=(u)_{B_{2r}\sm B_r}:=\fint\limits_{B_{2r}\sm B_r}{u\;dx}:=\frac{1}{\La^2(B_{2r}\sm B_r)}\int\limits_{B_{2r}\sm B_r}{u\;dx}.$ Then \eqref{eq:1.15} becomes
\begin{equation}
\int\limits_{B(x_0,r)}{|\grad u|^2\; dx}\le  \frac{C(\ga)}{r^2} \int\limits_{B(x_0,2r)\sm B(x_0,r)}{|u-a|^2\; dx}.
\label{eq:1.16}
\end{equation}
Recalling Poincaré's inequality for the annulus, stating, that for every $1\le p<\infty$ there exists a constant $c(p)>0,$ which only depends on $p,$ s.t. for all $v\in W^{1,p}(B_{2r}\sm B_r)$ it holds
\begin{equation}
\int\limits_{B(x_0,2r)\sm B(x_0,r)}{|v-(v)_{B_{2r}\sm B_r}|^p\; dx}\le c(p) r^p \int\limits_{B(x_0,2r)\sm B(x_0,r)}{|\grad v|^p\; dx}.
\label{eq:1.17}
\end{equation}
Applying, \eqref{eq:1.17} with $p=2,$ to \eqref{eq:1.16} yields 
\begin{equation}
\int\limits_{B(x_0,r)}{|\grad u|^2\; dx}\le C(\ga) \int\limits_{B(x_0,2r)\sm B(x_0,r)}{|\grad u|^2\; dx}.
\label{eq:1.18}
\end{equation}
Applying Widman's hole filling technique, see \cite{W71}, which means, adding $\displaystyle{C(\ga)\int\limits_{B(x_0,r)}{|\grad u|^2\; dx}}$ to both sides yields
\begin{equation*}
\int\limits_{B(x_0,r)}{|\grad u|^2\; dx}\le C_W(\ga) \int\limits_{B(x_0,2r)}{|\grad u|^2\; dx},
\label{eq:1.19}
\end{equation*}
where $C_W(\ga)=\frac{C(\ga)}{C(\ga)+1}<1,$ for all $0<\ga<\infty.$ 
Introducing the notation $\displaystyle{\phi(r):=\int\limits_{B(x_0,r)}{|\grad u|^2\; dx}}$ and rewriting the previous equation yields
\begin{equation*}
\phi(r)\le C_W(\ga) \phi(2r),\;\mbox{for all}\; 0\le 2r<\dist(x_0,\p\Om).
\label{eq:1.20}
\end{equation*}
By Lemma \ref{lem:iter}, which is stated below, there exist a ball $B(x_0,r')$  with $r'=r'(\dist(x_0,\p\Om))>0,$ an exponent $0<\alpha=\alpha(C_W(\ga))\le1,$ and a constant $L=L(\dist(x_0,\p\Om))>0$ s.t. 
\begin{equation*}
\phi(r)\le L^2 r^\alpha,\;\mbox{for all}\; 0<r<r'.
\label{eq:1.21}
\end{equation*}
Finally, the Dirichlet growth theorem, see Theorem \ref{thm:DGT} below, implies local Hölder-continuity. Indeed, noting first that for $n=p=2$ the singular set $\Sigma$ of Theorem \ref{thm:DGT} is empty. Hence, there exist $\mu=\min\{\frac{\alpha}{2},1\}\in(0,1]$ and a constant $c=c(n,p,\mu)>0$ s.t. for any $x,y\in B_n(x_0,r')$ s.t. $|x-y|\le \frac{\delta(x)}{2}$ with $\delta(x)=r'-|x-x_0|$ it holds
\begin{equation*}
|u(x)-u(y)|\le cL \delta(x)^{1-\frac{n}{p}-\mu}|x-y|^\mu.
\end{equation*}
This completes the proof.\vspace{1cm}

The section is completed by the results used in the above argument. We start by showing that $\rho'$ satisfies a monotonicity inequality.

\begin{lem} Let $\rho$ be defined as before. Then the following inequality holds
\begin{equation}
d\rho'(d)\ge0\; \forall\; d\in\R.
\label{eq:MC1}
\end{equation}
\label{lem:L2.2.1}
\end{lem}
\textbf{Proof:}\\
Let $d\le0.$ Then $\rho'(d)=0$ and \eqref{eq:MC1} holds. Moreover, $\rho_1'(0)=0$ and $\rho_1'(s)=\gamma$ for all $s\le s_0$ and $\rho$ convex on $\R$ implies that $\rho'$ is monotonically increasing on $(0,\infty).$ Hence, for all $d>0$ $\rho'(d)\ge0.$ Finally, $\rho'(d)d\ge0$ for all $d\in\R.$ \\

Next, we state the version of the Iteration Lemma we used in this context. 

\begin{lem}[Iteration Lemma]Let $\be\in(0,1),$ $\phi:[0,2r_1)\rightarrow[0,\infty)$ be a non-decreasing function satisfying
\begin{equation}
\phi(r)\le \be \phi(2r),\;\mbox{for all}\; 0< r<r_1.
\label{eq:1.22}
\end{equation}
Then there exist $\alpha=\alpha(\beta)>0,$ $r_2=r_2(r_1)>0$ and $c=c(r_1)>0$ s.t.
\begin{equation}
\phi(r)\le c r^\alpha ,\;\mbox{for all}\; 0<r<r_2.
\label{eq:1.23}
\end{equation}
\label{lem:iter}
\end{lem}
\textbf{Proof:}\\ This can be deduced from \cite{MG83}, Lemma 2.1, Chapter 3.\\

Lastly, the Dirichlet Growth Theorem is stated as given in Morrey's monograph \cite[Theorem 3.5.2]{CM66}. A proof is included for the convenience of the reader.

\begin{thm}[Dirichlet Growth Theorem] Let $B_n(x_0,R)\ss \Om\ss\R^n$ and $u\in W^{1,p}(B_n(x_0,R))$ with $1\le p\le n.$ Suppose, that the inequality 
\begin{equation}
\int\limits_{B_n(x,r)}{|\grad u|^p\;dx}\le L^p\left(\frac{r}{\delta(x)}\right)^{n-p+p\mu}
\label{eq:1.DGT1}
\end{equation}
holds for all $x\in B_n(x_0,R),$ all $0\le r\le\delta(x):=R-|x-x_0|$ and some constants $0<\mu\le1$ and $L>0.$\\

Then there exists a subset $\Sigma:=\Sigma(u,n,p)\ss B_n(x_0,R)$ s.t. $\dim_{\mathcal{H}}\Sigma\le n-p$ and $u\in C^{0,\mu}(B_n(x_0,r)\sm\Sigma)$ for all $0\le r\le \frac{\delta(x)}{2}.$ Moreover, there exists a constant $c=c(n,p,\mu)>0$ s.t. for any $x,y\in B_n(x_0,R)\sm\Sigma$ s.t. $|x-y|\le \frac{\delta(x)}{2}$ it holds
\begin{equation}
|u(x)-u(y)|\le c L\delta(x)^{1-\frac{n}{p}}\left(\frac{|x-y|}{\delta(x)}\right)^\mu.
\label{eq:1.DGT}
\end{equation}

\label{thm:DGT}
\end{thm}
\textbf{Proof:}\\
\textbf{Step 1: Approximation.}   
Initially, note that for any $u\in W^{1,p}(B_n(x_0,R))$ with $1\le p\le n$ one can assign vales $u(x)=\lim\limits_{\sigma\ra0}\fint\limits_{B(x,\sigma)}{u(z)\;dz}$ for any $x\in B_n(x_0,R)$ up to a subset $\Sigma:=\Sigma(u,n,p)\ss B_n(x_0,R)$ s.t. $\dim_{\mathcal{H}}\Sigma\le n-p,$ see \cite[\S3.8]{GIU03}. 
Now by density of $C^1(B_n(x_0,R))$ in $W^{1,p}(B_n(x_0,R))$ it is enough to show that \eqref{eq:1.DGT} is true for any $u\in C^1(B_n(x_0,R))$ satisfying \eqref{eq:1.DGT1} and $x,y\in B_n(x_0,R)$ s.t. $|x-y|\le \frac{\delta(x)}{2}$. Indeed, assume the later is true then we know that for any $u\in W^{1,p}(B_n(x_0,R))$ there exists a sequence $\{u_\ve\}\ss C^1(B_n(x_0,R))$ s.t. $\|u_\ve- u\|_{W^{1,p}(B_n(x_0,R))}\ra0$ as $\ve\ra 0.$  
Then note that for every $\ve>0$ we know that $u_\ve(x)=\lim\limits_{\sigma\ra0}\fint\limits_{B(x,\sigma)}{u_\ve(z)\;dz}$ for any $x\in B_n(x_0,R)$ up to a subset $\Sigma'\ss B_n(x_0,R)$ s.t. $\dim_{\mathcal{H}}\Sigma'\le n-p.$ Combining the two sets $\Sigma''=\Sigma\cap\Sigma'$ still satisfies $\dim_{\mathcal{H}}\Sigma''\le n-p,$ suppressing $\Sigma''$ for $\Sigma.$ Hence, we have  for any $x\in B_n(x_0,R)\sm\Sigma$ that
 \begin{align*}\lim\limits_{\ve\ra0}|u(x)-u_\ve(x)|\le\lim\limits_{\ve\ra0}\lim\limits_{\sigma\ra0}\fint\limits_{B(x,\sigma)}{|u(z)-u_\ve(z)|\;dz}=0.\end{align*}
The latter can be seen by exchanging limits and $\|u_\ve- u\|_{W^{1,p}(B_n(x_0,R))}\ra0$ as $\ve\ra 0.$ 
Then for any $x,y\in B_n(x_0,R)\sm\Sigma$ s.t. $|x-y|\le \frac{\delta(x)}{2}$ it holds
 \begin{align*}|u(x)-u(y)|=\lim\limits_{\ve\ra0}|u_\ve(x)-u_\ve(y)|\le c L\delta(x)^{1-\frac{n}{p}}\left(\frac{|x-y|}{\delta(x)}\right)^\mu,\end{align*}
completing the argument.\\

\textbf{Step 2: Showing, that \eqref{eq:1.DGT} is true for $C^1-$maps.} So assume from now on $u\in C^1(B_n(x_0,R))$ satisfying \eqref{eq:1.DGT1} and $x,y\in B_n(x_0,R)$ s.t. $|x-y|\le \frac{\delta(x)}{2}$. By the triangle inequality for an arbitrary $w\in B_n(x_0,R)$ it holds
\begin{align*}|u(x)-u(y)|\le|u(x)-u(w)|+|u(w)-u(y)|.\end{align*}
Denote the midpoint of $x$ and $y$ by $z,$ i.e. $z=\frac{x+y}{2}$ and the distance between $x$ and $z$ by $l:=\frac{|x-y|}{2}.$ Now we average the previous inequality over $B(z,l)$, i.e. \begin{align*}|u(x)-u(y)|\le\fint\limits_{B_n(z,l)}{|u(x)-u(w)|\;dw}+\fint\limits_{B_n(z,l)}{|u(w)-u(y)|\;dw}.\end{align*}
It is enough to consider one of these integrals, the other one can be treated similarly. By applying the mean value theorem we get
\begin{align*}
\fint\limits_{B_n(z,l)}{|u(w)-u(y)|\;dw}\le\frac{2l}{\La^n (B_n(z,l))}\int\limits_{B_n(z,l)}\int\limits_0^1{{|\grad u(y+t(w-y))|\;dt}\;dw}.
\end{align*}
Swapping the integrals and substituting $v=y+t(w-y)$ with $\bar{z}=y+t(z-y)$ we obtain
\begin{align*}
\frac{2l}{\La^n (B_n(z,l))}\int\limits_0^1\int\limits_{B_n(\bar{z},lt)}{{|\grad u(v)|\;dv}\;t^{-n}dt}.
\end{align*}
By Hölder's inequality, assumption \eqref{eq:1.DGT1} and the scaling property of the $n-$dimensional ball we get
\begin{align*}
\int\limits_{B_n(\bar{z},lt)}{|\grad u(v)|\;dv}&\le\La^n(B_n(\bar{z},lt))^{1-{\frac{1}{p}}}\left(\int\limits_{B_n(\bar{z},lt)}{|\grad u(v)|^p\;dv}\right)^{1/p}&\\
&\le \La^n(B_n(\bar{z},lt))^{1-{\frac{1}{p}}} L\left(tl/\d\right)^{\frac{n}{p}-1+\mu}&\\
&= \omega_n^{1-{\frac{1}{p}}} L(tl)^{n-1+\mu}\d^{1-\frac{n}{p}-\mu},&
\end{align*}
where we denoted the volume of the $n-$dimensional unit ball by $\omega_n.$ Together,
\begin{align*}
\fint\limits_{B_n(z,l)}{|u(w)-u(y)|\;dw}\le&\frac{\omega_n^{1-{\frac{1}{p}}}}{\La^n (B_n(z,l))}Ll^{n+\mu}\d^{1-\frac{n}{p}-\mu}\int\limits_0^1{t^{\mu-1}\;dt}&\\
\le&cL\d^{1-\frac{n}{p}-\mu}|x-y|^{\mu}&
\end{align*}
with $c=c(n,p,\mu):=\frac{\omega_n^{-{\frac{1}{p}}}}{n-1+\mu}>0.$ The integral can always be performed, since $\mu>0.$ A similar estimate can be established for $\fint\limits_{B_n(z,l)}{|u(x)-u(w)|\;dw}$ completing the proof.\\

\begin{re}\ \\ (i) In the case $p>n,$ the Sobolev embedding $W^{1,p}\hookrightarrow C^{0,\mu}$ implies Hölder continuity.\\
(ii) Note that if $\alpha\le2$ in Lemma \ref{lem:iter} one can choose the exponent $0<\mu\le1$ from the previous statement as $\mu=\frac{\al}{2}.$ If $\al>2$ one can easily reduce the exponent by the estimate $r^\al\le r^2,$ where we assumed $r<1,$ which does not affect the generality.
\end{re}

\section{From $W^{1,2}$ to $W_{loc}^{2,2}$ for any $0<\ga<1$ via Difference Quotients}
\label{sec:2.3}
In the last paragraph we have seen that all stationary points must be locally Hölder continuous. It is natural to ask whether the regularity can be improved even further. In this section we show that if $0<\ga<1,$ i.e. the integrand is uniformly convex, then stationary points of the functional $I,$ as defined in \eqref{eq:1.1}, are in $W_{loc}^{2,2}$. This agrees from what we would expect by standard theory. But one of the requirements, to apply that theory, is that $\grad_\xi^2W(\grad u(\cd))$ needs to be locally bounded from above. But, again, $\xi\mapsto\grad_\xi^2W(\xi)$ grows quadratically in $\xi.$\\ 
Why does the standard method fail precisely? Usually, one would like to test the ELE \eqref{eq:ELE1.1} with $\vp^s(x):=-\D^{-h,s}(\eta^2(x)\D^{h,s}u(x)),$ where $\D^{h,s}f(x)=h^{-1}(f(x+he_s)-f(x))$ and $\D^{-h,s}f(x)=h^{-1}(f(x)-f(x-he_s))$ for a.e. $x\in\Om_h:=\{x\in\Om:\dist(x,\p\Om)>h\},$ $h\in \R_0^+,$ $s=1,2,$ and where $\{e_1,e_2\}$ denotes the standard basis. This yields the equation 
\begin{eqnarray}
\int\limits_{\Om_h}{\eta^2|\grad \Delta^{h,s}u|^2+\eta^2\Delta^{h,s}(\rho'(d)\cof\grad u)\cdot\grad\Delta^{h,s}u\; dx}\nonumber\\
=-2\int\limits_{\Om_h}{\eta\Delta^{h,s}\grad_\xi W(\grad u)\cdot(\grad\eta\ot \Delta^{h,s}u) \; dx},
\end{eqnarray}
which is very similar to \eqref{eq:ELE1.2}. The first major difficulty is that the best bound from below on the 2nd term of the LHS, we are aware of is \eqref{eq:2.6}.
By applying this bound we get the inequality 
\begin{eqnarray}
(1-\ga)\int\limits_{\Om_h}{\eta^2|\grad \Delta^{h,s}u|^2\; dx}\le C\int\limits_{\Om_h}{|\eta||\Delta^{h,s}\grad_\xi W(\grad u)||\grad\eta||\Delta^{h,s}u| \; dx},
\label{eq:2.3.1001}
\end{eqnarray}
which is problematic, since the LHS becomes negative if $\ga\ge1.$ But even if $\ga<1$ we need, for the RHS to be bounded, that the quantity $\Delta^{h,s}\grad_\xi W(\grad u)\in L^2(\Om_h).$ But 
\begin{align}
\D^{h,s}D_\xi W(\grad u(x))=\int\limits_0^1{\frac{d}{dt}D_\xi W(\grad u_t(x))\;dt}=\int\limits_0^1{D_\xi^2 W(\grad u_t(x))\;dt}\D^{h,s}\grad u(x),
\label{eq:3.DQ1}
\end{align}
where $u_t(x):=tu(x+he_s)+(1-t)u(x)$ is the convex combination for all $x\in \Om_h,$ $t\in[0,1].$ So our quantity is in $L^2(\Om_h),$ if $\int\limits_0^1{D_\xi^2 W(\grad u_t(\cdot))\;dt}\in L^\infty(\Om_h),$ since $\D^{h,s}\grad u\in L^2(\Om_h).$ But unfortunately, because $|D_\xi^2 W(\grad u_t(\cd))|\sim|\grad u_t(\cd)|^2,$ $D_\xi^2 W(\grad u_t(\cdot))$ is in $L^1(\Om_h),$ but not necessarily in the better space $L^2(\Om_h).$ Hence, the RHS of \eqref{eq:2.3.1001} might blow up.\\

Situations, where the second derivative of a (uniformly convex) integrand might be unbounded, have been studied in \cite{DM09}. De Maria's genius idea is to play the problematic term against a slight variation of itself. For this consider the generic term
\begin{equation}\int\limits_{\Om_h}{(\eta^2(x+he_s)-\eta^2(x))f(\D^{h,s}D_\xi W(\grad u(x)))\;dx}, \label{eq:2.3.1000a}\end{equation}
where $f$ is just a generic function. 
Now for the most part of the set $\Om_h$ the term \eqref{eq:2.3.1000a} disappears. Even though the integrand $f$ might be large, there is only a small set, on which $f$ can contribute. As a consequence, this introduces another `smallness' to the integral, which might be enough to control its behaviour. However, we will see, that even on this smaller set the integral might blow up if $\ga\ge1.$\\
 
How does one need to modify the test function in order to get in a position s.t. the crucial terms appear in this regularised form, as described in \eqref{eq:2.3.1000a}? 
We will be using 
\[\vp(x)=-\D^{-h,s}((\eta^2(x+he_s)-(1-\al)\eta^2(x))\D^{h,s} u(x)),\]
a.e. in $\Om_h$ and for some $0<\al<1.$ Note, in \cite{DM09} a slightly simpler version of the test function is used. The reason is, that compared to the work of De Maria we do not explicitly use the strong ellipticity inequality, even though it is satisfied in the regime $0<\ga<1.$ We can avoid this problem by the concrete form of our integrand $W$ and by introducing an additional small parameter $\al.$ This has the advantage, that our method could be used for arguments when $\ga$ exceeds $1,$ even though it seems unlikely that Theorem \ref{thm:2.3.DM08} can be extended to the regime $\ga\ge 1.$ However, it might be possible to use this argument to show, for $\ga\ge 1,$ with the natural modifications, higher integrability for stationary point in $W_{loc}^{2,q}$ for some(or any)
$1\le q<2.$ But, at the present, this remains open. We want to emphasise another time, that the fact $0<\ga<1$ will be crucial in the following discussion.\\
 
We start with the main result of this section, stating, that all stationary points of the functional \eqref{eq:1.1} are of class $W_{loc}^{2,2}.$ 

\begin{thm} Let $0< \ga<1$ and $u\in \A_{u_0}(\Om,\R^2)$ be a stationary point of the functional \eqref{eq:1.1}. Then $u\in W_{loc}^{2,2}(\Om,\R^2).$
\label{thm:2.3.DM08}
\end{thm}
\textbf{Proof:}\\
Recall the ELE
\begin{equation}
\int\limits_{\Om}{(\grad u+\rho'(d)\cof \grad u) \cdot\grad\vp\; dx}=0\; {\mbox{for all}}\; \vp\in W_0^{1,2}(\Om,\R^2).
\end{equation}

Let $x_0\in \Om$ and $r>0$ s.t. $B(x_0,3r)\ss\ss \Om.$ We will use the notation $\xb:=x+he_s$ and $x^-:=x-he_s$ for every $x\in\Om_h:=\{x\in\Om:\dist(x,\p\Om)>h\}.$ Further, choose $h_0:=\frac{r}{10}$ and let $0<h<h_0,$ then by construction $\supp{\eta(\cdot)},\supp{\eta(\cd\pm he_s)}\ss\ss B(x_0,\frac{7r}{4}).$ As a test function, motivated by \cite{DM09}, we choose 
\begin{equation}\vp=-\D^{-h,s}(\ta^h(x)\D^{h,s} u),\;\mb{with}\; \ta^h(x):=\eta^2(\xb)-(1-\al)\eta^2(x) \label{eq:2.3.98}\end{equation}
where $0<\al<1$ and $\eta\in C_c^\infty(\Om)$ is a standard mollifier satisfying the properties $\eta\equiv1$ in $B_r=B(x_0,r),$ $\eta\equiv0$ in $\Om\sm B(x_0,\frac{3r}{2}),$ $0\le\eta\le1$ and there exists $c>0$ s.t. $|\grad\eta|\le\frac{c}{r}$ and $|\grad^2\eta|\le\frac{c}{r^2}.$
Differentiating \eqref{eq:2.3.98} gives
\[\grad\vp=-\D^{-h,s}(\ta^h\D^{h,s} \grad u)-\D^{-h,s}(\grad \ta^h\ot\D^{h,s} u),\]

Then the ELE becomes{\small{
\begin{equation}{
\small{-\int\limits_{B_{2r}}{(\grad u+\rho'(d)\cof \grad u)\cd\D^{-h,s}(\ta^h\D^{h,s} \grad u)\;dx}
=\int\limits_{B_{2r}}{(\grad u+\rho'(d)\cof \grad u)\cd\D^{-h,s}(\grad \ta^h\ot\D^{h,s} u)\;dx}.}}
\label{eq:2.3.99}
\end{equation}}}
First recall the integration by parts formula for difference quotients for two functions $v,w\in L^2(\Om)$ of which at least one has support in $\Om_h,$ then 
\[\int\limits_{\Om}{v\D^{h,s}w\;dx}=-\int\limits_{\Om}{(\D^{-h,s}v)w\;dx}.\]
Applying this to the LHS of \eqref{eq:2.3.99} gives
\begin{flalign*}
LHS=&\int\limits_{B_{2r}}{\ta^{h}\D^{h,s}(\grad u+\rho'(d)\cof \grad u)\cd\D^{h,s}\grad u\;dx}&\\
=&\int\limits_{B_{2r}}{\eta^2(\xb)\D^{h,s}\grad u(x)\cd\D^{h,s}\grad u\;dx}-(1-\al)\int\limits_{B_{2r}}{\eta^2(x)\D^{h,s}\grad u(x)\cd\D^{h,s}\grad u\;dx}&\\
&+\int\limits_{B_{2r}}{\ta^{h}\D^{h,s}(\rho'(d)\cof \grad u)(x)\cd\D^{h,s}\grad u\;dx},&
\end{flalign*}
where we used the definition of $\ta^h.$ 

Recalling the product rule for difference quotients, for $f:\Om_h\rightarrow \R$ and $g:\Om\rightarrow \R$ it holds that
\begin{equation}
\D^{h,s}(fg)(x)=\D^{h,s}(f)(x)g(x)+f(\bar{x})\D^{h,s}g(x).
\label{eq:2.3.100PRFD}
\end{equation}
We point out, that this formula remains true even if $f$ and $g$ are vector- or matrixvalued and for various types of products (for example, scalar- or tensorproducts). Applying the product rule to the RHS of \eqref{eq:2.3.99} yields
\begin{flalign*}
RHS=&\int\limits_{B_{2r}}{(\grad u+\rho'(d)\cof \grad u)\cd\grad\ta^{h}(x)\ot\D^{-h,s}\D^{h,s} u(x)\;dx}&\\
&+\int\limits_{B_{2r}}{(\grad u+\rho'(d)\cof \grad u)\cd\D^{-h,s}(\grad\ta^{h})(x)\ot\D^{h,s} u(x^-)\;dx}.&
\end{flalign*}
By plugging this back into the ELE \eqref{eq:2.3.99} and rearranging terms we get
\begin{flalign*}
\int\limits_{B_{2r}}{\eta^2(\xb)|\D^{h,s}\grad u|^2(x)\;dx}=&(1-\al)\int\limits_{B_{2r}}{\eta^2(x)|\D^{h,s}\grad u|^2(x)\;dx}&\\
&-\int\limits_{B_{2r}}{\ta^{h}\D^{h,s}(\rho'(d)\cof \grad u)(x)\cd\D^{h,s}\grad u\;dx}&\\
&+\int\limits_{B_{2r}}{(\grad u+\rho'(d)\cof \grad u)\cd\grad\ta^{h}(x)\ot\D^{-h,s}\D^{h,s} u(x)\;dx}&\\
&+\int\limits_{B_{2r}}{(\grad u+\rho'(d)\cof \grad u)\cd\D^{-h,s}(\grad\ta^{h})(x)\ot\D^{h,s} u(x^-)\;dx}&\\
=:&(i)+(ii)+(iii)+(iv).&
\end{flalign*}
\textbf{About $(i):$}
We start with the following estimate: for $\be\in\R,$ we have
\begin{flalign*}\be^2\eta^2(x)=|\be(\eta(x)-\eta(\xb))+\be\eta(\xb)|^2
=\be^2(\eta(x)-\eta(\xb))^2+\be^2\eta^2(\xb)+2\be^2(\eta(x)-\eta(\xb))\eta(\xb).
\end{flalign*}
Then by Cauchy's inequality we get, for some $c'>0,$
\begin{flalign*}2\be^2(\eta(x)-\eta(\xb))\eta(\xb)\le\frac{1}{c'}(\eta(x)-\eta(\xb))^2+c'\be^4\eta^2(\xb).
\end{flalign*}
Choosing $\be^2=(1-\al)$ yields,
\begin{flalign*}(1-\al)\eta^2(x)\le\left((1-\al)+\frac{1}{c'}\right)(\eta(x)-\eta(\xb))^2+(1-\al)(1+c'(1-\al))\eta^2(\xb).
\end{flalign*}

Then by the latter estimate we can estimate $(i)$ by
\begin{flalign}
(1-\al)\int\limits_{B_{2r}}{\eta^2(x)|\D^{h,s}\grad u|^2\;dx}
\le& \left((1-\al)+\frac{1}{c'}\right)\int\limits_{B_{2r}}{|\eta(\xb)-\eta(x)|^2|\D^{h,s}\grad u|^2\;dx}&\nonumber\\
&+(1-\al)(1+c'(1-\al))\int\limits_{B_{2r}}{\eta(\xb)^2|\D^{h,s}\grad u|^2\;dx}&\nonumber\\
=:&(i.a)+(i.b).&\label{eq:2.3.1000a+b}
\end{flalign}
Using the mean value theorem and $|\grad \eta|\le\frac{c}{r}$ we can control $(i.a)$ via
\begin{flalign}
(i.a)\le&h^2\left((1-\al)+\frac{1}{c'}\right)\int\limits_{B_{\frac{7r}{4}}}{|\grad\eta|^2|\D^{h,s}\grad u|^2\;dx}&\nonumber\\
\le&\frac{h^2c^2\left((1-\al)+\frac{1}{c'}\right)}{r^2}\int\limits_{B_{\frac{7r}{4}}}{|\D^{h,s}\grad u|^2\;dx}&\nonumber\\
\le&\frac{c^2\left((1-\al)+\frac{1}{c'}\right)}{r^2}\int\limits_{B_{\frac{7r}{4}}}{|\grad u(\xb)|^2+|\grad u(x)|^2\;dx}&\nonumber\\
\le&\frac{c^2\left((1-\al)+\frac{1}{c'}\right)}{r^2}\int\limits_{B_{2r}}{|\grad u(x)|^2\;dx},&
\label{eq:2.3.1000.a.2}
\end{flalign}
where we used
\[\int\limits_{B_{\frac{7r}{4}}}{|\grad u(\xb)|^2\;dx}\le\int\limits_{B_{2r}}{|\grad u(x)|^2\;dx}.\]
From \eqref{eq:2.3.1000a+b} together with \eqref{eq:2.3.1000.a.2} we get
\begin{equation}
(i)\le (1-\al)(1+c')\int\limits_{B_{2r}}{\eta^2(\xb)|\D^{h,s}\grad u|^2\;dx}+\frac{C}{r^2}\int\limits_{B_{2r}}{|\grad u(x)|^2\;dx},
\label{eq:2.3.100.i}
\end{equation}
where we simplified the coefficient in front of the first integral by using the estimate $(1+c'(1-\al))\le(1+c').$ It will be crucial that for every $0<\al<1$ it is possible to choose $c'>0$ so small that $(1-\al)(1+c')<1.$ Moreover, we introduced a generic constant $C>0,$ which is also allowed to depend on $\ga.$ We can be generous, with the constant in front of the rightmost integral, as long as it remains finite.\\

\textbf{About $(ii):$}\\
Rewriting,
$\ta^h(x)=\al\eta^2(\xb)+(1-\al)(\eta^2(\xb)-\eta^2(x)),$
yields,
\begin{flalign*}
(ii)=&-\int\limits_{B_{2r}}{\ta^h(x)\D^{h,s}(\rho'(d)\cof \grad u)(x)\cd\D^{h,s}\grad u\;dx}&\\
=&-\al\int\limits_{B_{2r}}{\eta^2(\xb)\D^{h,s}(\rho'(d)\cof \grad u)(x)\cd\D^{h,s}\grad u\;dx}&\\
&-(1-\al)\int\limits_{B_{2r}}{(\eta^2(\xb)-\eta^2(x))\D^{h,s}(\rho'(d)\cof \grad u)(x)\cd\D^{h,s}\grad u\;dx}&
\end{flalign*}

For the first term we use the lower bound of Lemma \ref{HR:LUBs}, which is given by
\begin{equation}
\D^{h,s}(\rho'(d)\cof \grad u)(x)\cd\D^{h,s}\grad u(x)\ge-\ga|\D^{h,s}\grad u(x)|^2\end{equation}
for the second one we use the upper bound, shown in the same lemma, given by
\begin{equation}
|\Delta^{h,s}(\rho'(d)\cof\grad u)(x)|\le\ga|\Delta^{h,s}\grad u(x)|+\frac{2\ga}{h}|\grad u(x)|
\label{eq:2.16}
\end{equation}
which hold a.e. in $\Om_{h}.$ Hence,
\begin{flalign*}
(ii)\le&\al\ga\int\limits_{B_{2r}}{\eta^2(\xb)|\D^{h,s}\grad u|^2\;dx}&\\
&+(1-\al)\ga\int\limits_{B_{2r}}{|\eta^2(\xb)-\eta^2(x)||\D^{h,s}\grad u|^2\;dx}&\\&+\frac{2(1-\al)\ga}{h}\int\limits_{B_{2r}}{|\eta^2(\xb)-\eta^2(x)||\grad u(x)||\D^{h,s}\grad u|\;dx}&\\
=:&(ii.a)+(ii.b)+(ii.c)&
\end{flalign*}
(ii.a) is already of the desired form. For (ii.b) we apply, again, the mean value theorem and $|\grad\eta|\le \frac{c}{r}$ to obtain

\begin{flalign}
(ii.b)\le&(1-\al)\ga\int\limits_{B_{2r}}{|\eta(x)-\eta(\xb)||\eta(x)+\eta(\xb)||\D^{h,s}\grad u|^2\;dx}&\nonumber\\
\le&(1-\al)\ga h\int\limits_{B_{\frac{7r}{4}}}{|\grad\eta||\eta(x)+\eta(\xb)||\D^{h,s}\grad u|^2\;dx}&\nonumber\\
\le&Ch^2\ga^2\int\limits_{B_{\frac{7r}{4}}}{|\grad\eta||\D^{h,s}\grad u(x)|^2\;dx}+\ve(1-\al)^2\int\limits_{B_{2r}}{|\eta(x)+\eta(\xb)|^2|\D^{h,s}\grad u(x)|^2\;dx}&\nonumber\\
\le&C(\ga) h^2\int\limits_{B_{\frac{7r}{4}}}{|\grad\eta|^2|\D^{h,s}\grad u(x)|^2\;dx}+2\ve(1-\al)^2\int\limits_{B_{2r}}{(\eta^2(x)+\eta^2(\xb))|\D^{h,s}\grad u(x)|^2\;dx}.&\label{eq:2.16.a}
\end{flalign}
Here the $\ve$ results from another application of Cauchy's inequality with the weight $\ve>0,$ where we can choose $\ve$ arbitrarily small. The counterpart $\ve^{-1}$ has been absorbed into $C.$\\
Then the first term of the latter line can be controlled, as before in \eqref{eq:2.3.1000.a.2}, by
\begin{flalign*}
C(\ga) h^2\int\limits_{B_{\frac{7r}{4}}}{|\grad\eta|^2|\D^{h,s}\grad u(x)|^2\;dx}\le \frac{C(\ga)}{r^2} \int\limits_{B_{\frac{7r}{4}}}{|\grad u(\xb)|^2+|\grad u(x)|^2\;dx}\le\frac{C(\ga)}{r^2} \int\limits_{B_{2r}}{|\grad u(x)|^2\;dx}.
\end{flalign*}
The latter term of \eqref{eq:2.16.a} can be estimated in a similar fashion as has been done for $(i):$
\begin{flalign*}
&2\ve(1-\al)^2\int\limits_{B_{2r}}{(\eta^2(x)+\eta^2(\xb))|\D^{h,s}\grad u(x)|^2\;dx}&\\
\le&(2\ve(1-\al)^2+2\ve(1-\al)^2(1+2c'\ve(1-\al)^2))\int\limits_{B_{2r}}{\eta^2(\xb)|\D^{h,s}\grad u|^2\;dx}+\frac{c}{r^2}\int\limits_{B_{2r}}{|\grad u(x)|^2\;dx}.&
\end{flalign*}
Together we get
\begin{equation}
(ii.b)\le (2\ve(1-\al)^2(2+2c''\ve(1-\al)^2))\int\limits_{B_{2r}}{\eta^2(\xb)|\D^{h,s}\grad u(x)|^2\;dx}+\frac{C(\ga)}{r^2} \int\limits_{B_{2r}}{|\grad u(x)|^2\;dx}.
\label{eq:2.3.100}
\end{equation}

A similar calculation can be done for (ii.c):
\begin{flalign}
(ii.c)=&+\frac{2(1-\al)\ga}{h}\int\limits_{B_{2r}}{|\eta^2(\xb)-\eta^2(x)||\grad u(x)||\D^{h,s}\grad u|\;dx}&\nonumber\\
\le&2(1-\al)\ga \int\limits_{B_{\frac{7r}{4}}}{|\grad\eta||\eta(x)+\eta(\xb)||\grad u(x)||\D^{h,s}\grad u|\;dx}&\nonumber\\
\le&C\ga^2\int\limits_{B_{\frac{7r}{4}}}{|\grad\eta|^2|\grad u(x)|^2\;dx}+\ve(1-\al)^2\int\limits_{B_{2r}}{|\eta(x)+\eta(\xb)|^2|\D^{h,s}\grad u(x)|^2\;dx}&\nonumber\\
\le&(2\ve(1-\al)^2(2+2c''\ve(1-\al)^2))\int\limits_{B_{2r}}{\eta^2(\xb)|\D^{h,s}\grad u(x)|^2\;dx}+\frac{C(\ga)}{r^2} \int\limits_{B_{2r}}{|\grad u(x)|^2\;dx}.&
\label{eq:2.3.101}
\end{flalign}

Collecting, (ii.a), \eqref{eq:2.3.100} and \eqref{eq:2.3.101} yields, 
\begin{equation}
(ii)\le (\al\ga+4\ve(1-\al)^2(2+2c''\ve(1-\al)^2))\int\limits_{B_{2r}}{\eta^2(\xb)|\D^{h,s}\grad u(x)|^2\;dx}+\frac{C(\ga)}{r^2} \int\limits_{B_{2r}}{|\grad u(x)|^2\;dx}.
\label{eq:2.3.100.ii}
\end{equation}

\textbf{About $(iii):$}\\
The third term can be estimated by
\begin{flalign*}
(iii)=&\int\limits_{B_{2r}}{(\grad u+\rho'(d)\cof \grad u)\cd\grad\ta^{h}(x)\ot\D^{-h,s}\D^{h,s} u(x)\;dx}&\\
\le&(1+\ga)\int\limits_{B_{2r}}{|\grad\ta^{h}(x)||\grad u||\D^{-h,s}\D^{h,s} u|\;dx}.&
\end{flalign*}
Next $\grad \ta^h$ is given by
\[
\grad\ta^{h}(x)=2\eta(\xb)\grad\eta(\xb)-2(1-\al)\eta(x)\grad\eta(x), \;{\mb{for all}}\; x\in\Om_h.
\]
By the latter expression and the triangle inequality we get
\begin{flalign*}
(iii)\le&2(1+\ga)\int\limits_{B_{2r}}{\eta(\xb)|\grad\eta(\xb)||\grad u||\D^{-h,s}\D^{h,s} u|\;dx}&\\
&+2(1+\ga)(1-\al)\int\limits_{B_{2r}}{\eta(x)|\grad\eta(x)||\grad u||\D^{-h,s}\D^{h,s} u|\;dx}&\\
=:&(iii.a)+(iii.b)&
\end{flalign*}
Both terms are very similar, let's start with the second one. By Young's inequality and $|\grad\eta|\le\frac{c}{r}$ we get
\begin{flalign*}
(iii.b)\le&\ve\int\limits_{B_{2r}}{\eta^2(x)|\D^{-h,s}\D^{h,s} u|^2\;dx}+\frac{(1+\ga)^2(1-\al)^2}{\ve}\int\limits_{B_{2r}}{|\grad\eta(x)|^2|\grad u(x)|^2\;dx}&\\
\le&\ve\int\limits_{B_{2r}}{\eta^2(x)|\D^{-h,s}\D^{h,s} u|^2\;dx}+\frac{c(1+\ga)^2}{\ve r^2}\int\limits_{B_{2r}}{|\grad u(x)|^2\;dx}&
\end{flalign*}
The first term on the RHS seems difficult. But to our advantage De Maria shows in \cite{DM09} how this term can be controlled. First we rewrite the integrand by
\begin{flalign}
&\ve\int\limits_{B_{2r}}{\eta^2(x)|\D^{-h,s}\D^{h,s} u(x)|^2\;dx}&\nonumber\\
=&\ve\int\limits_{B_{2r}}{\frac{1}{h^2}|\eta(x)\D^{h,s} u(x)-\eta(x)\D^{h,s} u(x^-)|^2\;dx}&\nonumber\\
=&\ve\int\limits_{B_{2r}}{\frac{1}{h^2}|\eta(x^-)\D^{h,s}u(x^-)-\eta(x)\D^{h,s}u(x)+(\eta(x)-\eta(x^-))\D^{h,s} u(x^-)|^2\;dx}&\nonumber\\
\le&2\ve\int\limits_{B_{2r}}{|\D^{-h,s}(\eta\D^{h,s}u)(x)|^2\;dx}+2\ve\int\limits_{B_{2r}}{|\D^{-h,s}\eta(x)|^2|\D^{h,s}u(x^-)|^2\;dx},&\label{HR:(iii).1}
\end{flalign}
where we used $(a+b)^2\le2(a^2+b^2)$ and the second term of \eqref{HR:(iii).1} can be covered as follows
\begin{flalign*}
2\ve\int\limits_{B_{2r}}{|\D^{-h,s}\eta(x)|^2|\D^{h,s}u(x^-)|^2\;dx}\le&\frac{2C\ve}{r^2}\int\limits_{B_{\frac{7r}{4}}}{|\D^{h,s}u(x^-)|^2\;dx}&\\
\le&\frac{2C\ve}{r^2}\int\limits_{B_{\frac{7r}{4}+h_0}}{|\D^{h,s}u(x)|^2\;dx}&\\
\le&\frac{2C\ve}{r^2}\int\limits_{B_{2r}}{|\grad u(x)|^2\;dx},&
\end{flalign*}
where we used Nirenberg's lemma, see \ref{A.Nir}, in the latter estimate. The first term of \eqref{HR:(iii).1} can be treated, again, using Nirenberg's lemma and the product rule, as follows
\begin{flalign*}
&2\ve\int\limits_{B_{2r}}{|\D^{-h,s}(\eta\D^{h,s}u)(x)|^2\;dx}&\\
\le&2\ve\int\limits_{B_{2r}}{|\grad(\eta\D^{h,s}u)(x)|^2\;dx}&\\
\le&4\ve\int\limits_{B_{2r}}{|\grad\eta(x)\D^{h,s}u(x)|^2\;dx}+4\ve\int\limits_{B_{2r}}{\eta^2(x)|\grad\D^{h,s}u(x)|^2\;dx}&\\
\le&\frac{c\ve}{r^2}\int\limits_{B_{\frac{7r}{4}}}{|\D^{h,s}u(x)|^2\;dx}+4\ve\int\limits_{B_{2r}}{\eta^2(x)|\grad\D^{h,s}u(x)|^2\;dx}&\\
\le&\frac{c\ve}{r^2}\int\limits_{B_{2r}}{|\grad u(x)|^2\;dx}+c\ve\int\limits_{B_{2r}}{\eta^2(\xb)|\grad\D^{h,s}u(x)|^2\;dx}.&
\end{flalign*}
In the last step we used Nirenberg's lemma to bound the first term, and step (i) in the rightmost term. Hence,
\begin{flalign*}
(iii.b)\le c\ve\int\limits_{B_{2r}}{\eta^2(\xb)|\D^{h,s}\grad u|^2\;dx}+\frac{C(\ga)}{r^2}\int\limits_{B_{2r}}{|\grad u(x)|^2\;dx}.
\end{flalign*}
The (iii.a)-term can be controlled analogously:
\begin{flalign*}
(iii.a)\le&\ve\int\limits_{B_{2r}}{\eta^2(\xb)|\D^{-h,s}\D^{h,s} u|^2\;dx}+\frac{(1+\ga)^2}{\ve}\int\limits_{B_{2r}}{|\grad\eta(\xb)|^2|\grad u(x)|^2\;dx}&\\
\le&\ve\int\limits_{B_{2r}}{\eta^2(\xb)|\D^{-h,s}\D^{h,s} u|^2\;dx}+\frac{c(1+\ga)^2}{\ve r^2}\int\limits_{B_{2r}}{|\grad u(x)|^2\;dx}&
\end{flalign*}
The first term can be covered by
\begin{flalign}
&\ve\int\limits_{B_{2r}}{\eta^2(\xb)|\D^{-h,s}\D^{h,s} u(x)|^2\;dx}&\nonumber\\
=&\ve\int\limits_{B_{2r}}{\frac{1}{h^2}|\eta(\xb)\D^{h,s} u(x)-\eta(\xb)\D^{h,s} u(x^-)|^2\;dx}&\nonumber\\
=&\ve\int\limits_{B_{2r}}{\frac{1}{h^2}|\eta(x)\D^{h,s}u(x^-)-\eta(\xb)\D^{h,s}u(x)+(\eta(\xb)-\eta(x))\D^{h,s} u(x^-)|^2\;dx}&\nonumber\\
\le&2\ve\int\limits_{B_{2r}}{|\D^{h,s}(\eta(x)\D^{h,s}u(x-he_s))|^2\;dx}+2\ve\int\limits_{B_{2r}}{|\D^{h,s}\eta(x)\D^{h,s}u(x^-)|^2\;dx}&\label{eq:2.3.102}
\end{flalign}
First term of \eqref{eq:2.3.102}:
\begin{flalign*}
&2\ve\int\limits_{B_{2r}}{|\D^{h,s}(\eta(x)\D^{h,s}u(x-he_s))|^2\;dx}&\\
\le&2\ve\int\limits_{B_{2r}}{|\grad(\eta(x)\D^{h,s}u(x^-))|^2\;dx}&\\
\le&4\ve\int\limits_{B_{2r}}{|\grad\eta(x)\D^{h,s}u(x^-)|^2\;dx}+4\ve\int\limits_{B_{2r}}{\eta^2(x)|\grad\D^{h,s}u(x^-)|^2\;dx}&\\
\le&\frac{c\ve}{r^2}\int\limits_{B_{2r}}{|\grad u(x)|^2\;dx}+c\ve\int\limits_{B_{2r}}{\eta^2(\xb)|\grad\D^{h,s}u(x)|^2\;dx}&
\end{flalign*}
Rightmost term of \eqref{eq:2.3.102}:
\begin{flalign*}
2\ve\int\limits_{B_{2r}}{|\D^{h,s}\eta(x)|^2|\D^{h,s}u(x^-)|^2\;dx}\le&\frac{c\ve}{r^2}\int\limits_{B_{\frac{7r}{4}}}{|\D^{h,s}u(x^-)|^2\;dx}&\\
\le&\frac{c}{r^2}\int\limits_{B_{\frac{7r}{4}+h_0}}{|\D^{h,s}u(x)|^2\;dx}&\\
\le&\frac{c}{r^2}\int\limits_{B_{2r}}{|\grad u(x)|^2\;dx}&
\end{flalign*}

Hence,
\[
(iii.a)\le c\ve\int\limits_{B_{2r}}{\eta^2(\xb)|\D^{h,s}\grad u|^2\;dx}+\frac{C(\ga)}{r^2}\int\limits_{B_{2r}}{|\grad u(x)|^2\;dx}.
\]

Combining the bounds for (iii.a) and (iii.b) yields
\begin{equation}
(iii)\le c\ve\int\limits_{B_{2r}}{\eta^2(\xb)|\D^{h,s}\grad u(x)|^2\;dx}+\frac{C(\ga)}{r^2} \int\limits_{B_{2r}}{|\grad u(x)|^2\;dx}.
\label{eq:2.3.100.iii}
\end{equation}

\textbf{About $(iv):$}\\
Starting, similarly to (iii) by
\begin{flalign*}
(iv)=&\int\limits_{B_{2r}}{(\grad u+\rho'(d)\cof \grad u)\cd\D^{-h,s}(\grad\ta^{h})(x)\ot\D^{h,s} u(x^-)\;dx}&\\
\le&(1+\ga)\int\limits_{B_{2r}}{|\grad u||\D^{-h,s}(\grad\ta^{h})(x)||\D^{h,s} u(x^-)|\;dx}.&
\end{flalign*}
We can write
\begin{flalign*}
\D^{-h,s}(\grad\ta^{h})(x)=&2\D^{-h,s}[\eta(\xb)\grad\eta(\xb)-(1-\al)\eta(x)\grad\eta(x)]&\\
=&2[\D^{-h,s}(\eta(\xb))\grad\eta(\xb)+\eta(\xb)\D^{-h,s}(\grad\eta(\xb))&\\
&-(1-\al)\D^{-h,s}(\eta(x))\grad\eta(x)-(1-\al)\eta(x)\D^{-h,s}(\grad\eta(x))].&
\end{flalign*}
Using this latter expansion we can estimate $(iv)$ by
\begin{flalign*}
(iv)=&\int\limits_{B_{2r}}{(\grad u+\rho'(d)\cof \grad u)\cd\D^{-h,s}(\grad\ta^{h})(x)\ot\D^{h,s} u(x^-)\;dx}&\\
\le&2(1+\ga)\int\limits_{B_{2r}}{|\grad u||\D^{-h,s}(\eta(\xb))||\grad\eta(x)||\D^{h,s} u(x^-)|\;dx}&\\
&+2(1+\ga)\int\limits_{B_{2r}}{|\grad u||\eta(\xb)||\D^{-h,s}(\grad\eta(\xb))||\D^{h,s} u(x^-)|\;dx}&\\
&+2(1-\al)(1+\ga)\int\limits_{B_{2r}}{|\grad u||\D^{-h,s}(\eta(x))||\grad\eta(x)||\D^{h,s} u(x^-)|\;dx}&\\
&+2(1-\al)(1+\ga)\int\limits_{B_{2r}}{|\grad u||\eta(x^-)||\D^{-h,s}(\grad\eta(x))||\D^{h,s} u(x^-)|\;dx}&\\
=:&(iv.a)+(iv.b)+(iv.c)+(iv.d).&
\end{flalign*}
We can treat the first part in the following way 
\begin{flalign*}
(iv.a)=&2(1+\ga)\int\limits_{B_{2r}}{|\grad u||\D^{-h,s}(\eta(\xb))||\grad\eta(x)||\D^{h,s} u(x^-)|\;dx}&\\
\le&C(\ga)\int\limits_{B_{2r}}{|\grad\eta(x)|^2|\grad u|^2\;dx}+\ve\int\limits_{B_{2r}}{|\D^{-h,s}(\eta(\xb))|^2|\D^{h,s} u(x^-)|^2\;dx}&\\
\le&\frac{C(\ga)}{r^2}\int\limits_{B_{2r}}{|\grad u|^2\;dx}+\ve\int\limits_{B_{\frac{7r}{4}}}{|\grad\eta|^2|\D^{h,s} u(x^-)|^2\;dx}&\\
\le&\frac{C(\ga)}{r^2}\int\limits_{B_{2r}}{|\grad u|^2\;dx}+\frac{\ve}{r^2}\int\limits_{B_{\frac{7r}{4}}}{|\D^{h,s} u(x^-)|^2\;dx}&\\
\le&\frac{C(\ga)}{r^2}\int\limits_{B_{2r}}{|\grad u|^2\;dx}+\frac{c\ve}{r^2}\int\limits_{B_{2r}}{|\grad u(x)|^2\;dx}&\\
\le&\frac{C(\ga)}{r^2}\int\limits_{B_{2r}}{|\grad u|^2\;dx},&
\end{flalign*}
where we used Young's inequality, the mean value theorem $|\grad\eta|<\frac{c}{r}$ and Nirenberg's lemma.
One can proceed quite similarly to part $(iv.c):$
\begin{flalign*}
(iv.c)=&2(1-\al)(1+\ga)\int\limits_{B_{2r}}{|\grad u||\D^{-h,s}(\eta(x))||\grad\eta(x)||\D^{h,s} u(x^-)|\;dx}&\\
\le&C(\ga)\int\limits_{B_{2r}}{|\grad\eta(x)|^2|\grad u|^2\;dx}+\ve\int\limits_{B_{2r}}{|\D^{-h,s}(\eta(x))|^2|\D^{h,s} u(x^-)|^2\;dx}&\\
\le&\frac{C(\ga)}{r^2}\int\limits_{B_{2r}}{|\grad u|^2\;dx}+\ve\int\limits_{B_{\frac{7r}{4}}}{|\grad\eta|^2|\D^{h,s} u(x^-)|^2\;dx}&\\
\le&\frac{C(\ga)}{r^2}\int\limits_{B_{2r}}{|\grad u|^2\;dx}+\frac{\ve}{r^2}\int\limits_{B_{\frac{7r}{4}}}{|\D^{h,s} u(x^-)|^2\;dx}&\\
\le&\frac{C(\ga)}{r^2}\int\limits_{B_{2r}}{|\grad u|^2\;dx}+\frac{c\ve}{r^2}\int\limits_{B_{2r}}{|\grad u(x)|^2\;dx}&\\
\le&\frac{C(\ga)}{r^2}\int\limits_{B_{2r}}{|\grad u|^2\;dx}.&
\end{flalign*}
For (iv.b) we use Young's inequality and the mean value theorem but this time we apply it to $\grad \eta,$ combined with $|\grad^2\eta|\le \frac{c}{r^2}$ we get
\begin{flalign*}
(iv.b)=&2(1+\ga)\int\limits_{B_{2r}}{|\grad u||\eta(\xb)||\D^{-h,s}(\grad\eta(\xb))||\D^{h,s} u(x^-)|\;dx}&\\
\le&C(\ga)\int\limits_{B_{2r}}{|\D^{-h,s}(\grad\eta(\xb))||\grad u|^2\;dx}+\ve\int\limits_{B_{2r}}{|\eta(\xb)||\D^{-h,s}(\grad\eta(\xb))||\D^{h,s} u(x^-)|^2\;dx}&\\
\le&\frac{C(\ga)}{r^2}\int\limits_{B_{2r}}{|\grad u|^2\;dx}+\frac{c\ve}{r^2}\int\limits_{B_{2r}}{|\eta(\xb)|^2|\D^{h,s} u(x^-)|^2\;dx}&\\
\le&\frac{C(\ga)}{r^2}\int\limits_{B_{2r}}{|\grad u|^2\;dx}+\frac{c\ve}{r^2}\int\limits_{B_{\frac{7r}{4}}}{|\D^{h,s} u(x^-)|^2\;dx}&\\
\le&\frac{C(\ga)}{r^2}\int\limits_{B_{2r}}{|\grad u|^2\;dx}&
\end{flalign*}
Similarly for (iv.d):
\begin{flalign*}
(iv.d)=&2(1-\al)(1+\ga)\int\limits_{B_{2r}}{|\grad u||\eta(x^-)||\D^{-h,s}(\grad\eta(x))||\D^{h,s} u(x^-)|\;dx}&\\
\le&C(\ga)\int\limits_{B_{2r}}{|\D^{-h,s}(\grad\eta(x))||\grad u|^2\;dx}+\ve\int\limits_{B_{2r}}{|\eta(\xb)||\D^{-h,s}(\grad\eta(x))||\D^{h,s} u(x^-)|^2\;dx}&\\
\le&\frac{C(\ga)}{r^2}\int\limits_{B_{2r}}{|\grad u|^2\;dx}+\frac{c\ve}{r^2}\int\limits_{B_{2r}}{|\eta(\xb)|^2|\D^{h,s} u(x^-)|^2\;dx}&\\
\le&\frac{C(\ga)}{r^2}\int\limits_{B_{2r}}{|\grad u|^2\;dx}+\frac{c\ve}{r^2}\int\limits_{B_{\frac{7r}{4}}}{|\D^{h,s} u(x^-)|^2\;dx}&\\
\le&\frac{C(\ga)}{r^2}\int\limits_{B_{2r}}{|\grad u|^2\;dx}.&
\end{flalign*}

Together,
\begin{equation}
(iv)\le \frac{C(\ga)}{r^2} \int\limits_{B_{2r}}{|\grad u(x)|^2\;dx}.
\label{eq:2.3.100.iv}
\end{equation}
\textbf{Conclusion:}\\
Collecting, \eqref{eq:2.3.100.i}, \eqref{eq:2.3.100.ii}, \eqref{eq:2.3.100.iii} and \eqref{eq:2.3.100.iv} yields the following inequality
\begin{flalign}
\int\limits_{B_{2r}}{\eta^2(\xb)|\D^{h,s}\grad u|^2(x)\;dx}\le&((1-\al)(1+c')+\al\ga+c\ve)\int\limits_{B_{2r}}{\eta^2(\xb)|\D^{h,s}\grad u(x)|^2\;dx}&\nonumber\\
+&\frac{C(\ga)}{r^2} \int\limits_{B_{2r}}{|\grad u(x)|^2\;dx}.\label{eq:2.3.100.v}
\end{flalign}
We claim that the prefactor of the first term of the RHS can be chosen to be strictly smaller than 1 if $0\le\ga<1$ and $\al,\ve$ small enough. To see this for any $0<\al,\ga<1$ choose $c'>0$ s.t. $1+c'=(1-\al(1-\ga))^{-1}.$ Then 
\[1-((1-\al)(1+c')+\al\ga)=\frac{\al^2\ga(1-\ga)}{1-\al(1-\ga)}>0\]
Then there is even enough space to fit in some small $c\ve>0$ for $\ve>0$ small enough, s.t.
\[(1-\al)(1+c')+\al\ga+c\ve<1.\]
This allows one to absorb the first term into the LHS of \eqref{eq:2.3.100.v} yielding
\begin{flalign}
\int\limits_{B_{2r}}{\eta^2(\xb)|\D^{h,s}\grad u|^2(x)\;dx}\le&\frac{C(\ga)}{r^2}\int\limits_{B_{2r}}{|\grad u(x)|^2\;dx},&
\end{flalign}
where the RHS is independent of $h.$ Finally, Nirenberg's lemma implies that $u\in W_{loc}^{2,2}(\Om,\R^2).$\vspace{1cm}

This section is completed by a technical lemma establishing the important upper and lower bounds of the difference quotients on quantities arising from the nonlinear part of $W.$

\begin{lem} Let $h\in\R^+,$ $s\in\{0,1\},$ $u\in W^{1,2}(\Om,\R^2)$ and $\rho$ be as before. Then\\
(i) the lower bound
\begin{equation}
\Delta^{h,s}(\rho'(d)\cof\grad u)(x)\cdot\Delta^{h,s}\grad u(x)\ge-\ga|\Delta^{h,s}\grad u(x)|^2
\label{eq:2.6}
\end{equation}
and\\
(ii) the upper bound
\begin{equation}
|\Delta^{h,s}(\rho'(d)\cof\grad u)(x)|\le\ga|\Delta^{h,s}\grad u(x)|+\frac{2\ga}{h}|\grad u(x)|
\label{eq:2.16}
\end{equation}
hold a.e. in $\Om_{h}.$
\label{HR:LUBs}
\end{lem}
{\bf{Proof:}}\\
(i) Let $x\in \Om_h,$ fix an $s\in \{0,1\}$ and define $\bar{x}:=x+he_s.$ We start by applying the product rule for difference quotients \eqref{eq:2.3.100PRFD} to obtain 
\begin{eqnarray}
\Delta^{h,s}(\rho'(d)\cof\grad u)(x)\cdot\Delta^{h,s}\grad u(x)&=&\Delta^{h,s}(\rho'(d))(x)(\cof\grad u)(x)\cdot\Delta^{h,s}\grad u(x)\nonumber\\
&&+\rho'(d)(\bar{x})\Delta^{h,s}(\cof\grad u)(x)\cdot\Delta^{h,s}\grad u(x)\nonumber\\
&=:&I+II.
\label{eq:MP,PR2}
\end{eqnarray}
By linearity of the gradient and the cofactor in two dimensions it holds that $\Delta^{h,s}\grad u(x)=\grad\Delta^{h,s} u(x)$ and $\Delta^{h,s}(\cof\grad u)(x)=\cof\Delta^{h,s}\grad u(x)$ for all $x\in \Om_h.$ Let us now introduce the following notation $A:=\grad u(\bar{x}),$ $B:=\grad u(x).$ Note that for every matrix $N\in \R^{2\ti2},$ the determinant satisfies the following scaling behaviour $d_{N/h}=h^{-2}d_N$ for $h\in \R\setminus\{0\}.$ Using these properties and definitions to rewrite the second term of \eqref{eq:MP,PR2} as 
\begin{eqnarray}
II&=&\rho'(d)(\bar{x})\Delta^{h,s}(\cof\grad u)(x)\cdot\Delta^{h,s}\grad u(x)\nonumber\\
&=&\rho'(d)(\bar{x})\cof\Delta^{h,s}\grad u(x)\cdot\Delta^{h,s}\grad u(x)\nonumber\\
&=&h^{-2}\rho'(d_A)\cof (A-B) \cdot(A-B)\nonumber\\
&=&2h^{-2}\rho'(d_A)d_{A-B}.
\label{eq:2.7}
\end{eqnarray}
The first term of (\ref{eq:MP,PR2}) takes the form
\begin{eqnarray}
I=h^{-1}\Delta^{h,s}(\rho'(d))(x)\cof B\cdot(A-B).
\label{eq:2.8}
\end{eqnarray}
For two matrices $N,N'\in\R^{2\times2}$ the determinant of $N+N'$ can be expanded in the following way
\begin{eqnarray}
d_{N+N'}=d_{N'}+d_{N}+\cof N\cdot N'.
\label{eq:2.9}
\end{eqnarray}
Now choose $N'=A-B$ and $N=B$ and we get 
\begin{eqnarray}
d_{A}=d_{A-B}+d_{B}+\cof B\cdot (A-B)
\label{eq:2.10}
\end{eqnarray}
or 
\begin{eqnarray}
\cof B\cdot (A-B)=d_{A}-d_{B}-d_{A-B}=h\Delta^{h,s}d(x)-d_{A-B}.
\label{eq:2.11}
\end{eqnarray}
Therefore,
\begin{eqnarray}
I=\Delta^{h,s}(\rho'(d))(x)\Delta^{h,s}d(x)-h^{-1}\Delta^{h,s}(\rho'(d))(x)d_{A-B}.
\label{eq:2.12}
\end{eqnarray}
The whole expression can now be written as
\begin{eqnarray}
I+II&=&\Delta^{h,s}(\rho'(d))(x)\Delta^{h,s}d(x)\nonumber\\
&&+h^{-2}(2\rho'(d_A)-(\rho'(d_A)-\rho'(d_B))d_{A-B}\nonumber\\
&=&\Delta^{h,s}(\rho'(d))(x)\Delta^{h,s}d(x)+h^{-2}(\rho'(d_A)+\rho'(d_B))d_{A-B}\nonumber\\
&=:&III+IV.
\label{eq:MP,PR3}
\end{eqnarray}
The latter term can be bounded from below through
\begin{eqnarray}
IV&\ge& -h^{-2}|\rho'(d_A)+\rho'(d_B)||d_{A-B}|\nonumber\\
&\ge&-h^{-2}\ga|A-B|^2,
\label{eq:2.13}
\end{eqnarray}
where we used $\rho'\le\ga$ and $|d_N|\le\frac{1}{2}|N|^2.$ Furthermore, realize that $III\ge0$  is always true, since $\rho$ is convex and $\rho'$ satisfies a monotonicity inequality. Indeed, consider first $\Delta^{h,s}d(x)\ge0$ then $d_A\ge d_B$ and by monotonicity of $\rho',$ $\Delta^{h,s}\rho'(d)(x)\ge0.$ Hence, $III\ge0.$ Secondly, let $\Delta^{h,s}d(x)<0$ then $d_A< d_B$ implying $\rho'(d_A)\le\rho'(d_B)$ and $\Delta^{h,s}\rho'(d)(x)\le0$ and again $III\ge0.$

This implies
\begin{equation}
\Delta^{h,s}(\rho'(d)\cof\grad u)\cdot\Delta^{h,s}\grad u\ge-\ga|\Delta^{h,s}\grad u|^2.
\label{eq:2.14}
\end{equation}
\\

(ii) The upper bound follows again from the product rule \eqref{eq:2.3.100PRFD} and the triangle inequality
\begin{eqnarray}
|\Delta^{h,s}(\rho'(d)\cof\grad u)(x)|&=&|\rho'(d)(\bar{x})\Delta^{h,s}(\cof\grad u)(x)+\Delta^{h,s}(\rho'(d))(x)(\cof\grad u)(x)|\nonumber\\
&\le&\ga|\Delta^{h,s}\grad u(x)|+\frac{2\ga}{h}|\grad u(x)|.
\end{eqnarray}
\vspace{1cm}

\begin{re}
1. A warning might be in order here. By Standard Theory one would expect that if the integrand  is uniformly convex(that is $0<\ga<1$), $u\in W_{loc}^{2,2}(\Om,\R^2)$ then in 2d one would expect that by Meyers' theorem the integrability is automatically improved to $u\in W_{loc}^{2,p}(\Om,\R^2)$ for some $p>2$ and from there Schauder Theory applies.  
However, recall that Meyers' is applied in the following way: By $u\in W_{loc}^{2,2}(\Om,\R^2)$ one knows that $\grad u\in W_{loc}^{1,2}(\Om,\R^4)$ solves the `linearized' elliptic PDE 
\[
-\div(A(x)\grad^2 u)=0\;\mb{in}\;D'(\Om,\R^4),
\]
 where $A(\cd):=\grad_\xi^2W(\grad u(\cd))$ is the coefficients matrix. Note, that $\xi\mapsto W(\xi)$ is smooth, hence the regularity of $A$ is determined by the regularity of $\grad u.$ Now to satisfy all the assumptions of Meyers' Theorem one needs $A\in L_{loc}^\infty.$  But, again, $\xi\mapsto\grad_\xi^2W(\xi)$ grows quadratically in $\xi.$ Hence, at the present stage, we can only guarantee $A\in L^2_{loc},$ but not necessarily $A\in L_{loc}^\infty.$ Hence, we need to show higher integrability by hand, this can be done for a general $\ga>0$. This is the subject of the next section.\\

2.  It remains an open question, if Theorem \ref{thm:2.3.DM08} can be extended to the range when $\ga\ge1?$ We don't expect this. Indeed, assume for a second that we could. Then it is very likely that the argument could be extended to general polyconvex functionals with $p-$growth. This clearly would include the stationary points, which are Lipschitz but not any better, constructed by Szekelyhidi, in the famous paper \cite{Sz04}, and together with the result in the next section would show that they actually have to be smooth, which would lead to a contradiction. However, we know by the seminal work of Kristensen and Mingione, in \cite{KM07}, that the Szekelyhidi construction is forbidden for global ($\omega-$)minimizers.\footnote{Recall that we call a map $u\in W_{u_0}^{1,q}(\Om,\R^m),$ for some $1\le q<\infty,$ $\omega-$ or almost minimizer, if there is a non-decreasing concave function $\om:\R^+\ra \R^+$ satisfying $\om(0)=0$ s.t. for any $B_R\ss \Om$ it holds
\[\int\limits_{B_R} {F(x,u(x),\grad u(x))\;dx}\le (1+\om(R^2))\int\limits_{B_R} {F(x,v(x),\grad v(x))\;dx}\] for any $u-v\in W_{0}^{1,q}(B_R,\R^m).$ See, for instance \cite[Definition 3]{KM07}.} Hence, there is hope to show that global ($\omega-$)minimizers are in $W_{loc}^{2,2}(\Om,\R^2).$ 
It might also be possible to show, that global ($\omega-$)minimizers or stationary points are of the class $W_{loc}^{2,q}(\Om,\R^2)$ for some $1\le q<2.$
We think that this might be possible to show the latter, by altering the proof of Theorem \ref{thm:2.3.DM08}.
\end{re}

\section{From $W_{loc}^{2,2}$ to $C_{loc}^{\infty}$ for any $0<\ga<\infty.$ Higher-order regularity and Schauder Theory}
\label{sec:2.4}

The information that $u$ is of class $W_{loc}^{2,2}$ allows one to take partial derivatives of $\grad_\xi W(\grad u):$ For $i,j,k\in\{1,2\}$ we get 
\begin{flalign*}
\p_k(\grad_\xi W(\grad u))_{ij}=&\p_{k}\p_iu_j+\rho'(d)\p_{\xi_{ab}}\p_{k}(\cof\grad u)_{ij}&\\
&+\rho''(d)(\cof\grad u)_{ij}(\cof\grad u)_{ab}\p_{k}\p_au_b&\\
=&(\grad\p_{k}u+\rho'(d)\cof\grad u,_k+\rho''(d)(\cof\grad u\cd\grad\p_{k}u)\cof\grad u)_{ij}&.
\end{flalign*}
A priori by considering the $L^p-$spaces the first and second term on the RHS are in $L^2.$
The rightmost one and therefore the quantity on the LHS, however, are in $L^q$ for all $1\le q<2$ but not necessarily in a better space. \\
In case we want to emphasise the quasilinear structure of this quantity, a different way of expressing it is possible via
\begin{flalign*}
\p_k(\grad_\xi W(\grad u))_{ij}=&\d_{ia}\d_{jb}\p_{k}\p_au_b+\rho'(d)\p_{\xi_{ab}}((\cof\xi)_{ij})\p_{k}\p_au_b&\\
&+\rho''(d)(\cof\grad u)_{ij}(\cof\grad u)_{ab}\p_{k}\p_au_b&\\
=&(\grad_{\xi}^2 W(\grad u))_{ijab}\p_{k}\p_au_b&
\end{flalign*}
where 
\[(\grad_{\xi}^2 W(\grad u))_{ijab}=\d_{ia}\d_{jb}+\rho'(d)\p_{\xi_{ab}}((\cof\xi)_{ij})+\rho''(d)(\cof\grad u)_{ij}(\cof\grad u)_{ab}.\]

Now for every $k\in \{1,2\}$ we can test against $\vp_k\in C_c^\infty(\Om,\R^{2})$ and get
\begin{equation}
\sum\limits_{i,j}\int\limits_\Om {\p_k(\grad_\xi W(\grad u))_{ij}\p_{i}\vp_{kj}\;dx}=0.
\label{eq: 2.HR.1}
\end{equation}
Instead of having a system of equations we can sum over all $k:$
\begin{equation}
\sum\limits_{i,j,k}\int\limits_\Om {\p_k(\grad_\xi W(\grad u))_{ij}\p_{i}\vp_{kj}\;dx}=0,
\label{eq: 2.HR.2}
\end{equation}
needs to be satisfied for arbitrary $\vp_k\in C_c^\infty(\Om,\R^{2}),$ $k\in \{1,2\}.$\footnote{Note, this is indeed the same. It's easy to see that \eqref{eq: 2.HR.2}  follows from \eqref{eq: 2.HR.1}. So assume that \eqref{eq: 2.HR.2} holds and for the sake of a contradiction, that \eqref{eq: 2.HR.1} is false. Wlog. there exists $\vp_1$ s.t.  $\int\limits_\Om {\p_1(\grad_\xi W(\grad u))_{ij}\p_{x_i}\vp_{1j}\;dx}\not=0.$ By \eqref{eq: 2.HR.2} $\int\limits_\Om {\p_2(\grad_\xi W(\grad u))_{ij}\p_{x_i}\vp_{2j}\;dx}\not=0$ for all $\vp_2\in C_c^\infty(\Om,\R^{2})$ in particular $\vp_2\equiv0,$ a contradiction.} As usual we suppress the sums from now on.\\

By introducing $\phi:=(\vp_1,\vp_2)\in C_c^\infty(\Om,\R^{2\ti2})$ we can write this PDE in a closed form 
\[\int\limits_{\Om}\grad_{\xi}^2 W(\grad u)\grad^2 u\cd\grad\phi=0,\;\mb{for all}\;\phi\in C_c^\infty(\Om,\R^{2\ti2})\]
The multiplication needs to be understood in the above sense.\\
We want to emphasize again that the quantity $\grad_\xi^2W(\grad u(\cd))\grad^2 u(\cd)$ is in $L^q,$ for all $1\le q<2,$ but not necessarily in $L^2.$ Therefore, we can only test with functions in $W_0^{1,q'}$, where $2<q'<\infty$ is the dual Hölder exponent.\\

Now we can state the main result of this section. We show that every stationary point of \eqref{eq:1.1} which is additionally in $W_{loc}^{2,2}$ is already $C_{loc}^{1,\al}.$ This is obtained by establishing a reverse Hölder inequality.\\

\begin{thm}
Let $0<\ga<\infty,$ $\Om\ss\R^2$ be open and bounded, let $p\in[\frac{8}{5},2)$ and the corresponding $\beta:=2(\frac{2}{p}-1)>0,$ let $u\in \A_{u_0}\cap W_{loc}^{2,2}(\Om,\R^2),$ let $x\mt d_{\grad u(x)}:=\det\grad u(x)\in C_{loc}^{0,\be}(\Om)$ and $u$ solves the following quasilinear PDE in a weak sense, i.e.
\begin{equation}\int\limits_{\Om}\grad_{\xi}^2 W(\grad u)\grad^2 u\cd\grad\phi=0\;\mb{for all}\;\phi\in {W_{0}^{1,q'}(\Om,\R^{2\ti2})} \;\mb{and any}\;2<q'<\infty. \label{quasiPDE1}\end{equation}

Then there is $\d>0$ s.t. $u\in W_{loc}^{2,2+\d}(\Om,\R^2).$ In particular, there are constants $\d=\d(\ga,\|\rho''\|_{L^\infty},p)>0$ and $C=C(\ga,\|\rho''\|_{L^\infty},p)>0$ s.t. the following Reverse Hölder Inequality is satisfied
\begin{equation}
\left(\fint\limits_{B_{r}}{|\grad^2 u|^{2+\d}\;dx}\right)^{\frac{1}{2+\d}}\le C\left(\fint\limits_{B_{2r}}{|\grad^2 u|^2\;dx}\right)^{\frac{p}{2}}\;\mb{for any}\; B(x,2r)\ss\ss B.
\end{equation}
Moreover, there exists $\al\in(0,1]$ s.t. $u\in C_{loc}^{1,\al}(\Om,\R^2).$
\end{thm}
\textbf{Proof:}\\
\textbf{Step 1: Approximating $u.$} First note that by Sobolev embedding we have higher integrability $u\in W_{loc}^{1,p},$ for every $1\le p<\infty.$\\
Define $\tilde{u}:=u\chi_{B(x_0,2r)}$ and further $u_\ve:=\eta_\ve*\tilde{u}.$\\
Note that $\supp\tilde{u}\ss B(x_0,2r)$ and by definition $\supp\eta_\ve\ss B(0,\ve).$  Recall the following convergence properties of the mollification: 
\begin{align*}
&u_\ve\ra u \;\mb{in}\; W^{1,p} \;\mb{for all}\;1\le p<\infty&\\
&\grad^2 u_\ve\ra \grad^2 u \;\mb{in}\; L^2&
\end{align*}
These properties can be found in \cite[Appendix C.4]{LE10}.\\
An easy consequence is that 
\[\det\grad u_\ve\ra\det\grad u \;\mb{in}\; L^p\]
for every $1\le p<\infty.$ Indeed, since
\begin{flalign*} 
\|\det\grad u_\ve-\det\grad u\|_{L^p}=&\frac{1}{2}\|\cof\grad u_\ve\cd\grad u_\ve-\cof\grad u\cd\grad u\|_{L^p}&\\
\le&\frac{1}{2}\|(\cof\grad u_\ve-\cof\grad u)\cd\grad u_\ve\|_{L^p}+\frac{1}{2}\|\cof\grad u\cd(\grad u_\ve-\grad u)\|_{L^p}&\\
\le&\frac{1}{2}\|\cof\grad u_\ve-\cof\grad u\|_{L^{2p}}\|\grad u_\ve\|_{L^{2p}}+\frac{1}{2}\|\cof\grad u\|_{L^{2p}}\|\grad u_\ve-\grad u\|_{L^{2p}}&
\end{flalign*}
Since $\|\grad u_\ve-\grad u\|_{L^{2p}}\ra0$ for all $1\le p<\infty$ the determinants converge strongly in $L^p.$
Additionally, due to Lipschitz continuity of each $\rho^{(k)}(\cd)$ we have $\rho^{(k)}(d_{\grad u_\ve})\ra\rho^{(k)}(d_{\grad u})$ strongly in $L^p$ for all $k\in\N$. \\

\textbf{Step 2: Testing the Equation.} For each $k=1,2$ we are testing the system by
\[\vp_k=\eta^2\p_k u_\ve.\] 
Note, $u_\ve$ is smooth with compact support, hence $\vp_k\in C_c^\infty$ can be used as a test function in \eqref{quasiPDE1} yielding
\begin{flalign*}
\int\limits_{B_{2r}}{\eta^2(x)\grad\p_k u\cd \grad\p_k u_\ve\;dx}=&-\int\limits_{B_{2r}}{\eta^2(x)\p_k(\rho'(d_{\grad u})\cof\grad u)\cd\grad\p_ku_\ve\;dx}&\\
&-2\int\limits_{B_{2r}}{\eta(x)\grad_\xi^2W(\grad u)\grad \p_ku\cd(\p_k u_\ve\ot\grad\eta)\;dx}&\\
=:&(I)_\ve+(II)_\ve.&
\end{flalign*}

\textbf{Step 3: Convergence of $(I)_\ve$:}\\

\textbf{Weak global bound on $(I)_\ve$:}
Notice, again, that $\p_k(\rho'(d_{\grad u})\cof\grad u)\in L^q$ for $1\le q<2.$ Hence, for every $\ve>0$ is in $\p_k(\rho'(d_{\grad u})\cof\grad u)\cd\grad\p_ku_\ve\in L^q$ for $1\le q<2,$ since $\grad\p_ku_\ve\in L^\infty.$\\

As a first step we give a crude estimate to $(I)_\ve$ for small enough $\ve>0:$
Applying Lemma \ref{lem:2.HR.2}.(ii) yields 
\begin{flalign}
(I)_\ve=&-\int\limits_{B_{2r}}{\eta^2\p_k(\rho'(d_{\grad u})\cof\grad u)\cd\grad\p_ku_\ve\;dx}&\nonumber\\
\le&\ga\int\limits_{B_{2r}}{\eta^2\max\left\{
|\grad\p_s u|^2,\frac{1}{2}|\grad\p_s u|^2+\frac{1}{2}|\grad\p_s u_\ve|^2\right\}\;dx}&\nonumber\\
\le&\ga\int\limits_{B_{2r}\cap\Om_1}{\eta^2|\grad\p_s u|^2\;dx}+\frac{\ga}{2}\int\limits_{B_{2r}\cap\Om_2}{\eta^2(|\grad\p_s u|^2+|\grad\p_s u_\ve|^2)\;dx}&\nonumber\\
\le&\ga\int\limits_{B_{2r}}{|\grad\p_k u|^2\;dx},&
\label{eq:Ch1.4.101}
\end{flalign}
for small enough $\ve>0$ and where $\Om_1,\Om_2\ss\Om$ are disjoint sets s.t. $|\grad\p_s u|^2\ge|\grad\p_s u_\ve|^2$ a.e. in $\Om_1$ and $|\grad\p_s u|^2<|\grad\p_s u_\ve|^2$ a.e. in $\Om_2.$ \\
In the last step we used the following property of the mollification:\\
Let $V\ss\ss W\ss\ss U$ and $V,W,U$ open and bounded sets. Assume $f\in L_{loc}^p(U)$ and $1\le p<\infty.$ If $\ve>0$ small enough then for $\|f_\ve\|_{L^p(V)}\le\|f\|_{L^p(W)},$ see again \cite[Appendix C.4, p.631]{LE10}.\\
Indeed, we can set $U=\Om,$ $W=B_{2r}$ and we can find an open set $V$ s.t. $\supp\eta\ss\ss V\ss\ss B_{2r}.$ Moreover, choosing $p=2,$ $f:=\grad\p_k u$ yields the estimate.\\

Note that the above inequality is good enough in the case that $0<\ga<1,$ see Remark \ref{re:2.3.1} below. However, it is too crude if $\ga\ge1.$ Hence, we need to give a more refined argument. The strategy will be as follows: Instead of the crude estimate given in \eqref{eq:Ch1.4.101} we would like to consider the $\limsup\limits_{\ve\ra0}(I)_\ve$ and then get a more refined bound on this limit.\\

\textbf{Existence of $\limsup\limits_{\ve\ra0}(I)_\ve$:}
To guarantee the existence of $\limsup\limits_{\ve\ra0}(I)_\ve$ we want to apply a version of the Reverse Fatou's Lemma. The standard Fatou's Lemma requires a pointwise majorant on the considered sequence. It is a subtle point, that $\ga |\grad\p_k u|^2$ turns out to be a global bound as shown in \eqref{eq:Ch1.4.101} but not a pointwise majorant for our sequence. Luckily there is a version of Fatou's  which only requires an pointwise individual integrable majorant for every member, as long as this sequence of majorants converges themself.
We are exactly in such a situation, to see this by Lemma \ref{lem:2.HR.2}.(ii) we have
\begin{flalign}
-\eta^2\left(\p_k(\rho'(d_{\grad u})\cof\grad u)\cd\grad\p_su_\ve\right)(x)\le&\ga\max\left\{
|\grad\p_k u(x)|^2,\frac{1}{2}|\grad\p_k u(x)|^2+\frac{1}{2}|\grad\p_k u_\ve(x)|^2\right\}&
\label{eq:Ch1.4.1011}
\end{flalign}
for a.e.\! $x\in B_{2r}$ and for every $\ve>0.$
Define the sequence 
\begin{align*}g_\ve:=\ga\max\left\{
|\grad\p_k u|^2,\frac{1}{2}|\grad\p_k u|^2+\frac{1}{2}|\grad\p_k u_\ve|^2\right\}\in L^1(B_{2r},[0,\infty])\end{align*}
for every $\ve>0.$
Obviously, $g_\ve\ra \ga |\grad\p_k u|^2=:g$ converges strongly in $L^1(B_{2r},[0,\infty])$  and $g\in L^1(B_{2r},[0,\infty]).$ \\
Furthermore, define the sequence of integrable functions
\begin{align*}f_\ve:=-\eta^2\left(\p_k(\rho'(d_{\grad u})\cof\grad u)\cd\grad\p_ku_\ve\right)\in  L^1(B_{2r},[0,\infty])\end{align*} for every $\ve>0.$  Moreover, note  \begin{align*}f_\ve\ra f:=-\eta^2\left(\p_k(\rho'(d_{\grad u})\cof\grad u)\cd\grad\p_ku\right)\;\mb{pw. a.e. on}\;B_{2r}. \end{align*}
This convergence follows from the fact that $\grad\p_ku_\ve\ra\grad\p_ku$ converges pw. a.e. Note, that while the $f_\ve$'s are integrable, for the limit $f$ we can only guarantee that $f$ is measurable and $f^+$ integrable. Indeed, Lemma \ref{lem:2.HR.2}.(ii) guarantees an upper bound on $\int\limits_{B_{2r}}{f^+\;dx}$ given by 
\begin{align*}\int\limits_{B_{2r}}{f\;dx}\le\int\limits_{B_{2r}}{f^+\;dx}\le \ga\int\limits_{B_{2r}}{|\grad\p_k u|^2\;dx}<\infty.
\end{align*}
However, we do not have a lower bound on $\int\limits_{B_{2r}}{f\;dx}$ so it might be $-\infty.$ Luckily, our version of the Reverse Fatou's Lemma is such that it can still be applied.\\

Now after all those preparations we are finally in the position to apply our version of the Reverse Fatou's Lemma. By applying Lemma \ref{lem:2.3.RF1} we get
\begin{align*}
\limsup\limits_{\ve\ra0}(I)_\ve=&\limsup\limits_{\ve\ra0}\int\limits_{B_{2r}}{f_\ve\;dx}&\\
\le&\int\limits_{B_{2r}}{f\;dx}&\\
=:&(I).&
\end{align*}

Now we want to bound $(I)$ from above. The first case we have to take care of is that $(I)$ might be $-\infty.$ In this case we can just estimate $(I)$ by $0.$\\

So assume from now on, that $(I)$ is finite. Making use of the product rule yields
\begin{align*}
(I)=&-\int\limits_{B_{2r}}{\eta^2\p_k(\rho'(d_{\grad u})\cof\grad u)\cd\grad\p_ku\;dx}&\\
=&-\int\limits_{B_{2r}}{\eta^2\rho'(d_{\grad u})\cof\grad u,_k\cd\grad u,_k\;dx}-\int\limits_{B_{2r}}{\eta^2\rho''(d_{\grad u})(\cof\grad u\cd\grad\p_ku)^2\;dx}.&
\end{align*}
Realising that the rightmost term is always non-positive implies
\begin{align*}
(I)\le&-\int\limits_{B_{2r}}{\eta^2\rho'(d_{\grad u})\cof\grad u,_k\cd\grad u,_k\;dx}.&
\end{align*}

So far we have seen, that the $\limsup\limits_{\ve\ra0}(I)_\ve$ is either $-\infty,$ which is harmless as argued, or it is finite and the explicit upper bound is given by 
\[
\limsup\limits_{\ve\ra0}(I)_\ve\le (I)\le(I.a):=-\int\limits_{B_{2r}}{\eta^2\rho'(d_{\grad u})\cof\grad u,_k\cd\grad u,_k\;dx}.\]
\vspace{0.5cm}

\textbf{Step 4: Bounding $(I.a)$ via compensated compactness:}\\

\textbf{Approximating $(I.a):$} For reasons, becoming obvious in a second, we introduce, again, an approximation:
\begin{align}
(I.a)=&-\int\limits_{B_{2r}}{\eta^2\rho'(d_{\grad u})\cof\grad u_{\ve,k}\cd\grad u_{\ve,k}\;dx}&\nonumber\\
&+\int\limits_{B_{2r}}{\eta^2\rho'(d_{\grad u})\cof(\grad u_{\ve,k}-\grad u_{,k})\cd\grad u_{\ve,k}\;dx}&\nonumber\\
&+\int\limits_{B_{2r}}{\eta^2\rho'(d_{\grad u})\cof\grad u_{,k}\cd(\grad u_{\ve,k}-\grad u_{,k})\;dx}&\nonumber\\
=:&(I.a.a)_\ve+(I.a.b)_\ve+(I.a.c)_\ve.&
\end{align}
The terms $(I.a.b)_\ve$ and $(I.a.c)_\ve$ vanish when $\ve\ra0.$ Indeed by Hölder's we get
\begin{align*}|(I.a.c)_\ve|\le&\ga\int\limits_{B_{2r}}{|\cof\grad u_{,k}\cd(\grad u_{\ve,k}-\grad u_{,k})|\;dx}\\
\le&\ga\|\grad u_{,k}\|_{L^2}\|\grad u_{\ve,k}-\grad u_{,k}\|_{L^2}\ra0
\end{align*}
when $\ve\ra0,$ since $\grad^2u\in L^2$ and $\grad^2u_\ve\ra\grad^2u$ strongly in $L^2.$ One can argue similarly for $(I.a.b)_\ve.$\\

\textbf{Strategy for estimating $(I.a.a)_\ve$:} Recall now that the goal is to establish a reverse Hölder inequality. For this sake, we need to estimate $(I.a.a)_\ve$ in terms of the $L^p-$norms of the 2nd derivative for some $p<2$.\\
The integrand of $(I.a.a)_\ve$ contains the product $\cof\grad u_{\ve,k}\cd\grad u_{\ve,k}.$ Obviously, these quantities satisfy $\div\cof\grad u_{\ve,k}=0$ and $\mb{curl} \grad u_{\ve,k}=0$ in a distributional sense. This prompts the idea of making use of compensated compactness results, in particular, the famous div-curl lemma. Some difficulties arise here. First, both factors contain second derivatives of $u,$ which we would like to control in an $L^p$ norm with $p<2,$ this is done below. For this sake, the div-curl Lemma \ref{lem:2.HR.1} needs to be applied with $s=t=p<2,$ which is possible although, the quantity $\cof\grad u_{\ve,k}\cd\grad u_{\ve,k}$ is only controlled in the weak space $\mathcal{H}^\al$ with $\al=\frac{p}{2}<1.$ Secondly, we need the dual space to such a space and the duality inequality, to bring the Hardy space into play in the first place.\\

\textbf{Extending the domain to $\R^2$:} For any map $v:B_{2r}\ra\R^2$ we will use $\ol{v}$ to denote the trivial extension, i.e. $\ol{v}:=v\chi_{B_{2r}},$ where $\chi$ is the characteristic function. Define $f:=\cof\grad u_{\ve,k}$ and  $g:=-\grad u_{\ve,k}.$ Then we can rewrite $(I.a.a)_\ve$ by
\begin{align}
(I.a.a)_\ve=&\int\limits_{\R^2}{\ol{\eta^2\rho'(d_{\grad {u}})}(\ol{f}\cd \ol{g})\;dx}.&\label{eq:Ch1.4.10221}
\end{align}
Note, that the extensions $\ol{f}$ and $\ol{g}$ are indeed smooth on the whole space. This is true because  $\grad u_{\ve,k}$ lives on a compact subset of $B_{2r}$ and disappears smoothly on the boundary of that compact set in such a way that it is smooth on the full set $\ol{B_{2r}}$ and remains zero up to the boundary of  $\ol{B_{2r}}.$ Moreover, they satisfy the cancelation conditions in $\ol{B_{2r}},$ hence, the extensions $\ol{f}$ and $\ol{g}$ are indeed smooth and satisfy the cancelation conditions $\div\ol{f}=0$ and $\mb{curl}\ol{f}=0$ on the full space.\\

As a next step we want to apply a Fefferman-Stein Type Duality Inequality, however, for $\mathcal{H}^\al$ when $0<\al<1$ then the corresponding dual spaces 
are given by $\dot{\Lambda}_\beta/\{constants\}(\R^2),$ with $\beta=2(\frac{1}{\al}-1).$ We will need $H:=\ol{\eta^2\rho'(d_{\grad {u}})}\in \dot{\Lambda}_\beta(\R^2).$ 
As a first step $H$ must be continuous. For this first note, that $x\mt\rho'(d_{\grad {u(x)}})$ is continuous because of the assumption that  $x\mt d_{\grad u(x)}$ is Hölder-continuous.
Note, however, that $u$ might be non zero on the boundary of $B_{2r},$ hence $H$ might jump on the boundary. Obviously, $H$ is continuous in $\supp\eta$ and $H=0$ on $\R^2\sm \ol{B_{2r}}.$
However, $\eta$ disappears on $\ol{B_{2r}}\sm \supp\eta$ guaranteeing the continuity of $H.$ That $H$ is actually in $\dot{\Lambda}_\beta(\R^2)$ is shown below.\\

\textbf{Compensated compactness:} 
Applying now the duality inequality \eqref{eq:Ch1.4.103} to the RHS of \eqref{eq:Ch1.4.10221}, with a for now free $0<\al<1$ and the corresponding $\beta=2(\frac{1}{\al}-1),$ gives
\begin{align*}
(I.a.a)_\ve \le&C(\al) \|H\|_{\dot{\Lambda}_\beta(\R^2)}\|\ol{f}\cd \ol{g}\|_{\mathcal{H}^\al(\R^2)}.&
\end{align*}
Making use of the div-curl Lemma \ref{lem:2.HR.1} with $\al=\frac{p}{2}\in(\frac{2}{3},1)$ and $s=t=p,$ and using \eqref{eq:2.4.HS.equiv} yields
\[\|\ol{f}\cd \ol{g}\|_{\mathcal{H}^\al(\R^2)}\le  C(p)\|\ol{f}\|_{L^p(\R^2)}\|\ol{g}\|_{L^p(\R^2)}\le C(p)\|\grad u_{\ve,k}\|_{L^p(B_{2r})}^2.\]

As a last step we need to show, that $H\in\dot{\Lambda}_\beta(\R^2)$ with $\beta=2(\frac{2}{p}-1)$ and bound the norm by a constant.\\

\textbf{Showing $H\in\dot{\Lambda}_\beta(\R^2)$ and an estimate on the norm}: 
If $x,y\in \R^2\sm\ol{B_{2r}}$ then $|H(x)-H(y)|=0.$ The case when $x\in \R^2\sm\ol{B_{2r}}$ and $y\in \ol{B_{2r}}$ or vice versa can be reduced to the case below. Indeed, letting $z$ be the point where the straight line connecting $x$ and $y$ hits $\partial B_{2r},$ then by the triangle inequality and $H(y)=H(z)=0$ we get 
\begin{align*}|H(x)-H(y)|\le|H(x)-H(z)|+|H(z)-H(y)|=|H(x)-H(z)|.\end{align*}
By the result below and the fact that $|x-z|\le|x-y|$ by choice we finally get
\begin{align*}|H(x)-H(y)|\le C|x-z|^\be\le C|x-y|^\be.\end{align*}

Let now $x,y\in \ol{B_{2r}}.$ Then we have
\begin{align}|H(x)-H(y)|\le \ga|\eta^2(x)-\eta^2(y)|+|\rho'(d_{\grad u})(x)-\rho'(d_{\grad u})(y)|.\label{eq:Ch1.4.10221.a}\end{align}
By $\rho\in C^\infty$ and the assumption $d_{\grad u}\in C_{loc}^{0,\be}$ for the rightmost term we get 
\begin{align*}|\rho'(d_{\grad u})(x)-\rho'(d_{\grad u})(y)|\le&\|\rho''\|_{L^\infty}|d_{\grad u}(x)-d_{\grad u}(y)|\le C \|\rho''\|_{L^\infty}|x-y|^\be.&\end{align*}

In order to establish an estimate on the first term of the RHS of \eqref{eq:Ch1.4.10221.a} we follow the strategy of Morrey's proof of the Dirichlet Growth Theorem.\\
Let $h:=\eta^2.$ Then by the triangle inequality for an arbitrary $w\in \R^2$ it holds
\[|h(x)-h(y)|\le|h(x)-h(w)|+|h(w)-h(y)|.\]
Denote the midpoint of $x$ and $y$ by $z,$ i.e. $z=\frac{x+y}{2}$ and the distance between $x$ and $z$ by $l:=\frac{|x-y|}{2}.$ Now we average the previous inequality over $B(z,l)$, i.e. \[|h(x)-h(y)|\le\fint\limits_{B(z,l)}{|h(x)-h(w)|\;dw}+\fint\limits_{B(z,l)}{|h(w)-h(y)|\;dw}.\]
It is enough to control one of them the other one is similar. For this sake, let $w_t:=w+t(x-w).$ By the Mean Value Theorem we see
\begin{flalign*}
\fint\limits_{B(z,l)}{|h(x)-h(w)|\;dw}=&\fint\limits_{B(z,l)}{\left|\int\limits_0^1\grad h(w_t)dt\cd(x-w)\right|\;dw}&\\
\le&\fint\limits_{B(z,l)}\int\limits_0^1{|\grad h(w_t)|\;dt\;dw}|x-y|&\\
=&\frac{|x-y|}{|B(z,l)|}\int\limits_0^1\int\limits_{B(\bar{z},lt)}{|\grad h(v)|\;dv\;t^{-2}dt},&
\end{flalign*}
where we exchanged integrals and substituted $v=y+t(w-y)$ and $\bar{z}=y+t(z-y).$
Expressing the function explicitly, yields
\begin{align}
\int\limits_{B(\bar{z},lt)}{|\grad h(v)|\;dv}=&\int\limits_{B(\bar{z},lt)}{|2\eta(v)\p_k\eta(v)|\;dv}&\label{eq:CH1.4.105}
\end{align}
Applying Hölder's inequality for any $1<q<\infty$ we get 
\begin{align*}
\int\limits_{B(\bar{z},lt)}{|2\eta(v)\p_k\eta(v)|\;dv}\le &\frac{C}{r}\int\limits_{B(\bar{z},lt)\cap B(x_0,2r)}{\;dv}&\\
\le& \frac{C}{r}|B(\bar{z},lt)|^{\frac{1}{q}}|B(x_0,2r)|^{\frac{(q-1)}{q}}&\\
=&Cr^{1-\frac{2}{q}}(lt)^{\frac{2}{q}},&
\end{align*}
where we used $\eta\le 1,$ $|\grad\eta|\le \frac{c}{r},$ $|\rho'|\le\ga$ and $|B(\bar{z},lt)|=\pi(lt)^2.$
Finally,
\begin{flalign*}
|h(x)-h(y)|\le&\frac{|x-y|}{|B(z,l)|}Cr^{1-\frac{2}{q}}l^{\frac{2}{q}}\int\limits_0^1{t^{\frac{2}{q}-2}\;dt}&\\
\le&C(q)r^{1-\frac{2}{q}}|x-y|^{{\frac{2}{q}}-1}.
\end{flalign*}
 It is important to match $\beta$ exactly, since on the full space $\R^2$ there is no obvious relation between (homogenous) Lipschitz/Hölder-spaces with different exponents. Consequently, ${{\frac{2}{q}}-1}=\beta=2(\frac{2}{p}-1)$ needs to be satisfied. This choice yields $\frac{2}{q}=1+\beta>1$ guaranteeing that the integrals above can always be performed.\\

Hence, we showed that $H\in\dot{\Lambda}_{\beta}$ and the norm is bounded via
\[\|H\|_{\dot{\Lambda}_{\beta}}\le C(\ga,p)r^{-\beta}.\]

\textbf{Conclusion of Step 3+4:}
Finally, combining the above we obtain
\begin{align}
(I.a.a)_\ve\le C(\ga,p)r^{-\beta}\|\grad u_{\ve,k}\|_{L^p(B_{2r})}^2,
\label{eq:Ch1.4.1023}
\end{align}
yielding
\[
-\infty\le\limsup\limits_{\ve\ra0}(I)_\ve\le(I.a)\le\limsup\limits_{\ve\ra0}(I.a.a)_\ve\le C(\ga,p)r^{-\beta}\|\grad u_{,k}\|_{L^p(B_{2r})}^2.
\]
\vspace{0.5cm}

\textbf{Step 5: Upper bound on $(II)_\ve$:}
Splitting up the integrals is fine since all the integrands are well behaved. Hence we get
\begin{flalign}
(II)_\ve=&-2\int\limits_{B_{2r}}{\eta(x)\grad_\xi^2W(\grad u)\grad \p_ku\cd(\p_k u_\ve\ot\grad\eta)\;dx}&\nonumber\\
=&-2\int\limits_{B_{2r}}{\eta(x)\grad \p_k u\cd(\p_k u_\ve\ot\grad\eta)\;dx}&\nonumber\\
&-2\int\limits_{B_{2r}}{\eta(x)\rho'(d_{\grad u})(\p_{\xi_{ab}}(\cof\xi)_{ij})\p_{k}\p_au_b\p_k (u_\ve)_i\p_j\eta\;dx}&\nonumber\\
&-2\int\limits_{B_{2r}}{\eta(x)\rho''(d)(\cof\grad u)_{ij}(\cof\grad u)_{ab}\p_{x_k}\p_au_b\p_k (u_\ve)_i\p_j\eta\;dx}&\nonumber\\
=:&(II.a)_\ve+(II.b)_\ve+(II.c)_\ve.&\label{eq:Ch1.4.99}
\end{flalign}

\textbf{Estimating $(II.a)_\ve+(II.b)_\ve$:} the first two terms $(II.a)_\ve$ and $(II.b)_\ve$ can be controlled by Hölder's inequality (for any $1<p<2$ and $p^{-1}+p'^{-1}=1$) and recalling the properties of $\eta$ to obtain
\begin{flalign}
(II.a)_\ve+(II.b)_\ve\le&C(\ga)\int\limits_{B_{2r}}{\eta|\grad\p_ku||\p_k u_\ve||\grad\eta|\;dx}&\nonumber\\
\le&C(\ga)\left(\int\limits_{B_{2r}}{|\eta|^p|\grad\p_ku|^p\;dx}\right)^{\frac{1}{p}}\left(\int\limits_{B_{2r}}{|\p_k u_\ve|^{p'}|\grad\eta|^{p'}\;dx}\right)^{\frac{1}{p'}}&\nonumber\\
\le&C(\ga)r^{-1}\left(\int\limits_{B_{2r}}{|\grad\p_ku|^p\;dx}\right)^{\frac{1}{p}}\left(\int\limits_{B_{2r}\sm B_r}{|\p_k u_\ve|^{p'}\;dx}\right)^{\frac{1}{p'}}.&\label{eq:Ch1.4.991}
\end{flalign}
Now the idea is to give an estimate of the form 
 \begin{align*}
\|\p_k u_\ve\|_{L^{p'}(B_{2r}\sm B_r)}\lesssim \|\grad\p_k u_\ve\|_{L^p(B_{2r})}.
 \end{align*}

For this sake, recall the Sobolev-Poincare inequality for balls\footnote{See, \cite[Theorem 1, §1.4.4.]{M11}.}: there is a constant $C(s^*)>0$ s.t. for any $v\in W_0^{1,s^*}$ it holds
 \begin{align}\|v\|_{L^{s}(B(x,R))}\le C(s^*)R^{1+\frac{2}{s}-\frac{2}{s^*}}\|\grad v\|_{L^{s^*}(B(x,R))},\label{eq:Ch1.4.992} \end{align}
where $s^*=\frac{2s}{2+s}.$ By applying the latter inequality with $R=2r,$ $s=p'$ and $s^*=(p')^*=\frac{2p'}{2+p'}$ we obtain 
 \begin{align*}
\|\p_k u_\ve\|_{L^{p'}(B_{2r})}\le& C(p)r^{1+\frac{2}{p'}-\frac{2}{(p')^*}}\|\grad\p_k u_\ve\|_{L^{(p')^*}(B_{2r})}.&
\end{align*}
Making use of Hölder's to go from $L^{(p')^*}$ to $L^{p}$ we get
 \begin{align*}
\|\p_k u_\ve\|_{L^{p'}(B_{2r})}\le C(p)r^{1+\frac{2}{p'}-\frac{2}{(p')^*}+\frac{2}{(p')^*}-\frac{2}{p}}\|\grad\p_k u_\ve\|_{L^p(B_{2r})}, \end{align*}
keeping in mind, that this can only be true, if the condition $(p')^*\le p$ is satisfied. Realise that the exponent in the latter inequality simplifies to 
 \begin{align*}\omega_1:=1+\frac{2}{p'}-\frac{2}{p}. \end{align*}
Hence, we have seen that
 \begin{align*}
\|\p_k u_\ve\|_{L^{p'}(B_{2r}\sm B_r)}\le\|\p_k u_\ve\|_{L^{p'}(B_{2r})}\le C(p)r^{\omega_1}\|\grad\p_k u_\ve\|_{L^p(B_{2r})}.
 \end{align*}
By combining the latter inequality with \eqref{eq:Ch1.4.991}  we obtain
 \begin{align}(II.a)_\ve+(II.b)_\ve\le C(\ga,p)r^{(\omega_1-1)}\|\grad\p_k u\|_{L^p(B_{2r})}\|\grad\p_k u_\ve\|_{L^p(B_{2r})}.\label{eq:Ch1.4.993} \end{align}
\vspace{0.25cm}

\textbf{Estimating $(II.c)_\ve$:} For the rightmost term in \eqref{eq:Ch1.4.99} we argue similarly as before, again by Hölder's inequality and the properties of $\eta$ one obtains
\begin{flalign*}
(II.c)_\ve=&-2\int\limits_{B_{2r}}{\eta(x)\rho''(d)(\cof\grad u)_{ij}(\cof\grad u)_{ab}\p_{k}\p_au_b\p_k (u_\ve)_i\p_j\eta\;dx}&\\
\le&C(\|\rho''\|_{L^\infty})\int\limits_{B_{2r}}{\eta(x)|\p_k u_\ve||\grad\p_ku||\grad u|^2|\grad\eta|\;dx}&\\
\le& C(\|\rho''\|_{L^\infty})r^{-1}\left(\int\limits_{B_{2r}}{|\grad \p_ku|^p\;dx}\right)^{\frac{1}{p}}\left(\int\limits_{B_{2r}\sm B_r}{|\p_k u_\ve|^{{p'}}|\grad u|^{{2p'}}\;dx}\right)^{\frac{1}{{p'}}}.&
\end{flalign*}
By applying Hölder's inequality with $\frac{2}{3}+\frac{1}{3}=1$ to the rightmost term of the latter line and by recalling the property $\|\p_k u_\ve\|_{L^{q}(B_{2r})}\le\|\p_k u\|_{L^{q}(B_{2r})}$ for any $1\le q<\infty$ and for small enough $\ve>0$ we get 
\begin{flalign*}
\left(\int\limits_{B_{2r}\sm B_r}{|\p_k u_\ve|^{{p'}}|\grad u|^{{2p'}}\;dx}\right)^{\frac{1}{{p'}}}\le& \left(\int\limits_{B_{2r}\sm B_r}{|\grad u|^{{3p'}}\;dx}\right)^{\frac{2}{{3p'}}}\left(\int\limits_{B_{2r}\sm B_r}{|\p_k u_\ve|^{{3p'}}\;dx}\right)^{\frac{1}{{3p'}}}&\\
\le&\left(\int\limits_{B_{2r}\sm B_r}{|\grad u|^{{3p'}}\;dx}\right)^{\frac{1}{{p'}}},&
\end{flalign*}
if $\ve>0$ is sufficiently small.\\

We can proceed like before. Applying again \eqref{eq:Ch1.4.992} yields 
\begin{align*}
\|\grad u\|_{L^{3p'}(B_{2r})}\le& C(p)r^{1+\frac{2}{3p'}-\frac{2}{(3p')^*}}\|\grad^2 u\|_{L^{(3p')^*}(B_{2r})}.&
\end{align*}
Assuming $(3p')^*\le p$ and applying Hölder's to go from $L^{(3p')^*}$ to $L^{p}$ we get
 \begin{align*}
\|\grad u\|_{L^{3p'}(B_{2r})}\le C(p)r^{1+\frac{2}{3p'}-\frac{2}{(3p')^*}+\frac{2}{(3p')^*}-\frac{2}{p}}\|\grad^2 u\|_{L^p(B_{2r})}. \end{align*}
 The exponent in the latter inequality simplifies to 
 \begin{align*}\omega_2:=1+\frac{2}{3p'}-\frac{2}{p}. \end{align*}
This shows
 \begin{align*}
\|\grad u\|_{L^{3p'}(B_{2r}\sm B_r)}^3\le\|\grad u\|_{L^{3p'}(B_{2r})}^3\le C(p)r^{3\omega_2}\|\grad^2 u\|_{L^p(B_{2r})}^3.
\end{align*}
Then
 \begin{equation}
 (II.c)_\ve\le C(\|\rho''\|_{L^\infty},p)r^{\left(3\omega_2-1\right)}\|\grad\p_k u\|_{L^p(B_{2r})}\|\grad^2 u\|_{L^p(B_{2r})}^3. 
 \label{eq:Ch1.4.994} 
 \end{equation}
\vspace{0.25cm}

\textbf{Estimating $(II)_\ve$:} 
By combining \eqref{eq:Ch1.4.993} and \eqref{eq:Ch1.4.994} we have
\[(II)_\ve\le C(\ga,\|\rho''\|_{L^\infty},p)\|\grad\p_k u\|_{L^p(B_{2r})}\left[r^{(\omega_1-1)}\|\grad\p_k u_\ve\|_{L^p(B_{2r})}+r^{(3\omega_2-1)}\|\grad^2 u\|_{L^p(B_{2r})}^3\right],\]
for any $1<p<2$ s.t. $(p')^*\le p$ and $(3p')^*\le p,$ where the latter condition is stricter and reduces the range to $p\in[\frac{8}{5},2).$ 
\vspace{0.5cm}

\textbf{Step 6: Conclusion:}\\

Taking the limsup $\ve\ra0$ on both sides, and collecting all expressions from above yields

\[
\int\limits_{B_{r}}{|\grad^2 u|^2\;dx}\le C(\ga,\|\rho''\|_{L^\infty},p)[(r^{-\be}+ r^{(\omega_1-1)})\|\grad^2u\|_{L^p}^2+r^{(3\omega_2-1)} \|\grad^2 u\|_{L^{p}}^4]
\]
Recalling $\beta=2(\frac{2}{p}-1),$  $\omega_1-1=-\be,$ and dividing the latter by $|B_{2r}|$ yields
\begin{equation}
\fint\limits_{B_{r}}{|\grad^2 u|^2\;dx}\le C(\ga,\|\rho''\|_{L^\infty},p)(1+r^{(3\omega_2-1+\beta)} \|\grad^2 u\|_{L^{p}(B_{2r})}^2)\left(\fint\limits_{B_{2r}}{|\grad^2 u|^p\;dx}\right)^{\frac{2}{p}}.
\label{eq:Ch1.4.99.99}
\end{equation}
Now, notice firstly $3\omega_2-1+\be=\frac{2}{p'}-\frac{2}{p}$ and secondly by Hölder's inequality
 \begin{align*}
r^{\left(\frac{2}{p'}-\frac{2}{p}\right)}\|\grad^2 u\|_{L^{p}(B_{2r})}^2\le C(p)r^{\left(\frac{2}{p'}-\frac{2}{p}+\be\right)}\|\grad^2 u\|_{L^2(B_{2r})}^2.
\end{align*}
Then the exponent $\frac{2}{p'}-\frac{2}{p}+\be=\frac{2}{p'}+\frac{2}{p}-2=0$ vanishes for any $p\in[\frac{8}{5},2).$ Furthermore, by Lebesgue's Absolute Continuity Theorem we can choose $r\le1$ so small s.t. $\|\grad^2 u\|_{L^2(B_{2r})}^2\le1$ and \eqref{eq:Ch1.4.99.99} becomes\footnote{A similar reasoning can be found in \cite[p.142]{MG83}.}
\begin{equation}
\fint\limits_{B_{r}}{|\grad^2 u|^2\;dx}\le C(\ga,\|\rho''\|_{L^\infty},p)\left(\fint\limits_{B_{2r}}{|\grad^2 u|^p\;dx}\right)^{q},
\label{eq:Ch1.4.100}
\end{equation}

where $q=\frac{2}{p}.$ By choosing $f:=|\grad^2 u|^p,$ Gehring's lemma, see \cite[Ch.5, Prop.1.1]{MG83} guarantees the existence of some $\d=\d(\ga,\|\rho''\|_{L^\infty},p)>0$ and $\tilde{C}=\tilde{C}(\ga,\|\rho''\|_{L^\infty},p)>0$ s.t.  
\[
\left(\fint\limits_{B_{r}}{|\grad^2 u|^{2+\d}\;dx}\right)^{\frac{1}{2+\d}}\le \tilde{C}\left(\fint\limits_{B_{2r}}{|\grad^2 u|^2\;dx}\right)^{\frac{p}{2}}
\]
and in particular $u\in W_{loc}^{2,2+\d}(\Om,\R^2)$ and by Sobolev embedding 
$u\in C_{loc}^{1,\al},$ for some $\al>0.$ \vspace{0.5cm}

\begin{re}
Note, that for the uniformly convex case $(0<\ga<1)$ the argument can be simplified and the assumption $x\mt d_{\grad u(x)}\in C_{loc}^{0,\be}(\Om)$ can be dropped, in particular, the difficult limit taking and compensated compactness arguments can be avoided. The estimate found in \eqref{eq:Ch1.4.101} is good enough if $(0<\ga<1).$ Indeed, combined with the bound on $(II)_\ve$ from above one gets the slightly weaker estimate
\[
\fint\limits_{B_{r}}{|\grad^2 u|^2\;dx}\le C(\ga,\|\rho''\|_{L^\infty},p)\left(\fint\limits_{B_{2r}}{|\grad^2 u|^p\;dx}\right)^{q}+\ga\fint\limits_{B_{2r}}{|\grad^2 u|^2\;dx},
\]
replacing \eqref{eq:Ch1.4.100}.
 However, Gehring's lemma still applies and one can conclude as above.
\label{re:2.3.1}
\end{re}\vspace{0.5cm}

It is well known fact that, sufficiently regular (for instance $C^{1,\al}$) solutions to some elliptic systems (with nice coefficients), are already smooth. The following corollary shows that this is indeed true for our system.

\begin{cor}[Maximal Smoothness] Let $0<\ga<\infty,$ $\Om\ss\R^2$ be open and bounded, assume $u\in\A_{u_0}\cap W_{loc}^{2,2}(\Om,\R^2)$ to be a stationary point of the functional, as given in \eqref{eq:1.1}.\\

Then $u$ is smooth.
\end{cor}
\textbf{Proof:}\\
A proof is standard and we will only outline the strategy of elliptic regularity theory in $2\times2$ dimensions:\\
All previous results apply and $u$ is of class $\A_{u_0}\cap W_{loc}^{2,2}(\Om,\R^2)\cap C_{loc}^{1,\al}(\Om,\R^2).$ Moreover, we know that $\grad u\in W_{loc}^{1,2}(\Om,\R^4)$ solves the `linearized' elliptic PDE 
\[
-\div(A(x)\grad^2 u)=0\;\mb{in}\;D'(\Om,\R^4)
\]
 where $A(\cd):=\grad_\xi^2W(\grad u(\cd))$ is the coefficients matrix. Note, that $\xi\mapsto W(\xi)$ is smooth, hence the regularity of $A$ is determined by the regularity of $\grad u.$ Before it was not known whether or not $A$ is in $L^2,$ making it impossible to test with $W_0^{1,2}$ functions. Since, $\grad u(\cd)$ is now of class $C_{loc}^{0,\al},$ so is $A$ which makes upgrading to $W_0^{1,2}$ functions possible and hence, $\grad u\in W_{loc}^{1,2}(\Om,\R^4)$ is a weak solution to the following system\footnote{For more details on step 1 consult \cite[p.46]{MG83}.}
\begin{equation}
-\div(A(x)\grad^2 u)=0\;\mb{in}\;W_0^{1,2}(\Om,\R^4).
\label{eq:2.3.104}
\end{equation}

From this point Schauder theory takes over. Since, $\grad u\in W_{loc}^{1,2}(\Om,\R^4)$ solves \eqref{eq:2.3.104} weakly, $\grad u$ is of class $ C_{loc}^{2,\al}$, see \cite[Chapter 3, Theorem 3.2]{MG83}. In turn, this improves the regularity of the coefficients matrix $A$ to $C_{loc}^{1,\al}$, differentiating the PDE again, and applying \cite[Chapter 3, Theorem 3.3]{MG83}., obtaining $u\in C_{loc}^{3,\al}.$ Repeating this argument over and over improves the regularity of $u$ to class $C_{loc}^{\infty}.$\\\vspace{0.5cm}

Next we show the important lower bound on the quantity $\p_s(\rho'(d_{\grad v})\cof\grad v)\cd\grad \p_sv_\ve$ we used in the proof above. This follows fairly straightforward from the lower bound we showed for the `differenced' version in Lemma \ref{HR:LUBs}.

\begin{lem}\label{lem:2.HR.2} Let $\Om\ss\R^2$ be open and bounded and $v\in W^{2,2}(\Om).$\\

(i). Then for $s=1,2$ it holds
\begin{equation}\p_s(\rho'(d_{\grad v})\cof\grad v)\cd\grad \p_sv\ge-\ga|\grad\p_s v|^2\end{equation}
a.e. in $\Om.$\\
(ii). For every $s=1,2$ and $\ve>0$ small enough it holds 
\begin{equation}\p_s(\rho'(d_{\grad v})\cof\grad v)\cd\grad \p_sv_\ve\ge-\ga\max\{|\grad\p_s v|^2,\frac{1}{2}|\grad\p_s v|^2+\frac{1}{2}|\grad\p_s v_\ve|^2\}.\end{equation}
a.e. in $\Om.$
\end{lem}
\textbf{Proof:}\\
(i). By Nirenberg's lemma \cite[p.45]{MG83} we know that $\D^{h,s}\grad v\ra\grad\p_sv$ in $L^2(\Om,\R^4),$ in particular,
 $\D^{h,s}\grad v\ra\grad\p_sv$ a.e. in $\Om.$ By the standard/discrete product rule it is straightforward  to see that
\[\p_s(\rho'(d_{\grad v})\cof\grad v)=\lim\limits_{h\ra0}\D^{h,s}(\rho'(d_{\grad v})\cof\grad v)\]
a.e. in $\Om.$ By Lemma \ref{HR:LUBs} we get
\[\p_s(\rho'(d_{\grad v})\cof\grad v)=\lim\limits_{h\ra0}\D^{h,s}(\rho'(d_{\grad v})\cof\grad v)\ge-\ga\lim\limits_{h\ra0}|\D^{h,s}\grad v|^2=-\ga|\grad\p_s v|^2\]
a.e. in $\Om.$\\ 
(ii). There are two possible cases. Assume first that the inequality in (i) is strict, i.e.
\[\p_s(\rho'(d_{\grad v}(x))\cof\grad v(x))\cd\grad \p_sv(x)>-\ga|\grad\p_s v(x)|^2,\]
for some $x\in\Om.$
Since 
\[\p_s(\rho'(d_{\grad v}(x))\cof\grad v(x))\cd\grad \p_sv_\ve(x)\ra\p_s(\rho'(d_{\grad v}(x))\cof\grad v(x))\cd\grad \p_sv(x)\]
Then the same holds true if $\ve$ is small enough:\footnote{Indeed, we have a converging sequence $f_\ve\ra f$ if $\ve\ra 0$ with $f>g.$ Then $f_\ve>g$ for small enough $\ve.$ Assume not, then there exists a sequence $(\ve_n)_{n\in\N}>0$ with $\ve_n\ra 0$ if $n\ra\infty$ s.t. $f_{\ve_n}\le g$ this contradicts immediately the convergence $f_\ve\ra f$.}
\[\p_s(\rho'(d_{\grad v}(x))\cof\grad v(x))\cd\grad \p_sv_\ve(x)>-\ga|\grad\p_s v(x)|^2.\]

Assume instead that for some $x\in \Om$ equality holds:
\[\p_s(\rho'(d_{\grad v}(x))\cof\grad v(x))\cd\grad \p_sv(x)=-\ga|\grad\p_s v(x)|^2.\]
Then since, the above inequality is a Cauchy-Schwarz inequality, equality only holds if 
\[\p_s(\rho'(d_{\grad v}(x))\cof\grad v(x))=-\ga\grad\p_s v(x).\]
But then for every $\ve>0$ we get 
\begin{align*}\p_s(\rho'(d_{\grad v}(x))\cof\grad v(x))\cd\grad \p_sv_\ve(x)=&-\ga\grad\p_s v(x)\cd\grad \p_sv_\ve(x)&\\
\ge&-\frac{\ga}{2}(|\grad\p_s v(x)|^2+|\grad\p_s v_\ve(x)|^2).&
\end{align*}
Putting both cases together yields the claim.\vspace{0.5cm}

The next statement is a generalisation to the the Reverse Fatou's Lemma. Instead of a pointwise integrable upper limit dominating every member of the sequence, it will be enough to bound each member individually from above, as long as this sequence of upper limits converges by itself. Our proof is inspired by \cite[Thm 2.1]{FKZ16}.

\begin{lem} [A version of the Reverse Fatou's Lemma]
Let $(X,\Sigma,\mu)$ be a measure space, with $\mu(X)<\infty.$ Let $g_n, g\in L^1(X,[-\infty,\infty],\mu)$ for any $n\in \N$ be a sequence, s.t.
\begin{align*}g_n \ra g \;\mb{strongly in} \; L^1(X,[-\infty,\infty],\mu).\end{align*}
Suppose further that there is another sequence $(f_n)_{n\in \N}\ss L^1(X,[-\infty,\infty],\mu)$ s.t. 
\begin{align*}f_n\le g_n \;\mb{a.e. for every} \;n\in\N\end{align*}
and s.t. there exists a measurable functions $f:X\ra[-\infty,\infty]$ 
s.t. \begin{align*}f_n\ra f\mb{ptw. a.e.}\end{align*} 
and either (a) $\int\limits_{X}{f\;d\mu}=+\infty$ or (b) $f^+$ is integrable.  \\

Then it holds
\begin{align}\limsup\limits_{n\ra\infty}\int\limits_{X}{f_n\;d\mu}\le \int\limits_{X}{f\;d\mu}.\end{align}
\label{lem:2.3.RF1}
\end{lem}
\textbf{Proof:}\\
(a) $\int\limits_{X}{f\;d\mu}=+\infty.$ Trivial, since $\limsup\limits_{n\ra\infty}\int\limits_{X}{f_n\;d\mu}\le+\infty$ is always true.\\
(b) Assume $f^+$ to be integrable. \\
Note that we prove the more classical $\liminf-$version. This is done be making the necessary changes:  $F:=-f,F_n:=-f_n, G_n=-g_n, G=-g.$ Then all the statements above are still valid except $F_n\ge G_n$ and $F^-$ is now integrable. Then it is enough to show
\begin{align}\liminf\limits_{n\ra\infty}\int\limits_{X}{F_n\;d\mu}\ge \int\limits_{X}{F\;d\mu}.\label{eq:2.3.RF1.1}\end{align}

We start by considering the following decomposition
\begin{align*}\int\limits_{X}{F_n\;d\mu}=\int\limits_{X\cap\{F_n<-K\}}{F_n\;d\mu}+\int\limits_{X\cap\{F_n\ge-K\}}{F_n\;d\mu}=:A_{n,K}+B_{n,K},\end{align*}
which holds for any $n\in\N$ and any $K\in \R.$
The proof is divided into two steps. First one needs to show that $A_{n,K}\ra0$ for large enough $n$ and  $K.$ The second one is to show that $B_{n,K}$ can be controlled from below by $\int\limits_{X}{F\;d\mu}$ for large enough $n$ and  $K.$ Together with Step 1 that implies \eqref{eq:2.3.RF1.1}. \\

Step 1: $\liminf\limits_{K\ra+\infty}\liminf\limits_{n\ra+\infty} A_{n,K}\ge0.$\\

Since, $F_n\ge G_n$ we have
\begin{align*}A_{n,K}\ge& \int\limits_{X}{G_n\chi_{\{F_n<-K\}} \;d\mu}&\\
\ge& \int\limits_{X}{G_n\chi_{\{G_n<-K\}} \;d\mu}&
\end{align*}
for any $n\in\N$ and $K\in \R$. 
We can now split up the latter integral as follows
\begin{align}\int\limits_{X}{G_n\chi_{\{G_n<-K\}} \;d\mu}=\int\limits_{X}{(G_n-G)\chi_{\{G_n<-K\}} \;d\mu}+\int\limits_{X}{G\chi_{\{G_n<-K\}} \;d\mu}.\label{eq:2.3.RF1.2}\end{align}
Then the first term on the RHS converges to $0$ when $n\ra\infty,$ since  
\begin{align*}\int\limits_{X}{|(G_n-G)\chi_{\{G_n<-K\}}| \;d\mu}\le\int\limits_{X}{|G_n-G|\;d\mu}\ra0.\end{align*}
We can treat the rightmost term of \eqref{eq:2.3.RF1.2} by splitting up the integral as follows
\begin{align}\int\limits_{X}{G\chi_{\{G_n<-K\}} \;d\mu}=\int\limits_{X}{G\chi_{\{G_n<-K,|G_n-G|\le1\}} \;d\mu}+\int\limits_{X}{G\chi_{\{G_n<-K,|G_n-G|>1\}} \;d\mu}.\label{eq:2.3.RF1.3}\end{align}
Since $G\in L^1$ the rightmost term disappears when $n\ra\infty.$ Indeed, we have \\
\begin{align*}\mu(\{G_n<-K,|G_n-G|>1\})\le\mu(\{|G_n-G|>1\})\le\int\limits_{\{|G_n-G|>1\}}{|G_n-G| \;d\mu}\ra0\end{align*} when $n\ra\infty.$ Since $G\in L^1$ we have
\begin{align*}\int\limits_{X}{G\chi_{\{G_n<-K,|G_n-G|>1\}} \;d\mu}\ra0\end{align*}
when $n\ra\infty.$ \\
We can treat the first term of \eqref{eq:2.3.RF1.3} as follows. Realize first that
\begin{align*}\{G_n<-K,|G_n-G|\le1\}\ss\{G<-K+1\}\end{align*}
and hence by the Dominated Convergence Theorem we have
\begin{align*}G\chi_{\{G_n<-K,|G_n-G|\le1\}}\ra0 \;\mb{strongly in}\; L^1(X,[-\infty,\infty],\mu).\end{align*}
Together this shows
\begin{align*}\liminf\limits_{K\ra+\infty}\liminf\limits_{n\ra+\infty} A_{n,K}\ge0\end{align*}
completing Step 1.\\
 
Step 2:  $\liminf\limits_{K\ra+\infty}\liminf\limits_{n\ra+\infty} B_{n,K}\ge\int\limits_{X}{F\;d\mu}.$ \\
Since, $\mu(X)<\infty$ by Egorov's Theorem we can find for any $\ve>0$ we can find $X^\ve\ss X$ s.t. $\mu(X\sm X^\ve)<\ve$ and for any $\eta>0$ there is $N_\eta>0$ s.t. for any $n\ge N_\eta$ we have  $|F_n-F|<\eta$ on $X^\ve.$ With this at hand we can split up the integrals again as follows
\begin{align*}B_{n,K}=\int\limits_{X^\ve\cap\{F_n\ge-K\}}{F_n\;d\mu}+\int\limits_{(X\sm X^\ve)\cap\{F_n\ge-K\}}{F_n\;d\mu}.\end{align*}
For the rightmost term we find the simple estimate
\begin{align*}\int\limits_{(X\sm X^\ve)\cap\{F_n\ge-K\}}{F_n\;d\mu}\ge-K\mu(X\sm X^\ve)\ge-K\ve\end{align*} for any $K\in\R.$ Then since we have $|F_n-F|<\eta$ on $X^\ve$ this implies
\begin{align*}B_{n,K}\ge\int\limits_{X^\ve\cap\{F_n\ge-K\}}{F_n\;d\mu}-K\ve
\ge\int\limits_{X^\ve\cap\{F_n\ge-K\}}{(F-\eta)\;d\mu}-K\ve
\end{align*}
for any $n\ge N_\eta.$ Now we take $\ve_K:=K^{-2}$ and $\eta:=K^{-1}$ and $n\ge N(K)$.
Define $\phi_{n,K}:=F\chi_{\{X^{(\ve_K)}\cap\{F_n\ge-K\}\}}.$
Elementary set theoretic considerations yield for arbitrary sets $U,V,W$  with $V\ss U$ it holds
\begin{align*}U\sm(V\cap W)=(U\sm V)\cup (V\sm W)\ss(U\sm V)\cup (U\sm W).\end{align*}
Then for any $n\in \N$ we have
\begin{align*}\mu(X\sm\{X^{(\ve_K)}\cap\{F_n\ge-K\}\})\le\mu(X\sm X^{(\ve_K)})+\mu(X\sm\{F_n\ge-K\})\ra0\end{align*}
when $K\ra+\infty.$\\

Since  $F_n\ge G_n$ for any $n\in \N$ the same holds for the limit $F\ge G$ and hence we have $\phi_{n,K}=F\chi_{\{X^{(\ve_K)}\cap\{F_n\ge-K\}\}}\ge G$ for any $n\in \N$ and $K\in\R.$
Then for any $n\in \N,$ $(\phi_{n,K})_{K\in \R}$ is a sequence of measurable functions, which is bounded from below and we can apply the standard version of Fatou's Lemma to obtain\footnote{see \cite[Cor 5.34]{WZ15}.}
\begin{align*}\liminf\limits_{K\ra+\infty}\int\limits_{X}{\phi_{n,K}\;d\mu}\ge\int\limits_{X}{F\;d\mu}\;\mb{for any}\; n\in \N.\end{align*}
By putting everything together we get
\begin{align*}\liminf\limits_{K\ra+\infty}\liminf\limits_{n\ra+\infty}B_{n,K}\ge \int\limits_{X}{F\;d\mu},
\end{align*}
completing the proof.\vspace{0.5cm}

\textbf{Hardy spaces and Duality}\\
Initially, we recall the definition of Hardy spaces, where we follow \cite{U01}.
\begin{de}[Hardy spaces $\mathcal{H}^r$]
Let $S(\R^n)$ be the Schwartz class and $S'(\R^n)$ be its dual. For every $f\in S'(\R^n),$ and for  
$\Phi\in S(\R^n)$ s.t. $\int{\Phi \;dx}\not=0,$ a positive real parameter $\zeta>0$ and $x\in\R^n$ we define
 \[N_{\Phi,\zeta}f:=\sup\{f*(\Phi)_t(y):|x-y|<\zeta t\}.\]
Then for $r>0$ we can define the Hardy space by
\[\mathcal{H}_{\Phi,\zeta}^r(\R^n):=\{f\in S'(\R^n)|N_{\Phi,\zeta}f\in L^r(\R^n)\}.\]
Next, define the `norm'  by
\[\|f\|_{\mathcal{H}^r(\R^n)}:=\|N_{\Phi,\zeta}f\|_{L^r(\R^n)}.\]
Obviously, this quantity only defines a norm if $r\ge1$ and $\mathcal{H}^r$ is a Banach space. However, when $0<r<1,$ $\mathcal{H}^r(\R^n)$ is just a metric space. Hence, the notion of a dual space does not make any sense, in the later case.\\
It is important to note that these spaces are independent of $\Phi>0$ and $\zeta>0,$ i.e. for any $(\Phi,\zeta)$ and $(\Phi',\zeta')$ it holds 
\[\mathcal{H}_{\Phi,\zeta}^r(\R^n)=\mathcal{H}_{\Phi',\zeta'}^r(\R^n),\]
with equivalent `norms',i.e. for any choices of $(\Phi,\zeta)$ and $(\Phi',\zeta')$  there are constants $c=c(\Phi,\zeta)>0$ and $C=C(\Phi,\zeta)>0$ s.t. it holds
\[c\|N_{\Phi',\zeta'}f\|_{L^r(\R^n)}\le\|N_{\Phi,\zeta}f\|_{L^r(\R^n)}\le C\|N_{\Phi',\zeta'}f\|_{L^r(\R^n)}.\]
As a consequence, we can pick one and we can suppress the $\Phi,\zeta-$ dependence.\\

Lastly, recall that for any $1<r<\infty$ it holds \[\mathcal{H}^r(\R^n)=L^r(\R^n)\] with the equivalence of the norms, i.e. for any $1<r<\infty$ we can find constants $c_r>0$ and $C_r>0$ s.t. it holds
\begin{equation}c_r\|f\|_{L^r(\R^n)}\le\|f\|_{\mathcal{H}^r(\R^n)}\le C_r\|f\|_{L^r(\R^n)}.\label{eq:2.4.HS.equiv}\end{equation}
\end{de}

Secondly, we introduce the notion of a homogenous Lipschitz-space as can be found in \cite[Section 6.3]{G14}.

\begin{de}[Homogenous Lipschitz-spaces $\dot{\Lambda}_\beta$, $0<\beta\le1$] Define the semi-norm
\[\|f\|_{\dot{\Lambda}_\beta(\R^n)}:=\sup\limits_{\substack{x,y\in \R^n\\ x\not=y}}\frac{|f(x)-f(y)|}{|x-y|^\beta}.\]
Then we call
\[\dot{\Lambda}_\beta(\R^n):=\{f\in C(\R^n)|\|f\|_{\dot{\Lambda}_\beta(\R^n)}<\infty\}\]
\end{de}
 homogenous Lipschitz-space.\footnote{Note, that $\|\cd\|_{\dot{\Lambda}_\beta(\R^n)}$ is only a semi-norm on $\dot{\Lambda}_\beta(\R^n).$ However, $\dot{\Lambda}_\beta(\R^n)$ can be turned, as usual, into a normed space by $\dot{\Lambda}_\beta(\R^n)/\{constants\}.$}\\

\begin{re}1. Both types of spaces can be generalised to the case if the considered maps are vector-valued, as usual, we write $\mathcal{H}^r(\R^m,\R^n)$ and $\dot{\Lambda}_{\beta}(\R^m,\R^n).$ \\

2. For more information on Lipschitz spaces, for instance, a more general definition if $\be>1,$ or for the notion of inhomogeneous Lipschitz spaces see \cite[Section 6.3]{G14}.
\end{re}

Next, the duality result is stated, followed by a discussion of the literature, and where to find the proof.
\begin{lem} [$\mathcal{H}^r-\dot{\Lambda}_\beta$-Duality]
Assume $r<1,$ and $g\in \mathcal{H}^r(\R^n),$ and $f\in\dot{\Lambda}_\beta(\R^n)$ with $\beta:=n(\frac{1}{r}-1).$ Then there exists a constant $c(r,\Phi,\zeta)>0$ s.t.
\begin{equation}\int\limits_{\R^n}{fg\;dx}\le c(r,\Phi,\zeta)\|f\|_{\dot{\Lambda}_\beta(\R^n)}\|g\|_{\mathcal{H}^r(\R^n)}.
\label{eq:Ch1.4.103}\end{equation}
\end{lem}
\begin{re}
It is tough to find this statement in the literature.\footnote{Alternatively to the source we follow see \cite{JSW84}.} That is why we outline that the above is indeed true. We will follow Uchiyama's presentation, which he is using in his book \cite{U01} to stick to a single notation. \\

1. It has been shown that every function $f\in \mathcal{H}^r(\R^n)$ for $r\le1$ posses an atomic decomposition and hence $\mathcal{H}^r(\R^n)=\mathcal{H}_{at}^r(\R^n).$
Moreover, the norms are equivalent, too, i.e. there exists constant $c=c(r,\Phi,\zeta)>0$ and $C=C(r,\Phi,\zeta)>0$ s.t.
\[c(r,\Phi,\zeta)\|f\|_{\mathcal{H}_{at}^r(\R^n)}\le\|f\|_{\mathcal{H}^r(\R^n)} \le C(r,\Phi,\zeta)\|f\|_{\mathcal{H}_{at}^r(\R^n)},\]
see \cite[Eq.(0.8)-(0.9), p.7]{U01} for the statement and references to the original literature.\\

2. On the other hand Coifman and Weiss obtain in \cite[Theorem B, p.593]{Co77} the duality between $\mathcal{H}^r$ and $\dot{\Lambda}_\beta.$ 
It is important to realise that Coifman and Weiss use an equivalent notion to define these spaces by using a measure 
\[|f(x)-f(y)|\le \mu(B)^\alpha,\]
where $B$ is a ball containing $x,y$, see \cite[Eq.(2.2), p.591]{Co77}. By doing so obviously $\mu(B)\sim |x-y|^n$ and we get the relation $\beta=n\alpha.$\\

3. In particular, in \cite[Eq (3.19), p.634]{Co77} they obtain the duality inequality
\[\int\limits_{\R^n}{fg\;dx}\le \|f\|_{\dot{\Lambda}_\beta(\R^n)}\|g\|_{\mathcal{H}_{at}^r(\R^n)}.\]
This is indeed the case, however, they use a slightly different notation compared to Uchiyama. Coifman and Weiss use $|g|_{r,1}^{1/r}$ while Uchiyama has already absorbed the the exponent $1/p$ in the definition of $\|g\|_{\mathcal{H}_{at}^r(\R^n)},$ yielding the same result. Moreover, with the equivalence of the norms in point 1 we get 
\[\int\limits_{\R^n}{fg\;dx}\le C(r,\d,\phi)\|f\|_{\dot{\Lambda}_\beta(\R^n)}\|g\|_{\mathcal{H}^r(\R^n)},\]
which agrees with the conjecture.
\end{re}
We end this section by stating the famous div-curl Lemma first shown in \cite{CLMS93}.
\begin{lem}[div-curl Lemma]\label{lem:2.HR.1}
Let $\frac{n}{n+1}<s,t<\infty$ s.t. $\frac{1}{s}+\frac{1}{t}<1+\frac{1}{n},$ and $s, t$ are such that at least one of them is strictly larger than 1, and assume $F\in \mathcal{H}^s(\R^n,\R^n)$  and $G\in \mathcal{H}^t(\R^n,\R^n)$ satisfy \[\mb{curl}\, F=0 \;\mb{in}\; D'(\R^n,\R^\sigma)\;\mb{and}\; \div G=0 \;\mb{in}\; D'(\R^n),\] where $\sigma:=\mb{dim}\!\mb{curl}\!(\mathcal{H}^s(\R^n,\R^n)).$ Then $F\cd G\in \mathcal{H}^r(\R^n)$ for $\frac{1}{r}=\frac{1}{s}+\frac{1}{t}$ with
\begin{equation}\|F\cd G\|_{\mathcal{H}^r(\R^n)}\le C\|F\|_{\mathcal{H}^s(\R^n,\R^n)}\|G\|_{\mathcal{H}^t(\R^n,\R^n)},\end{equation}
for some positive constant $C>0.$
\end{lem}
\begin{re} Recall, that it holds
\[\div G=0 \;\mb{in}\; D'(\R^n) \;\mb{iff} \;\int\limits_{\R^n}{G\cd\grad h\;dx}=0 \;\mb{for any}\; h\in C_c^\infty(\R^n)\]
for any dimension $n\in \N\sm\{0\}.$\\
The distributional version of 'curl' is more complicated. For the relevant dimension $n=2$ it holds
\[\mb{curl}\, F=0 \;\mb{in}\; D'(\R^2) \;\mb{iff} \;\int\limits_{\R^2}{(F_2h,_1-F_1h,_2)\;dx}=0 \;\mb{for any}\; h\in C_c^\infty(\R^2). \]
A similar version can be obtained for $n=3.$ In higher dimensions one might want to use differential forms to define a distributional version of the 'curl'. The details are omitted for simplicity.
\end{re}

\clearpage{\pagestyle{empty}\cleardoublepage} 

\chapter{Radially symmetric $M-$covering stationary points of a polyconvex functional}
\section{Introduction and basic results}
In this section we consider the set of admissible functions which coincide with a $M-$covering map on the boundary.
For this sake, we define a $M-$covering map $u_M:S^1\rightarrow S^1$ for $M\in\N,$ $M\ge1$ via its representative
\begin{eqnarray}
\tilde{u_M}:[0,2\pi)&\rightarrow& S^1,\\
\th&\mapsto& (\cos M\th,\sin M\th).
\end{eqnarray}
We will use the notation $e_R(M\th):=(\cos M\th,\sin M\th)^T$ and $e_\th(M\th):=(-\sin M\th,\cos M\th)^T$ for all $M\in\N\sm\{0\}.$\\

Here we are interested in minimizing the functional \eqref{eq:1.1} in the special case that $u$ agrees with the $M-$covering map $u_M$ on the boundary, i.e.
\begin{equation}
\min\limits_{u\in \A_{u_M}}I(u).
\label{eq:3.MPMC}
\end{equation}
We use the method devised by P. Bauman, N.C. Owen and D. Phillips in two striking papers \cite{BOP91, BOP91MS}. In particular, this section will follow very closely the latter one without us mentioning it all the time.\\

In this chapter we will restrict the set of possible maps to radially symmetric $M-$covering maps: 
\[\A_r^M=\{u\in W^{1,2}(B, \R^2)|\exists r:[0,1]\rightarrow\R\; \mbox{s.t.}\; u(x)=r(R)e_R(M\theta)\;\mbox{and}\;r(1)=1\}\ss \A_{u_M}.\]
Then we are able to show the existence similarly to Theorem \ref{thm:2.1.1Ex}.
\begin{thm}[Existence]
There exists a minimizer $u\in \A_r^M$ s.t.
\[I(u)=\min\limits_{v\in \A_r^M} I(v).
\label{eq:}
\]
\end{thm}
{\bf{Proof:}} This can be shown analogously to the proof of Theorem \ref{thm:2.1.1Ex} using the direct method of the Calculus of Variations.

Further, we will restrict the set of test functions, to maps of the form  
\begin{equation}
\vp(x)=g(R)e_R(N\th) \;\mbox{where}\; g\in C_c^\infty((0,1)).
\label{eq:TFNC}
\end{equation}
Initially, we will discuss some space dependence of the radial part $r$ which must be satisfied if $u\in \A_r^M.$ Moreover, we will show that if $u\in \A_r^M$ is a stationary point of $I$, i.e. satisfies the ELE \eqref{eq:ELE1.1} and the test function is of the form \eqref{eq:TFNC} (with $M=N$), then the ELE reduces to an ODE of $r.$ If $r$ is a solution to this ODE, then we will show some further properties on $r.$  

\begin{lem}[Elementary properties]Let $M\in\N,$ $M\ge1$ and $u\in \A_r^M.$
\begin{enumerate}[label=(\roman*)]
\item Then $u\in \A_r^M$ if and only if $r$ is absolutely continuous on each compact subset of $(0,1]$ and $r\in L^2((0,1), R^{-1}dR)$ and $\dor\in L^2((0,1), RdR).$\\ 
\item If $M\not=N$ then $W$ is a Null-Lagrangian in the class of $\A_r^M$ and $N-$covering test functions, i.e.
\[\int\limits_B{(\grad u+\rho'(d)\cof\grad u)\cd\grad\vp\;dx=0 }\]
holds for all $u\in\A_r^M$ and all test functions of the form $\vp(x)=g(R)e_R(N\th)$ where $g\in C_c^\infty((0,1)),$ if $M\not=N.$
\item If additionally, $u$ is a minimizer to (\ref{eq:3.MPMC}) and $M=N$ then $r$ satisfies the ODE
\begin{equation}
\left(\frac{M^2r}{R}+M\rho'(d)\dot{r}\right)=\left(R\dot{r}+M\rho'(d)r\right)^\cdot \;\mbox{in}\; D'((0,1)).
\label{eq:3.15}
\end{equation}
\item Assume $r$ solves \eqref{eq:3.15}. Then $r\in C([0,1])$ and $r(0)=0.$
\item Moreover, $r\in C^{\infty}((0,1]).$
\end{enumerate} 
\end{lem}
{\bf{Proof:}}\\
(i): It is straightforward to see that for every $u\in\A_r^M$ it holds that
\[\|u\|_{W^{1,2}}^2=\int\limits_B{|u|^2+|\grad u|^2\;dx}=2\pi\int\limits_0^1{r^2R+\frac{M^2r^2}{R}+\dor^2R\;dR}\]
Since the LHS is finite the RHS needs to be finite, too, implying $r\in L^2((0,1), R^{-1}dR)$ and $\dor\in L^2((0,1), RdR).$\\
Let $[a,b]\ss(0,1].$
Since $r\in L^2((0,1], R^{-1}dR)$ it follows
\begin{equation}
\int\limits_a^b{r^2\;dR}\le b\int\limits_a^b{r^2\;\frac{dR}{R}}<\infty. 
\label{eq:3.24}
\end{equation}
Hence, $r\in L_{loc}^2((0,1]).$ Similar $\dor\in L^2((0,1], RdR)$ leads to 
\begin{equation}
\int\limits_a^b{\dor^2\;dR}\le\frac{1}{a}\int\limits_a^b{\dor^2\;RdR}<\infty
\label{eq:3.25}
\end{equation}
which implies $\dor\in L_{loc}^2((0,1])$ and $r\in W_{loc}^{1,2}((0,1]).$ Then $r$ agrees up to a set of measure zero with a function $\tilde{r}$ on $(0,1],$ where $\tilde{r}$ is absolutely continuous on any compact subset of $(0,1]$ (as always we identify $r$ with $\tilde{r}$). The latter is a consequence of \eqref{eq:3.25} and the fundamental theorem of calculus for Sobolev functions, see \cite[U1.6, p.71-72]{A12}.\\

(ii)-(iii): First we calculate some important quantities:
\begin{align}
\grad u=&{u,}_{R}\ot e_R+u,_{\ta}\ot e_\th=\dot{r}(R)e_R(M\th)\ot e_R(\th)+\frac{Mr}{R}e_\th(M\th)\ot e_\th(\th)&\\
\grad \vp=&\dot{g}(R)e_R(N\th)\ot e_R(\th)+\frac{Ng}{R}e_\th(N\th)\ot e_\th(\th)&\\
\cof\grad u=&\frac{Mr}{R}e_R(M\th)\ot e_R(\th)+\dot{r}(R)e_\th(M\th)\ot e_\th(\th)&\\
d=&\det\grad u=\frac{1}{2}\grad u\cdot \cof\grad u=\frac{Mr\dot{r}}{R}&\\
\grad u\cdot\grad \vp=&\dot{g}\dot{r}e_R(M\th)\cdot e_R(N\th)+\frac{MNgr}{R^2}e_\th(M\th)\cdot e_\th(N\th)&\\
\cof\grad u\cdot\grad \vp=&\frac{M\dot{g}r}{R}e_R(M\th)\cdot e_R(N\th)+\frac{Ng\dot{r}}{R}e_\th(M\th)\cdot e_\th(N\th)&\\
e_R(M\th)\cdot e_R(N\th)=&e_\th(M\th)\cdot e_\th(N\th)=\cos(M\th)\cos(N\th)+\sin(M\th)\sin(N\th)=\cos((M-N)\th).&
\end{align}

 Hence, the ELE \eqref{eq:ELE1.1} becomes
\begin{equation}
\int\limits_{0}^1\int\limits_{0}^{2\pi}{\left(\dot{g}\dot{r}+\frac{Mr\dot{g}}{R}\rho'(d)+\frac{MNgr}{R^2}+\rho'(d)\frac{Ng\dot{r}}{R}\right)\cos((M-N)\th)\;d\th RdR}=0
\label{eq:3.21}
\end{equation}
which is automatically true for all $r$ if $M\not=N$ and for $M=N$ takes the form
\begin{equation}
2\pi\int\limits_{0}^1{\dot{g}\left(\dot{r}R+M\rho'(d)r\right)+g\left(\frac{M^2r}{R}+M\rho'(d)\dot{r}\right)\;dR}=0
\label{eq:3.22}
\end{equation}
for all $g\in C_c^{\infty}((0,1)).$\\

(iv): Now we show that for any $r$ with $r\in L^2((0,1), R^{-1}dR)$ and $\dor\in L^2((0,1), RdR),$ it must hold that $\lim\limits_{R\ra0}r(R)=0.$ Suppose not, then wlog.\!\! there exists a strictly montonic decreasing sequence $\{R_j\}_{j\in \N}$ s.t. $R_j\ra0$ for $j\ra\infty$ and $|r(R_j)|>2\ve$ for any $j\in\N.$
Then since $\dor\in L^2((0,1), RdR),$ we can find $N\in\N$ so large that for any $n\ge N$ it holds
$\int\limits_0^{R_n}\dor^2R\;dR<\ve^2.$ Using the latter together with the fundamental theorem of calculus, and H\"older's inequality, then for any $n\ge N$ and any $R\in[R_n/e_n,R_n]$ with $e_n:=1-1/n$ (Note, that wlog. we can assume $R_{n+1}\le \frac{R_n}{e_n}$ if not consider the sequence $\wt{R}_{n+1}:=\min\{R_{n+1},\frac{R_n}{e_n}\}$) we obtain
 \[|r(R)-r(R_n)|=\left|\int_{R}^{R_n}\dot r(R')\, dR'\right|\le 
 \left| \int_{R}^{R_n}\frac1{R'}\, dR'\right|^\frac12\cdot 
 \left| \int_{R}^{R_n} 
 \left|\dot r(R')\right|^2R'\, dR'
 \right|^\frac12\le |\ln e_n|^\frac12 \ve.
 \]      
 Therefore by the reverse triangle inequality we have $|r(R)|>\ve$ for any $R\in[R_n/e_n,R_n]$. Then
 \[
 \int_{0}^{1}{r^2\over R}\,d R\ge\sum\limits_{n=N}^\infty\int_{R_n/e_n}^{R_n}{r^2\over R}\,d R> \sum\limits_{n=N}^\infty\int_{R_n/e_n}^{R_n}{\ve^2\over R}\,d R=\ve^2\sum\limits_{n=N}^\infty\left|\ln \left(1-\frac1{n}\right)\right|=+\infty,
 \]
 contradicting the integrability of ${r^2\over R}$ and showing that $r\in C([0,1])$ with $r(0):=0.$\\

(v): Next we show that $\dor\in C((0,1]).$\\
Let
\begin{equation}
q(w,a,b)=bw+M\rho'\left(\frac{Mwa}{b}\right)a
\label{eq:3.23}
\end{equation}
with $1\ge b>0,a,w\in\R.$
Then $q(\cdot,a,b)$ is a homeomorphism from $\R$ to $\R$ for all $a\in \R$ and  $b\in(0,1].$ Indeed, for $a=0,$ $q(w,0,b)=bw$ is a homeomorphism from $\R$ to $\R.$ Now let $b\in(0,1]$ and $a>0.$ Then $w\mt bw$ is strictly monotonically increasing and maps $\R$ to $\R$ continuously. Moreover, $w\mt M\rho'\left(\frac{Mwa}{b}\right)a$ is monotonically increasing and continuous as well, hence, $w\mt q(w,a,b)$ is a homeomorphism on $\R.$ If $a<0$ then $w\mt\rho'\left(\frac{Mwa}{b}\right)$ decreases but $w\mt M\rho'\left(\frac{Mwa}{b}\right)a$ still increases, and we can argue as above.\\
Assume that there exists $R_0\in(0,1]$ and a sequence $R_j\ra R_0$ for $j\ra\infty,$ s.t. $\dor(R_j)\ra\pm\infty.$ But then
\[q(\dor(R_j),r(R_j),R_j)=R_j\dor(R_j)+M\rho'\left(\frac{M\dor(R_j)r(R_j)}{R_j}\right)r(R_j)\ra\pm\infty+c=\pm\infty\]
which is impossible, since $R\ra q(\dor(R),r(R),R)$ is continuous on $(0,1].$ Hence, $\dor\in C((0,1]).$\\
As a last step we improve the regularity to $C^\infty((0,1]).$\\
Using the ODE we can represent $q$ by
\begin{equation}
q(\dor(R),r(R),R)=R\dor(R)+M\rho'(d)r(R)=c-\int\limits_R^1{\frac{M^2r(R')}{R}+M\rho'(d)\dor(R')\;dR'}
\label{eq:3.23}
\end{equation}
Then by the analysis above, $q$ is actually $C^1((0,1]).$ Further the derivative w.r.t. $w$ is $\p_wq=b+M^2\rho''(\frac{Mwa}{b})\frac{a^2}{b}>0$ for all $1\ge b>0,a,w\in\R.$ Then the implicit function theorem gives full regularity $r\in C^\infty((0,1]).$\\\vspace{0.5cm}

Collecting these results, from now on, we will consider solutions $r\in C([0,1])\cap C^\infty((0,1])$ to the boundary value problem
\begin{align}\left\{\begin{array}{ccc}
\frac{M^2r}{R}=\dor+R\ddor+M\rho''(d)\dod r & {\mbox{in}}& (0,1),\\
r(0)=0,&r(1)=1.&
\end{array}\label{eq:3.0.1}
\right.\end{align}

\section{The simple case $M=1$}
\label{sec:3.1}
In the $M=1$ case the 1-covering map $e_R(\th)$ agrees with the identity hence, $u_1=\Id$ on the boundary of the unit ball in $\R^2.$ 
A short calculation shows that $r(R)=R$ and $r(R)=\frac{1}{R}$ are both solutions to the ODE in \eqref{eq:3.0.1} but only $r(R)=R$ matches both boundary conditions.
Indeed, $r(R)=\frac{1}{R}$ satisfies $r(1)=1$ but not $r(0)=0.$ The latter condition originates from the fact that $u$ belongs to $W^{1,2}(\Om,\R^2).$ Note that $u=\frac{1}{R}e_R(\th)=\frac{x}{|x|^2}$ is excluded by the considered function space. \\

In this case we will show that $r(R)=R$ is not only a solution to \eqref{eq:3.0.1} but also $u=\Id=Re_R(\th)$ is the unique global minimizer to \eqref{eq:1.1}.

\begin{pro}
Let $0<\ga<\infty.$ Then $u=\Id$ is the unique smooth global minimizer of the functional \eqref{eq:1.1}.
\label{prop:3.1.1}  
\end{pro}
\textbf{Proof:}
Let $u=\Id+\vp$ with $\vp\in C_c^\infty(B,\R^2).$ Then 
\begin{align*}
I(u)-I(Id)&=\frac{1}{2}\int\limits_B{|\grad(x+\vp(x))|^2-|\grad x|^2\;dx}+\int\limits_B{\rho(d_{\grad u})-\rho(d_{\grad x})\;dx}&\\
&\ge\int\limits_B{\frac{1}{2}|\grad\vp|^2+\frac{1}{2}|\Id|^2-\grad x\cd\grad\vp\;dx}+\int\limits_B{\rho'(d_{\grad x})(d_{\grad u}-d_{\grad x})\;dx}&\\
&\ge\int\limits_B{\frac{1}{2}|\grad\vp|^2\;dx}+\rho'(d_{\Id})\int\limits_B{(d_{\grad u}-d_{\Id})\;dx}&\\
&=\frac{1}{2}\int\limits_B{|\grad\vp|^2\;dx}\ge0,&
\end{align*}
where we used that $\vp\in C_c^\infty(B,\R^2)$ and the fact that the determinant is a Null-Lagrangian. Note that equality holds iff $u\equiv \Id$ $(\vp\equiv0),$ hence $\Id$ is the unique global minimizer independent of $\ga$ and $\rho$.

\begin{re} Proposition \ref{prop:3.1.1} amounts to the statement that the functional $I(u)=\int\limits_\Om\frac{1}{2}|\grad u|^2+\rho(\det\grad u)\;dx$ is (strictly) quasiconvex.
\end{re}

\section{The general case $M\ge1$ and classical BOP-Theory}
In this paragraph it is shown that functions are at least of class $C^1$ on the whole interval $[0,1]$. Again, we use the method invented in \cite{BOP91MS}. The plan is as follows: First we will consider two auxiliary functions (namely, $d$ the determinant and $z$ defined below) which depend on $r$ and $\dor.$ We will study their behaviour in particular, close to the origin. This will then reduce the number of possibilities how $d$ and $z$ can behave. Then we can discuss these cases one by one and study the behaviour of $r$ and $\dor$.\\

We start with the following observations on the determinant. 
 
\label{sec:3.2}
\begin{lem}Let $M\in\N,\; M\ge1$ and assume that $r$ solves the BVP \eqref{eq:3.0.1}. Then $d\in C^\infty((0,1])$ and $\dod\ge0$ in $(0,1].$ Moreover, it holds that $d\in C([0,1])$ and $d\ge0$ in $[0,1].$
\label{lem:Mdet}
\end{lem}
{\bf{Proof:}}\\
The smoothness of $d$ in $(0,1]$ follows by the smoothness of $r$ in the same interval.
\[
\dod=\left(\frac{Mr\dor}{R}\right)^\cdot=\frac{M\dor^2}{R}+\frac{Mr\ddor}{R}-\frac{Mr\dor}{R^2}
\label{eq:3.d1}
\]
Now by multiplying the strong form of the ODE \eqref{eq:3.15} by $\frac{r}{R}$ we can express the term $r\ddor$ as
\begin{eqnarray}
r\ddor=\frac{M^2r^2}{R^2}-\frac{r\dor}{R}-M\rho''(d)\frac{r^2}{R}\dod.
\label{eq:3.d2}
\end{eqnarray}
Hence,
\[\dod=\frac{M\dor^2}{R}-\frac{Mr\dor}{R^2}+\frac{M^3r^2}{R^3}-\frac{Mr\dor}{R^2}-M^2\rho''(d)\frac{r^2}{R^2}\dod\]
and rearranging the equation yields, 
\begin{eqnarray}
\dod&=&(1+M^2\rho''(d)\frac{r^2}{R^2})^{-1}\left[\frac{M\dor^2}{R}-\frac{2Mr\dor}{R^2}+\frac{M^3r^2}{R^3}\right] \nonumber\\
&=&\frac{M[(R\dor-r)^2+(M^2-1)r^2]}{R^3+M^2\rho''(d)r^2R}\ge0.
\label{eq:3.d3}
\end{eqnarray}

Assume now that $\lim\limits_{R\ra0}d(R)\in[-\infty,0).$ By the smoothness of $d$ in $(0,1],$ there exists $\d>0$ s.t. $d(R)<0$ for all $R\in(0,\d).$ This implies that for all $R\in(0,\d)$ either $r(R)>0$ and $\dor(R)<0$ or $r(R)<0$ and $\dor(R)>0.$ Consider the first case $r,-\dor>0$ on $(0,\d).$ By the mean value theorem we get that there exists a $\xi\in(0,\frac{\d}{2})$ s.t. $\dor(\xi)=\frac{2r(\frac{\d}{2})}{\d}>0,$ contradicting $\dor(R)<0$ for all $R\in(0,\d).$ (Analogously, for the other case). This proves $\lim\limits_{R\ra0}d(R)\in[0,+\infty].$ By $\dod\ge0$ and $r$ smooth in $(0,1]$ $d$ cannot attain $+\infty.$ Therefore, the limit $d(0):=\lim\limits_{R\ra0}d(R)\in[0,\infty)$ exists and is nonnegative. Again by $\dod\ge0,$ $d$ remains nonnegative throughout the whole interval $[0,1].$  \\\vspace{0.5cm}

As a consequence, the nonnegativity and the monotonic growth of $d$ are transferred on to $r.$

\begin{lem} Let $M\in\N,$ $M\ge1$ and $r$ solves the BVP \eqref{eq:3.0.1}. Then $r(R)\ge0$ for all $R\in[0,1]$ and $\dor(R)\ge0$ for all $R\in(0,1].$ 
\end{lem}
{\bf{Proof:}}\\
By Lemma \ref{lem:Mdet} we know that
\[d=\frac{Mr\dor}{R}=\frac{M(r^2)^\cdot}{2R}\ge0.\]
Hence, $(r^2)^\cdot\ge0$ for all $R\in(0,1].$ For the sake of a contradiction, assume that there exists $R_0\in(0,1]$ s.t. $r(R_0)<0.$ Then by the continuity of $r$ and since $r^2$ grows monotonically, $r$ remains negative up to the boundary, i.e. $r(R)\le r(R_0)<0$ for all $R\in[R_0,1].$ This is not compatible with the boundary condition $r(1)=1$ yielding $r\ge0$ in $[0,1].$\\
The second claim follows in a similar fashion. Again we make the assumption that there exists $R_0\in(0,1]$ s.t. $\dor(R_0)<0.$ By continuity, there even exists an interval $(R_1,R_2]$ with $0<R_1<R_2\le1$ s.t. $\dor(R)<0$ for all $R\in (R_1,R_2]$ then by monotonicity of $r^2$ we know that $0\le r(R_2)-r(R_1)$ but on the other hand by the fundamental theorem of calculus we have
\[0\le r(R_2)-r(R_1)=\int\limits_{R_1}^{R_2}{\dor(R)\;dR}<0\]
leading again to a contradiction. Hence, $\dor\ge0$ in $(0,1].$\\\vspace{0.5cm}

Next we introduce the function $z(x):=\frac{1}{2}|\grad u(x)|^2+f(\det\grad u(x))$ for all $x\in\ol{B},$ where $f(d):=d\rho'(d)-\rho(d)$ for all $d\in\R.$ In the following lemma it is shown that $z$ satisfies a maximum principle in $\ol{B}\sm\{0\}$. This follows closely Chapter $3$ of \cite{BOP91}, in particular Theorem $3.2-3.3.$

\begin{lem} Let $M\in\N,\; M\ge2$ and assume $r$ solves the BVP \eqref{eq:3.0.1}. Then $z$ satisfies the strong maximum principle in $(0,1].$
\label{lem:3.4}
\end{lem}
{\bf{Proof:}}\\
It is enough to show that $z$ is a subsolution to an elliptic equation, i.e. $\D z+c_M\rho''(d)\doz\ge0$ in $(0,1],$ where $c_M:=\frac{M}{2(M-1)}.$ This is indeed enough to apply the strong maximum principle, see \cite[\S6.4.2 Thm 3]{LE10}. \\
Initially, note that for $u\in\A_r^M,$  $z(x)=z(R),$ $z(R)=\frac{\dor^2}{2}+\frac{M^2r^2}{2R^2}+f(d),$ $f(d)=d\rho'(d)-\rho(d)$ and $\D z=\frac{\doz}{R}+\ddoz.$\\
Now in order to calculate $\D z$ we first need to calculate $\doz$ and
 $\ddoz.$ Taking the derivative of $z$ wrt. $R$ yields,
\begin{equation}
\doz=\dor\ddor+\frac{M^2r\dor}{R^2}-\frac{M^2r^2}{R^3}+\rho''(d)d\dod
\label{eq:3.z1}
\end{equation}
where we used $(f(d))^\cdot=\rho''(d)d\dod.$ The strong version of the ODE \eqref{eq:3.15} is given by
\begin{equation}
\frac{M^2r}{R}=\dor+R\ddor+M\rho''(d)\dod r,
\label{eq:3.z2}
\end{equation}
which, when multiplied by $\frac{\dor}{R}$ leads to 
\begin{equation}
\rho''(d)d\dod=\frac{M^2r\dor}{R^2}-\frac{\dor^2}{R}-\dor\ddor.
\label{eq:3.z3}
\end{equation}
Substituting, $\rho''(d)d\dod$ in \eqref{eq:3.z1} via \eqref{eq:3.z3} yields,
\begin{equation*}
\doz=-\frac{M^2r^2}{R^3}+\frac{2M^2r\dor}{R^2}-\frac{\dor^2}{R}.
\label{eq:3.z4}
\end{equation*}
The second derivative of $z$ is then given by
\[
\ddoz=-\frac{2M^2r\dor}{R^3}+\frac{3M^2r^2}{R^4}+\frac{2M^2\dor^2}{R^2}+\frac{2M^2r\ddor}{R^2}-\frac{4M^2r\dor}{R^3}-\frac{2\dor\ddor}{R}+\frac{\dor^2}{R^2}.
\label{eq:3.z5}
\]
and the Laplacian becomes
\[
\D z=-\frac{2M^2r^2}{R^4}+\frac{2M\dor^2}{R^4}+\frac{2Mr\ddor}{R^2}-\frac{4Mr\dor}{R^3}-\frac{2\dor\ddor}{R}.
\label{eq:3.z6}
\]
The equations (\ref{eq:3.z3}) and (\ref{eq:3.d2}) allow us to replace the terms with a second derivative $\ddor$
\begin{equation}
\D z=\frac{2}{R^4}[M^2r^2(1+M)+R^2\dor^2(1+M)-2MRr\dor]-\rho''(d)\dod\left(\frac{2M^2r^2}{R^3}-\frac{2d}{R}\right).
\label{eq:3.z7}
\end{equation}
Now define
\begin{equation}
s(R):=M^2r^2(1+M)+R^2\dor^2(1+M)-2RMr\dor.
\label{eq:3.z8}
\end{equation}
Completion of the square, yields
 \[
s(R)=(Mr-R\dor)^2 +M^3r^2+MR^2\dor^2\ge0.
\label{eq:3.z9}
\]
To deal with the `$\rho''$-terms' of (\ref{eq:3.z7}) we  use the form (\ref{eq:3.d3}) of $\dod$ to obtain

\begin{align*}
\left(\frac{2M^2r^2}{R^3}-\frac{2d}{R}\right)\dod=&\left(\frac{2M^2r^2}{R^3}-\frac{2d}{R}\right)\frac{M[(R\dor-r)^2+(M^2-1)r^2]}{R^3+M^2\rho''(d)r^2R}&\\
=&\frac{2M^2[M^3r^4-(2M+M^2)Rr^3\dor+(2+M)R^2r^2\dor^2-R^3r\dor^3]}{R^3(R^3+M^2\rho''(d)r^2R)}&
\label{eq:3.z10}
\end{align*}
Denoting the expression in the brackets by $t(R),$ i.e.
\[
t(R):=M^3r^4-(2M+M^2)Rr^3\dor+(2+M)R^2r^2\dor^2-R^3r\dor^3
\label{eq:3.z11}
\]
we get 
 \begin{eqnarray*}
\D z&=&\frac{2s}{R^4}-\rho''(d)\frac{2M^2t}{R^3(R^3+M^2\rho''(d)r^2R)}\\
&=&\frac{2}{R^4}\frac{[R^3s+\rho''(d)M^2R(r^2s-t)]}{R^3+M^2\rho''(d)r^2R}.
\label{eq:3.z12}
\end{eqnarray*}
Since $\rho''(d)\ge0$ the denominator is nonnegative and we know that $s\ge0,$ therefore, the first term is nonnegative. To complete the proof it is enough to show that
 \begin{eqnarray}
\frac{2}{R^4}\frac{M^2R(r^2s-t)}{R^3+M^2\rho''(d)r^2R}+c_M\doz^2\ge0.
\label{eq:3.z13}
\end{eqnarray}
In order to prove (\ref{eq:3.z13}) we show the slightly stronger statement
 \begin{eqnarray}
2M^2(r^2s-t)+c_M\doz^2R^6\ge0.
\label{eq:3.z14}
\end{eqnarray}
Assume that (\ref{eq:3.z14}) holds, adding the nonnegative term $c_M\doz^2R^3(M^2\rho''(d)r^2R)$ to the left hand side and dividing by $R^3(R^3+M^2\rho''(d)r^2R)$ yields
\[
\frac{2M^2(r^2s-t)+c_M\doz^2R^3(R^3+M^2\rho''(d)r^2R)}{R^3(R^3+M^2\rho''(d)r^2R)}\ge0,
\label{eq:3.z15}
\]
which agrees with (\ref{eq:3.z13}). First we compute the quantities 
\[
\doz R^6 =M^4r^4+R^4\dor^4+(4M^4+2M^2)R^2r^2\dor^2-4M^4Rr\dor
\label{eq:3.z16}
\]
and
\[
2M^2(r^2s-t) =2M^4r^4-2M^2 R^2r^2\dor^2+2M^2R^3r\dor^3+2M^4Rr^3\dor.
\label{eq:3.z17}
\]
Then \eqref{eq:3.z14} becomes
\begin{flalign}
&2M^2(r^2s-t)+c_M\doz R^6& \nonumber\\
&=c_MR^4\dor^4+(2+c_M)M^4r^4+(4M^4c_M+2M^2c_M-2M^2) R^2r^2\dor^2 \nonumber\\
&\hspace{0.5cm}+2M^2(1-2c_M)R^3r\dor^3+2M^4(1-2c_M)Rr^3\dor& \nonumber\\
&=\left(\frac{2c_M-1}{2}\right)(M^{5/4}r-M^{1/4}R\dor)^4+\left(c_M-\left(\frac{2c_M-1}{2}\right)M\right)R^4\dor^4& \nonumber\\
&\hspace{0.5cm}+\left((2+c_M)M^4-\left(\frac{2c_M-1}{2}\right)M^5\right)r^4& \nonumber\\
&\hspace{0.5cm}+\left(4M^4c_M+2M^2c_M-2M^2-6\left(\frac{2c_M-1}{2}\right)M^3\right) R^2r^2\dor^2.&
\label{eq:3.z18}
\end{flalign}
In the last step we completed the quartic form. Note that $\left(\frac{2c_M-1}{2}\right)=\frac{c_M}{M}=\frac{1}{2(M-1)}\ge0.$ Furthermore, 
the coefficients satisfy
\begin{align*}
c_M-\left(\frac{2c_M-1}{2}\right)M=&c_M-\left(\frac{c_M}{M}\right)M=0,&\\
(2+c_M)M^4-\left(\frac{2c_M-1}{2}\right)M^5=&M^4\left((2+c_M)-\left(\frac{c_M}{M}\right)M\right)=2M^4\ge0.&
\end{align*}
The last coefficient in \eqref{eq:3.z18} is nonnegative if it satisfies the following condition
\begin{align*}
4M^4c_M+2M^2c_M-2M^2-6\left(\frac{2c_M-1}{2}\right)M^3=&4M^4c_M+2M^2c_M-2M^2-6c_MM^2&\\
=&4M^4c_M-4M^2c_M-2M^2\ge 0&
\end{align*}
equivalently,
\begin{align*}
1\le(2M^2-2)c_M=M\left(\frac{2M^2-2}{2M-2}\right),
\end{align*}
which is true for all $M\in\N,$ $M\ge2.$ This finishes the proof since all terms in (\ref{eq:3.z18}) are nonnegative yielding $2M^2(r^2s-t)+c_M\doz R^6\ge0$ and with the discussion above  $\D z+(\frac{M}{2(M-1)})\rho''(d)\doz\ge0$ in $(0,1]$.\\\vspace{0.5cm}

The latter statement now allows to control the behaviour close to the origin. Indeed, $z$ is monotonic close to $0.$ This gives constraints on the quantity $\frac{R\dor}{r}$ which will be useful later.

\begin{lem}\label{lem:3.7} Let $M\in\N,$ $M\ge2.$ and $r$ solves the BVP \eqref{eq:3.0.1}. Then there exists $\d>0$ s.t. $z$ is monotone on $(0,\d).$\\
Additionally, assume $r>0$ in $(0,\d).$ Then one of the following conditions holds identically in $(0,\d):$
\begin{align}
\dot{z}\ge0 \;&\mbox{equivalently}\; 0<M^2-M\sqrt{M^2-1}\le\frac{R\dor}{r}\le M^2+M\sqrt{M^2-1},\label{eq:3.23a}&\\
\dot{z}\le0 \;&\mbox{and}\; \left(\frac{r}{R}\right)^\cdot> 0 \;\mbox{or}\label{eq:3.23b}&\\
\dot{z}\le0 \;&\mbox{and}\; \left(\frac{r}{R}\right)^\cdot< 0.&
\label{eq:3.23c}
\end{align} 
\label{lem:3.z}
\end{lem}

By the previous Lemma \ref{lem:3.4} we know either that $z(x)$ is constant, in which case the monotonicity is given, or $z$ does not attain a maximum in $B\sm\{0\}.$ It follows that $R\mapsto z(R)$ can only have one local minimum in $(0,1],$ and therefore $\doz$ can only change sign once. Hence, there exists a $\d>0$ s.t. $R\mapsto z(R)$ is monotone on $(0,\d).$\\   

Suppose now that $\d>0$ s.t.\! the above holds and $r(R)>0$ for all $R\in(0,\d)$ and recall $z=\frac{\dor^2}{2}+\frac{M^2r^2}{2R^2}+f(d).$ Then the derivative of $z$ is given by
\begin{flalign*}
\dot{z}=\dor\ddor+\frac{M^2r\dor}{R^2}-\frac{M^2r^2}{R^3}+(f(d))^\cdot
\label{eq:3.23d}
\end{flalign*} 
Using the definition of $f$ we obtain
\[
(f(d))^\cdot=(\rho'(d)d-\rho(d))^\cdot=\rho'(d)\dod+\rho''(d)d\dod-\rho'(d)\dod=\rho''(d)d\dod.
\label{eq:3.23e}
\]
On the other hand, since $r$ is a strong solution to the ODE \eqref{eq:3.15} in $(0,1],$ all derivatives exist in a strong sense. In particular,
\begin{align*}
\frac{M^2r}{R}+M\rho'(d)\dot{r}=&\left(R\dot{r}+M\rho'(d)r\right)^\cdot& \nonumber\\
=&\dor+R\ddor+M\rho'(d)\dot{r}+M\rho''(d)\dod r.&
\label{eq:3.23f}
\end{align*}
Rearranging the above equations and multiplying it by $\frac{\dor}{R}$ yields
\[
(f(d))^\cdot=\rho''(d)d\dod=\frac{M^2r\dor}{R^2}-\frac{\dor^2}{R}-\dor\ddor.
\label{eq:3.23g}
\]
Therefore,
\begin{align}
\dot{z}(R)=-\frac{r^2}{R^3}\left[\left(\frac{R\dor}{r}\right)^2-2M^2\left(\frac{R\dor}{r}\right)+M^2\right].
\label{3.23h}
\end{align}

This polynomial of variable $\frac{R\dor}{r}$ has roots at $\lambda_{\pm}=M^2\pm M\sqrt{M^2-1}.$ It is negative between these roots and positive otherwise. In particular, if $\doz\ge0$ then , from \eqref{3.23h}, it follows that $M^2-M\sqrt{M^2-1}\le\frac{R\dor}{r}\le M^2+M\sqrt{M^2-1},$ which agrees with \eqref{eq:3.23a}. Now for \eqref{eq:3.23b} and \eqref{eq:3.23c} from $\dot{z}\le0$ on $(0,\d)$ it follows $\frac{R\dor}{r}\le M^2- M\sqrt{M^2-1}<1$ or $\frac{R\dor}{r}\ge M^2+ M\sqrt{M^2-1}>1.$ Then the calculation
\begin{eqnarray*}
\frac{R\dor}{r}&\gtrless& 1 \nonumber\\
\Rightarrow\frac{\dor}{R}&\gtrless& \frac{r}{R^2} \nonumber\\
\Rightarrow\left(\frac{r}{R}\right)^\cdot&=&\frac{\dor}{R}-\frac{r}{R^2}\gtrless 0
\label{3.23i}
\end{eqnarray*}
shows (\ref{eq:3.23b}) and (\ref{eq:3.23c}).\\\vspace{0.5cm}

Up to this point we have narrowed the number of possibilities, how $d$ and $z$ can behave close to $0,$ enough so that there are only a few cases left which we now can discuss individually. The shape of $d$ and $z$ will demand certain conditions on $r,$ which can either be matched by $r$ or will lead to contradictions, excluding these cases. This will reduce the number of types even further and leaves only the following restrictive statement:   
\begin{lem}Let $M\in\N,$ $M\ge2$ and suppose $r$ solves the BVP \eqref{eq:3.0.1}. Then $d(0)=0.$
Moreover, one of the two situations occurs:\\
i) $r$ is lifting-off delayed, i.e. there exists, $0<\d<1$ s.t. $r\equiv0$ on $[0,\d].$ Then $r\in C^\infty([0,1]).$\\
ii) $r$ is lifting-off immediately, i.e. $r,\dor>0$ away from zero. Then $r\in C^1([0,1])$ and $\dor(0)=0.$
\label{lem:3.2.100}
\end{lem}
{\bf{Proof:}}\\
We show that only $d(0)=0$ is possible and either $r\equiv 0$ near zero or $r>0$ away from zero and $\doz\ge0.$ All other situations are excluded by contraposition.\\
   
Recall, that $\lim\limits_{R\ra0} d(R)$ exists and agrees with $d(0)$ due to continuity, which was proven in lemma \ref{lem:Mdet}. Also notice that lemma \ref{lem:3.7} guarantees that $\lim\limits_{R\ra0} z(R)$ makes sense, however the limit might be $+\infty.$ In particular, we know $0\le\lim\limits_{R\rightarrow0} z(R)\le+\infty$ by the definition of $z$ and the behaviour of $f.$ \\

1. Case: $\lim\limits_{R\rightarrow0} d(R)=:l\in(0,\infty).$\\
By continuity of $d$ on $(0,1]$ there exists $\delta>0$ s.t. $r\dor>0$ on $(0,\delta)$ implying $r,\dor>0$ on $(0,\d)$ (even on $(0,1]$, by monotonicity).\\
First, assume $\lim\limits_{R\rightarrow0} z(R)=+\infty.$ In this case only $\doz\le0$ near zero is possible implying $\frac{r}{R}$  to be monotone on $(0,\d).$ The mean value theorem guarantees the existence of the two sequences $R_j'\ra0$ and $R_j\in(0,R_j')$ for any $j\in \N$ s.t. 
\begin{equation}
\lim\limits_{j\rightarrow\infty} \dor(R_j)=\lim\limits_{j\rightarrow\infty} \frac{r(R_j')}{R_j'}=:m.
\label{eq:3.l.9.1}
\end{equation}
By lemma \ref{lem:3.7} we know that there are two different cases, either $\dor>\frac{r}{R}$ or $\dor<\frac{r}{R}$ on $(0,\d).$ In the first case by $\dor>\frac{r}{R}$ and \eqref{eq:3.l.9.1} we have
\[
l=\lim\limits_{j\ra\infty} d(R_j')=\lim\limits_{j\ra\infty} \frac{M\dor(R_j') r(R_j')}{R_j'}\ge\lim\limits_{j\ra\infty} M \frac{r^2(R_j')}{R_j'^2}=Mm^2.
\]
Hence, we know $m\le\sqrt{\frac{l}{M}}<+\infty$ and together with the property $\dor>\frac{r}{R}$ on $(0,\d)$ it holds
\[
\lim\limits_{j\rightarrow\infty} z(R_j)=\lim\limits_{j\rightarrow\infty}\left[\frac{\dor^2(R_j)}{2}+\frac{M^2r^2(R_j)}{2R_j^2}+f(d(R_j))\right]\le\left[\frac{1}{M}+M\right]\frac{l}{2}+f(d(l))<+\infty,\]
contradicting $\lim\limits_{R\rightarrow0} z(R)=+\infty.$

In the 2nd case, that is $\dor<\frac{r}{R}$ on $(0,\d),$ it holds
\[
l=\lim\limits_{j\ra\infty} d(R_j)=\lim\limits_{j\ra\infty} \frac{M\dor(R_j) r(R_j)}{R_j}\ge\lim\limits_{j\ra\infty} M\dor^2(R_j)=Mm^2.
\]
Again $m\le\sqrt{\frac{l}{M}}<+\infty$ and together with the property $\dor<\frac{r}{R}$ on $(0,\d)$ we have
\[
\lim\limits_{j\rightarrow\infty} z(R_j')=\lim\limits_{j\rightarrow\infty}\left[\frac{\dor^2(R_j')}{2}+\frac{M^2r^2(R_j')}{2R_j'^2}+f(d(R_j'))\right]\le\left[\frac{1}{M}+M\right]\frac{l}{2}+f(d(l))<+\infty,\]
contradicting $\lim\limits_{R\rightarrow0} z(R)=+\infty.$\\
\vspace{0.5cm}

Now assume the limit exists, i.e. $\lim\limits_{R\rightarrow0} z(R)=:n\in[0,\infty).$ For the sake of contradiction we show that the right limit of $\dor$ exists in $0$ and that it is nonzero.\\
Introduce the new variables $\nu_1(R):=\dor(R)$ and $\nu_2(R):=\frac{Mr(R)}{R}.$ Then we can interpret the functions $d$ and $z$ as functions depending on these new variables $d(\nu_1,\nu_2)=\nu_1\nu_2$ and $z(\nu_1,\nu_2)=\frac{\nu_1^2}{2}+\frac{\nu_2^2}{2}+f(\nu_1\nu_2)$ on the set
$\mathcal{V}=\{(\nu_1,\nu_2):\nu_1>0,\nu_2>0\}.$\\
Consider 
\[
K:=\{(\nu_1,\nu_2):\frac{l}{2}\le d\le 2l\}\cap\{(\nu_1,\nu_2):s-1\le z\le s+1\}
\]
here $l>0$ and $s>1$ are parameters, in particular, $s$ is not the function introduced in \eqref{eq:3.z8}. $K$ is a compact subset of $\mathcal{V}$ due to the continuity of $d$ and $z.$ In particular, $\{d=l\}\cap\{z=s\}$ consists of at most two points $\{(a,b),(b,a)\}\ss\mathcal{V}.$\\ Now again by the continuity of $r$ and $\dor$ we know that for all $\ve>0$ there exists $R_0(\ve)>0$ s.t.\@ for all $0<R<R_0,$ $(\nu_1(R),\nu_2(R))\in B_\ve(a,b)\cup B_\ve(b,a)\ss K.$ Now we need to show that $R\mapsto(\nu_1(R),\nu_2(R))$ remains in one of these balls for all $0<R<R_0.$ If $a=b$ this is immediately true. So suppose $a\not=b.$ Then we can choose $\ve>0$ so small that the balls become disjoint, i.e. $B_\ve(a,b)\cup B_\ve(b,a)=\emptyset.$ Recall that $r\in C^\infty((0,1])$ therefore, the curve $R\mapsto(\nu_1(R),\nu_2(R))=(\dor(R),\frac{Mr(R)}{R})$ is connected and remains in one ball, say $B_\ve(a,b).$ Since $\ve>0$ was arbitrary, we see that
\begin{eqnarray*}
\lim\limits_{R\rightarrow0}\dor(R)=a\in(0,+\infty)
\end{eqnarray*}
Hence, $r\in C^1([0,1])$ and $\dor(0)>0.$\\

Consider the rescaled function $r_\ve(R):=\ve^{-1}r(\ve R)$ for $0<\ve<1$ with the derivatives $\dor_\ve(R)=\dor(\ve R)$ and $\ddor_\ve(R)=\ve\ddor(\ve R)$ for all $R\in(0,1].$ Since $r$ solves the ODE (\ref{eq:3.15}) strongly in $(0,1)$ so does the rescaled version. Indeed note that 
\begin{align*}
d_{r_\ve}(R)=&\frac{Mr_\ve(R)\dor_\ve(R)}{R}=\frac{Mr(\ve R)\dor(\ve R)}{\ve R}=d_r(\ve R)\;\;\mbox{and}&\\
\dod_{r_\ve}(R)=&(d_r(\ve R))^\cdot=\ve \dod_r(\ve R).&
\end{align*}
Hence,
\begin{eqnarray*}
\dor_\ve(R)+R\ddor_\ve(R)+M\rho''(d_{r_\ve}(R))\dod_{r_\ve}(R)r_\ve(R)-\frac{M^2r_\ve(R)}{R}=\nonumber\\
=\dor(\ve R)+(\ve R)\ddor(\ve R)+M\rho''(d_r(\ve R))\ve \dod_r(\ve R)\ve^{-1}r(\ve R)-\frac{M^2r(\ve R)}{\ve R}=0.
\end{eqnarray*}
where the last equality holds since it agrees with the strong form of the ODE (\ref{eq:3.15}) evaluated at $\ve R.$\\
Now if $\ve\ra 0,$ then the rescaled function $r_\ve$ converges uniformly to the linear map $r_0(R):=aR$, i.e. $r_\ve\ra aR$ and $\dor_\ve\ra a$ uniformly in $[0,1].$ But then the weak form of the ODE of $r_\ve$ converges to the weak form of $r_0:$
\begin{eqnarray*}
0&=&\lim\limits_{\ve\ra0}\int\limits_{0}^1{\dot{g}\left(\dor_\ve R+M\rho'(d_{r_\ve})r_\ve\right)+g\left(\frac{M^2r_\ve}{R}+M\rho'(d_{r_\ve})\dor_\ve\right)\;dR}\nonumber\\
&=&\int\limits_{0}^1{\dot{g}\left(\dor_0(R)R+M\rho'(d_{r_0})r_0\right)+g\left(\frac{M^2r_0}{R}+M\rho'(d_{r_0})\dor_0\right)\;dR}.
\end{eqnarray*}
Therefore, $r_0$ is a weak solution to the ODE (\ref{eq:3.15}). Since $r_0$ is smooth it also needs to satisfy (\ref{eq:3.z2}), which is not true since plugging $r_0$ into (\ref{eq:3.z2}) yields
\[M^2a=a,\]
which is not satisfied since $M\ge2,\;a>0.$\\
 
2. Case: $\lim\limits_{R\rightarrow0} d(R)=0.$ \\
There are only two possible scenarios: Either $r\equiv 0$ in $[0,\d]$ or $r>0$ in $(0,\d).$

If $r\equiv0$ near zero then $\dor\equiv0$ in $(0,\d]$ and we can easily see $\dor(0)=0$ and $\dor\in C^1([0,1]).$ Moreover, this argument works for all derivatives of $r$ yielding $r\in C^\infty([0,1]).$\\

Assume instead that $r$ lifts-off immediately, i.e. there exists $\d>0$ s.t. $r>0$ on $(0,\d).$ Then lemma \ref{lem:3.z} holds, assume first $\doz\le0.$ Again, $\frac{r}{R}$ is monotone on $(0,\d)$ and the mean value theorem implies the existence of of the two sequences $R_j'\rightarrow0$ and $R_j\in(0,R_j')$ for any $j\in \N$ s.t. 
\[
\lim\limits_{j\rightarrow\infty} \dor(R_j)=\lim\limits_{j\rightarrow\infty} \frac{r(R_j')}{R_j'}=:m.
\]
Assume now $\dor>\frac{r}{R}$ on $(0,\d).$ Then \[0=\lim\limits_{R\rightarrow0} d(R)=\lim\limits_{j\rightarrow\infty}\frac{M\dor(R_j')r(R_j')}{R_j'}\ge\lim\limits_{j\rightarrow\infty}\frac{Mr^2(R_j')}{R_j'^2}=Mm^2.\] 
Hence, $m=0.$ Now, by $\dor>\frac{r}{R}$ on $(0,\d)$ we get $\dor(R_j)\ra0$ and $\frac{r(R_j)}{R_j}\ra0$ and therefore $\lim\limits_{R\rightarrow0} d(R_j)=0$ if $j\ra\infty.$ This yields
\[
\lim\limits_{j\rightarrow\infty} z(R_j)=\lim\limits_{j\rightarrow\infty}\left[\frac{\dor^2(R_j)}{2}+\frac{M^2r^2(R_j)}{2R_j^2}+f(d(R_j))\right]=0.\]
Since $z\ge0$ and by assumption $\doz\le0$ on $[0,\d)$ we have $z=0$ on $[0,\d).$ By the non-negativity of $d\ge0$ and $f(d)\ge0$ (following from $f(0)=0$ and recalling $(f(d))^\cd=\rho''(d)d\dod\ge0$) we finally know that $r=\dor=0$ on $[0,\d),$ 
contradicting the assumption that $r$ lifts off immediately. One can argue similarly in the case when $\dor<\frac{r}{R}$ on $(0,\d).$ \\

Finally, assume $\doz\ge0.$ Then $(\ref{eq:3.23a})$ holds, i.e. $\dor\sim\frac{r}{R}$ on $(0,\d).$ Since $d(R)\rightarrow0$ for $R\rightarrow0$ 
\[
\dor(0)=\lim\limits_{R\rightarrow0} \dor(R)=\lim\limits_{R\rightarrow0}\frac{r(R)}{R}=0
\label{eq:}
\]
and again $r\in C^1([0,1])$ with $\dor(0)=0.$\\\vspace{0.5cm}

We end this paragraph, by showing that the constructed maps $u=re_R(M\th)$ s.t. $r$ solves the BVP \eqref{eq:3.0.1} are stationary points of the functional \eqref{eq:1.1}.

\begin{lem} Let $u\in\A_r^M$ with $u=r e_R(M\th)$ s.t. $r$ solves the BVP \eqref{eq:3.0.1}. Then $u$ solves the ELE \eqref{eq:ELE1.1} weakly, i.e. 
\[\int\limits_B{\grad_\xi W(\grad u)\cd\grad \vp \;dx}=0 \;\mb{for all}\; \vp\in C_c^\infty(B,\R^2).\]
\end{lem}
{\bf{Proof:}}\\
Lemma 3.6 of \cite{BOP91MS} applies and shows that $u$ solves the ELE strongly in $B\sm\{0\},$ this can be reformulated in the following sense, $u$ satisfies
\[\int\limits_B{\grad_\xi W(\grad u)\cd\grad \vp \;dx}=0\]
for all $\vp\in C_c^\infty(B,\R^2)$ with $\vp\equiv0$ near the origin.\\ Now we can follow the strategy of Theorem 3.11 in \cite{BOP91MS} to upgrade this to arbitrary test functions. For this sake, take $\eta_\ve\in C^\infty(B)$ s.t. $\eta_\ve\equiv0$ on $B_\ve$ and $\eta_\ve\equiv1$ on $B\sm B_{2\ve},$ $0\le \eta_\ve \le1$ and there exists $c>0$ s.t. $|\grad \eta_\ve|\le \frac{c}{\ve}.$ Take an arbitrary test function $\psi\in C_c^\infty(B,\R^2)$ and set $\vp=\eta_\ve\psi$ then $\vp$ vanishes close to the origin. Hence,
\begin{align*}
0&=\int\limits_B{\grad_\xi W(\grad u)\cd\grad \vp \;dx}&\\
&=\int\limits_B{\grad_\xi W(\grad u)\cd\eta_\ve\grad \psi \;dx}+\int\limits_B{\grad_\xi W(\grad u)\cd(\grad\eta_\ve\ot \psi) \;dx}&
\end{align*}

Then the first integral converges:
\[
\lim\limits_{\ve\ra0}\int\limits_B{\grad_\xi W(\grad u)\cd\eta_\ve\grad \psi \;dx}=\int\limits_B{\grad_\xi W(\grad u)\cd\grad \psi \;dx}
\]
by dominated convergence. 
Since, $\grad u\in C^0$ there exists $C>0$ s.t. $\|\grad_\xi W(\grad u)\|_{C^0}\le C.$ Then
\[
\left|\int\limits_B{\grad_\xi W(\grad u)\cd(\grad\eta_\ve\ot \psi) \;dx}\right|\le C\|\psi\|_{C^0}\int\limits_B{|\grad\eta_\ve| \;dx}\le \frac{C}{\ve}\|\psi\|_{C^0}\La^2(B_{2\ve})
\]
Hence, the second integral vanishes for $\ve\ra0$ and 
\[\int\limits_B{\grad_\xi W(\grad u)\cd\grad \psi \;dx}=0\]
holds for every $\psi\in C_c^\infty(B,\R^2).$

\section{Advanced BOP-Theory}
\label{sec:3.3}
In the following, we want to investigate if these stationary points are even more regular. A first step in that direction is the next lemma. This has been observed for the BOP-Case by Yan and Bevan, see \cite[Lem 3.(i)]{BY07}  and \cite[Lem 1.(i)]{Y06}. 
\begin{lem} Let $M\in\N,$ $M\ge1$ and suppose $r$ solves the BVP \eqref{eq:3.0.1}. Then\\
(i) $\liminf\limits_{R\ra0}\ddor(R)\ge0$ and\\
(ii) $\lim_{R\ra0}\ddor(R)R=0.$
\end{lem}
\textbf{Proof:}\\
(i) Assume not, then $\liminf\limits_{R\ra0}\ddor(R)<0.$ But then there exits a small interval s.t. $\ddor<0$ on $(0,\d).$ By the mean value theorem it follows the existence of some $\xi\in(0,\d)$ , s.t.
 \[0>\ddor(\xi)=\frac{\dor(\d)}{\d}\ge0,\]
which is a contradiction.\\ 
(ii) Recall the ODE
\[
\ddor(R)R=\frac{M^2r}{R}-\dor-M\rho''(d)\dod r
\]
Then 
\[
0\le\liminf\limits_{R\ra0}{\ddor(R)R}\le\limsup\limits_{R\ra0}{\ddor(R)R}\le\limsup\limits_{R\ra0}{\frac{M^2r(R)}{R}}=0.
\]
This shows the second claim. \\\vspace{1cm}

The general boundary value problem we are considering is 
\begin{align}\left\{\begin{array}{ccc}
Lr=M\rho''(d)\dod r & {\mbox{in}}& (0,1),\\
r(0)=0,&r(1)=1,&
\end{array}\label{eq:3.3.1}
\right.\end{align}
where $L$ refers to the linear part of the considered ODE, i.e.
\[Lr(R):=\frac{M^2r}{R}-\dor-R\ddor,\]
which makes sense for all $R\in[0,1]$ by the previous lemma. \\

\subsection{Delayed lift-off solution for arbitrary $\rho$}
From now on, we will distinguish between the different shapes of the solutions to the BVP \eqref{eq:3.3.1}. For this recall that, according to Lemma \ref{lem:3.2.100}, $r$ is called immediate lift-off solution, if $r,\dor>0$ away from zero.  We will denote such solutions by $r_0,$ while we will call delayed lift-off solutions by $r_\d$ for $0<\d<1.$ As a reminder, $r_\d$ is a delayed lift-off solution, if there is $0<\d<1$ s.t. $r\equiv0$ on $[0,\d].$ The $\d$ indicates the point such that $r_\d\equiv0$ in $[0,\d],$ but also $r_\d(R)>0$ for all $R\in(\d,1].$\\
The first statement will be that an immediate lift-off solution $r_\d$ is zero up to $\d,$ and that it needs to solve the BVP \eqref{eq:3.3.3} below. This fact will be crucial for the uniqueness result, Lemma \ref{Lem:3.3.3}.
\begin{lem}Let $0<\ga<\infty.$ 
If for some $\d>0,$ $r_\d$ solves the BVP \eqref{eq:3.3.1} then \[r_{\d}(R)=\left\{\begin{array}{ccc}
0& {\mbox{in}}& (0,\d],\\
\tilde{r}_{\d,0}& {\mbox{in}}& (\d,1],
\end{array}
\right.\label{eq:3.3.2}\]
where $\tilde{r}_{\d,0}\in C^\infty((\d,1))\cap C^0([\d,1])$ is the unique solution of
\begin{align}\left\{\begin{array}{ccc}
Lr=M\rho''(d)\dod r & {\mbox{in}}& (\d,1),\\
r(\d)=0 , &r(1)=1.&
\end{array}
\right.\label{eq:3.3.3}\end{align}
\end{lem}
\textbf{Proof:}\\
Assume $r_\d$ is a solution to the BVP \eqref{eq:3.3.1}. Then, $r_\d=0$ on $[0,\d]$ by definition and $r_\d|_{[\d,1]}\in C^{\infty}([\d,1])$ is the unique solution to 
\[\left\{\begin{array}{ccc}
Lr=M\rho''(d)\dod r & {\mbox{in}}& (\d,1),\\
r^{(k)}(\d)=0\;\mb{for all}\;k\in\N , &r(1)=1&
\end{array}
\right.\]
Hence, $r_\d|_{[\d,1]}$ solves \eqref{eq:3.3.3}.\\\vspace{0.5cm}

Furthermore, we demonstrate that such a solution, if it exists, needs to be unique. 
\begin{lem}Let $0<\ga<\infty.$ Assume there exists a solution to the BVP \eqref{eq:3.3.1} of the form $r_\d$, $\d>0.$ Then there exists a unique $\d>0$ in the sense that
\[r_{\d}(R)=\left\{\begin{array}{ccc}
0& {\mbox{in}}& (0,\d],\\
\tilde{r}_{\d,0}& {\mbox{in}}& (\d,1]
\end{array}
\right.\]
and $\tilde{r}_{\d,0}$ lifts off immediately, i.e. $\tilde{r}_{\d,0}(R)>0$ for all $R\in(\d,1].$
\label{Lem:3.3.3}
\end{lem}
\textbf{Proof:}\\
We can always choose $\d$ to be maximal in the sense that $r_\d\equiv0$ on $[0,\d]$ and $\tilde{r}_{\d,0}$ lifts off immediately. For the uniqueness assume that $r_\d$ and $r_{\d'}$ are two solutions to \eqref{eq:3.3.1} for $0<\d\le\d'<1$. Then both need to satisfy  \[\left\{\begin{array}{ccc}
Lr=M\rho''(d)\dod r & {\mbox{in}}& (\d,1),\\
r(\d)=0\;, &r(1)=1.&
\end{array}
\right.\]
The Picard-Lindelöf Theorem\footnote{see, \cite{GT12} Theorem 2.5, Corollary 2.6.} yields, $r_\d\equiv r_{\d'}$ on $[\d,1]$, trivially on the complete interval $[0,1].$ Hence, $\d=\d'$ by maximality.\\ 

\textbf{Conclusion:}
For arbitrary behaviour of $\rho$ and $0<\ga<\infty$ the above discussion shows that (smooth) delayed lift-off solutions to the BVP \eqref{eq:3.3.1} either do not exist or there is at most one. The latter lemma shows that the representation and $\d$ may vary, but the solution remains the same.

\subsection{Delayed lift-off $\rho$}

Lets assume $\rho$ itself is a delayed lift-off function, i.e. there exists $\tilde{s}>0$ s.t.  $\rho\equiv0$ on $[0,\tilde{s}]$. Since, $d(0)=0$ and $d$ continuous this implies that there exists a small interval $[0,\d),$ s.t. $R\mt\rho(d(R))\equiv0.$ So the $\rho-$term of the considered functional \eqref{eq:1.1} vanishes for a short amount of time. Hence, the functional reduces to the Dirichlet energy, at least close to the origin. We know, that stationary points to the Dirichlet energy are harmonic functions, which are smooth. This gives hope, that in case of a delayed lift-off $\rho,$ all solutions are smooth close to the origin.\\
The following result states that this is exactly true. We only consider immediate lift-off solutions $r_0,$ since we already know, from the previous discussion, that delayed solutions are smooth.
     
\begin{lem}Let $0<\ga<\infty$ and assume there exists $\tilde{s}>0$ s.t. $\rho(s)=0$ for all $s\in[0,\tilde{s}]$ and $\rho(s)>0$ for all $s>\tilde{s}.$ Further suppose there exists a solution of the BVP \eqref{eq:3.3.1} of the form $r_0.$\\
If $r_0$ solves \eqref{eq:3.3.1} then there exists a unique $\d=\d(\tilde{s},r_0,\dor_0)>0$ and some $0<a<1$ s.t. 
\begin{align}r_{0}(R)=\left\{\begin{array}{ccc}
a\left(\frac{R}{\d}\right)^M& {\mbox{in}}& (0,\d],\\
\tilde{r}_{\d,a}& {\mbox{in}}& (\d,1]
\end{array}\label{eq:3.3.4}
\right.\end{align}
and $\tilde{r}_{\d,a}$ the unique smooth solution of
\begin{align}\left\{\begin{array}{ccc}
Lr=M\rho''(d)\dod r & {\mbox{in}}& (\d,1),\\
r(\d)=a , &r(1)=1.&
\label{eq:3.3.5}
\end{array}
\right.\end{align}
and $\tilde{r}_{\d,a}$ is not in the kernel of $L$ for at least a short period of time, i.e. there exists $\ve>0$ s.t. $L\tilde{r}_{\d,a}(R)>0$ for all $R\in(\d,\d+\ve].$ Moreover, $(a\left(\frac{R}{\d}\right)^M)^{(k)}(\d)=\tilde{r}_{\d,a}^{(k)}(\d)$ for all $k\in \N.$
\end{lem}
\textbf{Proof:}\\
Since $d(0)=0$ and $d\in C([0,1])$ and the delay of $\rho$ there exists a unique $\d=\d(\tilde{s},r_0,\dor_0)>0$ s.t. $R\mapsto\rho(d(R))\equiv0$ for all $R\in[0,\d]$ and $\rho(d(R))>0$ for all $R>\d.$ Then $d(R)>0$ for $\d<R\le1$ and there exists an $\ve>0$ s.t. $\dod(R)>0$ for $R\in(\d,\d+\ve].$ Hence, $L\tilde{r}_{\d,a}(R)=M\rho''(d)\dod r>0$ for $R\in(\d,\d+\ve].$
Then $r_0$ needs to solve the following ODE
\[\left\{\begin{array}{ccc}
Lr=0& {\mbox{in}}& (0,\d),\\
r(0)=0,&r(\d)=a&
\end{array}
\right.\]
for some $0<a<1.$ Indeed $a$ can not exceed $1$ $(a>1)$ since $\dor\ge0$ and $r(1)=1.$ If $a=1$ then $r\equiv1$ in $[\d,1]$ implying $d\equiv0$ in $[0,1],$ a contradiction. $a=0$ is excluded by the assumption that the solution $r_0$ is an immediate lift-off function.\\
Then $r_0|_{[\d,1]}$ solves \eqref{eq:3.3.5} uniquely and $r_0$ takes the form \eqref{eq:3.3.4}.
Since $r_0\in C^\infty((0,1])$ all derivatives of $r_0$ need to agree at $\d,$ i.e. $(a\left(\frac{R}{\d}\right)^M)^{(k)}(\d)=\tilde{r}_{\d,a}^{(k)}(\d)$ for all $k\in \N.$\\

\begin{re} By construction $\d$ depends not only on $\tilde{s}$ but also on $r_0,\dor_0$. For every $r_0$ the $\d$ may vary, destroying any chance for a uniqueness result similar to Lemma \ref{Lem:3.3.3}.
\end{re}

\textbf{Conclusion:}
We can not guarantee the existence of solutions to the ODE of the form $r_0.$ But if they exist, they have to be smooth.\\
Moreover, combining this with our knowledge of delayed-lift off solutions, we are able to conclude that in case of a delayed $\rho$ all stationary points in the class $\A_r^M$ are smooth.\\

\subsection{Immediate lift-off $\rho$}

If $\rho$ is an immediate lift-off function ($\rho(s)>0$ for all $s>0$), then we are not (yet) able to show, that $r_0$ solutions need to be any smoother then $C^1.$\\
However, in the next statement we give a necessary condition, \eqref{eq:3.3.200}, which needs to be satisfied if $r_0$ is of class $C^{1,\al}$ for some $\al\in(0,1).$ As a consequence, this fixes the limit of the quantity $\frac{R\dor_0}{r_0},$ if $R$ tends to $0,$ which was discussed in Lemma \ref{lem:3.7}.\\

\begin{lem} (Necessary condition)\label{lem:3.3.200} Let $0<\ga<\infty$ and $\rho(d)>0$ for all $d>0.$ Assume there exists a solution to the BVP \eqref{eq:3.3.1} of the form $r_0.$\\

If  $r_0\in C^{1,\al}([0,1])$ for some $\al\in(0,1)$ then \[D_M:=\lim\limits_{R\ra0}\frac{R\dor_0}{r_0}=M\] and there exists $\d>0$ s.t.
\begin{equation}M^2-C_{\al,M}\rho''(d_{r_0})R^{2\al}\le\frac{R\dor_0}{r_0}+\frac{R^2\ddor_0}{r_0}< M^2 \;\mb{for all}\; R\in (0,\d).
\label{eq:3.3.200}\end{equation}
\end{lem}
\textbf{Proof:}\\
In the following we will suppress $r_0$ in favour of $r.$\\ 
First we note that $\rho(d)>0$ for all $d>0$ implies that there exists a $\d>0$ s.t. $R\mapsto\rho''(d(R))>0$ for all  $R\in(0,\d).$ This yields
\begin{equation}0<M\rho''(d)=\frac{Lr}{\dod r},\;\mb{on}\;(0,\d). \label{eq:3.3.98}\end{equation} 
Since, $r,\dor>0$ in $(0,1]$ we can infer $\dod>0$ and $Lr>0$ on $(0,\d).$
Since, $r\in C^{1,\al}([0,1])$ there exists $c_\al>0$ and $\d>0$ s.t.  
\begin{equation}|\ddor(R)|\le c_\al R^{\al-1},\;\mb{for all}\;R\in(0,\d). \label{eq:3.3.99}\end{equation}
Assume not. Then for all $\d>0,$ $c>0$ there exists $R\in (0,\d)$ s.t.
\begin{equation}[|\ddor(R)|> c R^{\al-1}. \label{eq:3.3.100}\end{equation}
Fix $\d>0.$ Then for all $c>0$ take $R_c\in (0,\d)$ s.t. the latter inequality holds. By continuity of $\ddor$ in $(0,1],$ there exists an $\ve>0$ s.t. \eqref{eq:3.3.100} holds even for all $R'\in(R_c-\ve,R_c+\ve).$ By integration we get
\[|\dor(R_c+\ve)-\dor(R_c)|=\int\limits_{R_c}^{R_c+\ve}{|\ddor(R')|\;dR'}>\frac{c}{\al}((R_c+\ve)^\al-R_c^\al).\]
Hence, for all $c>0$ we can find $R_c\in (0,\d)$ s.t. $\dor$ is not Hölder continuous at $R_c$ with Hölder constant $\frac{c}{\al}.$ Since $c>0$ is arbitrary this contradicts $r\in C^{1,\al}([0,1]).$\\ 

Then \eqref{eq:3.3.99} implies
\[0<\dod r=M\left(\frac{r\dor^2}{R}+\frac{r^2\ddor}{R}-\frac{r^2\dor}{R^2}\right)\le M\frac{r}{R}\left(\dor^2+r|\ddor|\right)\le c_{\al}^2\al^{-2}M\frac{r}{R}R^{2\al}\]
for all $R\in(0,\d).$ Together with \eqref{eq:3.3.98} we get
\[0<\frac{R}{r}Lr\le C_{\al,M}\rho''(d)R^{2(k-1+\al)} \;\mb{on}\;(0,\d),\]
where $C_{\al,M}:=M^2c_{\al}^2\al^{-2}>0.$
Using the explicit form of $Lr>0$ yields the claimed inequalities
\[M^2-C_{\al,M}\rho''(d)R^{2\al}\le\frac{R\dor}{r}+\frac{R^2\ddor}{r}< M^2 \;\mb{near}\; 0.\]
Taking the limit $R\ra0$ yields 
\[\lim\limits_{R\ra0}\left(\frac{R\dor}{r}+\frac{R^2\ddor}{r}\right)=M^2.\]
With the notation $D_M:=\lim\limits_{R\ra0}\frac{R\dor}{r}$ and $E_M:=\lim\limits_{R\ra0}\frac{R\ddor}{\dor}$ we get by L'H\^{o}pital's rule
\[D_M=\lim\limits_{R\ra0}\frac{R\dor}{r}=\lim\limits_{R\ra0}\frac{R\ddor+\dor}{\dor}=1+E_M\] and
\begin{align*}M^2&=\lim\limits_{R\ra0}\left(\frac{R\dor}{r}+\frac{R^2\ddor}{r}\right)=D_M+\lim\limits_{R\ra0}\left(\frac{R\dor}{r}\frac{R\ddor}{\dor}\right)&\\
&=D_M+D_ME_M=D_M+D_M(D_M-1)=D_M^2,&
\end{align*}
hence $D_M=M.$

\begin{re}
Improving the lower bound, in \eqref{eq:3.3.200}, up to the point where it matches the upper one, would show that $r$ can't be of class $C^{1,\al}.$ However, we don't know how to prove it and if it is even true.
\end{re}

\textbf{Conclusion:}
In this case there are two possibilities: There is at most one smooth delayed lift-off solution $r_\d$ for some $\d>0.$\\
Moreover, there could be an immediate lift-off solution of the form $r_0.$ But $r_0$ could have regularity anywhere between $C^1\ldots C^{1,\al}\ldots C^2\ldots C^\infty.$\\

\clearpage{\pagestyle{empty}\cleardoublepage} 
\chapter{Uniqueness in finite elasticity}
As mentioned in the introduction, in this chapter, we are interested in, what kind of circumstances are necessary to ensure that a minimizer of some functional is actually unique. We start by describing a buckling situation first discussed in \cite{BeDe20} and discuss a first uniqueness result in that setting.

\section{Uniqueness of minimizers to the `buckling' functional for small parameter}
\label{sec:4.1}
Up to now we do not claim any originality, but we need to recall all this knowledge to set our result in context and to make this section fairly self contained. In the next section we show that for $1<\ve\le\sqrt{2},$ $v=\Id$ is a global minimizer in the full class $\A^c$ and it is even the unique one if $\ve\in(1,\sqrt{2})$. This is a stronger result, in the above paper it has been shown that  $v_\ve=\Id$ is the only stationary point in the class $\mathcal{C}^c$ in $(1,\sqrt{2}).$\\

\subsection{Setting}
\textbf{Functional:} For every $\ve\ge1,$ we define the integrand $\W_\ve$ as
\[\W_\ve(x,\xi)=\frac{1}{\ve}|\xi^T \hat{x}|^2 + \ve|(\adj\; \xi) \hat{x}|^2,\;\mb{for all}\;x\in B\sm\{0\},\; \xi\in\R^{2\ti2},\]
and the functional $\bb{D}_\ve$ by 
\begin{align}\label{de} \bb{D}_{\ve}(v) = \int_{B}\W_\ve(x,\grad v) \, dx \end{align}
for all functions $v\in W^{1,2}(B,\R^2).$ For every $x\in B\sm\{0\}$ for the 2d-polar coordinates we make use of the following notations $\{\hat{x},\hat{x}^\perp\}=\{e_R,e_\th\}=\{e_R(\th),e_\th(\th)\}:=\{(\cos(\th),\sin(\th)),(-\sin(\th),\cos(\th))\}.$ \\

We recall the following facts about this functional:\\
(i) For any $\ve\in[1,\infty),$ $\bb{D}_{\ve}$ is uniformly convex with $\bb{D}_{\ve}\ge\frac{1}{\ve}\bb{D}_{1}.$\\
(ii) To explain the origin of the integrand, we start by noting that it holds \[|\xi|^2=|\xi^T \hat{x}|^2 + |(\adj\; \xi) \hat{x}|^2,\;\mb{for all}\;x\in B\sm\{0\},\; \xi\in\R^{2\ti2}.\] This explains that in case of $\ve=1,$ $\bb{D}_{1}$ agrees with the Dirichlet energy. So for a general $\ve\in(1,\infty),$ $\W_{\ve}$ can be thought of as a weighted perturbation of the Dirichlet functional. Lastly, while in case of $\ve=1$ the integrand does not explicitly depend on $x$ (this can be seen from the above relation, that the dependence on $x$ is clearly not present on the LHS), but this relation gets perturbed if $\ve>1$ hence $\W_{\ve}$ depends explicitly on $x,$ smoothly in $B\sm\{0\}$ with a discontinuity at the origin.\\

\textbf{Admissible sets:} Recall the definition of the unconstrained set
\[\A_{\footnotesize{\Id}}=\{v\in W^{1,2}(B,\R^2)| v_{|\p B}=\Id\},\]
and define the set of measure preserving maps as 
\[K=\{v\in W^{1,2}(B,\R^2)|\det\grad v=1 \;\mb{a.e.}\}.\]
We will concentrate on the set of admissible maps
\[\A_{\footnotesize{\Id}}^c=K\cap \A_{\footnotesize{\Id}},\]
suppressing the boundary condition as usual.\\

This situation has first been studied by J. Bevan and J. Deane in \cite{BeDe20}. They focused on stationary points in the subclass of so called twist maps defined by
\[
\mathcal{C}=\{v\in\A|v=Re_R(\th+k(R))\; and\; k\in \tilde{\mathcal{C}}\},
\]

where
\[
\tilde{\mathcal{C}}=\{k\in W^{1,1}((0, 1))|\; R^{1/2}k(R)\in L^2(0, 1),R^{3/2}k'(R) \in L^2(0, 1), k(1) = 0\}.
\]
As before we define the constrained version by
\[\mathcal{C}^c=K\cap \mathcal{C}.\]
Notice that $\Id\in \mathcal{C}$ and $\mathcal{C}^c$ by setting, $\Id=Re_R(\th+k)$ with $k\equiv0$ for any $\ve>1.$\\

\textbf{ELE and pressure:} Before we collect their main results in the next statement we need some preliminaries. Let $u$ be a stationary point to $\bb{D}_{\ve}$ in $\A^c.$ Then we can find a corresponding pressure $\la$ s.t. $(u,\la)$ satisfy the following ELE
\begin{equation}\int\limits_B{\grad_\xi\W_\ve(x,\grad u)\cd\grad\eta\;dx} =-2\int\limits_B{\la(x)\cof\grad u\cd\grad\eta \, dx}, \;\mb{for all} \;\eta\in C_c^\infty(B,\R^2).\label{eq:Uni.Comparison.0}\end{equation}
It is crucial to note that our pressure differs from the one in \cite{BeDe20}. The relation to their pressure function $P$ is given by $\la(x)=-\frac{P(x)}{2}.$ \\
If $u\in\mathcal{C}^c$ is a stationary point with pressure $\la$ it is shown that $\la$ only depends on $R,$ i.e. $\la(x)=\la(R),$ and with $p_\ve=\ve-\frac{1}{\ve}$ the ELE in $\tilde{\mathcal{C}}$ becomes 
\begin{equation}\la'(R)R\left(\frac{1}{\ve}+p_\ve\cos^2k(R)\right)=p_\ve\left(\frac{\sin^2k(R)}{\ve}-\ve\cos^2k(R)\right)-R^2k'^2(R)\label{eq:Uni.Comparison.1}\end{equation}
a.e. $R\in(0,1).$\\

Now to the statement:

\begin{thm}[Bevan, Deane; 2020]\footnote{The fact, that $\Id$ is a stationary point, can be deduced by realising that $k\equiv0$ is always a solution to the system (2.22)-(2.23) of \cite{BeDe20} independent of $\ve.$ The equation for the pressure steams mainly from equation (3.46). Together, this is Point $(1).$ Point $(2)$ can be found in \cite[Prop 2.16.(i)]{BeDe20}, Point $(3)$ corresponds to Theorem 1.1 and Prop 2.16.(ii) in \cite{BeDe20}.}
\begin{enumerate}
\item For any $\ve\in(1,\infty),$  $\Id$ is a stationary point of $\bb{D}_{\ve}$ and the corresponding pressure $\la_\ve$ satisfies $\la_\ve'(R)R=-p_\ve,$ for a.e. $R\in(0,1).$
\item  For any $\ve\in(1,\sqrt{2}],$ $\Id$ is the only stationary point of $\bb{D}_{\ve}$ in the class $\mathcal{C}^c.$
\item For any $\ve\in(\sqrt{2},\infty],$ there are infinitely many pairs of $C^1(B,\R^2)$ twist maps $v_{\pm,j}^\ve\in\mathcal{C}^c,$ where $k_j\in C^\infty,$ such that for each $j:$
\begin{enumerate}
\item $\bb{D}_{\ve}(v_{+,j}^\ve)=\bb{D}_{\ve}(v_{-,j}^\ve),$
\item $v_{\pm,j}^\ve$ is a stationary point of $\bb{D}_{\ve}$ in $\A^c,$
\item $\grad v_{\pm,j}^\ve(0)=\pm J,$ where $J=\begin{pmatrix}0&-1\\1&0\end{pmatrix}.$
\item $\bb{D}_{\ve}(v_{\pm,j}^\ve)<\bb{D}_{\ve}(v_{\pm,j+1}^\ve),$ for all $j\in\N.$
\end{enumerate}
\end{enumerate}
\label{Thm:Uni.Glob.1.1}
\end{thm}

Here we give a small interpretation: In the regime $(\sqrt{2},\infty],$ there are pairs of equal-energy minimizers in the class $\mathcal{C}^c$ and, moreover, they are stationary points in $\A^c.$ This already provides examples of non-unique minimizers in a subclass and non-unique stationary points wrt. the full class, contributing partial answers to a famous problem raised by Ball in \cite{BallOP}. This strongly suggests (even though they are not yet able to make this rigorous) that these pairs might even be global minimizers wrt. the full class $\A^c,$ which, if true, would generate a non-unique situation for global minimizers, giving a complete answer to the problem of Ball in incompressible elasticity.\\

\subsection{The identity as a global minimizer to $\bb{D}_{\ve}$ in the regime $\ve\in(1,\sqrt{2}].$} 
The main result of this section is that, in the regime $1<\ve<\sqrt{2},$ the identity is the unique global minimizer of $\bb{D}_{\ve}$ in the full class $\A^c.$ In the threshold case $\ve=\sqrt{2}$ the identity is still a global minimizer, however the uniqueness might be lost.\\

\begin{thm} In the regime $1<\ve<\sqrt{2}$ the identity is the unique global minimizer of $\bb{D}_\ve$ in the class $\A^c.$ Moreover, in the limiting case $\ve=\sqrt{2}$ the identity is a (not necessarily unique) global minimizer. 
\label{lem:Glob.1.5}
\end{thm}

\textbf{Proof.}\\
For any $\ve\in(1,\sqrt{2}),$ let $u\in \A^c$ be arbitrary, and set $\eta:=u-\Id\in W_0^{1,2}(B,\R^2).$\\

In order to show uniqueness, we start with the expansion 
\[\bb{D}_\ve(u)=\bb{D}_\ve(\Id)+\bb{D}_\ve(\eta)+\bb{H}_\ve(\Id,\eta),\]
where
\[\bb{H}_\ve(\Id,\eta):=2\int\limits_B{\frac{1}{\ve}e_r\cd \grad\eta^T e_r+\ve e_\th\cd \grad\eta^T e_\th\;dx}\]
denotes the mixed terms. Now $\Id$ is a global minimizer if 
\[\bb{D}_\ve(\eta)+\bb{H}_\ve(\Id,\eta)\ge0\]
holds for all $\eta=u-\Id\in W_{0}^{1,2}(B,\R^2).$ In view of this, we rewrite $\bb{H}_\ve(\Id,\eta)$ as follows. Noting that $\bb{H}_\ve(\Id,\eta)=\grad_\xi\W_\ve(x,\grad \eta)$ agrees with part of the ELE except the term, originating from the Lagrange-multiplier, we get
\[\bb{H}_\ve(\Id,\eta)=-2\int\limits_B{\la(R)\Id\cd \grad \eta\;dx}.\]
Expanding the Jacobian of $\eta$ and exploiting the fact that both $u$ and $\Id$ satisfy $\det\grad u=\det\Id =1$ a.e. yields
\[\det\grad\eta=-\Id\cd \grad \eta \;\mb{a.e.}.\]
By entering this in the equation above we get 
\begin{equation}\bb{H}_\ve(\Id,\eta)=2\int\limits_B{\la(R)\det \grad \eta\;dx}.\label{eq:Uni.Pf.Uni.1}\end{equation}
For the further argument we must ensure that the latter expression is independent of the $0-$mode of $\eta.$ The Fourier-representation for any $\eta\in C^\infty(B,\R^2)$ (For members of Sobolev- spaces one might approximate) is given by
\begin{align*}
\eta(x)=\sum\limits_{j\ge0}\eta^{(j)}(x),\;\mb{where}\;\eta^{(0)}(x)=\frac{1}{2}A_0(R), \;A_0(R)=\frac{1}{2\pi}\int\limits_{0}^{2\pi}{\eta(R,\th)\;d\th}\end{align*} 
and for any $ j\ge1$ we have
\begin{align*}
 \eta^{(j)}(x)=A_j(R)\cos(j\th)+B_j(R)\sin(j\th),
 \end{align*}
  where
\begin{align*}
A_j(R)=\frac{1}{2\pi}\int\limits_{0}^{2\pi}{\eta(R,\th)\cos(j\th)\;d\th}\;\mb{and}\;
B_j(R)=\frac{1}{2\pi}\int\limits_{0}^{2\pi}{\eta(R,\th)\sin(j\th)\;d\th}.
\end{align*} 
\footnote{Note, that this is true for an arbitrary $\eta\in C^\infty(B,\R^2).$ It is not necessary that $\eta$ separates the variables $R$ and $\th.$}
Now Lemma 3.2.(ii)+(iii) of \cite{JB14} indeed guarantee that $\eta^{(0)}$ does not contribute to the expression in \eqref{eq:Uni.Pf.Uni.1} hence
\[\bb{H}_\ve(\Id,\eta)=2\int\limits_B{\la(R)\det \grad \tilde{\eta}\;dx},\]
where $\tilde{\eta}=\eta-\eta^{(0)}.$ Applying Lemma 3.2.(iv) of \cite{JB14} yields
\[\bb{H}_\ve(\Id,\eta)=\int\limits_B{\la'(R)\tilde{\eta}\cd J \tilde{\eta}_{,\th}\;dRd\th},\]
where $J=\begin{pmatrix}0&-1\\1&0\end{pmatrix}.$
Making use of the fact that by Theorem \ref{Thm:Uni.Glob.1.1}.(a) the corresponding pressure to $\Id$ satisfies $\la'(R)R=-p_\ve$ implies
\begin{align*}
\bb{H}_\ve(\Id,\eta)\ge&-p_\ve \int\limits_B{R^{-2}|\tilde{\eta}||\tilde{\eta}_{,\th}|\;dx}.&
\end{align*}
Furthermore, we get 
\begin{align*}
\bb{H}_\ve(\Id,\eta)\ge&-p_\ve \left(\int\limits_B{R^{-2}|\tilde{\eta}|^2\;dx}\right)^{1/2}\left(\int\limits_B{R^{-2}|\tilde{\eta}_{,\th}|^2\;dx}\right)^{1/2}&\\
\ge&-p_\ve \int\limits_B{R^{-2}|\tilde{\eta}_{,\th}|^2\;dx}&\\
\ge&-p_\ve \int\limits_B{|\grad \tilde{\eta}|^2\;dx}&\\
\ge&-p_\ve \int\limits_B{|\grad \eta|^2\;dx},&
\end{align*}
The $\sim$ could be dropped since $ \int\limits_B{|\grad \tilde{\eta}|^2\;dx}\le \int\limits_B{|\grad \eta|^2\;dx}.$ Moreover, we made use of the Cauchy-Schwarz inequality and the following estimate
\begin{equation}\int\limits_B{R^{-2}|\tilde{\eta}_{,\th}|^2\;dx}\ge\int\limits_B{R^{-2}|\tilde{\eta}|^2\;dx}.\label{eq:Uni.Buckling.1}\end{equation}
The latter one can be seen, by computing the LHS of the latter inequality mode-by-mode. Indeed, for any $j\in\N\sm\{0\}$ we have
\[\int\limits_B{R^{-2}|\eta_{,\th}^{(j)}|^2\;dx}=j^2\int\limits_B{R^{-2}|\eta^{(j)}|^2\;dx}.\]
By summing over all modes $j\ge1,$ we obtain the claimed inequality.\\

As previously mentioned it holds $\bb{D}_\ve\ge \frac{1}{\ve}\bb{D}_1.$ Together with the above derivation we reach the condition 
\begin{equation}\bb{D}_\ve(\eta)+\bb{H}_\ve(\Id,\eta)\ge\left(\frac{2}{\ve}-\ve\right)\bb{D}_1(\eta)\ge0,
\label{eq:Uni.1.1}
\end{equation}
which needs to be satisfied for any $\eta\in W_{0}^{1,2}(B,\R^2).$ The prefactor $\frac{2}{\ve}-\ve$ is non-negative on $(1,\sqrt{2}],$ yielding that $\Id$ is indeed a global minimizer in the class $\A^c$. Moreover, in the open interval $(1,\sqrt{2})$ the prefactor is strictly positive, implying, since $\bb{D}_1(\eta)=0$ iff $\eta\equiv0$ (in $W_0^{1,2}$) holds, that $\Id$ is the unique global minimizer, completing the proof.

\subsection{Conclusion to the unconstrained minimization problem}
We can make use of our knowledge of the constrained problem to conclude that the identity can neither be a stationary point nor a minimizer of the unconstrained one. Moreover, we get an upper bound on the minimal energy. The remark below discusses the relevance of this result.\\

\begin{cor}(Unconstrained problem) Let $1<\ve<\infty.$ Then the identity is not a stationary point/minimizer of the unconstrained minimization problem $\bb{D}_\ve$  in either of the classes $\mathcal{C}$ and $\A.$ Moreover, it holds that
\[\min\limits_{\bar{v}\in\mathcal{C},\mathcal{A}}{\bb{D}_\ve(\bar{v})}<\bb{D}_\ve(\Id)=\pi\left(\ve+\frac{1}{\ve}\right).\]
\end{cor}
\textbf{Proof:}\\  
Clearly, $\min\limits_{\bar{v}\in\mathcal{C},\mathcal{A}}{\bb{D}_\ve(\bar{v})}\le\bb{D}_\ve(\Id)$ is always true, since $\Id\in \mathcal{C},\mathcal{A}$ and the minimal energy is always smaller or equal than any of its competing energies. However, we claim $\le$ can be replaced by $<,$ which is a minor improvement.\\

Now recall that $\Id$ is a stationary point of the constrained problem in the class $\mathcal{C}^c.$ From Theorem \ref{Thm:Uni.Glob.1.1}.(a) we know that $\Id$ satisfies \eqref{eq:Uni.Comparison.1} with $R\la'(R)=-p_\ve.$ For the sake of a contradiction, assume now that $\Id$ is a stationary point to the unconstrained problem in the class $\mathcal{C}.$ Then it satisfies again \eqref{eq:Uni.Comparison.1} but with the LHS being equal to $0.$ This implies $p_\ve=0$ but we know $p_\ve>0$ for all $1<\ve<\infty,$ hence, we have reached a contradiction.\\

\begin{re}The above needs to be viewed in sharp contrast to the following well-known result:
Let $\Om\ss\R^n$ and assume $f:\R^{n\ti m}\ra\R$ only depends on the gradient, is (strictly-) quasiconvex and satisfies $v_{|\p \Om}=Ax+b$ for some constant matrix $A\in \R^{n\ti m}$ and some $b\in\R^n$ on the boundary. Then we have
\[\int\limits_{\Om} {f(A)\;dx}\le \int\limits_{\Om} {f(A+\grad\psi(x))\;dx}\;\mb{for all}\;\psi\in W_0^{1,\infty}(\Om,\R^{n}).\]

This shows, in particular, that $v\equiv Ax+b$ is a global minimizer and it is unique if $f$ is strictly quasiconvex.\\

This theorem, if applicable in our case, would imply that the identity is the unique global minimizer of the unconstrained problem. However, our integrand $W(x,\xi)$ depends additionally on $x$ in a discontinuous way, and hence the above result does not apply. This shows that if the integrand depends on $x$ or $v,$ in general one can no longer expect the above statement to be true. This can also be seen from the quasiconvexity condition for functions of the form $f(x,s,\xi),$ under suitable assumptions on the dependence of $x,s$ and $\xi$, which reads:\\
\[\int\limits_{\Om} {f(x,s, \xi)\;dy}\le \int\limits_{\Om} {f(x,s, \xi+\grad\psi(y))\;dy}\;\mb{for all}\;\psi\in W_0^{1,\infty}(\Om,\R^{n}),\]
see for example \cite[Eq.(0.2), p.126]{AF84}. However, the above is not of the correct form. In order to get an analogous result, for this more general integrands, one needs that $x$ is integrated over, but that is not the case.
\label{re:Uni.1.10}
\end{re}

\section{Uniqueness in incompressible elasticity}
\label{sec:4.2}
\subsection{Uniqueness in small pressure situations}
In the previous paragraph we discussed the question: Given a functional and the set of admissible maps, the constraint and the identity on the boundary, is the problem uniquely minimized by the identity?\\
Here we will take a different stand: Given the set of admissible maps, the constraint, and a specified map on the boundary, if we want the affine map to be the (unique) minimizer, what should the functional look like, and what criteria must it satisfy? We even answer this question in slightly more generality, including the case of N-covering maps.\\

Consider the following problem: For all $u\in W^{1,2}(B,\R^2)$ define the energy by
\begin{equation}E(u)=\int\limits_B{f(x,\grad u)\;dx},  \label{eq:USPS:1.0}\end{equation}
where the integrand is of the quadratic form
\[f(x,\xi)=M(x)\xi\cd\xi, \;\mb{for any} \;x\in B, \;\xi\in\R^{2\ti2},\]
where $M\in L^\infty(B,\R^{16})$ and $M$ symmetric.\footnote{$M\in L^\infty$ is a minimal requirement in order to guarantee $E(u)<\infty$ for any $u\in W^{1,2}.$} Furthermore, we want $f$ to be uniformly convex, hence we require that there exists a constant $\nu>0$ s.t. 
\begin{equation}\nu\Id\le M(x) \;\mb{for a.e.}\; x\in B.\label{eq:USPS:1.0.1}\end{equation}

Now define the constrained minimization problem as follows

\begin{equation}
\min\limits_{u\in\A} E(u)
\label{eq:USPS:1.1}
\end{equation}
with the set of admissible maps being given by
\begin{align*}
\A^c=K\cap\{u\in W^{1,2}(B,\R^2): \; u_{|{\p B}}=g(x)+b, b\in\R^2, g \in W^{1,2}(\p B,\R^2)\},
\end{align*}
where $(R,\th)$ refer to the standard plane polar coordinates.\\

\begin{re}\emph{Note, that the above includes the following specific situations:\\
a) $g(\th)=Ax$ and any constant matrix in $A\in SL^+(2,\R)$ and\footnote{For a discussion why, in case of affine boundary data, it is enough to consider matrices in the class $SL^+(2,\R),$ see Remark \ref{re: Uni.2.1}.}\\
b) $g(\th)=\frac{1}{\sqrt{N}}e_R(N\th)$ the N-covering map on the boundary. }\\
\end{re}
\begin{re}\emph{Note, that in incompressible elasticity, for quadratic functionals in $2\ti2$ dimensions, convexity and rank-1 convexity agree.
It is well known, that in $2\ti2d$ for quadratic functionals, poly-, quasi- and rank-1 convexity agree, see \cite[Thm 5.25]{D08}.
 Moreover, we know, that because of the $2\ti2d-$incompressibility, convexity and polyconvexity agree. Hence, convexity and rank-1 convexity do agree. In this sense, our result is optimal, there is no larger class of semiconvexities to consider, except if one either allows for more general functionals or goes to higher dimensions.}
\end{re}

Before stating the main theorem, we need the following two technical lemmas:\\

\begin{lem} Let $\la\in W^{1,1}(B)$ and assume that $\|R\grad\la\|_{L^\infty(B,\R^2,\frac{dx}{R})}<\infty.$ \\

Then the following statements are true:\\

i) $\la\in BMO(B).$\\
ii) If $\vp_n\ra\vp\in W^{1,2}(B,\R^2)$ then $\int\limits_B{\la(x)\det\grad\vp_n\;dx}\ra\int\limits_B{\la(x)\det\grad\vp\;dx}.$\\
iii) It holds $\int\limits_B{|\grad\vp|^2\;dx}=\sum\limits_{j\ge0}\int\limits_B{|\grad\vp^{(j)}|^2\;dx}$ for any $\vp\in W^{1,2}(B,\R^2).$\\ 
iv) $\det\grad\vp^{(0)}=0$ for any $\vp\in W^{1,2}(B,\R^2).$\\
v)  $\int\limits_B{\la(x)\det\grad\vp\;dx}=-\frac{1}{2}\int\limits_B{((\cof\grad\vp)\grad\la(x))\cd\vp\;dx}$ for any $\vp\in W_0^{1,2}(B,\R^2).$\\
vi)$\int\limits_B{\la(x)\det\grad\vp\;dx}=-\frac{1}{2}\int\limits_B{((\cof\grad\vp^{(0)})\grad\la(x))\cd\tilde{\vp}\;dx}-\frac{1}{2}\int\limits_B{((\cof\grad\vp)\grad\la(x))\cd\tilde{\vp}\;dx}$ for any $\vp\in W_0^{1,2}(B,\R^2).$
\label{Lem:Uni.Tech.1}
\end{lem}
\textbf{Proof:}
\begin{itemize}
\item [i),ii)] For a proof of the two points, see \cite[Prop 3.2]{JB14}. The argument given there still applies if $\la$ depends on $x$ instead of $R.$\\
\item [iii),iv)] For a proof of these two points, see \cite[Lem 3.2]{JB14}. \\
\item [v)] This statement differs slightly, from \cite[Lem 3.2.(iv)]{JB14}, hence, we include a proof here.\\

Assuming $\vp\in C_c^\infty$ then a simple computation shows: 
\begin{align*}\int\limits_B{\la(x)\det\grad\vp\;dx}=&\int\limits_B{\la(x)J\vp,_R\cd\vp,_\th\;\frac{dx}{R}}&\\
=&-\int\limits_B{(\la(x)J\vp,_R),_\th\cd\vp\;\frac{dx}{R}}&\\
=&-\int\limits_B{\la(x),_\th J\vp,_R\cd\vp\;\frac{dx}{R}}-\int\limits_B{\la(x)J\vp,_{R\th}\cd\vp\;\frac{dx}{R}}&\\
=&-\int\limits_B{\la(x),_\th J\vp,_R\cd\vp\;\frac{dx}{R}}+\int\limits_B{\la(x),_RJ\vp,_{\th}\cd\vp\;\frac{dx}{R}}+\int\limits_B{(\la(x)J\vp,_{\th})\cd\vp,_R\;\frac{dx}{R}}&\\
=&-\int\limits_B{((\cof\grad\vp)\grad\la(x))\cd\vp\;dx}-\int\limits_B{(\la(x)J\vp,_{R})\cd\vp,_\th\;\frac{dx}{R}}&
\end{align*}
The result follows by bringing the rightmost term to the LHS and dividing by two.\\

Note, as a last step one needs to upgrade the above equation to hold, not just for $\vp\in C_c^\infty$ but instead for all $\vp\in W_0^{1,2}.$ This is slightly delicate, because of the weak spaces involved, for a proof, see \cite[Lem 3.2.(iv)]{JB14}.\\

\item [vi)] This is a version of (v) where it is emphasised in which way does the above expression depend on $\vp^{(0)}$.
Again, we assume $\vp\in C_c^\infty,$ and we start by noting $\vp,_\th=\tilde{\vp},_\th,$ hence,
\begin{align*}\int\limits_B{\la(x)\det\grad\vp\;dx}=&\int\limits_B{\la(x)J\vp,_R\cd\tilde{\vp},_\th\;\frac{dx}{R}}&\\
=&-\int\limits_B{((\cof\grad\vp)\grad\la(x))\cd\tilde{\vp}\;dx}+\int\limits_B{(\la(x)J\tilde{\vp},_{\th})\cd\tilde{\vp},_R\;\frac{dx}{R}}&
\end{align*}
then the rightmost term is just
\begin{align*}\int\limits_B{(\la(x)J\tilde{\vp},_{\th})\cd\tilde{\vp},_R\;\frac{dx}{R}}=-\int\limits_B{(\la(x)J\tilde{\vp},_R)\cd\tilde{\vp},_{\th}\;\frac{dx}{R}}=\frac{1}{2}\int\limits_B{((\cof\grad\tilde{\vp})\grad\la(x))\cd\tilde{\vp}\;dx},
\end{align*}
together with the above we get
\begin{align*}\int\limits_B{\la(x)\det\grad\vp\;dx}=&-\frac{1}{2}\int\limits_B{((\cof\grad\vp^{(0)})\grad\la(x))\cd\tilde{\vp}\;dx}-\frac{1}{2}\int\limits_B{((\cof\grad\vp)\grad \la(x))\cd\tilde{\vp}\;dx},
\end{align*}
completing the proof.

\end{itemize}
\vspace{0.5cm}
The uniqueness condition will be of the form $\|\grad \la(x)R\|_{L^\infty(B,\R^2,\frac{dx}{R})}\le C,$ for some constant $C>0$ and where $\la$ will be the corresponding pressure to some stationary point.  A priori, the condition only guarantees the existence of $\grad \la(x)R$ in a suitable space. In the next lemma we show that $\la$ and $\grad\la$ themselves exist in a suitable space. This, in particular, allows one to make use of the technical lemma above.\\

\begin{lem}Let $\mu:B\ra\R$ be a function satisfying 
\[\|\grad \mu(x)R\|_{L^\infty(B,\R^2,\frac{dx}{R})}<\infty.\]
Then $\mu\in W^{1,p}(B,\R^2,dx)$ for any $1\le p<2.$\\
\end{lem}
\textbf{Proof:}\\
A proof is straightforward. Indeed, it holds
\[\int\limits_B{|\grad \mu|^p\;dx}\le\|\grad \mu R\|_{L^\infty(\frac{dx}{R})}^p\int\limits_B{R^{1-p}\;\frac{dx}{R}},\]
where the latter integrand is integrable for all $1\le p<2.$ \vspace{0.5cm}

We are ready to give the main statement of this section. It is shown that under suitable boundary conditions, and if $u$ is a stationary point of some energy $E$ and the corresponding pressure satisfies the assumption that $\|\grad \la(x)R\|_{L^\infty(B,\R^2,\frac{dx}{R})}\le \frac{\sqrt{3}\nu}{2\sqrt{2}},$ then $u$ is a global minimizer to \eqref{eq:USPS:1.1}. Notice, that for a vector-valued, measurable function $f=f_Re_R+f_{\th}e_{\th}:B\ra\R^2,$ we define its $L^\infty-$norm via $\|f\|_{L^\infty(B,\R^2,\mu)}:=\max\{\|f_R\|_{L^\infty(B,\mu)},\|f_{\th}\|_{L^\infty(B,\mu)}\}.$

\begin{thm}[Uniqueness under small pressure] Assume $M\in L^{\infty}(B,\R^{16})$ to be symmetric and to satisfy \eqref{eq:USPS:1.0.1} for some $\nu>0$. Let $b\in \R^2,$ $g \in W^{1,2}(S^1,\R^2),$ and let $u\in \A^c$ be an arbitrary stationary point of the energy $E$ defined in \eqref{eq:USPS:1.0}.
Furthermore, assume that for the map $u$ the corresponding pressure $\la$ exists and satisfies 
\begin{equation}\|\grad \la(x)R\|_{L^\infty(B,\R^2,\frac{dx}{R})}\le \frac{\sqrt{3}\nu}{2\sqrt{2}}.\label{eq:Uni.SPC.101}\end{equation}\vspace{0.25cm}

Then problem \eqref{eq:USPS:1.1} is globally minimized by the map $u.$\\

Moreover, if the inequality is strict, i.e. $|\grad \la(x)R|< \frac{\sqrt{3}\nu}{2\sqrt{2}}$ on a non-trivial measurable set $U\ss B$ with respect to the measure $\frac{dx}{R},$ then the above map is the unique minimizer of E in $\A^c$.
\label{Thm:Uni.SPC.1} 
\end{thm}
\textbf{Proof:}\\
Let $u\in\A^c$ be a stationary point with pressure $\la$ and let $v\in \A^c$ be arbitrary and set $\eta:=v-u\in W_0^{1,2}(B,\R^2).$\\

Again, we start expanding the energy via 
\[E(v)=E(u)+E(\eta)+H(u,\eta),\]
where
\[H(u,\eta):=2\int\limits_B{M(x)\grad u\cd\grad\eta\;dx}\]
denotes the mixed terms.\\
Now, expanding the Jacobian of $\eta$ and exploiting the fact that both $u$ and $v$ satisfy $\det\grad u=\det\grad v =1$ a.e. yields
\[\det\grad\eta=-\cof\grad u \cd \grad \eta \;\mb{a.e.}\]
By the later identity and realising that $(u,\la)$ satisfy the ELE given in \eqref{eq:Uni.RP.2}, $H$ can be written as
\begin{equation}H(u,\eta)=2\int\limits_B{\la(x)\det \grad \eta\;dx}.\label{eq:Uni.SPC.1}\end{equation}
By Lemma \ref{Lem:Uni.Tech.1}.(vi) we have
\[H(u,\eta)=-\int\limits_B{(\cof\grad\eta^{(0)}\grad\la(x))\cd\tilde{\eta}\;dx}-\int\limits_B{(\cof\grad\eta\grad\la(x))\cd\tilde{\eta}\;dx}=:(I)+(II).\]
Now by noting that the $0-$mode is only a function of $R,$ we get
\begin{align*}
(\cof\grad\eta^{(0)}\grad\la(x))\cd\tilde{\eta}=\frac{\la,_\th}{R}(\eta_{1,R}^{(0)}\tilde{\eta}_2-\eta_{2,R}^{(0)}\tilde{\eta}_1).
\end{align*}
Instead of just $\la,_\th$ on the right hand side of the latter equation we would like to have the full gradient of $\la.$ This can be achieved by using the basic relations $e_\th\cd e_\th=1$ and  $e_R\cd e_\th=0$ to obtain 
\begin{align*}
(\cof\grad\eta^{(0)}\grad\la(x))\cd\tilde{\eta}=(\la,_RR e_R+\la,_\th e_\th)\cd(\eta_{1,R}^{(0)}\tilde{\eta}_2-\eta_{2,R}^{(0)}\tilde{\eta}_1)\frac{e_\th}{R}.
\end{align*}
Arguing similarly for (II), and a short computation shows
\begin{align}
H(u,\eta)=&-\int\limits_B(\la,_RR e_R+\la,_\th e_\th)\cd\left[(\tilde{\eta}_1\tilde{\eta}_{2,\th}-\tilde{\eta}_2\tilde{\eta}_{1,\th})\frac{e_R}{R}\right.&\nonumber\\
&+\left. (\tilde{\eta}_2(\eta_{1,R}^{(0)}+\eta_{1,R})-\tilde{\eta}_1(\eta_{2,R}^{(0)}+\eta_{2,R}))e_\th\right]\;\frac{dx}{R}.&
\label{eq:Uni.SPC.2}\end{align}
Now for our choice of the norm of $f$ via $\|f\|_{L^\infty(B,\R^2,\mu)}:=\max\{\|f_R\|_{L^\infty(B,\mu)},\|f_{\th}\|_{L^\infty(B,\mu)}\}$ the corresponding Hölder estimate is given by
\begin{align*}\int\limits_B f\cd g\;d\mu=\int\limits_B f_R g_R+f_{\th} g_{\th}\;d\mu\le \|f\|_{L^\infty(B,\R^2,\mu)}\int\limits_B (| g_R|+| g_\th|)\;d\mu.\end{align*}
Then by applying the latter in \eqref{eq:Uni.SPC.2} we obtain 
\begin{align*} 
H(u,\eta)\ge&-\|\grad \la(x)R\|_{L^\infty(B,\R^2,\frac{dx}{R})}\int\limits_B\left[\left|\tilde{\eta}_1\tilde{\eta}_{2,\th}-\tilde{\eta}_2\tilde{\eta}_{1,\th}\right|\frac{1}{R}\right.&\\
&+\left.\left|\tilde{\eta}_2(\eta_{1,R}^{(0)}+\eta_{1,R})-\tilde{\eta}_1(\eta_{2,R}^{(0)}+\eta_{2,R})\right|\right]\;\frac{dx}{R}.&
\end{align*}
By $\|\grad\la(x) R\|_{L^\infty(\frac{dx}{R})}\le\frac{\sqrt{3}\nu}{2\sqrt{2}}$ and a weighted Cauchy-Schwarz Inequality, we see
\begin{align*}
H(u,\eta)\ge&-\frac{\nu\sqrt{3}}{4\sqrt{2}}\left[2a\|\tilde{\eta}_1\|_{L^2(dx/R^2)}^2+2a\|\tilde{\eta}_2\|_{L^2(dx/R^2)}^2\right.&\\
&\left.+\frac{1}{a}\int\limits_B{\left[\frac{\tilde{\eta}_{2,\th}^2}{R^2}+(\eta_{2,R}^{(0)}+\eta_{2,R})^2+(\eta_{1,R}^{(0)}+\eta_{1,R})^2+\frac{\tilde{\eta}_{1,\th}^2}{R^2}\right]\;dx}\right].&
\end{align*}
Now again, by using the Cauchy-Schwarz inequality, Fourier-estimate \eqref{eq:Uni.Buckling.1} and combining some of the norms yields
\begin{align*}
H(u,\eta)\ge&-\frac{\nu\sqrt{3}}{4\sqrt{2}}\left[(2a+\frac{1}{a})\|\tilde{\eta}_1,_\th\|_{L^2(dx/R^2)}^2+(2a+\frac{1}{a})\|\tilde{\eta}_2,_\th\|_{L^2(dx/R^2)}^2\right.&\\
&\left.+\frac{2}{a}\|\eta_1,_R^{(0)}\|_{L^2(dx)}^2+\frac{2}{a}\|\eta_1,_R\|_{L^2(dx)}^2+\frac{2}{a}\|\eta_2,_R^{(0)}\|_{L^2(dx)}^2+\frac{2}{a}\|\eta_2,_R\|_{L^2(dx)}^2\right]&\\
\ge&-\frac{\nu\sqrt{3}}{4\sqrt{2}}\left[(2a+\frac{1}{a})\|\tilde{\eta},_\th\|_{L^2(dx/R^2)}^2+\frac{2}{a}\|\eta,_R^{(0)}\|_{L^2(dx)}^2+\frac{2}{a}\|\eta,_R\|_{L^2(dx)}^2\right].&
\end{align*}
Making use of $\tilde{\eta},_\th=\eta,_\th,$ which is true since the 0-mode does not depend on $\th,$ and $\|\eta,_R^{(0)}\|_{L^2(dx)}^2\le\|\eta,_R\|_{L^2(dx)}^2$ to obtain
\begin{align*}
H(u,\eta)\ge&-\frac{\nu\sqrt{3}}{4\sqrt{2}}\left[(2a+\frac{1}{a})\|\eta,_\th\|_{L^2(dx/R^2)}^2+\frac{4}{a}\|\eta,_R\|_{L^2(dx)}^2\right].&
\end{align*}
Choosing $a=\frac{\sqrt{3}}{\sqrt{2}}$ and again combining norms gives
\begin{align*}
H(u,\eta)\ge&-\frac{\nu\sqrt{3}}{4\sqrt{2}}\left[\frac{4\sqrt{2}}{\sqrt{3}}(\|\eta,_\th\|_{L^2(dx/R^2)}^2+\|\eta,_R\|_{L^2(dx)}^2)\right]&\\
=&-\nu D(\eta),&
\end{align*}
where $D(\eta):=\|\grad\eta\|_{L^2(dx)}^2$ denotes the Dirichlet energy.
This yields,
\[E(\eta)+H(u,\eta)\ge E(\eta)-\nu D(\eta)\ge0.\]
Recalling, that $\nu\Id\le M$ holds completes the proof.\\\vspace{0.5cm}

\begin{re} The prefactor $\frac{\sqrt{3}}{2\sqrt{2}}$ is the best we can reach at the moment. It remains an open question if it can be improved for general $\la.$
\end{re}
Even if we can not improve the prefactor in the general case, we can still improve it if $\la$ depends only on one of $R$ and $\th.$

\begin{cor} If either $\la(x)=\la(R)$ or $\la(x)=\la(\th)$ for all $x\in B$ then the prefactor in \eqref{eq:Uni.SPC.101} can be chosen as $1.$ 
\label{cor:Uni.SPC.Cor.1}
\end{cor}
\textbf{Proof:}\\
\textbf{(i) $\la(x)=\la(R):$} 
This case is significantly simpler and one can argue more along the lines of the proof to Theorem \ref{lem:Glob.1.5}.
The reason is that in this case it still holds that
\[H(v,\eta)=2\int\limits_B{\la(R)\det \grad \eta\;dx}=2\int\limits_B{\la(R)\det \grad \tilde{\eta}\;dx},\]
where $\tilde{\eta}=\eta-\eta^{(0)},$ eliminating the $0-$mode.\\ Then applying Lemma 3.2.(iv) of \cite{JB14} as before yields,
\[H(v,\eta)=\int\limits_B{\la'(R)R\tilde{\eta}\cd J \tilde{\eta}_{,\th}\;\frac{dx}{R^{2}}}.\]
Following the argument of Theorem \ref{lem:Glob.1.5} and using $\|\la'(R)R\|_{L^\infty(\frac{dx}{R})}\le\nu,$ we get
\begin{align*}
H(v,\eta)\ge&-\|\la'(R)R\|_{L^\infty(\frac{dx}{R})} \int\limits_B{|\tilde{\eta}||\tilde{\eta}_{,\th}|\;\frac{dx}{R^{2}}},&\\
\ge&-\nu \left(\int\limits_B{|\tilde{\eta}|^2\;\frac{dx}{R^{2}}}\right)^{1/2}\left(\int\limits_B{|\tilde{\eta}_{,\th}|^2\;\frac{dx}{R^{2}}}\right)^{1/2}&\\
\ge&-\nu  \int\limits_B{|\tilde{\eta}_{,\th}|^2\;\frac{dx}{R^{2}}}&\\
\ge&-\nu  \int\limits_B{|\grad \eta|^2\;dx}.&
\end{align*}
Note, as before the $\sim$ could be dropped, since $ \int\limits_B{|\grad \tilde{\eta}|^2\;dx}\le \int\limits_B{|\grad \eta|^2\;dx}.$ \vspace{0.5cm}

\textbf{(ii) $\la(x)=\la(\th):$}
Here we start with \eqref{eq:Uni.SPC.2} which simplifies to
\begin{align*}
H(u,\eta)=&-\int\limits_B{\la,_\th(\th)[(\eta_2,_R^{(0)}+\eta_2,_R)\tilde{\eta}_1+(\eta_1,_R^{(0)}+\eta_1,_R)\tilde{\eta}_2]\;\frac{dx}{R}}&
\end{align*}
By Hölder's inequality, Inequality \eqref{eq:Uni.Buckling.1} and  $\|\la,_\th\|_{L^\infty}(\frac{dx}{R})\le\nu$ we get
\begin{align*}
H(u,\eta)\ge&-\|\la,_\th\|_{L^\infty(\frac{dx}{R})}\int\limits_B{|(\eta_2,_R^{(0)}+\eta_2,_R)\tilde{\eta}_1+(\eta_1,_R^{(0)}+\eta_1,_R)\tilde{\eta}_2|\;\frac{dx}{R}}&\\
\ge&-\frac{\nu}{2}[2\|\tilde{\eta}_1,_\th\|_{L^2(dx/R^2)}^2+2\|\tilde{\eta}_2,_\th\|_{L^2(dx/R^2)}^2+\|\eta,_R^{(0)}\|_{L^2(dx)}^2+\|\eta,_R\|_{L^2(dx)}^2]&
\end{align*}
Using $\tilde{\eta},_\th=\eta,_\th$ and $\|\eta,_R^{(0)}\|_{L^2(dx)}^2\le\|\eta,_R\|_{L^2(dx)}^2$ we get
\begin{align*}
H(u,\eta)\ge&-\frac{\nu}{2}[2\|\eta,_\th\|_{L^2(dx/R^2)}^2+2\|\eta,_R\|_{L^2(dx)}^2]&\\
=&-\nu D(\eta).&
\end{align*}
Completing the proof.\\

\begin{re}[Relaxation of the assumptions]
For Theorem \ref{Thm:Uni.SPC.1} to hold it is enough to assume that $f$ is asymptotically convex, i.e. there exists $\nu\in L^\infty(B)$ s.t. $0\le \nu(x)\Id\le f(x,\xi)=M(x)\xi\cd\xi$ for a.e. $x\in B$, and any $\xi\in\R^2$ and $M\in L^\infty.$ Here $\nu$ is allowed to be $0.$\footnote{We do not have to worry about existence, for Theorem \ref{Thm:Uni.SPC.1} to hold, since if there is no stationary point the result simply does not apply.} Additionally, we have to assume $\nu(x)=\nu(R),$ since we do not know if the Fourier-estimate \eqref{eq:Uni.Buckling.1} is still true if $\nu$ depends on $\th.$ \\

In this case the small pressure condition can be relaxed to a pointwise estimate:
\[|\grad \la(x)R|\le \frac{\sqrt{3}\nu(R)}{2\sqrt{2}}\;\mb{for a.e.}\;x\in B\;\mb{wrt.}\; \frac{dx}{R},\]
with uniqueness if the above condition is strictly satisfied on a non-trivial set.
\end{re}

\subsection{Computational method of the pressure}
In this paragraph we develop a way of how to compute the quantity $\grad \lambda(x)R$ in situations when the corresponding stationary point $u$  separates the variables $R$ and $\th.$ Obviously, the class of maps 
\[\A_s^c:=\A^c\cap\{u=r(R)g(\th)+b|r\in L^2((0,1],R^{-1}\;dR)\wedge r'\in L^2((0,1],R\;dR),r(1)=1\}\]
where the variables are separated will play a key role. \\

As a first step, we see that the incompressibility condition implies that $\A_s^c$ is empty for many choices of $g$. If however, $g$ satisfies $\det\grad g=1$ a.e. in $[0,2\pi),$ then there is exactly one possible map $u=Rg(\th)+b\in \A_s^c.$\\

\begin{lem}[Rigidity of the class $\A_s^c$] Let $b\in\R^2,$ $g\in W^{1,2}(S^1,\R^2)$ and $u\in \A_s^c.$ Then the Jacobian of $u$ can be represented by
\[\det\grad u=\left(\frac{r(R)r'(R)}{R}\right)(Jg(\th)\cd g'(\th)), \;\mb{a.e. in B.}\]

Moreover, $u$ obeys $\det\grad u=1$ a.e. only if $u=Rg(\th)+b,$ where $\det\grad g=Jg(\th)\cd g'(\th)=1$ for almost any $\th\in[0,2\pi).$ 

\end{lem}
\textbf{Proof:}\\
The argument is straight forward. By the Incompressibility for any $u\in\A_s^c$ the following equation needs to hold for a.e. $x\in B$
\[1=\det\grad u=\left(\frac{r(R)r'(R)}{R}\right)(Jg(\th)\cd g'(\th)).\]
By assumption $(Jg(\th)\cd g'(\th))=1$ a.e. we get
\[\left(\frac{r(R)r'(R)}{R}\right)=1, \;\mb{with}\; r(1)=1.\]
which is uniquely solved by $r(R)=R,$ for any $R\in(0,1],$ finishing the proof.
\\\vspace{0.5cm}

\begin{re}\emph{
The above statement has two consequences. The first one is that for any $b\in \R^2$ and normalized $g,$ i.e. $g$ obeys $\det\grad g=1$ a.e. in $[0,2\pi),$ there is exactly one $u\in\A_s^c$ namely the one homogenous map $u=Rg(\th)+b,$ i.e.
\[\A_s^c=\{u=Rg(\th)+b\}.\]   
The second one is that for any non-normalized $g,$ there is no function $u\in \A_s^c$ at all.}\\
\end{re}

Secondly, for stationary points that are in the class $\A_s^c$ and therefore 1-homogenous we compute the corresponding quantity $\grad \la(x)R.$ This will be an important ingredient to actually apply our theory. The next lemma shows that $\grad \la(x)R$ needs to satisfy a PDE system weakly and that $\grad \la(x)R$ is completely determined by $M$ (corresponding to the functional) and the boundary conditions $g.$\\

\textbf{Notation:} Recall the notation for 2d polar coordinates $\{\hat{x},\hat{x}^\perp\}=\{e_R,e_\th\}=\{e_R(\th),e_\th(\th)\}:=\{(\cos(\th),\sin(\th)),(-\sin(\th),\cos(\th))\}.$ Additionally, we will use $\{e_{NR},e_{N\th}\}=\{e_R(N\th),e_\th(N\th)\}:=\{(\cos(N\th),\sin(N\th)),(-\sin(N\th),\cos(N\th))\},$ for any $N\in \N.$ Moreover, we will use the notation $M_{ijgk}=(M(e_i\ot e_j))\cd(g\ot e_k),$ for any combination of $i,j,k\in\{R,\th\}$ and any map $g \in \R^{2}.$ Especially, if $g=e_{Nl}$ for some $l\in\{R,\th\}$ we will use $M_{ij(Nl)k}$ for short.\\

\begin{lem}[Representation of the pressure] Let $1\le p\le\infty,$ $k\in\N\sm\{0\}$ and assume $M\in L^{\infty}\cap W^{k,p}(B,\R^{16},dx),$ $g \in W^{k+1,p}(S^1,\R^2)$  where $g$ obeys $\det\grad g=1$ a.e. in $[0,2\pi),$ $b\in \R^2,$ and let $u=Rg(\th)+b$ be a stationary point of the energy $E.$\\

Then there exists a corresponding pressure $\la\in W^{k,p}(B,\R,dx)$ and it satisfies the following system of equations a.e. wrt. $\frac{dx}{R}$ in $B:$
\begin{align}
\la(x),_\th(Jg\cd e_R)-\la(x),_RR(Jg'\cd e_R)=&-[M_{R\th (g+g'')\th}&\nonumber\\
&+((M,_\th)_{R\th g R}+(M,_\th)_{R \th g'\th})&\nonumber\\
&+ R((M,_R)_{R R g R}+(M,_R)_{R R g' \th})]&\nonumber\\
=:&h_1(M,g)&\nonumber\\
\la(x),_\th(Jg\cd e_\th)-\la(x),_RR(Jg'\cd e_\th)=&-[M_{\th \th (g+g'')\th}&\nonumber\\
&+((M,_\th)_{\th \th g R}+(M,_\th)_{\th \th g'\th})&\nonumber\\
&+ R((M,_R)_{\th R g R}+(M,_R)_{\th R g' \th})]&\nonumber\\
=:&h_2(M,g)&
\label{eq:Uni.RP.1}
\end{align}
\end{lem}
 \textbf{Proof:}\\
Let $u\in \A_s^c$ be a stationary point. If there exists a corresponding pressure $\la\in W^{1,p}$ then $u$ is a solution to
\begin{equation}\int\limits_B{M(x)\grad u\cd\grad\eta\;dx}=-\int\limits_B{\la(x)\cof\grad u\cd\grad\eta\;dx}\;\mb{for any} \eta\in C_c^\infty(B,\R^2).\label{eq:Uni.RP.2}\end{equation}

The strategy is as follows. For now we will assume that $\la\in W^{1,p}.$ In order to derive the above system of equations, we will start by entering the explicit form of $u$ and the representation $\eta=(\eta\cd e_R) e_R+(\eta\cd e_{\th}) e_{\th}$ into the latter equation. By some further computations, which are mainly integrations by parts, one obtains \eqref{eq:Uni.RP.1}. The last step discusses the existence of $\la\in W^{1,p}.$\\

\textbf{Step 1:} Computation of (LHS):\\
The derivative and the cofactor of the map $u=Rg(\th)+b$ are given by
\[\grad u=g\ot e_R+g'\ot e_\th,\]
\[\cof\grad u=Jg\ot e_\th-Jg'\ot e_R.\]

Plugging the above into the left-hand side of \eqref{eq:Uni.RP.2} and integrating by parts yields
\begin{align*}(LHS)=&\int\limits_B{M(x)(g\ot e_R+g'\ot e_\th)\cd\left(\eta,_R\ot e_R+\frac{1}{R}\eta,_\th\ot e_\th\right)\;dx}&\\
=&-\int\limits_B{R(M(x),_R)(g\ot e_R+g'\ot e_\th)\cd(\eta\ot e_R)\;\frac{dx}{R}}&\\
&-\int\limits_B{[M(x)((g+g'')\ot e_\th)+M(x),_\th(g\ot e_R+g'\ot e_\th)]\cd(\eta\ot e_\th)\;\frac{dx}{R}}&
\end{align*}
Now by expanding  $\eta=\al e_R+\be e_{\th}$ with $\al=(\eta\cd e_R)$ and $\be=(\eta\cd e_{\th})$ and the shorthand introduced above we get
\begin{align*}(LHS)=&-\int\limits_B{R[(M,_R)_{RRgR}+(M,_R)_{RRg'\th}]\al+R[(M,_R)_{\th RgR}+(M,_R)_{\th Rg'\th}]\be\;\frac{dx}{R}}&\\
&-\int\limits_B{M_{R\th(g+g'')\th}\al+[(M,_\th)_{R\th gR}+(M,_\th)_{R\th g'\th}]\al}&\\
&{+M_{\th\th(g+g'')\th}\be+[(M,_\th)_{\th\th gR}+(M,_\th)_{\th\th g'\th}]\be\;\frac{dx}{R}}&\\
=&\int\limits_B{h_1\al+h_2\be\;\frac{dx}{R}}.&
\end{align*}

\textbf{Step 2:} Computation of (RHS):  \\
Now by again using the explicit form of $\cof\grad u,$ and integration by parts we get
\begin{align*}
(RHS)=&-\int\limits_B{(\la(x)(Jg\cd\eta,_\th)-\la(x)R(Jg'\cd\eta,_R))\;\frac{dx}{R}}&\\
=&\int\limits_B{\la,_\th(x)(Jg\cd\eta)-\la,_R(x)R(Jg'\cd\eta)\;\frac{dx}{R}}.&
\end{align*}
Further we use the expression $\eta=\al e_R+\be e_{\th}$ with the notation $\al=(\eta\cd e_R)$ and $\be=(\eta\cd e_{\th})$ to derive
\begin{align*}
(RHS)=&\int\limits_B{(\la,_\th(x)(Jg\cd e_R)-\la,_R(x)R(Jg'\cd e_R))\al}&\\
&{+(\la,_\th(x)(Jg\cd e_\th)-\la,_R(x)R(Jg'\cd e_\th))\be\;\frac{dx}{R}}.&
\end{align*}
Together with Step 1 and the realization that in the above $\al,\be\in C_c^\infty(B)$ are arbitrary, the claimed equations need to be true a.e. in $B$.\\

\textbf{Step 3:} Existence of the pressure $\la\in W^{1,p}:$\\
The above can be rewritten as 
\[\int\limits_B{\div(\la\cof\grad u)\cd\eta\;dx}=\int\limits_B{h(M,g)\cd\eta\;\frac{dx}{R}},\]
where $h=(h_1,h_2).$ From the above we know that $\div(\la\cof\grad u)\in L^p(dx)$ iff $h(M,g)\in L^p(\frac{dx}{R}).$ Now consider $h_1(M,g)$ (similar for $h_2$) and we define
\begin{align*}
h_{11}:=&-M_{R\th (g+g'')\th}&\\
h_{12}:=&[((M,_\th)_{R\th g R}+(M,_\th)_{R \th g'\th})&\\
&+ R((M,_R)_{R R g R}+(M,_R)_{R R g' \th})].&
\end{align*}
Then for $h_{11}\in L^p(\frac{dx}{R})$ we need $M\in L^\infty(dx)$ and $g,g''\in L^p,$ which is true by assumption. Now, by Sobolev imbedding we have $W^{2,p}\hookrightarrow W^{1,\infty}([0,2\pi),\R^2),$ and hence, in order for $h_{12}\in L^p(\frac{dx}{R}),$ since, $g,g'\in L^\infty,$ it is enough to require that $\grad M\in L^{p}(dx).$
This is exactly how we chose the classes for $M$ and $g.$ This guarantees the existence of $\div(\la\cof\grad u)=(\cof\grad u)\grad\la\in L^p(dx).$ By further noting that $g,g'\in L^\infty$ implies $\grad u\in L^\infty(dx).$ Now by the relation $A^T\cof A=(\det A)\Id,$ which is true for any $A\in \R^{2\ti2},$ we have $\grad\la=\grad u^T(\cof\grad u)\grad\la\in L^p(dx).$ This guarantees the existence of $\grad\la\in L^p(dx)$ as desired. One can argue similarly for the higher integrability.\\

\begin{re} An analogous representation of the pressure is possible in any type of coordinates, under suitable assumptions on $M.$ It is important to realize that while the representation of the pressure may vary from coordinates to coordinates, all representations have to agree since the pressure is unique up to a constant.  
\end{re}

\subsection{Small pressure condition under affine boundary data}
To get an understanding of the behaviour of the pressure and the small pressure condition we study them in the easiest of all situations: the case of affine boundary conditions.\\

\begin{exa}[Affine boundary data] Let $1\le p\le\infty,$ and $M\in L^{\infty}\cap W^{1,p}(B,\R^{16})$ and $g=Ae_R,$ with $A\in SL^+(2,\R)$ a constant matrix and $u=Ax+b\in \A$ a stationary point of the energy $E.$\\ 

Then by Lemma \ref{eq:Uni.RP.1} there exists a corresponding pressure $\la\in W^{1,p}$ and it satisfies the following system of equations a.e. wrt. $\frac{dx}{R}$ in $B:$

\begin{align*}
\la(x),_\th(JAe_R\cd e_R)-\la(x),_RR(JAe_\th\cd e_R)=&-[((M,_\th)_{R\th (AR) R}+(M,_\th)_{R \th (A\th) \th})&\\
&+ R((M,_R)_{R R (AR) R}+(M,_R)_{R R (A\th) \th})]&\\
=:&h_1(M,Ae_R)&\\
\la(x),_\th(JAe_R\cd e_\th)-\la(x),_RR(JAe_\th\cd e_\th)=&-[((M,_\th)_{\th \th (AR) R}+(M,_\th)_{\th \th (A\th) \th})&\\
&+ R((M,_R)_{\th R (AR) R}+(M,_R)_{\th R (A\th) \th})]&\\
=:&h_2(M,Ae_R).&
\end{align*}
\end{exa} 
(i) $M=const$ then $\grad\la=0.$

An immediate consequence is that if $M=const$ then $\grad\la$ vanishes and the pressure term has no effect. Hence, stationary points to the constrained and unconstrained problem agree and the affine map $u=Ax+b$ is always the global minimizer. There is a much simpler way of seeing that, see Remark \ref{re: Uni.2.1}.\\

(ii) Let $A=\Id,$ $M(x)=diag(M_{R R R R},M_{R\th R \th},M_{\th R \th R},M_{\th \th \th \th})$ with $\nu>0$ and
$M_{R R R R}(x)=f(R)\ge\nu,$ $M_{R\th R \th}(x),M_{\th R \th R}(x)\ge\nu$ arbitrary, $M_{\th \th \th \th}=F(R)\ge\nu.$ 
Then the system becomes
\begin{align*}
\la(x),_RR=&-R(M,_R)_{R R R R}&\\
\la(x),_\th=&-(M,_\th)_{\th \th \th \th}&
\end{align*}
and the small pressure condition reads 
\[\sup_{x\in B}|\la(x),_RRe_R+\la(x),_\th e_\th|=\sup_{x\in B}|-R(M,_R)_{R R R R}(x)-(M,_\th)_{\th \th \th \th}(x)|\le \frac{\sqrt{3}\nu}{2\sqrt{2}}\]
Since, $M_{\th \th \th \th}(x)=F(R),$ and the fact that all off-diagonal elements are $0,$ $(M,_\th)_{\th \th \th \th}=0$ and so the condition becomes
\[|-Rf'(R)|\le \frac{\sqrt{3}\nu}{2\sqrt{2}}\]
Hence, $f$ needs to satisfy
\[f(1)+\frac{\sqrt{3}\nu}{2\sqrt{2}}\ln R\le f(R)\le f(1)-\frac{\sqrt{3}\nu}{2\sqrt{2}}\ln R \;\mb{for all}\;R\in(0,1].\]
In particular, this means that in order for $f$ to satisfy the small pressure condition it can not grow any faster than $c-\frac{\sqrt{3}\nu}{2\sqrt{2}}\ln R,$ for some constant $c,$ close to the origin.\\
 
\begin{re}\label{re: Uni.2.1}
1. ($A\in SL^+(2,\R),$ $M$ constant): Consider the case $g(\th)=Ae_R(\th),$ where $A$ is constant and $A\in SL^+(2,\R),$ i.e. $\det A=1.$ Assume further that $M(x)=M$ is constant. In this situation we know that the affine map $u=Ax+b$ is the unique global minimizer of the unconstrained problem. However, the setting is such that the affine map satisfies the constraint $\det\grad u= \det A=1.$ Now, noting that
\[E(u)=\min\limits_{uncon.} E(v)\le\min\limits_{con.} E(v),\]
the affine map needs to be the unique global minimizer to the constrained problem, as well.\\

2. ($A\in \R^{2\ti2}\sm SL^+(2,\R),$ $M$ arbitrary): Let $g(\th)=Ae_R(\th),$ $A\in \R^{2\ti2}\sm SL^+(2,\R)$ s.t. $\det A\not=1$ i.e. $u|_{\p B}=Ax+b$ and $M(x)$ arbitrary. Then there can not be any minimizer at all. This can be seen from the fact that for any $u\in\A^c$ because of the constraint $\det\grad u=1 \; a.e.$ it holds
\begin{eqnarray*}
\int\limits_B\det\grad u\;dx=\La^2(B).
\end{eqnarray*}
But on the other hand we know that the Jacobian is a Null-Lagrangian hence,
\begin{eqnarray*}
\int\limits_B\det\grad u\;dx=\int\limits_B\det\grad v\;dx
\end{eqnarray*}
for any $v\in\A.$ Choosing $v=Ax+b\in\A$ yields
\begin{eqnarray*}
\int\limits_B\det\grad u\;dx=\int\limits_B\det A\;dx=\det A\La^2(B)
\end{eqnarray*} 
for any $u\in\A.$ Since, $\det A\not=1$ the set $\A^c$ in this problem is empty.
\end{re}\vspace{0.5cm}

\subsection{Small pressure condition under N-covering boundary data and a counterexample to regularity}
\label{sec:4.2.3}
As in the paragraph before, we study the behavior of the pressure and the small pressure criteria with the N-covering boundary conditions. However, this leads naturally to the construction of a counterexample to regularity. In order to understand the idea, realize that $u=\frac{R}{\sqrt{N}}e_{NR}(\th)+b\in C^{0,1}(B,\R^2)\sm C^{1}(B,\R^2),$ for any $N\in\N\sm\{0,1\}$, since $\grad u$ is discontinuous at $0.$\footnote{From now on we make use of the notation  $e_{NR}(\th)=e_{R}(N\th)$ and $e_{N\th}(\th)=e_{\th}(N\th)$ allowing to suppress the argument for shortness.} If we are able to construct a concrete $M$ s.t. the corresponding pressure $\la$ satisfies the small pressure criteria, then we could guarantee that the above map has to be a minimizer of the corresponding energy $E.$\\

For this sake, we start by computing the pressure in the case of $M$ depending only on $\th$ and being diagonal wrt. the basis of polar coordinates.\\
 
\begin{lem}[Representation of the pressure, N-cover, M(\th)=diag] For $N\in\N\sm\{0,1\}$ let $g=\frac{1}{\sqrt{N}}e_{NR}$ and assume $M\in L^{\infty}\cap W^{k,p}(B,\R^{16})$ for some $1\le p\le\infty$ and $k\in\N\sm\{0\},$ where $M$ is of the specific form
\[M(x)=diag(M_{R R R R},M_{R\th R \th},M_{\th R \th R},M_{\th \th \th \th})=diag(\al(\th),\be(\th),\ga(\th),\d(\th))\] with $\nu>0$ and $\al,\be,\ga,\d\ge\nu$ for any $\th\in[0,2\pi).$ Furthermore, suppose $u\in \A_s^c$ is a stationary point of the energy $E.$\\

Then there exists a corresponding pressure $\la\in W^{k,p}(B)$ and it satisfies the following system of equations a.e. wrt. $\frac{dx}{R}$ in $B:$
\begin{align}
-\la(x),_\th\frac{1}{\sqrt{N}}\sin(\th_{N-1})+\la(x),_RR\sqrt{N}\cos(\th_{N-1})=&\sqrt{N}\be'\sin(\th_{N-1})&\nonumber\\
+\left[\sqrt{N}(N-1)\be+\sqrt{N}\d-\right.&\left.\frac{\al}{\sqrt{N}}\right] \cos(\th_{N-1})&\nonumber\\
=:&h_1& \label{eq:Uni.NC.1}\\
\la(x),_\th\frac{1}{\sqrt{N}}\cos(\th_{N-1})+\la(x),_RR\sqrt{N}\sin(\th_{N-1})=&-\sqrt{N}\d'\cos(\th_{N-1})&\nonumber\\
+\left[\sqrt{N}\be+\sqrt{N}(N-1)\d-\right.&\left.\frac{\ga}{\sqrt{N}} \right] \sin(\th_{N-1})&\nonumber\\
=:&h_2&
\label{eq:Uni.NC.2}
\end{align}
where we used the shorthand $\th_{k}:=k\th$ for any $k\in\R.$ \\
\label{Lem:Uni.NC.1}
\end{lem}
\textbf{Proof:}\\
By Lemma \ref{eq:Uni.RP.1} we know that the pressure $\la$ exists and system \eqref{eq:Uni.RP.1} is satisfied. Now we just have to verify that \eqref{eq:Uni.RP.1} agrees with the claimed system given by \eqref{eq:Uni.NC.1}\ and \eqref{eq:Uni.NC.2}. We start by verifying the first one.\\

\textbf{Step 1:}\\
The LHS of \eqref{eq:Uni.NC.1} is just a simple computation. Hence, we focus on the corresponding RHS, which we named $h_1.$ First note, that $M$ only depends on $\th$ hence, we are only left with
\begin{equation}h_1=-\left[M_{R\th (g+g'')\th}+(M,_\th)_{R\th g R}+(M,_\th)_{R \th g'\th}\right].\label{eq:Uni.NC.3}\end{equation}
We have $g=\frac{1}{\sqrt{N}}e_{NR}, g'=\sqrt{N}e_{N\th},g''=-N\sqrt{N}e_{NR}$ and hence
\[M_{R\th (g+g'')\th}=\left(\frac{1}{\sqrt{N}}-\sqrt{N}N\right)M_{R\th (NR)\th}=\left(\frac{1}{\sqrt{N}}-\sqrt{N}N\right)(M_{R\th R\th}(e_{NR}\cd e_R)+M_{R\th \th \th}(e_{NR}\cd e_\th)).\]
Using that $M_{R\th \th \th}=0$ and $M_{R\th R\th}=\be$ yields
\[M_{R\th (g+g'')\th}=\left(\frac{1}{\sqrt{N}}-\sqrt{N}N\right)\be\cos(\th_{N-1}).\]
For the second term of \eqref{eq:Uni.NC.3} consider
\[M_{R\th g R,\th}=(M,_\th)_{R\th g R}+M_{\th \th g R}-M_{R R g R}+M_{R\th g' R}+M_{R\th g \th}\]
Rearranging and a small computation yields,
\begin{align*}
(M,_\th)_{R\th g R}=&M_{R\th g R,\th}-M_{\th \th g R}+M_{R R g R}-M_{R\th g' R}-M_{R\th g \th}&\\
=&\frac{1}{\sqrt{N}}[\al-\be]\cos(\th_{N-1}).&
\end{align*}
Similarly, for the rightmost term of \eqref{eq:Uni.NC.3} we get
\begin{align*}
(M,_\th)_{R\th g' \th}=&M_{R\th g' \th,\th}-M_{\th \th g' \th}+M_{R R g' \th}-M_{R\th g'' \th}+M_{R\th g' R}&\\
=&-\sqrt{N}\be'\sin(\th_{N-1})+\sqrt{N}[\be-\d]\cos(\th_{N-1}).&
\end{align*}
Together,
\[h_1=\sqrt{N}\be'\sin(\th_{N-1})+\left[\sqrt{N}\d-\frac{\al}{\sqrt{N}}+\sqrt{N}(N-1)\be\right]\cos(\th_{N-1}).\]

\textbf{Step 2:}\\
We can argue similar as above. For $h_2$ we have
\[h_2=-\left[M_{\th \th (g+g'')\th}+(M,_\th)_{\th \th g R}+(M,_\th)_{\th \th g'\th}\right].\]
Then 
\begin{align*}
M_{\th \th (g+g'')\th}=&\left(\frac{1}{\sqrt{N}}-\sqrt{N}N\right)\d\sin(\th_{N-1})&\\
(M,_\th)_{\th \th g R}=&M_{\th \th g R,\th}+M_{R \th g R}+M_{\th R g R}-M_{\th \th g' R}-M_{\th \th g \th}&\\
=&\frac{1}{\sqrt{N}}[\ga-\d]\sin(\th_{N-1})&\\
(M,_\th)_{\th \th g' \th}=&M_{\th \th g' \th,\th}+M_{R \th g' \th}+M_{\th R g' \th}-M_{\th \th g'' \th}+M_{\th \th g' R}&\\
=&\sqrt{N}\d'\cos(\th_{N-1})+\sqrt{N}[\d-\be]\sin(\th_{N-1})&
\end{align*}
and finally
\[h_2=-\sqrt{N}\d'\cos(\th_{N-1})+\left[\sqrt{N}\be+\sqrt{N}(N-1)\d-\frac{\ga}{\sqrt{N}}\right]\sin(\th_{N-1}),\]
completing the proof.\vspace{1cm}

Next we compute the small pressure criteria in the same situation. Moreover, we will provide concrete solutions.\\

\begin{lem}[Small pressure condition] Let the assumptions be as above.\\

For any $N\in\N\sm\{0,1\},$ let $M=(a,1,a,1)\nu,$ where we pick $a$ to be constant and in the range
\[1\le N^2-N< a < N^2+N.\]
Then for this $M$ the corresponding pressure $\la$ is given by
\begin{align*}
\la(x)=c+\left[N-\frac{a}{N}\right]\ln(R) \;\mb{for any}\; x\in B
\end{align*}
for any real constant $c\in \R,$ which is independent of $R$ and $\th.$ Moreover, 
$\la\in W^{1,q}(B) \;\mb{for any}\;1\le q<2$ and $\la$ satisfies condition \eqref{eq:Uni.SPC.101} strictly.
\label{Lem:Uni.NC.2}
\end{lem}
\textbf{Proof:}\\
Define first
\[H_1=\left[\sqrt{N}(N-1)\be+\sqrt{N}\d-\frac{\al}{\sqrt{N}}\right]\;\mb{and}\;H_2=\left[\sqrt{N}\be+\sqrt{N}(N-1)\d-\frac{\ga}{\sqrt{N}}\right].\]
By solving the system \eqref{eq:Uni.NC.1} and \eqref{eq:Uni.NC.2} we obtain
\begin{align*}
\la,_RR=&(\be'-\d')\frac{\sin(2\th_{N-1})}{2}+\frac{1}{\sqrt{N}}(H_1\cos^2(\th_{N-1})+H_2\sin^2(\th_{N-1}))&\\
\la,_\th=&\sqrt{N}(H_2-H_1)\frac{\sin(2\th_{N-1})}{2}-N(\be'\sin^2(\th_{N-1})+\d'\cos^2(\th_{N-1})).&
\end{align*}
For the specific case of $M=(a,1,a,1)\nu$ they become
\begin{align*}
\la,_RR=\left[N-\frac{a}{N}\right] \;\mb{and}\;\la,_\th=0
\end{align*}
showing, in particular, that $\la$ depends only on $R$ i.e. $\la(x)=\la(R).$
Indeed, the pressure is then given by
\begin{align*}
\la(x)=c+\left[N-\frac{a}{N}\right]\ln(R) \;\mb{for any}\; x\in B
\end{align*}
and for any real constant $c\in \R,$ which is independend of $R$ and $\th.$
The small pressure condition \eqref{eq:Uni.SPC.101}, where we are allowed to use the improved constant by Corollary \ref{cor:Uni.SPC.Cor.1}, then becomes
\begin{align*}
\left|N-\frac{a}{N}\right|<1.
\end{align*}
Solving this inequality by case distinction yields the claimed bounds on $a.$ The integrability is then easily deduced, completing the proof.\\\vspace{0.5cm}

To make it more accessible for the reader we collect what we have shown so far in the following corollary.

\begin{cor}[Counterexample to regularity]\label{cor:4.2.Counter.1} Let $\nu>0,$ $b\in\R^2,$ assume $u_0=g+b$ with $g=\frac{1}{\sqrt{N}}e_{NR}$ and $N\in\N\sm\{0,1\}$ on the boundary, and let the energy be given by
\begin{equation}E(u)=\int\limits_B{f(x, \grad u)\;dx}\label{eq:Uni.CTR.1.1}\end{equation}
for any $u\in W^{1,2}(B,\R^2),$ where the integrand $f$ is of the quadratic form
\[f(x,\xi)=\nu\left[a(e_{R}^T\xi e_{R})^2+(e_{R}^T\xi e_{\th})^2+a(e_{\th}^T\xi e_{R})^2+(e_{\th}^T\xi e_{\th})^2\right],\]
for any $x\in B\sm\{0\},$ $\xi\in \R^{2\ti2}$ and some $a\in\left(N^2-N,N^2+N\right).$\\

Then the following statements are true:

\begin{enumerate} [label=(\roman*)]
\item Then we can find an $M\in L^\infty(B,\R^{16})$ s.t.
\[f(x,\xi)=\nu M(x)\xi\cd\xi\]
for any $x\in B\sm\{0\},$ $\xi\in \R^{2\ti2}$ and where $M$ is of the explicit form\footnote{Here the multiplication is understood through its action on $\xi\in \R^{2\ti2}$ which is given by
\[((a\ot b)(c\ot d))\xi\cd\xi=(a\ot b)_{ij}\xi_{ij} (c\ot d)_{kl}\xi_{kl}\; \mb{for}\;i,j,k,l\in\{1,2\}.\]
}
\begin{align*}M(x)=&a((e_R\ot e_R)(e_R\ot e_R))+((e_R\ot e_\th)(e_R\ot e_\th))&\\
&+a((e_\th\ot e_R)(e_\th\ot e_R))+((e_\th\ot e_\th)(e_\th\ot e_\th))&\end{align*}
and satisfies 
$\nu\Id\le\nu M(x)$ for any $x\in B\sm\{0\}$ and any $N\in\N\sm\{0,1\}.$
\item The maps $x\mapsto M(x)$ and $x\mapsto f(x,\xi),$ for any $\xi\in \R^{2\ti2}\sm\{0\},$ are discontinuous at $0.$
\item The maps $x\mapsto M(x)$ and $x\mapsto f(x,\xi),$ for any $\xi\in \R^{2\ti2}\sm\{0\},$ belong to
\[W^{1,q}\sm W^{1,2} \;\mb{for any}\;1\le q<2\] with the spaces $(B,\R^{16})$ and $(B)$ respectively.
\item The map \begin{equation}u=\frac{R}{\sqrt{N}}e_{NR}+b\in C^{0,1}(B,\R^2)\sm C^1(B,\R^2)\label{eq:Uni.CTR.1.2}\end{equation}
is a stationary point of $E,$ as defined in \eqref{eq:Uni.CTR.1.1}, and the corresponding pressure $\la$ exists and satisfies $\la\in W^{1,q}(B) \;\mb{for any}\;1\le q<2.$ 
\item Moreover, the map $u,$ as given in \eqref{eq:Uni.CTR.1.2}, is the unique global minimizer of $E,$ as defined in \eqref{eq:Uni.CTR.1.1}, in the class $\A^c.$
\item The minimal energy is given by
\[\min\limits_{v\in\A^c}E(v)=\frac{\nu\pi}{2}(1+a)\left(\frac{1}{N}+N\right).\]
\end{enumerate}
\end{cor}
\textbf{Proof:}
\begin{enumerate} [label=(\roman*)]
\item\!\!-(ii) trivial.
\item[(iii)]It is enough to show this point for $M.$ Note, that $M$ only depends on $\th,$ i.e. $M(x)=M(\th).$ Hence, the gradient is given by
\begin{equation}\grad M=\frac{1}{R}\p_\th M(\th)\ot e_\th \;\mb{for any}\; x\in B\sm\{0\}.\label{eq:Uni.CTR.1.3}\end{equation}
First realise that the derivative wrt. $\th$ only replaces $e_R$ with $e_\th$ (up to sign) and vice versa and therefore one can still bound the modulus of $\|\p_\th M(\th)\|_{L^\infty(B,\R^{16},\frac{dx}{R})}\le C(a)$ via some real constant $C(a)>0,$ which does depend only on $a.$ Then integrating $|\grad M|^q $ wrt.\! $dx$ using \eqref{eq:Uni.CTR.1.3} and by the latter discussion the claim follows.

\item[(iv)] As a consequence of $g\in C^\infty$ and point (iii), Lemma \ref{Lem:Uni.NC.1} guarantees that $u$ is a stationary point and the existence of $\la$ in the right spaces.

\item[(v)] By Lemma \ref{Lem:Uni.NC.2} we know that $\la$ satisfies the small pressure criteria strictly. Together, with Theorem \ref{Thm:Uni.SPC.1} this implies that the given $u=\frac{R}{\sqrt{N}}e_{NR}+b$ is indeed the unique global minimizer to the energy $E.$
\item[(vi)]  trivial.
\end{enumerate}

\begin{re} We want to emphasise the meaning of this again:\\
For the full ball $B\ss\R^2$ and smooth boundary conditions (however, with a topological change), we have constructed a uniformly convex functional, which depends discontinuously on $x,$ but depends smoothly on $\grad u,$ s.t. the corresponding energy is uniquely globally minimised by a map, which is everywhere Lipschitz but not any better.\\ 
The fact that such a simple counterexample exists, shows the rigidity of the incompressible case. In other words, admissible maps satisfying the constraint seem extremely rare.  
\end{re}

\textbf{Higher integrability and optimality}\\
Instead of discussing the counterexample on a level of continuity and differentiability one can also raise the question at the level of Sobolev-spaces.
 For this sake, we first need to know how integrable the function $u=\frac{R}{\sqrt{N}}e_{NR}+b$ is:\footnote{Recently this question has attracted some attention in the compressible case, of course that is extraordinary difficult, since the counterexamples are generated by means of convex integration. First results on higher integrability are available in \cite{RZZ18,RZZ20}.}

\begin{lem} (Higher integrability) For any $N\in\N\sm\{0,1\}$ and $b\in\R^2$ we have
 \[u=\frac{R}{\sqrt{N}}e_{NR}+b\in W^{2,q}(B,\R^2)\sm W^{2,2}(B,\R^2), \;\mb{for any}\;1\le q<2.\]
\end{lem}
\textbf{Proof:}\\
It is enough to provide the first and second derivatives:
\[D u=\frac{1}{\sqrt{N}}e_{NR}\ot e_{R}+\sqrt{N}e_{N\th}\ot e_{\th}\]
\[D^2 u=\frac{1}{R}\left(\frac{1}{\sqrt{N}}-\sqrt{N}N\right)e_{NR}\ot e_{\th}\ot e_{\th}\]
in components that is 
\[\p_{ij} u_k=\frac{1}{R}\left(\frac{1}{\sqrt{N}}-\sqrt{N}N\right)(e_{NR})_k (e_{\th})_i (e_{\th})_j,\; \mb{for}\;i,j,k\in\{1,2\}.\]
From here the claim easily follows.\\

\begin{re} (Optimality of the counterexample)\\
1. In elliptic regularity theory one is often able to establish a higher-order regularity result, i.e. there comes a point when the candidate is in such a good space that one can conclude it actually has to be maximally smooth. It is then natural to ask what is the weakest space your candidate needs to be in, so that one can still make the previous conclusion? 
A counterexample like ours gives a limit on such a possible space.\\

The best high-order regularity result available in our case, we are aware of, can be found in \cite{BOP92}:\\
They showed that for the special case of the Dirichlet functional and $u\in W^{2,q}(B,\R^2)$ with $q>2$ being a stationary point satisfying $\det\grad u=1$ a.e., then $u\in C^\infty(B,\R^2).$
It is very likely that a similar result can be established for a general p-growth functional with the necessary changes in $q$). Lets assume for a second that such a result is indeed possible. Intriguingly, this seems to leave a gap about $q=2.$\\

2. Clearly one could refine the level of spaces even further, by considering it on a Fractional/Besov-space level (for example $u\in W^{1+\al,2}$ for some $\al\in(0,1)$), or one could ask whether $\grad u\in BMO,$ or even in some van Schaftingen-space, which range between BMO and $W^{1,2}$ (introduced in \cite{VS06}). \\

3. The singular set $\Sigma$ of our examples provided, is of course just the origin $\Sigma=\{0\}.$ It remains open if there can be situations in incompressible elasticity where the stationary points/minimizers possess a richer $\Sigma.$
\end{re}

\subsection{High pressure situations and high frequency uniqueness}
\label{sec:4.2.4}
In general, we can not expect global uniqueness if the pressure is high, see \cite{BeDe20}. However, as long as the quantity $\grad \la(x)R$ is still essentially bounded we can at least show existence for variations that consist only of large enough Fourier-modes. This can be seen as a generalisation of \cite[Prop 3.3]{JB14}.

\begin{thm} [High frequency uniqueness] Assume $M\in L^{\infty}(B,\R^{16})$ to be symmetric and to satisfy \eqref{eq:USPS:1.0.1} for some $\nu>0$. Let $b\in \R^2,$ $g \in W^{1,2}(S^1,\R^2),$ and let $u\in \A^c$ be an arbitrary stationary point of the energy $E$ defined in \eqref{eq:USPS:1.0}.\\

i) \textbf{(purely high modes.)} Suppose the corresponding pressure $\la,$ exists and satisfies
\[\|\grad \la(x)R\|_{L^\infty(B,\R^2,\frac{dx}{R})}\le n\nu\]
for some $n\in\N.$ 
Then the map $u$ minimizes $E$ in the subclass 
\[\F_n^c=\left\{v\in \A^c|\;\eta=v-u\in W_0^{1,2}(B,\R^2) \;\mb{and}\; \eta=\sum\limits_{j\ge n}\eta^{(j)} \right\}.\]
Moreover, if the inequality is strict, i.e. $|\grad \la(x)R|< n\nu$ on a non-trivial measurable set $U\ss B$ wrt. $\frac{dx}{R},$ then the above map is the unique minimizer in $\F_n^c$.\\
  
ii) \textbf{($0-$mode and high modes.)} Suppose instead that the pressure $\la$ exists and satisfies
\[\|\grad \la(x)R\|_{L^\infty(B,\R^2,\frac{dx}{R})}\le \frac{\sqrt{3}m\nu}{2\sqrt{2}}\]
for some $m\in\N.$ 
Then the map $u$ minimizes $E$ in the subclass 
\[\F_{0,m}^c=\left\{v\in \A^c|\;\eta=v-u\in W_0^{1,2}(B,\R^2) \;\mb{and}\; \eta=\eta^{(0)}+\sum\limits_{j\ge m}\eta^{(j)} \right\}.\]
Moreover, if the inequality is strict, i.e. $|\grad \la(x)R|< \frac{\sqrt{3}m\nu}{2\sqrt{2}}$ on a non-trivial measurable set $U\ss B$ wrt. $\frac{dx}{R},$ then the above map is the unique minimizer in $\F_{0,m}^c$.  

\end{thm}
\textbf{Proof:}\\
i) The proof is a combination of the proofs of the Theorems \ref{lem:Glob.1.5} and \ref{Thm:Uni.SPC.1}.\\
Let $u\in\A^c$ be a stationary point with pressure $\la,$ let $v\in \F_n^c$ be arbitrary and set $\eta:=v-u\in W_0^{1,2}(B,\R^2).$ Then $\eta$ is of the form $\eta=\sum\limits_{j\ge n}\eta^{(j)}.$ \\

Recall from \eqref{eq:Uni.SPC.1} and Lemma \ref{Lem:Uni.Tech.1}.(v) that
\[H(u,\eta)=2\int\limits_B{\la(x)\det \grad \eta\;dx}=\int\limits_B{R(\cof\grad\eta\grad\la)\cd\eta\;\frac{dx}{R}}.\]
Then similarly to the argument given in the proof of Theorem \ref{lem:Glob.1.5} we can conclude
\begin{align*}
H(u,\eta)\ge&-\|\grad \la(x)R\|_{L^\infty(\frac{dx}{R})}\int\limits_B{|\grad\eta||\eta|\;\frac{dx}{R}}&\\
\ge&-n\nu\left(\int\limits_B{|\grad\eta|^2\;dx}\right)^{\frac{1}{2}}\left(\int\limits_B{|\eta|^2\;\frac{dx}{R^2}}\right)^{\frac{1}{2}}&\\
\ge&-\nu\left(\int\limits_B{|\grad\eta|^2\;dx}\right).&
\end{align*}
Where we used again the estimate
\begin{equation}\int\limits_B{R^{-2}|\eta_{,\th}|^2\;dx}\ge n^2\int\limits_B{R^{-2}|\eta|^2\;dx}.\label{eq:Uni.HFU.1}\end{equation}
which uses the fact that $\eta$ only contains Fourier-modes $n$ or higher.\\

ii) Following the proof of Theorem \ref{Thm:Uni.SPC.1} and using \eqref{eq:Uni.HFU.1} instead of \eqref{eq:Uni.Buckling.1} one arrives at

\begin{align*}
H(u,\eta)\ge&-\frac{\nu\sqrt{3}m}{4\sqrt{2}}\left[\left(\frac{2a}{m^2}+\frac{1}{a}\right)\|\eta,_\th\|_{L^2(dx/R^2)}^2+\frac{4}{a}\|\eta,_R\|_{L^2(dx)}^2\right].&
\end{align*}
Choosing $a=\frac{\sqrt{3}m}{\sqrt{2}}$ gives
\begin{align*}
H(u,\eta)\ge&-\frac{\nu\sqrt{3}m}{4\sqrt{2}}\left[\frac{4\sqrt{2}}{\sqrt{3}m}\left(\|\eta,_\th\|_{L^2(dx/R^2)}^2+\|\eta,_R\|_{L^2(dx)}^2\right)\right]&\\
=&-\nu D(\eta).&
\end{align*}

\begin{re}
1. Obviously, $\F_n^c\ss\F_{0,n}^c,$ for any $n\in\N.$\\
2. Also if $n=0$ or $m\in\{0,1\},$ one recovers Theorem \ref{Thm:Uni.SPC.1}.\\  
3. Clearly, there is a connection between $n$ and $m.$ In order to see this assume $\|\grad \la(x)R\|_{L^\infty(\frac{dx}{R})}\le n\nu$ for some $n\in\N.$ Then for $m=\lceil\frac{2\sqrt{2}n}{\sqrt{3}}\rceil$  it holds $\|\grad \la(x)R\|_{L^\infty(\frac{dx}{R})}\le\frac{\sqrt{3}m\nu}{2\sqrt{2}}$ and part $(ii)$ applies. Roughly speaking, one can trade in some of the high modes in order to include the 0-mode. 
\end{re}

\subsection{High frequency uniqueness for the $p-$Dirichlet functional}
\label{sec:4.2.5}
For $u_0\in L^p(B,\R^2)$ with $2\le p<\infty$ we define
\[
\A_{u_0}^p:=\{u\in W^{1,p}(B,\R^2): \; u_{|{\p B}}=u_0\}\;\mb{and}\;\A_{u_0}^{p,c}:=\A^p\cap K,
\]
and let 
\begin{equation}E(u)=\int\limits_B{f(x,\grad u)\;dx}\label{eq:HFU.I.2.1}\end{equation}
where $f$ is a version of the $p-$Dirichlet functional, i.e. $f$ is of the form
\[f(x,\xi)=\frac{\nu(x)}{p}|\xi|^p \]
for a.e. $x\in B$ and $\xi\in \R^{2\ti2}.$
Moreover,  $\nu\in L^\infty(B)$ is supposed to satisfy $\nu(x)\ge0$ a.e.\! in $B.$ Since $\nu$ is allowed to take on the value $0,$ the integrand could indeed disappear for some $x \in B.$ Since, additionally $f$ is convex in its 2nd variable we say $f$ is asymptotically convex. \vspace{0.5cm}

Now we can give an analogous result for these types of integrands.

\begin{thm} [High frequency uniqueness] Let $2\le p<\infty,$ assume $u_0\in L^2(B,\R^2)$ on the boundary and let $u\in \A_{u_0}^{p,c}$ to be a stationary point of $E,$ as given in \eqref{eq:HFU.I.2.1}. 
Furthermore, let $\sigma(x):=\nu(x)|\grad u(x)|^{\frac{p-4}{2}}\grad u(x)\in L^{\frac{2p}{p-2}}(B,\R^{2\ti2})$ and assume that there exists $l\in \N$ s.t.
\begin{equation}|\sigma,_\th|\le l|\sigma| \mb{for a.e.}\;x\in B\label{eq:HFU.C.G.01} \end{equation}
holds.\\

Then the following statements are true:

i) \textbf{(purely high modes.)}  Suppose the corresponding pressure $\la$ exists and satisfies
\begin{equation}|\grad \la(x)R|\le n\nu(x)|\grad u|^{p-2}\mb{for a.e.}\;x\in B \;\mb{wrt.}\;\frac{dx}{R}\label{eq:HFU.IC.G.02}\end{equation}
for some $n\in\N.$\\

Then $u$ is a minimizer of $E$ in the subclass 
\[\F_{n_*}^{p,\sigma,c}=\left\{v\in \A_{u_0}^{p,c}|\eta=v-u\in W_0^{1,p}(B,\R^2) \;\mb{and}\; \sigma\eta=\sum\limits_{j\ge n}(\sigma\eta)^{(j)} \right\},\]
where $n_*^2=n^2+l^2.$
Moreover, if the inequality is strict on a non-trivial set wrt.\! $\frac{dx}{R}$ then $u$ is the unique minimizer in $\F_{n_*}^{p,\sigma,c}$.\\

ii) \textbf{($0-$mode and high modes.)} Suppose the corresponding pressure $\la$ exists and satisfies
\begin{equation}|\grad \la(x)R|\le \frac{\sqrt{3}m\nu(x)|\grad u|^{p-2}}{2\sqrt{2}}\mb{for a.e.}\;x\in B \;\mb{wrt.}\;\frac{dx}{R}\label{eq:HFU.IC.G.03}\end{equation}
for some $m\in\N.$\\
Then $u$ is a minimizer of $E$ in the subclass 
\[\F_{0,m_*}^{p,\sigma,c}=\left\{v\in \A_{u_0}^{p,c}|\eta=v-u\in W_0^{1,p}(B,\R^2) \;\mb{and}\; \sigma\eta=(\sigma\eta)^{(0)}+\sum\limits_{j\ge m_*}(\sigma\eta)^{(j)} \right\},\]
where $m_*^2=m^2+l^2.$
Moreover, if the inequality is strict on a non-trivial set wrt.\! $\frac{dx}{R}$ then $u$ is the unique minimizer in $\F_{0,m_*}^{p,\sigma,c}$.
\end{thm}
\textbf{Proof:}\\
i) We start by the standard expansion
\begin{align}
E(v)-E(u)=&\int\limits_B{\frac{\nu(x)}{p}(|\grad u+\grad \eta|^{p}-|\grad u|^{p})\;dx}&\nonumber\\
\ge&\int\limits_B{\frac{\nu(x)}{2}|\grad u|^{p-2}|\grad\eta|^2+\nu(x)|\grad u|^{p-2}\grad u\cd\grad\eta\;dx}.&
\label{eq:HFU.IC.G.04}
\end{align}
where we used the following inequality\footnote{see, \cite{SS18}, Prop A.1, with $\sigma=0.$}

\[\frac{1}{p}|b|^{p}\ge\frac{1}{p}|a|^{p}+|a|^{p-2}a(b-a)+\frac{1}{2}|a|^{p-2}|b-a|^{2}.\] 
The ELE is given by
\begin{align*}
\int\limits_B{\nu(x)|\grad u|^{p-2}\grad u\cd\grad\eta\;dx}=&-\int\limits_B{2\la \cof\grad u\cd \grad\eta\;dx}\;\mb{for all}\;\eta\in C_c^\infty(B,\R^2).
\end{align*}
By the latter equation, Lemma \ref{Lem:Uni.Tech.1}.(v), \eqref{eq:HFU.IC.G.02} and Hölder's inequality we can estimate the rightmost term in \eqref{eq:HFU.IC.G.04} from below by
\begin{align}
\int\limits_B{\nu(x)|\grad u|^{p-2}\grad u\cd\grad\eta\;dx}=&-\int\limits_B{2\la \cof\grad u\cd \grad\eta\;dx}&\nonumber\\
=&\int\limits_B{R((\cof\grad\eta)\grad\la)\cd\eta\;\frac{dx}{R}}&\label{eq:HFU.IC.G.04a}\\
\ge&-\frac{n}{2}\int\limits_B{\nu(x)|\grad u|^{p-2}|\grad\eta||\eta|\;\frac{dx}{R}}&\nonumber\\
\ge&-\frac{n}{2}\left(\int\limits_B{\nu(x)|\grad u|^{p-2}|\grad\eta|^2\;dx}\right)^{\frac{1}{2}}\left(\int\limits_B{\nu(x)|\grad u|^{p-2}|\eta|^2\;\frac{dx}{R^2}}\right)^{\frac{1}{2}}&\nonumber\\
\ge&-\frac{1}{2}\int\limits_B{\nu(x)|\grad u|^{p-2}|\grad\eta|^2\;dx}.&\nonumber
\end{align}
For the last step we made used of the following version of the Fourier-estimate given by
\begin{equation} n^2\int\limits_B{|\sigma|^2|\eta|^2\;dx} \le \int\limits_B{|\sigma|^2|\eta_{,\th}|^2\;dx}.\label{eq:HFU.FE.2}\end{equation}
This is indeed true, for this sake, first assume that  $\sigma\in C^\infty(B,\R^{2\ti2}).$ An application of \eqref{eq:Uni.HFU.1}, the product rule and \eqref{eq:HFU.C.G.01} yields,   
\begin{align} n_*^2\int\limits_B{|\sigma\eta|^2\;dx} \le& \int\limits_B{|(\sigma\eta)_{,\th}|^2\;dx}&\nonumber\\
\le& \int\limits_B{|\sigma|^2|\eta_{,\th}|^2\;dx}+\int\limits_B{|\sigma_{,\th}|^2|\eta|^2\;dx}&\nonumber\\
\le& \int\limits_B{|\sigma|^2|\eta_{,\th}|^2\;dx}+l^2\int\limits_B{|\sigma|^2|\eta|^2\;dx}.&
\label{eq:HFU.FE.2a}
\end{align}
Absorbing the rightmost term of the latter expression into the LHS yields \eqref{eq:HFU.FE.2}. Additionally, \eqref{eq:HFU.FE.2a} justifies its own upgrade. It holds for any $\sigma\in L^2(B,\R^{2\ti2}),$ satisfying \eqref{eq:HFU.C.G.01}. The latter is true since $2\le p\le\infty$ and $\sigma\in L^{\frac{2p}{p-2}}(B,\R^{2\ti2}).$\\

ii) Similarly.\\

\begin{re} 
By thinking some more about the inequality by Sivaloganathan and Spector, one might be able to allow more general functionals of the form
\[f(x,\xi)=\frac{\nu(x)}{p}|\xi|^{p-2}M(x)\xi\cd\xi,\;\mb{for all}\; \xi\in\R^{2\ti2}.\]
\end{re}
\begin{re} 1. Notice that \eqref{eq:HFU.C.G.01} is especially satisfied if $p=2$ and $\nu(x)=\nu(R).$ So condition \eqref{eq:HFU.C.G.01} 
can be thought off as a natural extension of this fact to the case, where $p$ might be arbitrary and $\sigma(x)$ depends on $x$ instead of $R.$\\
2. It is worth mentioning, that despite the fact that the sets $F_{n_*}^{p,\sigma,c}$ and $F_{0,m_*}^{p,\sigma,c}$ depend on $\sigma$ it remains true that if $n_*=0$ or $m_*\in\{0,1\},$ one gets uniqueness in the full class $\A^{p,c}.$ Indeed, there are two cases to consider. Firstly, let $n_*=0$ or $m_*\le1$ s.t. $n=0$ or $m=0.$ Then $\grad\la\equiv0$ and by  \eqref{eq:HFU.IC.G.04} and \eqref{eq:HFU.IC.G.04a} one obtains
\begin{align*}
E(v)-E(u)\ge\int\limits_B{\frac{\nu(x)}{2}|\grad u|^{p-2}|\grad\eta|^2\;dx}.\end{align*}
This implies that $u$ is a global minimizer. Moreover, realizing that $|\grad u|^{p-2}>0$ a.e. one can conclude that it has to be the unique one, if there exists a nontrivial set, where $\nu>0$.\\
In the case when $m_*=1$ and $m=1,l=0,$ then  $\sigma_{,\th}=0$ and \eqref{eq:HFU.FE.2} can be applied with $m_*=m=1$ yielding the proof. \\
3. The result from above could be generalized to sets $\Om\ss \R^2,$ which are open, bounded, with $\p\Om\in C^{0,1}$ and homeomorphic to the unit ball with a Bilipschitz map connecting the two. Then the constants in \eqref{eq:HFU.IC.G.02} and \eqref{eq:HFU.IC.G.03} will additionally depend on the Lipschitz constants. 
\end{re}

\section{Uniqueness in compressible elasticity}
\label{sec:4.3}
\subsection{A High frequency uniqueness result for a polyconvex functional}
\label{sec:4.3.1}
It is natural to ask wether high frequency uniqueness can also be shown in compressible elasticity. We were able to give such a result for the polyconvex functional we were studying in the  Chapters $2$ and $3.$\\
 
Firstwe recall that for any $u\in W^{1,2}(B,\R^2)$ the functional was given by
\begin{equation}
I(u)=\int\limits_{B}{\frac{1}{2}|\grad u|^2+\rho(\det \grad u)\; dx}.
\label{eq:HFU.UC.1.1}
\end{equation}
Recall further that $\rho\in C^{k}(\R)$ was defined by
\begin{equation}
\rho(s)=\left\{\begin{array}{ccc}
0& {\mbox{if}}& s\le0,\\
\rho_1(s)& {\mbox{if}}& 0\le s\le s_0,\\
\gamma s+\ka&{\mbox{if}}& s_0\le s,
\end{array}
\right.
\label{eq:HFU.UC.1.2}
\end{equation}
for some constants $\gamma>0,$ $\ka<0$ and $s_0\ge0.$ 
Moreover, $\rho_1:[0,s_0]\rightarrow \R$ must be a convex $C^k-$function on $[0,s_0],$ satisfying the boundary conditions $\rho_1(0)=0$ and $\rho_1(s_0)=\gamma s_0+\ka$ and the connections need to be in such a way that $\rho$ is $C^k-$everywhere. Then $\rho$ is convex and $C^k$ on $\R$.\\

After this small repetition we can state the result.\\

\begin{thm} [High frequency uniqueness] Let $0<\ga<\infty,$ $\rho\in C^2(\R)$ satisfying \eqref{eq:HFU.UC.1.2}, assume $u_0\in L^2(B,\R^2)$ on the boundary and let $u\in \A_{u_0}$ satisfying
\[\|\grad^2 u(x)R\|_{L^\infty(B,\R^{8},\frac{dx}{R})}<\infty,\]
 be a stationary point of the energy $I$ as defined in \eqref{eq:HFU.UC.1.1}.\\
 
Then the following statements are true:\\
i) \textbf{(purely high modes.)} There exists $n\in\N$ s.t.
\begin{equation}\|\grad(\rho'(d_{\grad u}))R\|_{L^\infty(B,\R^2,\frac{dx}{R})}\le n,\label{eq:HFU.UC.1.10}\end{equation}
and $u$ is a minimizer of $I$ in the subclass
\[\F_n=\left\{v\in \A_{u_0}|\;\eta=v-u\in W_0^{1,2}\left(B,\R^2,\frac{dx}{R}\right) \;\mb{and}\; \eta=\sum\limits_{j\ge n}\eta^{(j)} \right\}.\]
Moreover, if the inequality is strict, i.e. $|\grad(\rho'(d_{\grad u}))R|< n$ on a non-trivial measurable set $U\ss B$ wrt. $\frac{dx}{R},$ then $u$ is the unique minimizer in $\F_n$.\\

ii) \textbf{($0-$mode and high modes.)} There exists $m\in\N$ s.t.
\begin{equation}\|\grad(\rho'(d_{\grad u}))R\|_{L^\infty(B,\R^2,\frac{dx}{R})}\le \frac{\sqrt{3}m\nu}{2\sqrt{2}},\label{eq:HFU.UC.1.11}\end{equation}
and $u$ is a minimizer of $I$ in the subclass 
\[\F_{0,m}=\left\{v\in \A|\;\eta=v-u\in W_0^{1,2}\left(B,\R^2,\frac{dx}{R}\right) \;\mb{and}\; \eta=\eta^{(0)}+\sum\limits_{j\ge m}\eta^{(j)} \right\}.\]
Moreover, if the inequality is strict, i.e. $|\grad(\rho'(d_{\grad u}))R|< \frac{\sqrt{3}m\nu}{2\sqrt{2}}$ on a non-trivial measurable set $U\ss B$ wrt. $\frac{dx}{R},$ then $u$ is the unique minimizer in $\F_{0,m}$.
\end{thm}

\textbf{Proof:}\\
i) Note, that if $u\in \A_{u_0}$ is defined like above then we can control the quantity, we are interested in, via
\[\|\grad(\rho'(d_{\grad u}))R\|_{L^\infty(B,\R^2,\frac{dx}{R})}\le\|\rho''\|_{L^\infty(\R)}\|\grad^2u(x)R\|_{L^\infty(B,\R^{8},\frac{dx}{R})} \|\grad u\|_{L^\infty(B,\R^4,\frac{dx}{R})}<\infty,\]
where $\|\rho''\|_{L^\infty(\R)}<\infty$ is satisfied, since $\rho''\in C^0(\R)$ and has compact support, which can be deduced from it's definition \eqref{eq:HFU.UC.1.2}.\\

This shows that there has to be an $n\in\N$ s.t. 
\[\|\grad(\rho'(d_{\grad u}))R\|_{L^\infty(\frac{dx}{R})}\le n.\]
Let $v\in\F_n$ and $\eta=v-u\in W_0^{1,2}(B,\R^2)$ and $\eta=\sum\limits_{j\ge n}\eta^{(j)}.$
By using the subdifferential inequality, the expansion of the determinant and the assumption that $u$ is a stationary point of $I,$ we get
\begin{align*}
I(v)-I(u)\ge & \frac{1}{2}D(\eta)+\int\limits_B{\rho'(d_{\grad u})d_{\grad \eta}\;dx}.&
\end{align*}
Then like above we get
\begin{align*}I(v)-I(u)=&\frac{1}{2}D(\eta)-\frac{1}{2}\int\limits_B{R((\cof\grad\eta)\grad(\rho'(d_{\grad u})))\cd\eta\;\frac{dx}{R}}&\\
\ge&\frac{1}{2}D(\eta)-\frac{1}{2}\|\grad(\rho'(d_{\grad u}))R\|_{L^\infty(\frac{dx}{R})}\int\limits_B{|\grad\eta||\eta|\;\frac{dx}{R}}&\\
\ge&0,&
\end{align*}
where we made use of the fact that $\eta$ only contains of Fourier-modes $\ge n$. \\

One can argue analogously for $(ii).$ However, one needs to follow, again, the proof of Theorem \ref{Thm: Uni.SPC.1}.\\

\vspace{0.5cm}

The next lemma presents two simple conditions on the determinant which, when satisfied, guarantee a unique global minimizer.

\begin{lem} [Uniqueness conditions] Let $0<\ga<\infty,$ $\rho\in C^1(\R)$ satisfying \eqref{eq:HFU.UC.1.2}, $u_0\in L^2(B,\R^2)$ on the boundary and assume $u\in \A_{u_0}$ be a stationary point of the energy $I$ as defined in \eqref{eq:HFU.UC.1.1}.\\

Assume, additionally, that\\
(i) $d_{\grad u}\ge s_0$ a.e. in $B,$ or\\ 
(ii) there exists a constant $c\in\R$ s.t. $d_{\grad u}=c$ a.e. in $B.$\\

Then $u$ is the unique global minimizer of $I.$ 

\end{lem}

\textbf{Proof:}\\
(i): Note that if $d_{\grad u}\ge s_0$ a.e. then $\rho'(d_{\grad u})=\ga$ a.e. in $B$ and hence
\begin{align*}
I(v)-I(u)\ge & \frac{1}{2}D(\eta)+\int\limits_B{\rho'(d_{\grad u})d_{\grad \eta}\;dx}&\\
= &\frac{1}{2} D(\eta)+\ga\int\limits_B{d_{\grad \eta}\;dx}= \frac{1}{2}D(\eta),&
\end{align*}
where we used in the last line that $d_{\grad \eta}$ is a Null-Lagrangian with $0$ boundary data.\\

Similarly for (ii).\\

\begin{re}
Realize, that condition (ii) is satisfied if $u=u_0\equiv Ax+b$ a.e. in $B.$ Sadly, in case of the BOP map, which we constructed in Chapter $3$, for the case $N\ge2,$ $\det \grad u$ starts at the origin with $\det \grad u(0)=0$ and increases montonically wrt. an increasing radial component $R,$ however, not linearly if $N\ge2$. Hence neither (i) nor (ii) is satisfied and the above lemma does not apply to the BOP map. 
\end{re}\vspace{0.25cm}

\subsection{Generalisation to polyconvex $p-$growth functionals and revisiting a result of Sivaloganathan and Spector}
\label{sec:4.3.2}
For $u_0\in L^p(B,\R^2)$ with $2\le p\le\infty$ define 
\[
\A^p=\{u\in W^{1,p}(B,\R^2): \; u_{|{\p B}}=u_0\}.
\]
Define
\begin{equation}E(u)=\int\limits_B{\Phi(x,\grad u)\;dx}\label{eq:HFU.UC.G.0}\end{equation}

where $\Phi$ is of the specific form
\[\Phi(x,\xi)=\frac{\nu(x)}{p}|\xi|^p+\Psi(x,\xi,\det\xi),\]
and where we want $(\xi,d)\mapsto \Psi(x,\xi,d),$ to be convex for a.e.\! $x \in B,$ making $\Psi(x,\cd)$ a convex representative of a polyconvex function a.e.  in $B$. \\
Moreover,  $\nu\in L^\infty(B)$ is supposed to satisfy $\nu(x)\ge0$ a.e.\! in $B.$ We want the function $\nu$ to be optimal, which means there can be no term of the form  $a(x)|\xi|^p$ in $\Psi.$ However, it is not necessary for $\Psi$ to be non-negative, as the example below shows.\\
We want $\Phi$ to be of $p-$growth, so we suppose that there exits $C\in L^\infty(B)$ with $C(x)\ge0$ a.e. in $B$ s.t.
\[0\le \Phi(x,\xi)\le \frac{C(x)}{p}(1+|\xi|^p)\;\mb{for all} \;\xi\in \R^{2\ti2}\;\mb{and a.e.} \; x\in B. \]
All of the above combined guarantees that $x\mapsto \Phi(x,\grad u(x))\in L^1(B)$ for any $u\in W^{1,p}(B,\R^2)$ and hence the corresponding energy is finite. Furthermore, since $\nu$ is allowed to take on the value $0,$ the integrand could indeed disappear for some $x \in B.$ Therefore, $\Phi$ is asymptotically polyconvex. \vspace{0.5cm}

 Now we can state the main theorem of this paragraph.
 
\begin{thm} [High frequency uniqueness] Let $2\le p\le\infty,$ assume $u_0\in L^2(B,\R^2)$ on the boundary and let $u\in \A_{u_0}^p$ to be a stationary point of $E$ as given in \eqref{eq:HFU.UC.G.0}. Furthermore, let $\sigma(x):=\nu(x)|\grad u(x)|^{\frac{p-4}{2}}\grad u(x)\in L^{\frac{2p}{p-2}}(B,\R^{2\ti2})$ and assume that there exists $l\in \N$ s.t.
\begin{equation}|\sigma,_\th|\le l|\sigma| \mb{for a.e.}\;x\in B\label{eq:HFU.UC.G.01} \end{equation}
holds.\\

Then the following statements are true:

i) \textbf{(purely high modes.)}  Assume there exists $n\in\N$ s.t.
\begin{equation}|\grad_x \p_d\Phi(x,\grad u,d_{\grad u})R|\le n\nu(x)|\grad u|^{p-2}\mb{for a.e.}\;x\in B \;\mb{wrt.}\;\frac{dx}{R}.\label{eq:HFU.UC.G.1} \end{equation}

Then $u$ is a minimizer of $E$ in the subclass 
\[\F_{n_*}^{p,\sigma}=\left\{v\in \A_{u_0}^p|\;\eta=v-u\in W_0^{1,p}(B,\R^2) \;\mb{and}\; \sigma\eta=\sum\limits_{j\ge {n_*}}(\sigma\eta)^{(j)} \right\},\]
where $n_*^2:=n^2+l^2.$
Moreover, if the inequality is strict on a non-trivial set wrt. $\frac{dx}{R},$ then $u$ is the unique minimizer in $\F_{n_*}^{p,\sigma}$.\\

ii) \textbf{($0-$mode and high modes.)} Assume there exists $m\in\N$ s.t.
\begin{equation}|\grad_x \p_d\Phi(x,\grad u,d_{\grad u})R|\le \frac{\sqrt{3}m\nu(x)|\grad u|^{p-2}}{2\sqrt{2}}\mb{for a.e.}\;x\in B \;\mb{wrt.}\;\frac{dx}{R}.\label{eq:HFU.UC.G.2}\end{equation}
Then $u$ is a minimizer of $E$ in the subclass 
\[\F_{0,m_*}^{p,\sigma}=\left\{v\in \A_{u_0}^p|\;\eta=v-u\in W_0^{1,p}(B,\R^2) \;\mb{and}\; \sigma\eta=(\sigma\eta)^{(0)}+\sum\limits_{j\ge m_*}(\sigma\eta)^{(j)} \right\},\]
where $m_*^2:=m^2+l^2.$
Moreover, if the inequality is strict on a non-trivial set wrt. $\frac{dx}{R},$ then $u$ is the unique minimizer in $\F_{0,m_*}^{p,\sigma}$.
\label{thm:HFU.UC.G.1}
\end{thm}
\textbf{Proof:}\\
i) We start again by the standard expansion
\begin{align}
E(v)-E(u)=&\int\limits_B{\frac{\nu(x)}{p}(|\grad u+\grad \eta|^{p}-|\grad u|^{p})}&\nonumber\\
&+{\Psi(x,\grad u+\grad \eta,\det\grad u+\grad \eta)-\Psi(x,\grad u,\det\grad u)\;dx}&\nonumber\\
\ge&\int\limits_B{\frac{\nu(x)}{2}|\grad u|^{p-2}|\grad\eta|^2+\nu(x)|\grad u|^{p-2}\grad u\cd\grad\eta}&\nonumber\\
&+{\p_\xi\Psi(x,\grad u,\det\grad u)\cd\grad\eta+\p_d\Psi(x,\grad u,\det\grad u)(d_{\grad\eta}+\cof\grad u\cd\grad\eta)\;dx},&\nonumber
\end{align}
where we used the subdifferential inequality for $\Psi$, and the following inequality\footnote{see, \cite{SS18}, Prop A.1, with $\sigma=0.$}
\[\frac{1}{p}|b|^{p}\ge\frac{1}{p}|a|^{p}+|a|^{p-2}a(b-a)+\frac{1}{2}|a|^{p-2}|b-a|^{2}.\] 

The ELE is given by
\begin{align*}
0=&\int\limits_B{\grad_\xi\Phi\cd\grad\eta\;dx}&\\
=&\int\limits_B{\nu(x)|\grad u|^{p-2}\grad u\cd\grad\eta+\p_\xi\Psi(x,\grad u,\det\grad u)\cd\grad\eta}&\\
&+{\p_d\Psi(x,\grad u,\det\grad u)\cof\grad u\cd\grad\eta\;dx}\;\mb{for all}\;\eta\in C_c^\infty(B,\R^2).&
\end{align*}
Hence, we have
\begin{align}
E(v)-E(u)=&\int\limits_B{\frac{\nu(x)}{2}|\grad u|^{p-2}|\grad\eta|^2+\p_d\Psi(x,\grad u,\det\grad u)d_{\grad\eta}\;dx}.&
\label{eq:HFU.UC.G.21}
\end{align}
The proof is completed, as before, by making use of Lemma \ref{Lem:Uni.Tech.1}.(v), \eqref{eq:HFU.UC.G.1}, Hölder's inequality and \eqref{eq:HFU.FE.2} to estimate the rightmost term in \eqref{eq:HFU.UC.G.21} by
\begin{align*}
\int\limits_B{\p_d\Psi d_{\grad\eta}\;dx}=&-\frac{1}{2}\int\limits_B{R((\cof\grad\eta)\grad_x\p_d\Psi)\cd\eta\;\frac{dx}{R}}&\\
\ge&-\frac{n}{2}\int\limits_B{\nu(x)|\grad u|^{p-2}|\grad\eta||\eta|\;\frac{dx}{R}}&\\
\ge&-\frac{1}{2}\int\limits_B{\nu(x)|\grad u|^{p-2}|\grad\eta|^2\;dx}.&\\
\end{align*}

\begin{re}
1. We can argue similarly to the incompressible case to see that, the cases $n_*=0$ or $m_*\in\{0,1\},$ indeed, imply uniqueness in the full class $\A^{p}.$ \\
2. Again a generalization to sets $\Om\ss \R^2,$ which are open, bounded, with $\p\Om\in C^{0,1}$ and homeomorphic to the unit ball with a Bilipschitz map connecting the two is possible.\\ 
Then the constants in \eqref{eq:HFU.UC.G.1} and \eqref{eq:HFU.UC.G.2} will additionally depend on the Lipschitz constants. 
\end{re}

How could the conditions \eqref{eq:HFU.UC.G.1} and \eqref{eq:HFU.UC.G.2} be verified in practice?\\

\textbf{Example [Alibert-Dacorogna-Marcellini (ADM) integrand]:}\\
Here we want to apply the above result to the ADM-integrand as given in \eqref{eq:Uni.ADM.1}. It is well known that the behaviour of this integrand depends on $\al,$ in particular, for $0\le\al<\frac{1}{2}$ it is uniformly convex while for $\frac{1}{2}\le\al\le 1$ it is genuinely polyconvex. In this case, the above condition reduces to a Reverse-Poincare-type estimate.\footnote{Here Reverse-Poincare-type estimate has the loose meaning of it beeing any inequality that bounds a higher derivative by a smaller derivative. Here the estimates are taken in a pointwise fashion, classically expressions like (Reverse-) Poincare estimates refer to inequalities which compare the $L^p-$norms rather than the pointwise quantities.}\\ 

\begin{cor}Let $u_0\in L^2(B,\R^2)$ on the boundary and let $u\in \A_{u_0}^4$ to be a stationary point of $E$ as given in \eqref{eq:HFU.UC.G.0}, where
\begin{equation}\Phi_\al(\xi)=\frac{1}{4}|\xi|^4-\frac{\al}{2}|\xi|^2\det\xi\;\mb{for all}\;\xi\in\R^{2\ti2}, \label{eq:Uni.ADM.1}\end{equation}
with $\al\in(0,1].$ Furthermore, $\nu=\frac{1}{4}$ and let $\sigma(x):=\grad u(x)\in L^4(B,\R^{2\ti2})$ and assume that there exists some $l\in \N$ s.t.
\begin{equation}|\grad u,_\th(x)|\le l|\grad u(x)| \mb{for a.e.}\;x\in B.\label{eq:HFU.UC.G.011} \end{equation}

Then the following statements are true:\\

i) \textbf{(purely high modes.)}  Assume there exists $n\in\N$ s.t.
\begin{equation}|\grad^2u(x)R|\le \frac{n}{4\al}|\grad u(x)|\;\mb{for a.e.}\;x\in B \;\mb{wrt.}\;\frac{dx}{R}.\label{eq:HFU.UC.G.3} \end{equation}

Then $u$ is a minimizer of $E$ in the subclass $\F_{n_*}^{4,\sigma},$ where $n_*^2=n^2+l^2.$
Moreover, if the inequality is strict on a non-trivial set wrt. $\frac{dx}{R},$  then $u$ is the unique minimizer in $\F_{n_*}^{4,\sigma}$.\\

ii) \textbf{($0-$mode and high modes.)} Assume there exists $m\in\N$ s.t.
\begin{equation}|\grad^2u(x)R|\le \frac{\sqrt{3}m}{8\sqrt{2}\al}|\grad u(x)| \;\mb{for a.e.}\;x\in B \;\mb{wrt.}\;\frac{dx}{R}.\label{eq:HFU.UC.G.4}\end{equation}
Then $u$ is a minimizer of $E$ in the subclass $\F_{0,m_*}^{4,\sigma},$ where $m_*^2=m^2+l^2.$
Moreover, if the inequality is strict on a non-trivial set wrt.\! $\frac{dx}{R},$ then $u$ is the unique minimizer in $\F_{0,m_*}^{4,\sigma}$.
\end{cor}

\textbf{Comparison to a result by Sivaloganathan and Spector}

A uniqueness result, similarly to Theorem \ref{thm:HFU.UC.G.1}, has been recently established in the non-linear elasticity(NLE) setting by Sivaloganathan and Spector.\\
 It might be of interest for comparison reasons to include their result, given in \cite{SS18}. We will give a version of this statement, however, we will not give the most general one and we will skip some of the details, to keep it simple.\\

 \begin{thm}[\cite{SS18}, Theorem 4.2]
Let $\Om\ss\R^2$ be open, bounded, connected and $\p\Om\in C^{0,1},$  assume $u_0\in C(\ol{\Om}, \R^2)$ on the boundary and let the set of admissible maps be defined by \[
\A=\{u\in W^{1,2}(\Om,\R^2)\cap C(\ol{\Om}, \R^2): \;  \det\grad u>0 \;\mb{a.e. and}\; \; u_{|{\p \Om}}=u_0\}.
\] 
The energy $E$ is given by \eqref{eq:HFU.UC.G.0}. The integrand is assumed to be uniformly polyconvex and given by
\[\Phi(x,\xi)=\frac{\nu(x)}{p}|\xi|^p+\Psi(x,\xi,\det\xi),\]
for some $2\le p<\infty$ a.e.\! $x \in \Om$ and any $\xi\in \R^{2\ti 2}.$ We want $(\xi,d)\mapsto \Psi(x,\xi,d),$ to be convex for a.e.\! $x \in \Om,$ making $\Psi(x,\cd)$ a convex representative of a polyconvex function a.e.  in $\Om$. \\
Moreover,  $\nu\in L^\infty(\Om)$ is supposed to satisfy $\nu(x)>\nu_0>0$ a.e.\! in $\Om.$
Furthermore, we assume $\Psi$  and $\Phi$ to be smooth enough and $\Phi$ to be frame indifferent.\\
 
Finally, let $u\in \A$ be a weak solution to the equilibrium equation and assume $u$ satisfies

\begin{equation}|\p_d\Phi(x,\grad u,d_{\grad u})R|\le \nu(x)|\grad u|^{p-2}\mb{for a.e.}\;x\in \Om.
\label{eq:Uni.SivSP.1}
\end{equation}

Then $u$ is a global minimizer of $E$ in $\A.$\\
Moreover, if the inequality is strictly satisfied on the full set $\Om,$ then $u$ is the unique minimizer in $\A.$\\

\end{thm}

\begin{re}
While the condition given by Sivaloganathan and Spector proves useful in the NLE-setting as they were able to show on various examples, the condition seems rather limiting in the finite elasticity setting. Indeed, consider, again the ADM-example and assume \eqref{eq:Uni.SivSP.1} instead. Then the condition reduces to $\al\le\frac{1}{2},$ yielding uniqueness in the regime $0<\al<\frac{1}{2},$ where the integrand is uniformly convex and uniqueness is known anyway. \eqref{eq:Uni.SivSP.1} seems too restrictive to allow for a genuinely polyconvex functional.\\
 
One could argue that our assumption seems unnatural. From what we learned while studying the ADM-example, one should expect a condition of the form
\begin{equation}|\grad_x\p_d\Phi R|\lesssim|\grad u|^{p-4}|\grad^2 u|.\end{equation}
But this would again just reduce, similarly to the previous case, to $\al\le\frac{1}{4}.$ \\

So, for us, the condition  
\begin{equation}|\grad_x\p_d\Phi R|\le\frac{\sqrt{3}}{2\sqrt{2}}|\grad u|^{p-2}\label{eq:HFU.UC.1.12}\end{equation}
seems most natural and clearly reduces to our condition, we discussed, in the previous paragraph. Indeed, for $p=2,$ consider $\Psi(d)=\rho(d)$ as defined in \eqref{eq:HFU.UC.1.2}. Then condition \eqref{eq:HFU.UC.1.12} becomes $|\grad_x(\rho'(d_{\grad_u})) R|\le\frac{\sqrt{3}}{2\sqrt{2}},$ revisiting \eqref{eq:HFU.UC.1.11} with $m=1$.  As we have seen the latter condition reduces in the ADM-example to a Reverse-Poincare-type inequality. It might be difficult to find solutions satisfying the Reverse-Poincare inequality, even though it is a fairly natural inequality in elliptic situations. Despite all that at least it applies to any $\al\in(0,1]$.\\

\end{re}

\clearpage{\pagestyle{empty}\cleardoublepage} 
\chapter{Conclusions and Outlook}

\section{Conclusions and Outlook for the chapters 2 and 3.}
\textbf{Conclusions for chapter 2:}\\
As a first result we have seen in section \ref{sec:2.2} that all stationary points need to be everywhere locally Hölder continuous for general boundary conditions. A similar result for the nonlinear elasticity case is not known (maybe not even true). Such a result in the NLE case is only known for special classes of functions, minimizers rather than stationary points and particular boundary conditions, like the positive twist maps considered in \cite{JB17} or the scenario explored in \cite{BK19}.\\

In paragraph \ref{sec:2.3} we have seen that as long as the integrand is uniformly convex $0<\ga<1$ then the stationary points need to be in $W_{loc}^{2,2}.$ This argument relied heavily on De Maria's seminal work.
It was a subtle fact, that for arbitrary $0<\ga<\infty$ Meyers' Theorem could not be applied to gain slightly higher integrability going to $W_{loc}^{2,2+\d}$ for some $\d>0.$ \\ 

Overcoming this type of issue is the content of section \ref{sec:2.4} For this sake, under the additional assumptions that stationary points are in $W_{loc}^{2,2}$ and possess Hölder-continuous Jacobians, a Reverse Hölder Inequality has been established by making use of advanced measure theoretic concepts and exotic Hardy Space Theory. This showed that for any $0<\ga<\infty$ stationary points need to be in $W_{loc}^{2,2+\d}$  for some $\d>0.$ This turned out to be enough to start the well-known `bootstrapping' argument, guaranteeing maximal smoothness to those stationary points.\\

\textbf{Conclusions for chapter 3:}\\
In section \ref{sec:3.1} we have seen that if $\Om=B$ and $u=\Id$ on $S^1$ then $u=\Id$ in $\ol{B}$ is the unique global minimizer, for any $0<\ga<\infty$.\\
However, it remains open if in addition to $u=\Id$ there are other stationary points (local minimizers). There could even be another radially symmetric maps $u(x)=r(R)e_R(\th)$ s.t. $u\not=\Id$ is a solution to the BVP \eqref{eq:3.0.1} for $M=1$.\\

In Section \ref{sec:3.2} and \ref{sec:3.3} we guaranteed, for any $0<\ga<\infty$, the existence of rsMc. stationary points that are at least of class $C^1.$ In the case of $\rho$ lifting-off delayed we showed that they need to be $C^\infty.$ For $0<\ga<1$ one of the rsMc. stationary points we discussed needs to be the unique global minimizers wrt. $\A_M$ and $u\in C^\infty.$  If $\ga\ge 1$ there might be additional stationary points, which might have a rough profile. We can not stress the meaning of this enough. This means for our integrand $I$ and M-covering bc. we can always guarantee a $C^1$ stationary point and even one that it is of class $C^\infty$ if $\rho$ lifts-off delayed. These are absolutely non-trivial results.\\

These results were obtained by implementing the BOP technique in our situation. We could see that some of the arguments have simplified significantly compared to the works by Baumann, Owen, Phillips, Yan and Bevan since we deal with an FE model and one can more directly make use of the ELE and does not need to 'tiptoe' around it using the weaker energy momentum equation. For this reason many calculations were much more concrete. Nevertheless, we also experienced, that our integrand imposes less strict conditions on the minimizers, allowing them, in the M-covering case to be in a wide range of possible classes somewhere between $C^1$ and $C^\infty.$\\

\textbf{Open questions and possible future research:}\\
\begin{enumerate}
\item Questions related to the general regularity theory:
\begin{enumerate}
\item Rather obvious but non the less intriguing questions, are can one generalise our results to p-growth functionals and to 3 or higher dimensions?

\item Can De Maria's method or any other method for that matter be used to show that $u$ is either in  $W_{loc}^{2,q}$ for some $1\le q<2$ or in a Fractional space $W_{loc}^{1+s,2}$ for some $s\in (0,1)$ or in a relevant Besov space in the regime $\ga\ge1$? If not in general then at least under some additional assumption. For instance assuming $u\in C^1\cap W^{1,2}$ with this additional information one might be able to gain some higher regularity or integrability. Here one might even reach $W_{loc}^{2,2}.$ Since the radial symmetric $M-$covering stationary points we discussed are of class $ C^1\cap W^{1,2}$ such a result would improve the regularity of these points.

\item Can one make use of some of the methods we used to gain insights in other situations like the NLE models? For instance, J. Bevan introduced the concept of positive twist maps in \cite{JB17} and studied them in the NLE case. He showed that minimizers in that class need to be locally Hölder continuous. The question now is can we add to that and show that these types of minimizers need to be even more regular? Furthermore, is it possible to establish a partial regularity for these types of maps? This would be intriguing since partial regularity is completely unknown in the NLE setting.

\item Another pressing question for further research might be: what other possible uses are there for the concept of positive twist in the Calculus of Variations?\\
For instance, is there a positive twist property in three dimensions and is one able once again to obtain some regularity for this class of functions, similarly to the work presented in \cite{JB17}?
Is it possible to identify flows, which generate positive twist maps? This could link the concept in a unique sense to Fluid Mechanics.

\end{enumerate}
\item Questions related to BOP-Theory:
\begin{enumerate}
\item We have seen that the radially symmetric MC stationary points can posses various possible shapes (immediate/delayed lift-off). We would like to determine, which of these (is) are the minimizer(s) at least w.r.t. the class $\A_r^M$? Does the minimizer change it's shape when $\ga$ passes through $1,$ from the regime $\ga<1$ to the other $\ga\ge1$? Note, for an immediate lift-off $\rho$ this might also imply a change in regularity!\\

\item Furthermore, we would like to narrow down the regularity of $r_\d$ in the case, where $\rho$ and $r_d$ are both lifting-off immediately. However, it seems to be very difficult (or it might even be impossible) to exclude non-smooth $C^1$ (but no better) solutions in the $M-$covering case. Indeed, the BOP Theory seems exhausted at this stage.To go beyond the presented work, one needs to consider even higher derivatives of the quantities discussed in \ref{sec:3.2} and \ref{sec:3.3}, which becomes more and more complicated.
\end{enumerate}
\end{enumerate}

\section{Conclusions and Outlook for chapter 4}
\textbf{Conclusions:}\\
We have introduced conditions in finite elasticity, which when satisfied can guarantee uniqueness of global minimizers to certain energies. These conditions turned out to be analogously to the ones introduced by Sivaloganathan and Spector in the NLE-setting. These conditions have surprisingly turned out to be at the heart of the matter in various ways. They connected and generalised many things, which were, until now, only known for much more specific situations. The central idea, following in the tradition of works by J. Bevan, of this method is actually computing the pressure (or some of its norms) and making use of it. Usually only more abstract results are obtained, like guaranteeing the existence of the pressure in some space, but there is no further usage of the pressure. \\

In particular, we were able to provide a counterexample to regularity, which also can be seen as a contribution to the understanding of the double-covering problem. Indeed, we have constructed a functional, which, roughly speaking, is fairly close to the Dirichtlet energy, however,  s.t. the map $u=\frac{R}{\sqrt{N}}e_{NR}+b$ is the unique global minimizer. It also answers the question why the Double Covering Problem turns out to be so difficult, because it is a high pressure problem and our method only guarantees high frequency uniqueness rather than actual uniqueness.\\ Uniqueness questions in high pressure situations remain completely open, a positive answer to the Double Covering Problem would be a first step in the understanding of such situations.\\ 

\textbf{Open questions and possible future research:}
\begin{enumerate}
\item Questions related to the small pressure condition:
\begin{enumerate}
\item Is the prefactor $\frac{\sqrt{3}}{2\sqrt{2}}$ which shows up in \eqref{eq:Uni.SPC.101} optimal for general situations?\\
\item If the latter is true, can this be used to construct a situation, where the energy and the  boundary conditions are such that there are multiple stationary points such that the corresponding pressures satisfy $\la=\frac{\sqrt{3}}{2\sqrt{2}},$ respectively, yielding that all these stationary points are actually global minimizer, constructing an example of non-unique global minimizers in the incompressible case. This did not fully work in the paper by \cite{BeDe20}, because of the continuity of the energy, wrt.\@ some parameter (Making use of the continuity one is able to exclude non-uniqueness in the threshold). However, one might get lucky by setting up a functional, without any dependence on a parameter. 

\item Consider the following situation: $\Om=B\ss \R^2,$ and fix $E=\bb{D}.$ Then taking $u_0=\Id$ on the boundary and assume $u=u_0=\Id$ is a stationary point and it  then $\grad\la_{\Id}=0.$ Now again $E=\bb{D},$ however, we take the double-cover $u_0=u_2$ on the boundary. Then $u=u_2$ is a stationary point and  $\|\grad \la_{u_0}(x)R\|_{L^\infty(B,\R^2,\frac{dx}{R})}=3/2.$ From this and our full analysis, we got the impression that the size of $\grad \la R$ seems to increase with the topological degree of the map imposing the boundary conditions. This inspires the following conjecture:
Given $\Om=B\ss \R^2,$ some energy $E,$ and fix some boundary conditions $u_0$ and assume $u\in \A_{u_0}$ to be a stationary point of $E$. Then the corresponding pressure $\la_{u}$ satisfies\footnote{To make this unmistakable clear, this is a conjecture we do not claim anything here. A relation between these quantities seems plausible. However, the dimensions could be different, or maybe one is only able to show an upper but no lower bound. We don't have any more knowledge, however, it seems worse to us mentioning it, since a result like this could turn out to be very useful.}
\[\|\grad \la_{u}(x)R\|_{L^\infty(B,\R^2,\frac{dx}{R})}\sim  E(u)\deg (u_0).\]
\end{enumerate}
\item High frequency uniqueness:
\begin{enumerate}
\item One can raise the question, what would be the analogs of these in higher dimensions, say, $\mb{dim}\!\!\!=\!3.$ This is, in particular, interesting in case of the compressible elasticity case since the major difference between the 3D and the 2D case is that polyconvex functionals are allowed to depend additionally on $\cof\grad u.$ It is an intriguing question how this would affect, for instance, the conditions given in Theorem \ref{thm:HFU.UC.G.1}. \\

\item To get an understanding of the novel notion of high frequency uniqueness and how rare is it for a system to satisfy such a result, it is necessary to get much more knowledge on this notion. This can be achieved by trying to establish high frequency uniqueness results in various situations, for instance, for quasiconvex functionals or for PDE's like the Navier-Stokes equation or the Monge-Ampere equation.

\end{enumerate} 
\item The counterexample we provided leaves a few open questions:\\
\begin{enumerate}
\item Recall that in our counterexample the integrand $f(x,\xi)$ depends discontinuously on $x.$ \\
It is natural to ask if one can give an example, where $f$ depends smoothly on $x$ or has no explicit $x$ dependency at all.\\
A first step in the direction of constructing such an example, could be to provide generalised versions of the Lemmas \ref{Lem:Uni.NC.1} and \ref{Lem:Uni.NC.2} allowing $\al, \be, \ga, \d$ to depend additionally on $R.$ By this method it might be possible to get an $f$ depending smoothly on $x$. For an integrand with no explicit $x$ dependency, one might need to use a completely different construction, maybe in euclidian coordinates rather than polar coordinates, so one does not have to deal with the in the latter coordinates naturally arising discontinuity at the origin. \\

\item Can one find an incompressible situation, where the boundary condition has no topological twist, i.e. is injective, s.t. one still can find an energy and a corresponding singular minimizer?
 
\item Is it possible to give a counterexample whose singular set has a larger Hausdorff-dimension. (The size of the later is clearly limited by the partial regularity result obtained in \cite{EVGA99}.) On the other hand, can one obtain a partial regularity result in the incompressible case without the degeneracy assumption given in \cite{EVGA99}. Moreover, can one provide an estimate on the Hausdorff-dimension in the incompressible case, which seems completely open.

\item Can one give a counterexample in the compressible case, using our method? This seems difficult, however, it would resolve a long standing open question to show that there is a compressible situation with a singular global minimizer.
\end{enumerate}
\end{enumerate}
\appendices
\chapter{Important mathematical tools}
\section{Trace Theorem}
\label{Ap:A.1}
We first state the classical trace theorem for domains with Lipschitz boundaries.
\begin{thm}[Trace theorem]Let $\Omega\ss\R^n$ be open and bounded with $\p \Om\in C^{0,1}$ and $1\le p<\infty.$ Then there exists an operator
\begin{equation}
T: W^{1,p}(\Omega)\ra L^{p}(\p\Omega)
\label{eq:A1.1}
\end{equation}
such that $T$ satisfies 
\[
Tu=u|_{\p \Om} \;\mb{if}\; u\in W^{1,p}(\Om)\cap C(\bar{\Om}) 
\]
and there exists a constant $C=C(p,\Om)>0$ s.t. 
 \[
\|Tu\|_{L^p(\p\Om)}\le C \|u\|_{W^{1,p}(\Om)}  \;\mb{for all}\; u\in W^{1,p}(\Om). 
\]
\label{thm:Trace theorem}
\end{thm}
\textbf{Proof:} For a proof of this statement, see \cite{A12}, A6.6, p.279-281 or for a proof of the weaker situation of $C^1$ boundaries, see \cite{LE10}, Section 5.5, Theorem 1, p.258.\\

The above theorem has been generalized in various ways, most recently by Ding, in \cite{D96}, who showed that $T: W^{s,2}(\Om)\ra W^{s-\frac{1}{2},2}(\p\Om)$ is a linear bounded operator as long as $\frac{1}{2}<s<\frac{3}{2}.$
If $s>\frac{3}{2}$ then $T: W^{s,2}(\Om)\ra W^{1,2}(\p\Om).$
Intriguingly, (as far as we know) the case $s=\frac{3}{2}$ remains open.\\

\begin{thm}Let $\Omega\ss\R^n$ be open and bounded with $\p \Om\in C^{0,1}$ and $1\le p<\infty.$ Then 
\[
Tu=0 \;\mb{on}\; \p\Om \;\mb{if and only if}\; u\in W_0^{1,p}(\Om).
\]
\label{thm:Trace theorem2}
\end{thm}
\textbf{Proof:} See, \cite{A12}, Lemma A6.10, p.284-285.\\

For $T(u-u_0)=0$ on $\p\Om$ we say $u=u_0$ on $\p\Om$ in the `trace sense'. By the latter theorem, $u=u_0$ on $\p\Om$ in the `trace sense' if and only if 
$u-u_0\in W_0^{1,p}(\Om).$

\section{Weak convergence of the determinant}
In the following lemma we will show weak convergence of the determinant in the $2\times 2-$case. For an extensive discussion, see \cite{D08}, Chapter 8.3.
\begin{lem}Let $\Omega\ss\R^2$ be open and bounded and $2\le p<\infty.$
Assume that
\begin{equation}
u_k\rightharpoonup u\in W^{1,p}(\Omega,\R^2).
\label{eq:A2.2}
\end{equation}
Then
\begin{equation}
\det\grad u_k\rightharpoonup \det\grad u\in D'(\Omega).
\label{eq:A2.3}
\end{equation}
Moreover, if $p>2$ then
\begin{equation}
\det\grad u_k\rightharpoonup \det\grad u\in L^{p/2}(\Omega).
\label{eq:A2.4}
\end{equation}
\label{lem:WCDet}
\end{lem}
\textbf{Proof:}
Multiplying the Jacobian by $\vp\in C_c^\infty(\Om)$ and integrating by parts yields, 
\[
\int\limits_{\Om}{(\det\grad u)\vp\,dx}=\int\limits_{\Om}{(u_1J\grad u_2)\cd\grad \vp\,dx},
\]
where 
$
J=\begin{pmatrix}
0&-1\\
1&0
\end{pmatrix}.
$

Rellich's Embedding Theorem\footnote{see, \cite{A02}, A6.4, p. 274.} implies that $u_k\ra u$ in $L^{p}(\Omega,\R^2)$ and $\grad u_k\rightharpoonup \grad u$ in $L^{p}(\Omega,\R^2).$ Therefore,
\begin{equation}
(u_{k1}J\grad u_{k2})\rightharpoonup (u_1J\grad u_2) \;\mb{in}\;D'(\Omega,\R^2).
\end{equation}
Hence,
\begin{equation}
\det\grad u_k\rightharpoonup \det\grad u \;\mb{in}\; D'(\Omega).
\end{equation}
Moreover, $(u_1J\grad u_2)\in L^{p/2}$ and if $p>2$ then the dual space is $L^{p'}$ s.t. $\frac{2}{p}+\frac{1}{p'}=1,$ in particular $p'<\infty$ and $C_c^\infty$ is dense in $L^{p'}$ and one can upgrade the weak convergence in $D'$ to weak convergence in $L^{p/2}.$ 

\begin{re}Note that the latter step is not possible if $p=2.$ Then the dual of $L^{1}$ is $L^\infty$ and $C_c^\infty$ is not dense in this space.
More generally, in the borderline case of $p=n$ and $u_k\rightharpoonup u$ in $W^{1,n}(\Omega, \R^n),$ $\Om\ss\R^n$, does not necessarily imply $\det\grad u_k\rightharpoonup\det\grad u$ in $L^1(\Omega)$, see \normalfont\cite{D08}\textit{Remark $8.(iii)$ and Example $8.6$ for a counterexample (The candidate is a highly oscillating function). Compare also } \cite{BM84}\textit{, Counterexample 7.3, p.247. Ball and Murat show in addition that this statement is not even true if one restricts to radially symmetric maps, see Counterexample 7.1, p. 245-246.}
\end{re}
The following result omits the missing weak $L^1-$convergence of the determinants, proving that the weaker notion of convergence in $D'$ is sufficient to obtain weak lower semicontinuity if, in addition, we can ensure that all determinants belong to $L^1.$
   
\begin{thm}[Proposition A.3, \cite{BM84}]Let $\Omega\ss\R^n$ be open and bounded, $f:\R\rightarrow\bar{\R}$ convex, lower semicontinuous and bounded below. Further, assume $v_j,v\in L^1(\Omega)$ for all $k\in\N$ and $v_j\rightharpoonup v$ in $D'(\Omega).$ Then 
\begin{equation}
\liminf\limits_{k\rightarrow\infty}\int\limits_\Omega{f(v_k)\;dx}\ge \int\limits_\Omega{f(v)\;dx}.
\label{eq:A2.4}
\end{equation}
\label{thm:A.2.1}
\end{thm}
\textbf{Proof:}\\ For a proof consult \cite{BM84}, Proposition A.3, p. 251-253.

\section{Nirenberg's lemma}

\begin{thm}[Nirenberg's lemma\footnote{Our proof follows closely the one given in \cite{GT} lemma 7.23-7.24. The statement can also be found in \cite{MG83}, p.45.}]
Let $\Omega\subset \R^n$ be an open domain.\\
i) Suppose $u\in W^{1,p}(\Omega,\R^N)$ with $1\le p<\infty.$ Then for all $\Om'\ss\ss\Om,$ $\Delta^{h,s}u\in L^p(\Om',\R^N)$ for all $0<h<\dist(\Om',\p\Om)$ and it holds
\begin{equation}
\|\Delta^{h,s}u\|_{L^p(\Omega',\R^N)}\le\|\p_s u\|_{L^p(\Omega,\R^N)}.
\label{eq:A2.22}
\end{equation}

ii) If $u\in L^p(\Omega,\R^N),$ $1<p<\infty,$ and there exists an $h_0>0$ and a constant $K>0,$ s.t. 
\[\sup_{0<h<h_0}\|\D^{h,s}u\|_{L^p(\Omega',\R^N)}\le K.
\label{eq:A2.24}\]
where the constant $K=K(\Om,u,h_0)$  is independent of $h$ but may depend on $\Om,\;u,\;h_0.$\\

Then $\p_s u\in L^p(\Om,\R^N)$ and
\[\|\p_s u\|_{L^p(\Om,\R^N)}\le K.\]
\label{A.Nir}
Moreover,
\begin{equation}
\Delta^{h,s}u\longrightarrow\p_s u\;\mbox{in}\; L^p(\Omega,\R^N).
\label{eq:A2.23}
\end{equation}
\end{thm}

\textbf{Proof.}
i) Fix $\Om'\ss\ss\Om$ and $0<h<\dist(\Om',\p\Om).$ By density, we can restrict to functions in the class $C^\infty(\Om,\R^N)\cap L^p(\Om,\R^N)$ (Note that $1\le p<\infty$). Then 
\begin{align}
\int\limits_{\Om'}{|\Delta^{h,s}u|^p\;dx}=\int\limits_{\Om'}{\left|\frac{1}{h}\int\limits_0^h {\p_su(x_1,\ldots,x_{s-1},x_{s}+t,x_{s+1},\ldots,x_{n})\;dt}\right|^p\;dx}\\
\le\int\limits_{\Om'}{\frac{1}{h}\int\limits_0^h {|\p_su(x_1,\ldots,x_{s-1},x_{s}+t,x_{s+1},\ldots,x_{n})|^p\;dt}\;dx},
\end{align}
where we used Jensen's inequality in the last line. By Fubini's theorem we are allowed to interchange the integrals. Furthermore, for all $0\le t\le h$ it holds
\[\int\limits_{\Om'}{|\p_su(x+te_s)|^p\;dx}=\int\limits_{\Om'+te_i}{|\p_su(x')|^p\;dx'}.\]
Then $\Om'+te_i\ss\Om$ for all $0\le t\le h$ and we get the estimate
\[\sup\limits_{0\le t\le h}\int\limits_{\Om'+te_i}{|\p_su(x')|^p\;dx'}\le\int\limits_{\Om}{|\p_su|^p\;dx}.\]
This integral does not depend on $t$ anymore. Hence,
\[\int\limits_{\Om'}{|\Delta^{h,s}u|^p\;dx}\le\frac{1}{h}\int\limits_0^h {\int\limits_{\Om'+te_i}{|\p_su(x'+te_s)|^p\;dx'}\;dt}\le\int\limits_{\Om}{|\p_su|^p\;dx}.\]
\\\vspace{0.5cm}

ii) By (\ref{eq:A2.24}), $\D^{h,s}u$ is a bounded sequence in $L^p(\Om,\R^N).$ By weak compactness in separable, reflexive Banach spaces (see, \cite{A12}, 6.11.(2), p.247) there exists a subsequence $h_k\rightarrow0$ for $k\ra\infty$ and $v\in L^p(\Om,\R^N)$ s.t. 
\[\D^{h_k,s}u\rhu v\;\mbox{in}\; L^p(\Om,\R^N)\;\mb{for}\;k\ra\infty.\]
For all $\vp\in C_c^\infty(\Om,\R^N),$ $\D^{-h,s}\vp$ converges uniformly to $\p_s \vp$ in $\Om$ and we have
\[\int\limits_{\Om}{v\cdot\vp\;dx}=\lim\limits_{k\ra\infty}\int\limits_{\Om'}{\D^{h_k,s}u\cdot\vp\;dx}=\lim\limits_{k\ra\infty}-\int\limits_{\Om'}{u(x)\cd\D^{-h_k,s}\vp(x)\;dx}=-\int\limits_{\Om}{u\cdot\p_s\vp\;dx}.\]
Then $v$ satisfies the definition of the weak derivative of $u,$ hence, $\p_su:=v\in L^p(\Om,\R^N).$ The estimate follows by the wlsc. of the norm,\footnote{see, again \cite{A12}, 6.2(4), p.239-240.}  i.e.
\[\|\p_s u\|_{L^p(\Om,\R^N)}\le\liminf\limits_{k\ra\infty}\|\D^{h_k,s} u\|_{L^p(\Om',\R^N)}\le K.\]
Combining this, with the mean value theorem we get
\[\|\p_su\|_{L^p(\Om,\R^N)}\le\liminf\limits_{h\ra0}\|\D^{h,s}u\|_{L^p(\Om',\R^N)}\le\limsup\limits_{h\ra0}\|\D^{h,s}u\|_{L^p(\Om',\R^N)}\le\|\p_su\|_{L^p(\Om,\R^N)}\]
showing $\|\D^{h,s}u\|_{L^p(\Om',\R^N)}\ra\|\p_su\|_{L^p(\Om,\R^N)}$ if $h\ra0$. Recalling that for $1<p<\infty,$\footnote{See, \cite{A12}, Ü6.5-Ü6.6, p.263-265.}
\[f_k\ra f\;\mb{in}\;L^p\Llra f_k\rhu f\;\mb{in}\;L^p\; \mb{and}\; \|f_k\|_{L^p}\ra\|f\|_{L^p}\]
if $k\ra\infty,$ implies $\D^{h,s}u\ra \p_su$ in $L^p(\Om,\R^N).$

\section{Behaviour of the integrand for small $\gamma<1.$}
\label{Ap:A.3}
Here we validate that our integrand is indeed uniformly convex if $\ga<1.$ This can be seen, by showing that $\grad W$ satisfies a monotonicity inequality.
\begin{lem} Let $W(A):=\frac{1}{2}|A|^2+\rho(d_A)$ for $A\in \R^{2\ti2}.$ Then for an arbitrary $\ga>0$ it holds 
\begin{equation}
(\grad_A W(A)-\grad_A W(B))\cdot(A-B)\ge(1-\ga)|A-B|^2 \;\;\mbox{for all}\;A,B\in \R^{2\times2}.
\label{eq:2.21}
\end{equation}
\end{lem}
\textbf{Proof:} We calculate
\begin{eqnarray*}
(\grad_A W(A)-\grad_A W(B))\cdot(A-B)&=&|A-B|^2\\
&&+(\rho'(d_A)\cof A-\rho'(d_B)\cof B)\cdot(A-B)\\
&=&|A-B|^2+\rho'(d_A)((d_A-d_B)+d_{A-B})\\
&&-\rho'(d_B)((d_{A}-d_{B})-d_{A-B})\\
&=&|A-B|^2+(\rho'(d_A)-\rho'(d_B))(d_A-d_B)\hspace{1cm}\\
&&+(\rho'(d_A)+\rho'(d_B))d_{A-B}\\
&\ge&(1-\ga)|A-B|^2
\label{eq:2.22}
\end{eqnarray*}
where we used again
\begin{eqnarray*}
\cof B\cdot (A-B)=(d_{A}-d_{B})-d_{A-B}
\label{eq:2.23}
\end{eqnarray*}
and
\begin{eqnarray*}
d_B=d_{-B}=d_{(A-B)-A}=d_A+d_{A-B}-\cof A\cdot(A-B)
\label{eq:2.24}
\end{eqnarray*}
which can be written as
\begin{eqnarray*}
\cof A\cdot(A-B)=(d_A-d_B)+d_{A-B}.
\label{eq:2.25}
\end{eqnarray*}
If $\ga<1$ the constant $1-\ga$ in the above inequality becomes strictly positive, yielding uniform convexity of the integrand. 

\section{A Reverse Fatou's Lemma for non-negative functions}
In a preliminary version we made use of the version of the Reverse Fatou's Lemma as given below. It only applies to sequences of non-negative functions, however a proof is elementary. We state it here, since it could be useful to some.\footnote{I found this version and the proof on the famous website stackexchange.com, see \cite{stack} for details. I was not able to find this statement else were. Firstly, we thank user "Kaa1el" for providing this statement and giving a detailed proof. Secondly, we include a detailed proof here for the convenience of the reader and to make sure everything is indeed sound. Again, I do not claim any originality whatsoever.}

\begin{lem} [A Reverse Fatou's Lemma for non-negative functions]
Let $(X,\Sigma,\mu)$ be a measure space. Let $(g_n)_{n\in \N}\ss L^1(X,[0,\infty],\mu)$ be a sequence, which converges strongly in $L^1(X,[0,\infty],\mu)$ to a function $g\in L^1(X,[0,\infty],\mu).$
Suppose further that there is another sequence of measurable functions $f_n: X\ra[0,\infty]$ s.t. $f_n\le g_n$ a.e. for every $n\in\N.$\\

Then it holds
\begin{align}\limsup\limits_{n\ra\infty}\int\limits_{X}{f_n\;d\mu}\le \int\limits_{X}{\limsup\limits_{n\ra\infty} f_n\;d\mu}.\end{align}
\end{lem}
\textbf{Proof:}\\
We start by realising first that for every $n\in\N$ it holds
\begin{align*}0\le g=g-g_n+g_n\le |g-g_n|+g_n\end{align*}
and
\begin{align*}0\le g_n-f_n=g_n-g+g-f_n\le |g-g_n|+g-f_n=:h_n \;\mb{a.e.}\end{align*}

Then $(h_n)_{n\in \N}$ is a sequence of non-negative measurable functions and we can apply the standard version of Fatou's Lemma to obtain\footnote{see \cite[Cor 5.34]{WZ15}.}
\begin{align}\int\limits_{X}{\liminf\limits_{n\ra\infty} h_n\;d\mu}\le \liminf\limits_{n\ra\infty}\int\limits_{X}{h_n\;d\mu}.\label{eq:2.3.RF1.1}\end{align}
Using the superadditivity of the limes inferior and the fact that $\liminf\limits_{n\ra\infty}(-f_n)=-\limsup\limits_{n\ra\infty}(f_n)$ we can estimate the LHS of \eqref{eq:2.3.RF1.1} from below by
\begin{align} \int\limits_{X}{\liminf\limits_{n\ra\infty} |g-g_n|\;d\mu}+\int\limits_{X}{g\;d\mu}-\int\limits_{X}{\limsup\limits_{n\ra\infty}(f_n)\;d\mu}\le \int\limits_{X}{\liminf\limits_{n\ra\infty} h_n\;d\mu}.\label{eq:2.3.RF1.2}\end{align}
Now by applying Fatou's Lemma to $|g-g_n|$ we get
\begin{align*} 0\le\int\limits_{X}{\liminf\limits_{n\ra\infty} |g-g_n|\;d\mu}\le\liminf\limits_{n\ra\infty} \int\limits_{X}{|g-g_n|\;d\mu}=0,\end{align*}
eliminating the leftmost term of \eqref{eq:2.3.RF1.2}.\\

On the other hand, recall that the superadditivity of the limes inferior becomes an additivity if one of the sequences converges and we know $|g-g_n|+g\ra g$ strongly in $L^1.$
Then the RHS of \eqref{eq:2.3.RF1.1}  becomes
\begin{align*}  \liminf\limits_{n\ra\infty}\int\limits_{X}{h_n\;d\mu}=\int\limits_{X}{g\;d\mu}+\liminf\limits_{n\ra\infty}\int\limits_{X}{(-f_n)\;d\mu}
=\int\limits_{X}{g\;d\mu}-\limsup\limits_{n\ra\infty}\int\limits_{X}{f_n\;d\mu}.\end{align*}

Together this yields
\begin{align*}\limsup\limits_{n\ra\infty}\int\limits_{X}{f_n\;d\mu}\le \int\limits_{X}{\limsup\limits_{n\ra\infty} f_n\;d\mu},\end{align*}
completing the proof.\vspace{0.5cm}

\clearpage{\pagestyle{empty}\cleardoublepage} 

\clearemptydoublepage






\fancyhead[RO,LE]{\thepage}
\fancyhead[CO,CE]{}
\fancyhead[LO,RE]{References}
\fancyfoot{}
\bibliography{LiteraturePhDMD}



\clearemptydoublepage




\clearemptydoublepage


\end{document}